\documentclass[12pt,final,dvips]{my_nmeauth}

\usepackage{etex} 
\usepackage{graphicx,subfigure,color,colordvi}
\usepackage{latexsym,amsbsy,epsfig,fancybox}
\usepackage{amstext}
\usepackage{amsfonts,textgreek}
\usepackage{mathrsfs}
\usepackage{url}
\usepackage{wrapfig,bm}
\usepackage{caption}
\usepackage{psfrag}
\usepackage{upgreek}
\usepackage{isomath}
\usepackage[dvips]{geometry}
\usepackage{latexsym}
\usepackage{graphicx}
\usepackage{array}
\usepackage{multirow}
\usepackage{booktabs}
\usepackage{colortbl}
\usepackage{verbatim}
\usepackage{hypernat}
\usepackage[colorlinks=true,citecolor=blue]{hyperref}
\usepackage{epstopdf,overpic}
\usepackage{enumerate}
\usepackage{amsmath,amssymb,xspace}

\usepackage[algoruled]{algorithm2e}

\newcommand{\eref}[1]{(\ref{#1})}
\newcommand{\fref}[1]{Figure~\ref{#1}}

\newcommand{\vm}[1]{\bm{#1}} 
\newcommand{\bsym}[1]{\bm{#1}}

\renewcommand{\Re}{{\rm{I\!R}}}
\newcommand{\vx}{\vm{x}}
\newcommand{\mat}[1]{\bm{#1}} 
\newcommand{\diffx}{d\vm{x}}
\newcommand{\diffs}{ds}
\newcommand{\transpose}{\mathsf{T}}
\newcommand{\tangent}{\mathsf{T}}
\newcommand{\smat}[2][ccccccccccccccccccccccccccccccccccccccccccccccccccc]{\left
[\begin{array}{#1}#2 \\ \end{array} \right]}

\newcommand{\cons}{\mathrm{c}}
\newcommand{\stab}{\mathrm{s}}

\renewcommand{\arraystretch}{0.7} 

\definecolor{forestgreen}{RGB}{34, 139, 34}

\begin{document}


\def\baselinestretch{1} 
\NME{1}{00}{0}{00}{00}

\runningheads{Silva-Valenzuela et al.}{A nodal integration scheme for meshfree Galerkin methods}

\noreceived{} \norevised{} \noaccepted{}

\title{A nodal integration scheme for meshfree Galerkin methods using the virtual element decomposition}

\author{R. Silva-Valenzuela\affil{1,2,3},
        A. Ortiz-Bernardin\affil{1,2}\comma\corrauth,
        N. Sukumar\affil{4},
        E. Artioli\affil{5} and
        N.~Hitschfeld-Kahler\affil{6,7}}

\address{\affilnum{1}\ Department of Mechanical Engineering,
                       Universidad de Chile,
                       Av.\ Beauchef 851, Santiago 8370456, Chile\\
\affilnum{2}\          Computational and Applied Mechanics Laboratory,
                       Center for Modern Computational Engineering,
                       Facultad de Ciencias F\'isicas y Matem\'aticas,
                       Universidad de Chile,
                       Av. Beauchef 851, Santiago 8370456, Chile\\
\affilnum{3}\          Department of Mechanical Engineering,
                       Universidad de La Serena,
                       Av. Benavente 980, La Serena 1720170, Chile\\                       
\affilnum{4}\          Department of Civil and Environmental Engineering,
                       University of California,
                       Davis, CA 95616, USA\\
\affilnum{5}\          Department of Civil Engineering and Computer Science,
                       University of Rome Tor Vergata,
                       Via del Politecnico 1, 00133 Rome, Italy\\
\affilnum{6}\          Department of Computer Science,
                       Universidad de Chile,
                       Av. Beauchef 851, Santiago 8370456, Chile\\
\affilnum{7}\          Meshing for Applied Science Laboratory,
                       Center for Modern Computational Engineering,
                       Facultad de Ciencias F\'isicas y Matem\'aticas,
                       Universidad de Chile,
                       Av. Beauchef 851, Santiago 8370456, Chile}
\corraddr{A. Ortiz-Bernardin,
          Department of Mechanical Engineering, Universidad de Chile,
          Av. Beauchef 851, Santiago 8370456, Chile.  E-mail:
          \url{aortizb@uchile.cl}}

\begin{abstract}
In this paper, we present a novel nodal integration scheme for meshfree 
Galerkin methods that draws on the mathematical framework of the virtual element 
method. We adopt linear maximum-entropy basis functions for the discretization 
of field variables, although the proposed scheme is applicable 
to any linear meshfree approximant. In our approach, the weak form integrals
are nodally integrated using nodal representative cells that carry the nodal 
displacements and state variables such as strains and stresses. The nodal 
integration is performed using the virtual element decomposition, wherein 
the bilinear form is decomposed into a consistency part and a stability part 
that ensure consistency and stability of the method. The performance 
of the proposed nodal integration scheme is assessed through 
benchmark problems in linear and nonlinear analyses of solids
for small displacements and small-strain kinematics. Numerical results 
are presented for linear elastostatics and linear elastodynamics, 
and viscoelasticity. We demonstrate that the proposed nodally 
integrated meshfree method is accurate, converges optimally, 
and is more reliable and robust than a standard cell-based 
Gauss integrated meshfree method.
\end{abstract}


\section{INTRODUCTION} \label{sec:intro}

Meshfree Galerkin methods are 
based on the weak form where the field 
variables are discretized using basis functions associated
with a set of scattered nodes that partition the domain
of analysis. Since the inception of the element-free Galerkin
method~\cite{belytschko:1994:EFG}, efficient and accurate numerical integration 
of the weak form integrals has attracted 
broad attention in meshfree Galerkin methods.
The different numerical integration techniques that exist 
for meshfree Galerkin methods can be grouped into two main
approaches: cell integration and nodal integration techniques. Nodal integration 
in meshfree methods is attractive since state variables (strains, stresses, and 
internal variables) can be stored at the nodes, thereby avoiding the need for 
remapping algorithms at integration points that arise in traditional Lagrangian
finite element large deformation simulations with remeshing.  In the optimal 
transportation meshfree method~\cite{li:OTM:2010}, remeshing is avoided by using 
the integration points as the carriers of material state information in tandem 
with changes in the support of the basis functions to enable Lagrangian finite 
deformation simulations. In the material point method~\cite{sulsky:MPM:1994}, 
integration points of an initial grid (usually, a Cartesian grid) carry the material 
state variables and the solution stage occurs in three phases: (a) information is 
mapped from particles to grid nodes; (b) equations of motion are solved on the 
grid nodes, and the updated information is mapped back to the particles 
to update their positions and velocities; and (c) the grid is reset. 
In this process, the mesh remains fixed throughout the simulation, and the 
deformation causes the material points to move within other cells in the grid.
In this paper, we use the virtual element decomposition~\cite{BeiraoDaVeiga-Brezzi-Cangiani-Manzini-Marini-Russo:2013} 
and follow Reference~\cite{ortiz:CSMGMVEM:2017} to devise an accurate and stable nodally integrated 
meshfree method for linear elastostatics and linear elastodynamics, as well as
nonlinear viscoelasticity.

The most simple cell integration technique requires the construction 
of nonoverlapping cells on which standard Gauss integration 
is performed. However, this integration scheme is inexact due to the 
following two properties of meshfree basis functions~\cite{dolbow:1999:NIO}: 
(1) They are nonpolynomial functions; 
and (2)	in general, the region that is defined by the intersecting supports 
of two overlapping nodal basis functions does not
coincide with the integration cell. These are two issues that often lead to consistency
errors (patch test is not passed) and stability problems due to underintegration.
Various approaches have been put forth in meshfree Galerkin 
methods to address integration errors on cells. For instance,
higher-order tensor-product Gauss quadrature is adopted in 
the element-free Galerkin method~\cite{belytschko:1994:EFG,dolbow:1999:NIO},
whereas the support of the nodal basis functions is used 
as the domain of integration in the meshless local 
Petrov-Galerkin method~\cite{atluri:1999:ACS} and also 
in the method of finite spheres~\cite{de:2001:FSI}.  More details on 
cell-based integration schemes in meshfree Galerkin methods can be 
found in Ortiz-Bernardin et al.~\cite{ortiz:CSMGMVEM:2017}.

Nodal integration techniques perform the integration of the weak
form integrals by sampling them at the nodes. This approach
requires the construction of nodal cells
that represent the volume of the integrals being sampled 
at the nodes. Nodal integration techniques are also
prone to integration errors and require to be stabilized.
A direct nodal integration (1-point) scheme leads to 
rank instabilities because meshfree basis functions 
have zero or nearly zero derivatives at the nodes. 
In an effort to stabilize the rank instability, 
Beissel and Belytschko~\cite{beissel:1996:NIO} 
introduced a least-squares residual-based method where the second-order 
derivatives stabilize the rank instability. With the aim of computing 
nodal derivatives away from the nodes to avoid rank instabilities, 
Chen et al.~\cite{chen:2001:ASC} devised a strain smoothing procedure 
known as stabilized conforming nodal integration (SCNI) approach, 
which is the cornerstone of various 
meshfree nodal-based~\cite{chen:2002:NLV,puso:2008:MAF,chen:2013:VCI} 
and cell-based~\cite{duan:2012:SOI,duan:2014:CEF,duan:2014:FPI} integration methods, 
and even smoothed finite element methods~\cite{liu:SFEM:2007}. Another method to 
compute derivatives away from the nodes is the stress-point method, which 
was first introduced in the smoothed particle hydrodynamics (SPH) meshfree 
method~\cite{dyka:SPTI:1997}. Although the SCNI approach~\cite{chen:2001:ASC} 
suppresses the rank instability and provides patch test satisfaction, 
instabilities due to nonzero low-energy modes are still 
encountered~\cite{puso:2008:MAF}. Puso et al.~\cite{puso:2008:MAF} 
proposed a penalty stabilization to stabilize the SCNI 
method. Hillman and Chen~\cite{hillman:ACSNI:2016}
combined the arbitrary order variationally consistent meshfree 
nodal integration framework~\cite{chen:2013:VCI} with a 
stabilization method devised from an implicit gradient expansion 
of the strains at the nodes. The resulting 
method is first-order variationally consistent and stable, and 
unlike the stabilization method of Puso et al.~\cite{puso:2008:MAF}, 
it is devoid of tunable parameters.

Recently, the virtual element
method~\cite{BeiraoDaVeiga-Brezzi-Cangiani-Manzini-Marini-Russo:2013,BeiraodaVeiga-Brezzi-Marini:2013}
(VEM) has been proposed, where an exact algebraic construction of
the stiffness matrix is realized without the explicit construction of
basis functions (basis functions are deemed as \textit{virtual}). In the VEM, the
stiffness matrix is decomposed into two parts: a consistent term
that reproduces a given polynomial space and a correction term
that provides stability. Such a decomposition (herein referred to as
the \textit{virtual element decomposition}) is formulated in the
spirit of the Lax equivalence theorem (consistency $+$ stability
$\to$ convergence) for finite-difference schemes and is sufficient
for the method to pass the patch test.  In polygonal 
and polyhedral finite elements, Talischi and Paulino~\cite{Talischi-Paulino:2014}, 
Gain et al.~\cite{Gain-Talischi-Paulino:2014}
and Manzini et al.~\cite{Manzini-Russo-Sukumar:2014} have
used the virtual element decomposition to pass the patch test. 
The VEM can be viewed as a stabilized Galerkin method on 
polytopal meshes~\cite{cangiani:2015}.
For meshfree Galerkin methods with cell-based integration,
Ortiz-Bernardin et al.~\cite{ortiz:CSMGMVEM:2017} used the virtual element 
decomposition to develop a method for linear elastostatics that is consistent and stable.

In this paper, on using the virtual element decomposition, a novel meshfree 
nodal integration technique for linear and nonlinear analyses of solids
for small displacements and small-strain kinematics is presented. In particular,
linear elastostatics, linear elastodynamics and viscoelasticity 
are considered. We use the acronym NIVED to refer to this method.
The distinctions and contributions of the present work vis-{\`a}-vis 
our prior work on meshfree Galerkin methods
using the virtual element decomposition~\cite{ortiz:CSMGMVEM:2017} are as follows:
\begin{itemize}
  \item focus herein is on nodal integration, whereas Gauss-like integration scheme is used
        in Ortiz-Bernardin et al.~\cite{ortiz:CSMGMVEM:2017};
  \item we develop the method for linear (static and dynamic) and nonlinear 
        analyses of solids for small displacements and small-strain 
        kinematics, whereas only linear elastostatics is treated in the prior work; and
  \item unlike the prior work, the stabilization of the stiffness matrix       in the nodal integration
        method does not require a tunable (scaling) parameter.
\end{itemize}
We consider maximum-entropy
basis functions (Section~\ref{sec:maxent}),
although the formulation is applicable to any linear meshfree approximant. 
The governing equations for linear elastostatics are described
in Section~\ref{sec:goveqn}. The main ingredients of the virtual
element framework and the development of
our proposed nodal integration technique are presented
in Section~\ref{sec:nived} for linear elasticity and 
in Section~\ref{sec:nonlinear_nived} for nonlinear analysis
of viscoelastic solids. The demonstration of the patch 
test satisfaction is provided in Section~\ref{sec:patch_test}.
The accuracy and convergence of the devised nodal integration method 
are assessed in Section~\ref{sec:num_experiments} 
through several examples in linear elastostatics, linear elastodynamics 
and viscoelasticity. A summary of our main findings 
with some final remarks is provided in Section~\ref{sec:conclusions}.

\section{MAXIMUM-ENTROPY BASIS FUNCTIONS} \label{sec:maxent}
Consider a convex domain represented by a set of $N$
scattered nodes and a prior (weight) function $w_{a}(\vm{x})$
associated with each node $a$.  We can write down the
approximation for a vector-valued 
function $\vm{v}(\vm{x})$ in the form:
\begin{equation}\label{eq:trial}
\vm{v}^h(\vm{x})=\sum_{a=1}^m\phi_a(\vm{x}) \vm{v}_a,
\end{equation}
where $\vm{v}_a:=\vm{v}(\vm{x}_a)$ are nodal 
coefficients, $\phi_a(\vm{x})$ is the meshfree basis function associated
with node $a$ and $m\leq N$ represents the number of nodes whose basis 
functions take a nonzero value at the point $\vm{x}$.

On using the Shannon-Jaynes (or relative) entropy
functional, the maximum-entropy basis functions $\{\phi_{a}(\vm{x})\geq 0\}_{a=1}^m$
are obtained via the solution of the following convex optimization
problem~\cite{sukumar:2007:OAC}:
\begin{subequations}\label{eq:maxent_problem}
\begin{align}
\min_{\boldsymbol{\phi} \in {\Re}_{+}^{m}} \sum_{a=1}^{m} &
\phi_{a}(\vm{x})\ln\left(\frac{\phi_{a}(\vm{x})}{w_{a}(\vm{x})}\right)
\\ \intertext{subject to the linear reproducing conditions:}
\sum_{a=1}^{m} \phi_{a}(\vm{x}) = 1 & \quad \sum_{a=1}^{m}
\phi_{a}(\vm{x})\,\vm{c}_a  = \vm{0},
\end{align}
\end{subequations}
where $\vm{c}_a=\vx_a - \vx$ are shifted nodal coordinates
and $\Re_{+}^{m}$ is the nonnegative orthant.  In this paper,
we use as the prior weight function either the Gaussian radial 
basis function given by~\cite{arroyo:2006:LME}
\begin{equation*}\label{eq:priors}
w_{a}(\vm{x}) = \exp\left(-\frac{\gamma}{h_{a}^2}\|\vm{c}_a\|^2\right),
\end{equation*}
where $\gamma$ is a parameter that controls the support 
size of the basis function and $h_a$
is a characteristic nodal spacing associated with node $a$, 
or the $C^2$ quartic polynomials given by~\cite{yau:2009:MCR}
\begin{equation}
w_a(q) = \left\{\begin{array}{lr}
    1-6q^{2}+8q^{3}-3q^{4} & \quad 0\leq q\leq 1 \\
    0 & \quad q > 1
    \end{array}
\right.,
\end{equation}
where $q=\|\vm{c}_a\|/(\gamma h_{a})$.

On using the method of Lagrange multipliers, the
solution to~\eref{eq:maxent_problem} is given
by~\cite{sukumar:2007:OAC}
\begin{equation}\label{eq:maxent_bfun}
\phi_{a}(\vm{x},\boldsymbol{\lambda}) =
\frac{w_{a}(\vm{x}) \exp (- \boldsymbol{\lambda}(\vm{x})\cdot
\vm{c}_a(\vm{x}))}{Z(\vm{x},\boldsymbol{\lambda}(\vm{x})) },
\quad  Z(\vm{x},\boldsymbol{\lambda}(\vm{x}))=\sum_{b=1}^m
w_{b}(\vm{x}) \exp (- \boldsymbol{\lambda}(\vm{x})\cdot
\vm{c}_b(\vm{x})),
\end{equation}
where the Lagrange multiplier vector $\boldsymbol{\lambda}(\vm{x})$ is
obtained as the minimizer of the dual optimization problem
($\vm{x}$ is fixed):
\begin{equation*}\label{eq:maxent_dual}
\boldsymbol{\lambda}^*(\vm{x}) = \arg \min_{\boldsymbol{\lambda}
\in {\Re}^{d}} \ln Z(\vm{x},\boldsymbol{\lambda}),
\end{equation*}
where $\boldsymbol{\lambda}^*$ is the converged solution that
gives the basis functions as $\phi_a(\vm{x})=\phi_{a}(\vm{x},\boldsymbol{\lambda}^*)$
for $a=1,\ldots,m$.

\section{GOVERNING EQUATIONS} \label{sec:goveqn}

\subsection{Strong form} \label{sec:strongform}

Consider an elastic body that occupies the
open domain $\Omega \subset \Re^2$ and is
bounded by the one-dimensional surface $\Gamma$ whose
unit outward normal is $\vm{n}_{\Gamma}$. The boundary is assumed
to admit decompositions $\Gamma=\Gamma_u\cup\Gamma_t$ and
$\emptyset=\Gamma_u\cap\Gamma_t$, where $\Gamma_u$ is the Dirichlet
boundary and $\Gamma_t$ is the Neumann boundary. The closure of
the domain is $\overline{\Omega}=\Omega\cup\Gamma$. Let
$\vm{u}(\vm{x}) : \overline{\Omega} \rightarrow \Re^2$ be
the displacement field at a point $\vm{x}$ of the elastic
body when the body is subjected to external tractions
$\bar{\vm{t}}(\vm{x}):\Gamma_t\rightarrow \Re^2$ and body forces $\vm{b}(\vm{x}):\Omega\rightarrow\Re^2$.
The imposed Dirichlet (essential) boundary conditions are
$\bar{\vm{u}}(\vm{x}):\Gamma_u\rightarrow \Re^2$. The boundary-value
problem for linear elastostatics is: find
$\vm{u}(\vm{x}): \overline{\Omega} \rightarrow \Re^2$
such that
\begin{subequations}\label{eq:strongform}
\begin{align}
\bsym{\nabla} \cdot \bsym{\sigma} + \vm{b} &= 0 \quad \mathrm{in}\,\,\Omega, \\ 
\vm{u} &= \bar{\vm{u}} \quad \mathrm{on}\,\,\Gamma_u, \\ 
\bsym{\sigma} \cdot \vm{n}_{\Gamma} &= \bar{\vm{t}} \quad \mathrm{on}\,\,\Gamma_t,
\end{align}
\end{subequations}
where $\bsym{\sigma}$ is the Cauchy stress
tensor.

\subsection{Weak form} \label{sec:weakform}

The Galerkin weak formulation reads:
find $\vm{u}(\vm{x})\in \mathscr{V}$ such that
\begin{equation}\label{eq:weakform}
\begin{split}
a(\vm{u},\vm{v}) & = \ell(\vm{v}) \quad \forall \vm{v}(\vm{x})\in \mathscr{W}, \\
a(\vm{u},\vm{v})=\int_{\Omega}\bsym{\sigma}(\vm{u}):\bsym{\varepsilon}(\vm{v})\,\diffx,
& \quad
\ell(\vm{v}) = \int_{\Omega}\vm{b}\cdot\vm{v}\,\diffx + \int_{\Gamma_t}\bar{\vm{t}}\cdot\vm{v}\,\diffs,
\end{split}
\end{equation}
where $\mathscr{V}$ and $\mathscr{W}$ are the 
continuous displacement trial and test spaces:
\begin{align*}
\mathscr{V} :=
\left\{\vm{u}(\vm{x}): \vm{u} \in [H^{1}(\Omega)]^2, \ \vm{u} = \bar{\vm{u}}
\ \textrm{on } \Gamma_u \right\}, \quad
\mathscr{W} := \left\{\vm{v}(\vm{x}): \vm{v} \in [H^{1}(\Omega)]^2, \ \vm{v} = \vm{0} \ \textrm{on } \Gamma_u
\right\},
\end{align*}
and $\bsym{\varepsilon}(\vm{v})$ is the small strain tensor 
that is given by
\begin{equation}
\bsym{\varepsilon}(\vm{v}) = \frac{1}{2}\left(\bsym{\nabla}\vm{v}
+\bsym{\nabla}^\transpose\vm{v}\right).
\label{eq:symstrain}
\end{equation}

\begin{figure}[!tbh]
  \centering
  \psfrag{E}{\Large $E$} 
  \psfrag{x}{\Large $\vm{x}_E$}
  \mbox{
  \subfigure[]{\label{fig:nodalpartition_a}
  \epsfig{file = 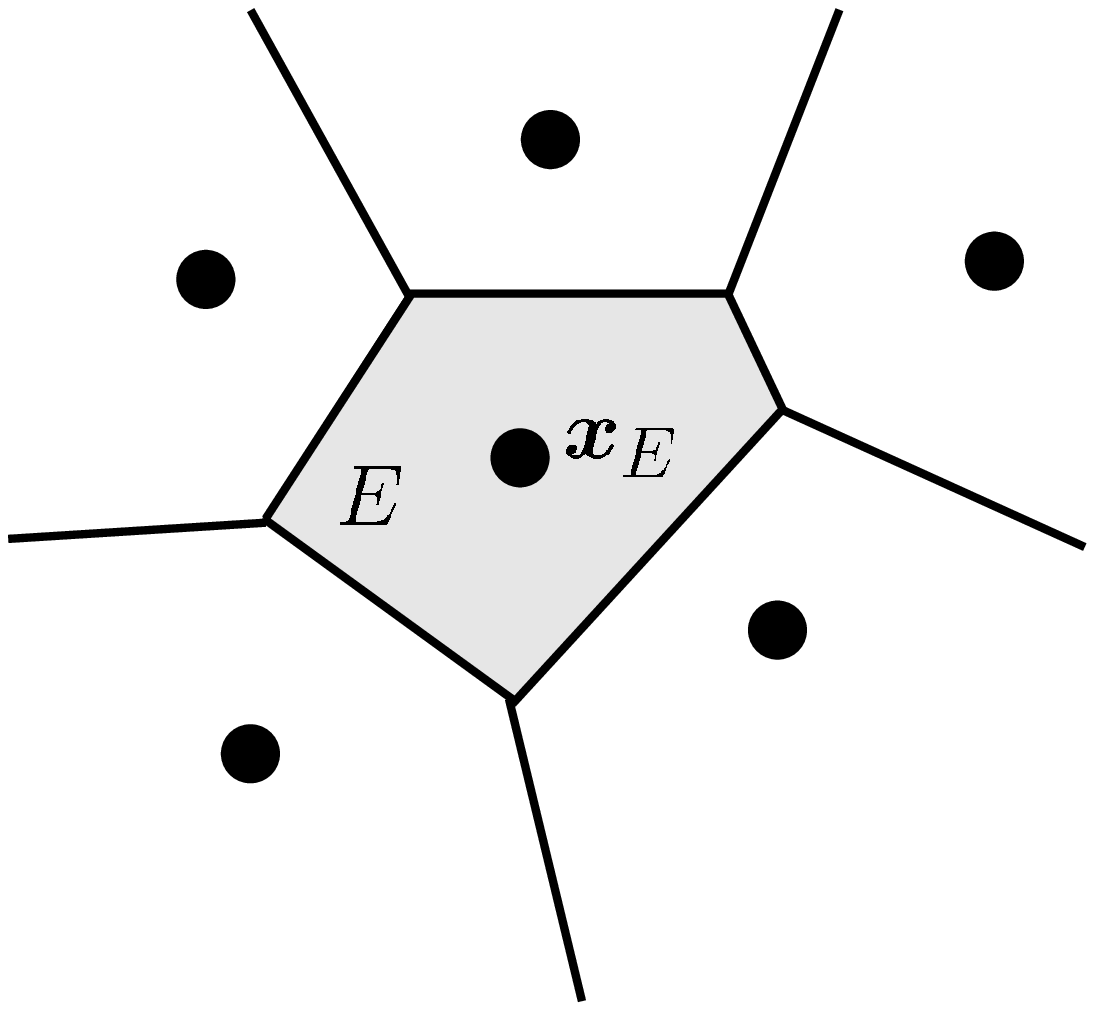,
  width = 0.35\textwidth}}
  \hspace{0.1\textwidth}
  \subfigure[]{\label{fig:nodalpartition_b}
  \epsfig{file = 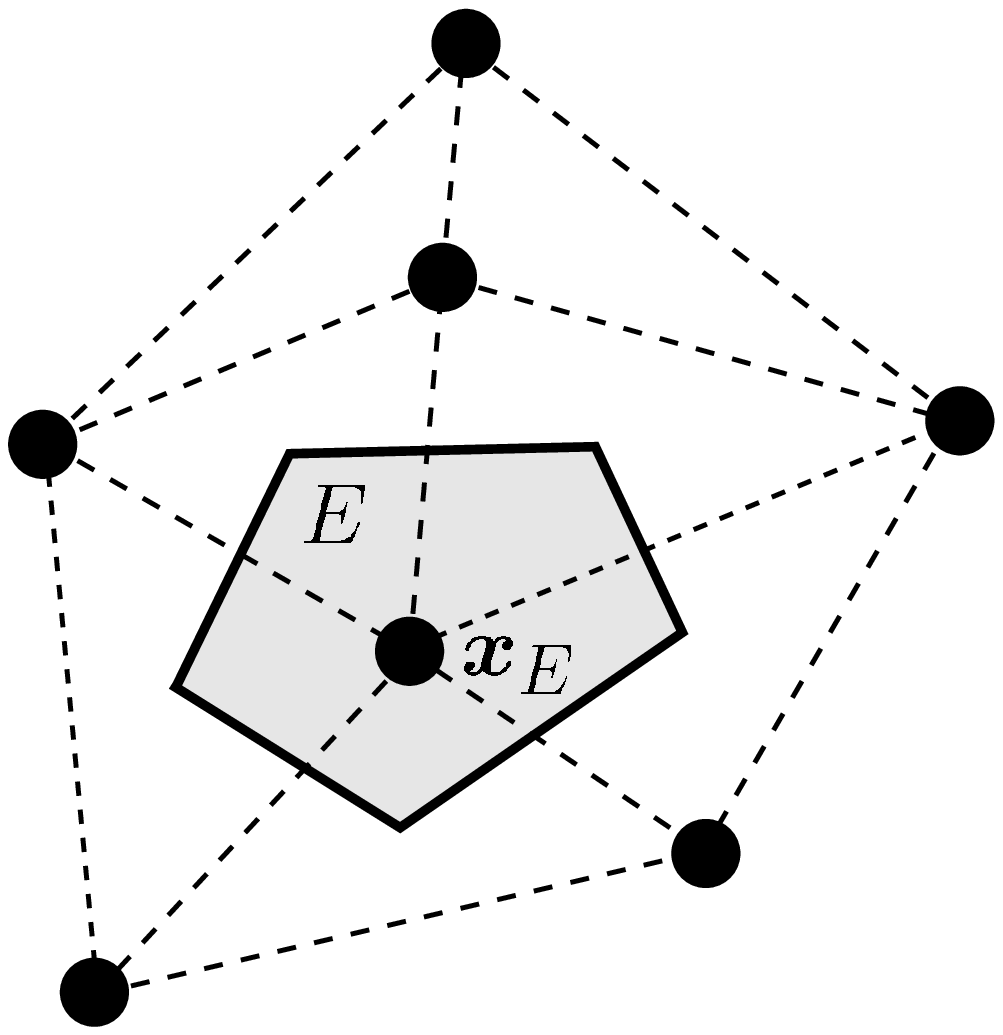,
  width = 0.3\textwidth}}
  }
  \caption{Construction of a nodal representative domain (shaded area) 
   using (a) a Voronoi diagram and (b) a triangular mesh. The node
   and its representative cell are denoted by $E$ and $\vm{x}_E$ are
   the nodal coordinates.}
  \label{fig:nodalpartition}
\end{figure}

\begin{figure}[!tbh]
\centering 
\mbox{
\subfigure[]{\label{fig:neumann_nodalpartition_a}
\psfrag{g}{\large $\Gamma_t$} 
\psfrag{s}{\large $|S|$} 
\psfrag{x}{\large $\vm{x}_S$}
\epsfig{file = 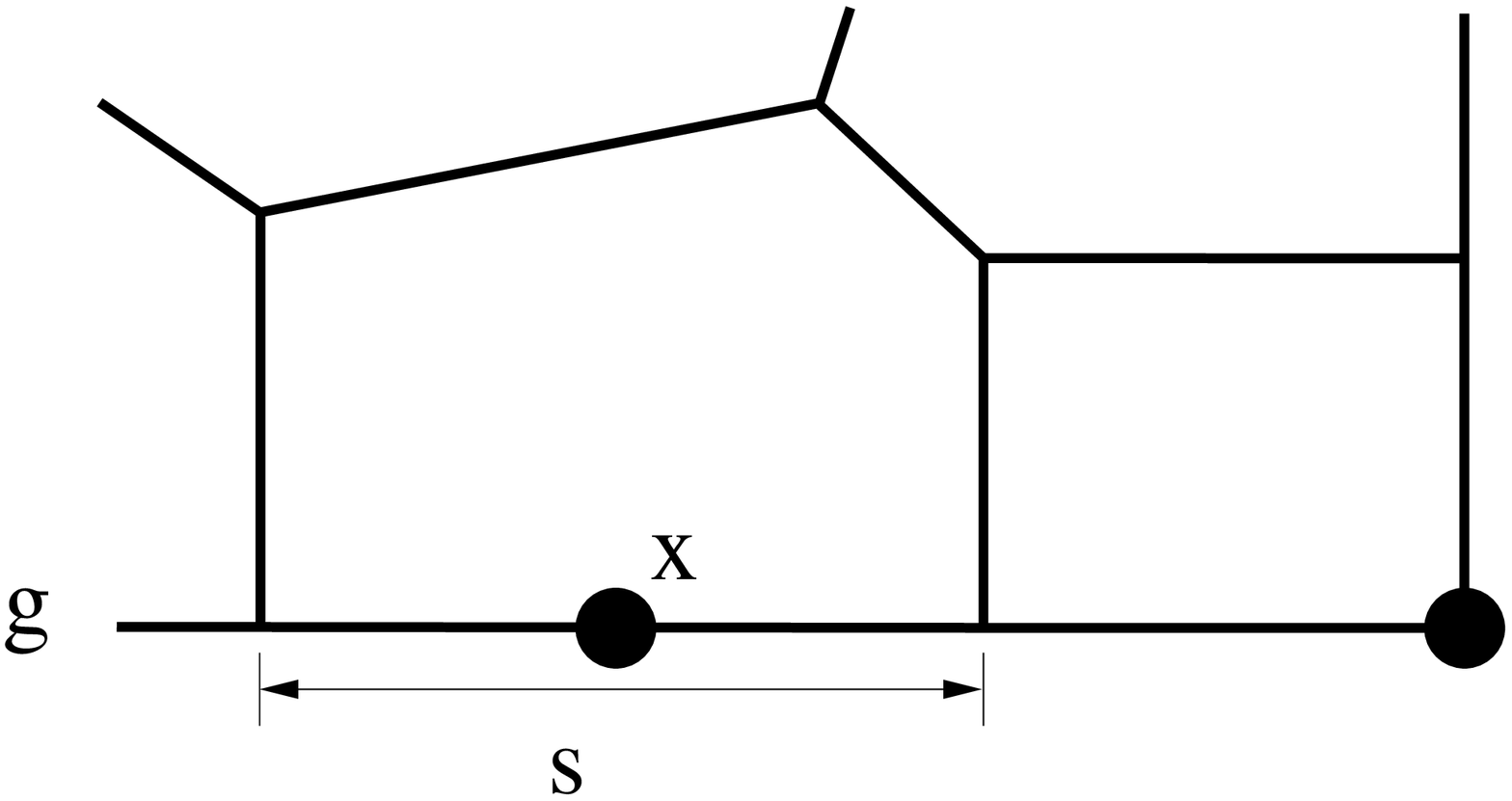,
width = 0.45\textwidth}}
\subfigure[]{\label{fig:neumann_nodalpartition_b}
\psfrag{g}{\large $\Gamma_t$} 
\psfrag{s}{\large $|S|$} 
\psfrag{x}{\large $\vm{x}_S$}
\epsfig{file = 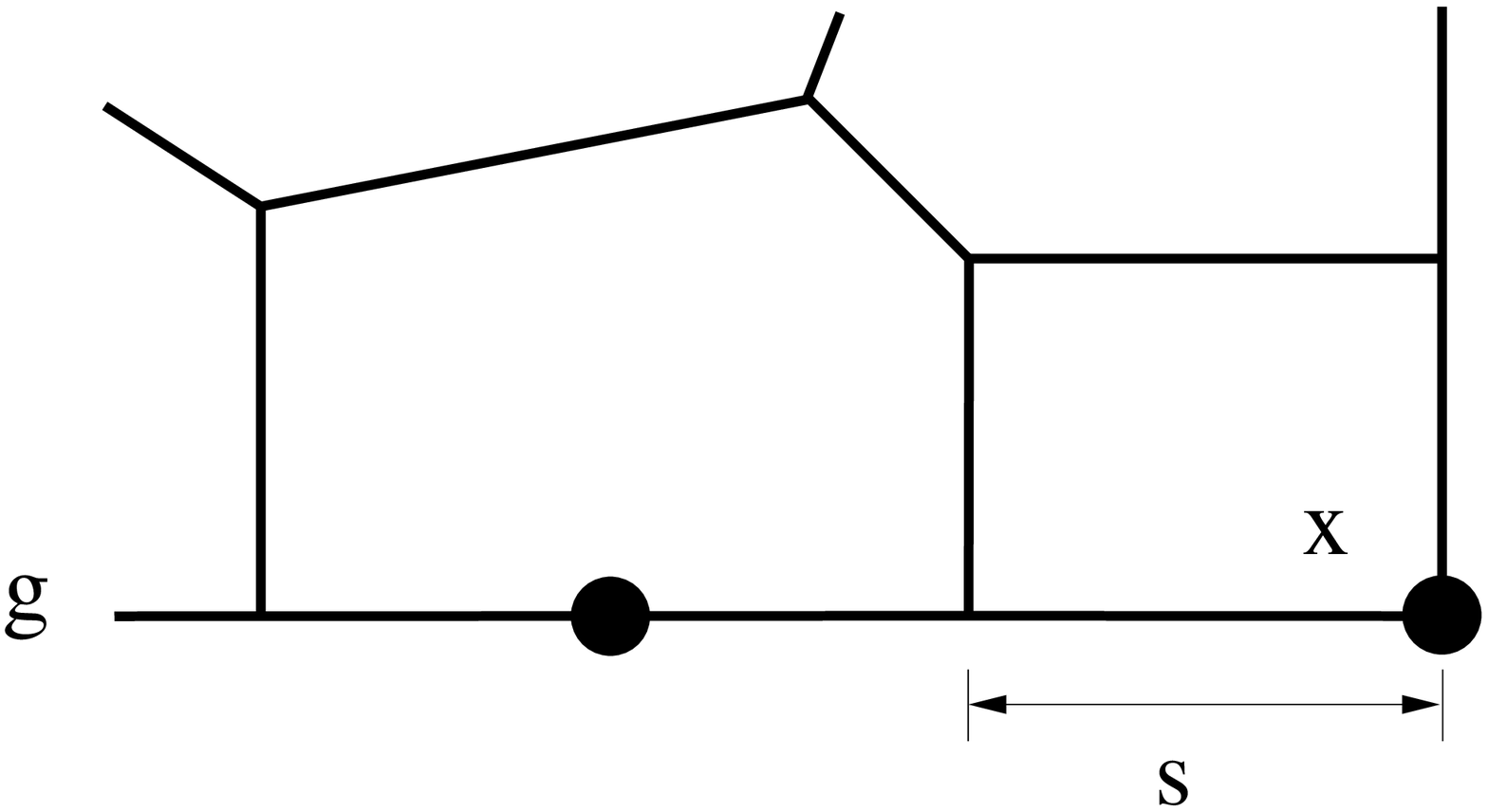,
width = 0.45\textwidth}}
}
\caption{Schematic representation of a node located on $\Gamma_t$ with coordinate $\vm{x}_S$ 
and associated length of influence $|S|$. (a) A node located inside a boundary edge and
(b) a node located at a boundary corner.} 
\label{fig:neumann_nodalpartition}
\end{figure}

We consider a set of scattered nodes that discretizes $\bar{\Omega}$.
Using these nodes, the domain $\Omega$ is 
partitioned into nonoverlapping nodal representative polygonal cells. 
This can be achieved using a Voronoi diagram (\fref{fig:nodalpartition_a}) 
or a triangular mesh where the centroids of triangles surrounding a given 
node are connected to form a polygon (\fref{fig:nodalpartition_b}).
We denote by $E$ a node and its associated cell. The node $E$ has 
coordinates $\vm{x}_E$ in the Cartesian coordinate system and 
the area of its associated cell is $|E|$. 
A node on $\Gamma_t$ 
is denoted by $S$ and its coordinates in the Cartesian coordinate system 
by $\vm{x}_S$. The length of influence of the node $S$ is represented by $|S|$.
\fref{fig:neumann_nodalpartition} presents two typical cases of a node
located on $\Gamma_t$.
The partition formed by 
all the nodes and their associated cells lying on $\Omega$ is denoted 
by $\mathcal{T}^h$, where $h$ is the maximum diameter 
of any cell in the partition. The one-dimensional partition formed
by all the nodes and their associated lengths of influence lying on $\Gamma_t$
is denoted by $\mathcal{E}^h$. 

Following a standard Galerkin procedure, we define the following
discrete local spaces at the nodal cell level:
\begin{equation*}
\mathscr{V}^h|_E :=
\left\{\vm{u}^h(\vm{x}): \vm{u}^h \in [ H^{1}(E)]^2\right\},\quad
\mathscr{W}^h|_E := \mathscr{V}^h|_E,
\end{equation*}
where $\vm{u}^h(\vm{x})$ is given in the form~\eref{eq:trial} with $\vm{x}\in E$. These
discrete local spaces are assembled to form the following discrete global
spaces:
\begin{align*}
\mathscr{V}^h &:=
\left\{\vm{u}(\vm{x})\in\mathscr{V}: \vm{u}|_E \in \mathscr{V}^h|_E 
\quad\forall E \in \mathcal{T}^h\right\},\\
\mathscr{W}^h &:=
\left\{\vm{v}(\vm{x})\in\mathscr{W}: \vm{v}|_E \in \mathscr{V}^h|_E 
\quad \forall E \in \mathcal{T}^h\right\}.
\end{align*}

Owing to the definition of the discrete global spaces, we can evaluate 
the weak form~\eref{eq:weakform} by sampling it locally at 
each node in $\mathcal{T}^h$ and summing through all of them, as follows:
\begin{subequations}\label{eq:weakform_disc}
\begin{equation}
a(\vm{u}^h,\vm{v}^h) =\sum_{E\in\mathcal{T}^h}a_E(\vm{u}^h,\vm{v}^h)=\sum_{E\in\mathcal{T}^h}|E|\,\bsym{\varepsilon}_E(\vm{v}^h):\mathcal{D}:\bsym{\varepsilon}_E(\vm{u}^h)\quad\forall\vm{u}^h,\,\forall\vm{v}^h\in\mathscr{V}^h,
\label{eq:weakform_disc_a}
\end{equation}
\begin{equation}
\ell(\vm{v}^h) = \sum_{E\in\mathcal{T}^h}\ell_{b,E}(\vm{v}^h)+\sum_{S\in\mathcal{E}^h}\ell_{t,S}(\vm{v}^h)=\sum_{E\in\mathcal{T}^h}|E|\,\vm{b}_E\cdot\vm{v}^h + \sum_{S\in\mathcal{E}^h}|S|\,\bar{\vm{t}}_S\cdot\vm{v}^h\quad\forall\vm{v}^h\in\mathscr{V}^h,
\label{eq:weakform_disc_b}
\end{equation}
\end{subequations}
where $\mathcal{D}$ is the material moduli tensor that defines
the constitutive relation $\bsym{\sigma}=\mathcal{D}:\bsym{\varepsilon}(\vm{u})$;
$\bsym{\varepsilon}_E := \bsym{\varepsilon}(\vm{x}_E)$, $\vm{b}_E := \vm{b}(\vm{x}_E)$ 
and $\bar{\vm{t}}_S := \bar{\vm{t}}(\vm{x}_S)$ are nodal
quantities. Equation~\eref{eq:weakform_disc} is the result of direct nodal 
integration (1-point) of the weak form integrals.


\section{NODAL INTEGRATION USING THE VIRTUAL ELEMENT DECOMPOSITION} \label{sec:nived}

As discussed in Section~\ref{sec:intro}, the direct nodal integration of 
the bilinear form as given in~\eref{eq:weakform_disc_a} is 
not viable due to stability issues. As a remedy, nodal integration techniques
add a stabilization term to the bilinear form~\cite{puso:2008:MAF}. 
Here, we take a different route and develop a nodal integration scheme 
for meshfree methods using the mathematical framework of the 
virtual element method. In this approach, the consistency and 
stability of the nodally integrated meshfree Galerkin method
is ensured by construction.

\subsection{Nodal integration cell}\label{sec:intmesh}

\begin{figure}[!tbh]
  \centering 
  \psfrag{n}{\Large $\vm{n}$}
  \psfrag{E}{\Large $E$} 
  \psfrag{x}{\Large $\vm{x}_E$}
  \epsfig{file = 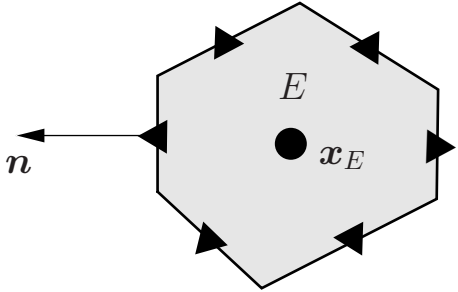,
  width = 0.4\textwidth} \caption{Schematic representation of a nodal 
  integration cell. The integration node and its associated nodal cell 
  is denoted by $E$, and the area of the nodal cell by $|E|$. The 
  coordinates of the integration node are $\vm{x}_E$ and the unit outward
  normal to the nodal cell's boundary is $\vm{n}$. The $\blacktriangle$ 
  stands for the 1-point Gauss rule over an edge of the nodal cell.} 
  \label{fig:nodalcell}
\end{figure}

The sampling of the weak form at an integration node in the partition
$\mathcal{T}^h$, where the partition is constructed
by means of one of the methods shown in~\fref{fig:nodalpartition}, is 
performed using 1-point Gauss rule over the edges of its nodal 
cell. It is stressed that in the nodal integration scheme,
the integration node and its associated nodal cell are interchangeable. 
This means that any quantity that is evaluated on the nodal cell 
is also a quantity evaluated at the node. Considering the Cartesian 
coordinate system, $\vm{x}_E$ are the coordinates of the integration 
node and $\vm{n}$ is the unit outward normal to the nodal 
cell's boundary. A sample nodal integration cell is shown in~\fref{fig:nodalcell}.

\subsection{Nodal contribution}\label{sec:nodalcontribution}

\begin{figure}[!tbh]
  \centering 
  \epsfig{file = 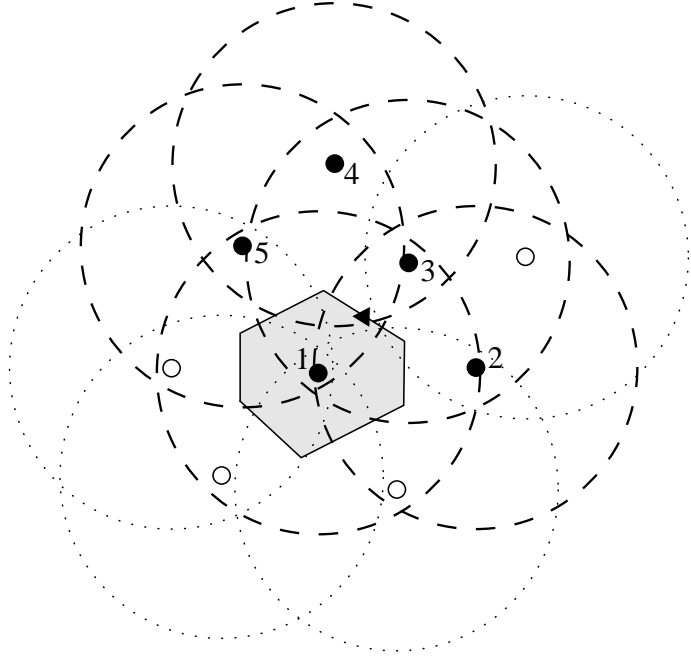,
  width = 0.6\textwidth} \caption{Construction of the local nodal contribution list. 
  Nodal basis functions are evaluated at the integration point $\blacktriangle$. 
  A circle centered at a node represents the nodal basis function support.
  The contribution list is formed by the nodes labeled 
  from $1$ to $5$ ($m=5$) since their supports (dashed circles) take
  a nonzero value at the integration point (i.e., contain the integration point). 
  Dotted circles do not contain the integration point and thus their 
  nodes are not part of the list.}
\label{fig:nodalcontribution}
\end{figure}

Since the NIVED method uses meshfree basis functions for the discretization 
of the field variables, we have to consider the contributions from the neighboring
nodes when performing Gauss integration over the edges of the nodal cell. The global 
nodal contribution list at an integration point is defined as the global indices 
of the nodes whose basis functions take a nonzero value at the integration point. 
We label these global indices from 1 to $m$ and construct a local nodal 
contribution list with them. We always keep track of the correspondence 
between the local and global nodal contribution lists as this correspondence
is used later in the assembly of the global stiffness matrix and global
force vector. The construction of the local nodal contribution list is 
schematically shown in~\fref{fig:nodalcontribution} for the evaluation 
of meshfree basis functions at an integration point on the edge of a nodal cell.

\subsection{Virtual element decomposition}\label{sec:vem_decomposition}

Due to the nonpolynomial nature of the linearly precise maximum-entropy
meshfree basis functions, the approximation of the displacement field using 
these functions contains a linear polynomial part plus some additional 
nonpolynomial terms. Let $[\mathcal{P}(E)]^2$ represent the
space of linear displacements over the nodal cell $E$.

Following the
standard VEM literature (see for instance, Reference~\cite{BeiraodaVeiga-Brezzi-Marini-Russo:2014})
the following projection operator onto the linear displacement space is defined:
\begin{equation}
\Pi: \mathscr{V}^h|_E  \to [\mathcal{P}(E)]^2, \,
\Pi\vm{p}=\vm{p} \ \forall \vm{p}\in [\mathcal{P}(E)]^2,
\label{eq:projection_operator}
\end{equation}
which allows the splitting of the meshfree approximation of the fields 
into their linear polynomial part and their nonpolynomial terms, respectively, 
as follows:
\begin{subequations}\label{eq:displacement_splitting}
\begin{align}
\vm{u}^h &=\Pi\vm{u}^h+(\vm{u}^h-\Pi\vm{u}^h),\\
\vm{v}^h &=\Pi\vm{v}^h+(\vm{v}^h-\Pi\vm{v}^h).
\end{align}
\end{subequations}
The projection $\Pi$ is required to satisfy the following orthogonality condition:
\begin{equation}
a_E(\vm{p}, \vm{v}^h-\Pi\vm{v}^h) = 0 \quad
\forall\vm{p}\in [\mathcal{P}(E)]^2, \ \ \vm{v}^h\in \mathscr{V}^h|_E.
\label{eq:ortho_condition}
\end{equation}
Using~\eref{eq:displacement_splitting} and~\eref{eq:ortho_condition},
we obtain the following decomposition of the local bilinear form:
\begin{equation}
a_E(\vm{u}^h,\vm{v}^h) = a_E(\Pi\vm{u}^h,\Pi\vm{v}^h)+ a_E(\vm{u}^h-\Pi\vm{u}^h,\vm{v}^h-\Pi\vm{v}^h).
\label{eq:vem_decomp1}
\end{equation}

The first term on the right-hand side of~\eref{eq:vem_decomp1} is computable
as it depends on the linear fields. However, the second term is
noncomputable as it depends on the nonpolynomial terms. In the framework of the VEM,
the second term is approximated by a bilinear form that can be conveniently
computed adopting the form of a stabilization term. We denote this stability bilinear
form by $s_E$ and rewrite~\eref{eq:vem_decomp1} as follows:
\begin{equation}
a^h_E(\vm{u}^h,\vm{v}^h) = a_E(\Pi\vm{u}^h,\Pi\vm{v}^h)+ s_E(\vm{u}^h-\Pi\vm{u}^h,\vm{v}^h-\Pi\vm{v}^h).
\label{eq:vem_decomp2}
\end{equation}
We refer to~\eref{eq:vem_decomp2} as the \textit{virtual element decomposition}.

The virtual element decomposition is endowed with the following crucial properties
for establishing the convergence of the 
VEM~\cite{BeiraoDaVeiga-Brezzi-Cangiani-Manzini-Marini-Russo:2013,BeiraodaVeiga-Brezzi-Marini:2013}:

For all $h$ and for all $E$ in $\mathcal{T}^h$
\begin{itemize}
\item \textit{Consistency}: For all $\vm{p} \in [\mathcal{P}(E)]^2$ and for all $\vm{v}^h\in \mathscr{V}^h|_E$
\begin{equation}\label{eq:consistency_cond}
a_E^h(\vm{p},\vm{v}^h)=a_E(\vm{p},\vm{v}^h).
\end{equation}
\item \textit{Stability}: There exists two constants $\alpha_*>0$ and $\alpha^*>0$, independent of $h$ and of $E$, such that
\begin{equation}\label{eq:stability_cond}
\forall\vm{v}^h\in \mathscr{V}^h|_E, \quad \alpha_*a_E(\vm{v}^h,\vm{v}^h)\leq a_E^h(\vm{v}^h,\vm{v}^h)\leq \alpha^*a_E(\vm{v}^h,\vm{v}^h).
\end{equation}
\end{itemize}

In view of the preceding properties, it is straightforward to recognize that the first
term on the right-hand side of~\eref{eq:vem_decomp2} provides
consistency (i.e., ensures patch test satisfaction) and the second term
lends stability. The stability property~\eref{eq:stability_cond} reveals
the necessary conditions that $s_E$ must possess: it must be symmetric and
positive definite on the kernel of $\Pi$ so that
property~\eref{eq:stability_cond} holds without violating~\eref{eq:consistency_cond}.

\subsection{Projection operator}\label{sec:projection_operator}

The explicit form of the projection operator $\Pi$ is obtained from the orthogonality
condition~\eref{eq:ortho_condition}.
Let the strain tensor~\eref{eq:symstrain} be written as
\begin{equation}
\bsym{\varepsilon}(\vm{v}) = \bsym{\nabla}\vm{v}-\bsym{\omega}(\vm{v}),
\label{eq:altsymstrain}
\end{equation}
where $\bsym{\omega}(\vm{v})$ is the skew-symmetric tensor given by
\begin{equation}
\bsym{\omega}(\vm{v}) = \frac{1}{2}\left(\bsym{\nabla}\vm{v}
-\bsym{\nabla}^\transpose\vm{v}\right).
\label{eq:skewstrain}
\end{equation}
Using~\eref{eq:altsymstrain}, noting that $\bsym{\sigma}(\vm{p})$
and $\bsym{\nabla}\Pi\vm{v}^h$ are constant fields over $E$ since
$\vm{p},\,\Pi\vm{v}^h \in [\mathcal{P}(E)]^2$, 
and that $\bsym{\sigma}:\bsym{\omega}=0$, 
we can write the orthogonality condition~\eref{eq:ortho_condition} as
\begin{align}\label{eq:projectionsforms}
a_E(\vm{p},\vm{v}^h-\Pi\vm{v}^h) &= \int_E\bsym{\sigma}(\vm{p}):\left[\bsym{\nabla}(\vm{v}^h-\Pi\vm{v}^h)
-\bsym{\omega}(\vm{v}^h-\Pi\vm{v}^h)\right]\,\diffx \nonumber\\
&= \bsym{\sigma}(\vm{p}):\left[\int_E\bsym{\nabla}\vm{v}^h\,\diffx-\bsym{\nabla}\Pi\vm{v}^h\int_E\,\diffx\right]\nonumber\\
&= \bsym{\sigma}(\vm{p}):\left[\int_E\bsym{\nabla}\vm{v}^h\,\diffx-|E|\bsym{\nabla}\Pi\vm{v}^h\right] \\
&= 0,
\end{align}
which leads to
\begin{equation}\label{eq:grad_pip}
\bsym{\nabla}\Pi\vm{v}^h=\frac{1}{|E|}\int_E \bsym{\nabla}\vm{v}^h\,\diffx.
\end{equation}

Note that~\eref{eq:grad_pip} defines $\bsym{\nabla}\Pi\vm{v}^h$ as the average value of
$\bsym{\nabla}\vm{v}^h$ over the cell $E$. On using~\eref{eq:altsymstrain}, \eref{eq:grad_pip}
can be rewritten as
\begin{equation}\label{eq:grad_pip_alt}
\bsym{\nabla}\Pi\vm{v}^h= \bsym{\varepsilon}(\Pi\vm{v}^h) + \bsym{\omega}(\Pi\vm{v}^h) 
=\frac{1}{|E|}\int_E \bsym{\varepsilon}(\vm{v}^h)\,\diffx+\frac{1}{|E|}\int_E \bsym{\omega}(\vm{v}^h)\,\diffx,
\end{equation}
which on integrating yields
\begin{equation}\label{eq:pip_projection_form}
\Pi\vm{v}^h=\left(\frac{1}{|E|}\int_E\bsym{\varepsilon}(\vm{v}^h)\,\diffx\right)\cdot\vm{x}+\left(\frac{1}{|E|}\int_E\bsym{\omega}(\vm{v}^h)\,\diffx\right)\cdot\vm{x}+\vm{a}_0.
\end{equation}
To determine $\vm{a}_0$ we need a projection operator onto 
constants $P_0: \mathscr{V}^h|_E \to \Re^2$ such that
\begin{equation}\label{eq:pipo}
P_0(\Pi\vm{v}^h)=P_0\vm{v}^h.
\end{equation}
Given that the field variables computed at the integration node
become the representative field variables of the cell, we define the projection 
operator onto constant as
\begin{equation}\label{eq:pipo_def}
P_0\vm{v}^h=\vm{v}^h(\vm{x}_E)=\vm{v}_E,
\end{equation}
where $\vm{x}_E$ are the coordinates of the integration node (see~\fref{fig:nodalcell}).

Applying~\eref{eq:pipo} to \eref{eq:pip_projection_form} gives
\begin{align}\label{eq:pipo_pipv}
P_0(\Pi\vm{v}^h)&=\left(\frac{1}{|E|}\int_E\bsym{\varepsilon}(\vm{v}^h)\,\diffx\right)\cdot P_0\vm{x}+\left(\frac{1}{|E|}\int_E\bsym{\omega}(\vm{v}^h)\,\diffx\right)\cdot P_0\vm{x}+ P_0 \vm{a}_0= P_0\vm{v}^h\nonumber\\
&=\left(\frac{1}{|E|}\int_E\bsym{\varepsilon}(\vm{v}^h)\,\diffx\right)\cdot\vm{x}_E+\left(\frac{1}{|E|}\int_E\bsym{\omega}(\vm{v}^h)\,\diffx\right)\cdot\vm{x}_E+\vm{a}_0=\vm{v}_E.
\end{align}
On solving for $\vm{a}_0$ in~\eref{eq:pipo_pipv}, we obtain
\begin{equation}\label{eq:constant_c}
\vm{a}_0=\vm{v}_E-\left(\frac{1}{|E|}\int_E\bsym{\varepsilon}(\vm{v}^h)\,\diffx\right)\cdot\vm{x}_E-\left(\frac{1}{|E|}\int_E\bsym{\omega}(\vm{v}^h)\,\diffx\right)\cdot\vm{x}_E.
\end{equation}

We introduce the definitions $\bsym{\varepsilon}_E(\vm{v}^h):=\frac{1}{|E|}\int_E\bsym{\varepsilon}(\vm{v}^h)\,\diffx$
and $\bsym{\omega}_E(\vm{v}^h) := \frac{1}{|E|}\int_E\bsym{\omega}(\vm{v}^h)\,\diffx$ since both are constant tensors
over the cell and thus they can be associated with the integration node. Finally, substituting~\eref{eq:constant_c}
into~\eref{eq:pip_projection_form} yields the projection operator onto the linear displacements, as follows:
\begin{equation}\label{eq:pi_operator}
\Pi\vm{v}^h=\bsym{\varepsilon}_E(\vm{v}^h)\cdot(\vm{x}-\vm{x}_E)+\bsym{\omega}_E(\vm{v}^h)\cdot(\vm{x}-\vm{x}_E)+\vm{v}_E.
\end{equation}

In~\eref{eq:pi_operator}, the cell averages  $\bsym{\varepsilon}_E(\vm{v}^h)$ and $\bsym{\omega}_E(\vm{v}^h)$
are evaluated on the boundary of $E$ by invoking the divergence theorem, which gives 
the following nodal measures:
\begin{equation}\label{eq:average_symstrain}
\bsym{\varepsilon}_E(\vm{v}^h)=\frac{1}{|E|}\int_E \bsym{\varepsilon}(\vm{v}^h)\,\diffx = \frac{1}{2|E|}\int_{\partial
E}\left(\vm{v}^h\otimes\vm{n}+\vm{n}\otimes\vm{v}^h\right)\,\diffs
\end{equation}
and
\begin{equation}\label{eq:average_skewstrain}
\bsym{\omega}_E(\vm{v}^h)=\frac{1}{|E|}\int_E \bsym{\omega}(\vm{v}^h)\,\diffx=\frac{1}{2|E|}\int_{\partial
E}\left(\vm{v}^h\otimes\vm{n}-\vm{n}\otimes\vm{v}^h\right)\,\diffs,
\end{equation}
respectively.

\subsection{Projection matrix}\label{sec:projection_matrix}
After some algebraic manipulations, \eref{eq:pi_operator} can be written as
\begin{equation}
\Pi\vm{v}^h = \mat{h}(\vm{x})\bsym{\upvarepsilon}_E(\vm{v}^h)+\mat{g}(\vm{x})\mat{r}_E(\vm{v}^h),
\label{eq:projmat1}
\end{equation}
where
\begin{equation}
\mat{h}(\vm{x}) = \smat{(x_1-x_{1E}) & 0 & \frac{1}{2}(x_2-x_{2E})\\ 0 & (x_2-x_{2E}) & \frac{1}{2}(x_1-x_{1E})},\quad
\mat{g}(\vm{x}) = \smat{1 & 0 & \frac{1}{2}(x_2-x_{2E})\\ 0 & 1 & -\frac{1}{2}(x_1-x_{1E})},
\label{eq:mtrxhg}
\end{equation}
\begin{equation}
\bsym{\upvarepsilon}_E(\vm{v}^h) = \smat{(\bsym{\varepsilon}_E)_{11} & (\bsym{\varepsilon}_E)_{22} & 2(\bsym{\varepsilon}_E)_{12}}^\transpose,\quad
\mat{r}_E(\vm{v}^h) = \smat{v_{1E} & v_{2E} & 2(\bsym{\omega}_E)_{12}}^\transpose.
\label{eq:mtrxer}
\end{equation}
Using meshfree basis functions leads to the following discrete representation of~\eref{eq:mtrxhg}:
\begin{equation}
\mat{h}^h(\vm{x}) = \sum_{a=1}^m\phi_a(\vm{x})\mat{h}(\vm{x}_a)=\mat{N}\mat{H}_E, \quad
\mat{g}^h(\vm{x}) = \sum_{a=1}^m\phi_a(\vm{x})\mat{g}(\vm{x}_a)=\mat{N}\mat{G}_E, \quad
\label{eq:disc_mtrxhg}
\end{equation}
where

\begin{equation}
\mat{H}_E=\smat{
(\mat{H}_E)_1\\ \vdots\\ (\mat{H}_E)_a\\ \vdots\\ (\mat{H}_E)_m}, \quad
(\mat{H}_E)_a = \smat{(x_{1a}-x_{1E}) & 0 & \frac{1}{2}(x_{2a}-x_{2E}) \\ 0 & (x_{2a}-x_{2E}) & \frac{1}{2}(x_{1a}-x_{1E})},
\label{eq:matrix_H}
\end{equation}

\begin{equation}
\mat{G}_E=\smat{
(\mat{G}_E)_1\\ \vdots\\ (\mat{G}_E)_a\\ \vdots\\ (\mat{G}_E)_m}, \quad
(\mat{G}_E)_a = \smat{1 & 0 & \frac{1}{2}(x_{2a}-x_{2E}) \\ 0 & 1 & -\frac{1}{2}(x_{1a}-x_{1E})},
\label{eq:matrix_G}
\end{equation}

\begin{equation}
\mat{N} =
\smat{\mat{N}_1 & \cdots & \mat{N}_a & \cdots & \mat{N}_m}, \quad \mat{N}_a=\smat{\phi_a & 0 \\ 0 &
\phi_a},
\label{eq:matrix_N}
\end{equation}
and explicitly replacing $\vm{v}^h$ in the form~\eref{eq:trial} into~\eref{eq:mtrxer} gives
\begin{equation}
\bsym{\upvarepsilon}_E^h=\bsym{\upvarepsilon}_E\left(\sum_{a=1}^m\phi_a(\vm{x})\vm{v}_a\right)=\mat{W}_E^\transpose\mat{q}, \quad
\mat{r}_E^h=\mat{r}_E\left(\sum_{a=1}^m\phi_a(\vm{x})\vm{v}_a\right)=\mat{R}_E^\transpose\mat{q},
\label{eq:disc_mtrxer}
\end{equation}
where
\begin{equation}
\mat{q}= \smat{\vm{v}_1^\transpose & \cdots & \vm{v}_a^\transpose & \cdots & \vm{v}_m^\transpose}^\transpose,\quad
\vm{v}_a=\smat{v_{1a} & v_{2a}}^\transpose,
\label{eq:matrix_q}
\end{equation}
\begin{equation}
\mat{W}_E=\smat{
(\mat{W}_E)_1\\ \vdots\\ (\mat{W}_E)_a\\ \vdots\\ (\mat{W}_E)_m}, \quad
(\mat{W}_E)_a = \smat{q_{1a} & 0 & q_{2a} \\ 0 & q_{2a} & q_{1a}},\quad
q_{ia}=\frac{1}{|E|}\int_{\partial E}\phi_a(\vm{x})n_i\,\diffs,\quad i=1,2
\label{eq:matrix_W}
\end{equation}
and 
\begin{equation}
\mat{R}_E=\smat{
(\mat{R}_E)_1\\ \vdots\\ (\mat{R}_E)_a\\ \vdots\\ (\mat{R}_E)_m}, \quad
(\mat{R}_E)_a = \smat{\phi_a^E & 0 & q_{2a}\\ 0 & \phi_a^E & -q_{1a}},\quad
\label{eq:matrix_R}
\end{equation}
where $\phi_a^E:=\phi_a(\vm{x}_E)$ is the value of the $a$-th basis function
at the integration node $E$. To evaluate $q_{ia}$ in~\eref{eq:matrix_W}, 
a 1-point Gauss rule is used on each edge 
of the nodal cell (see~\fref{fig:nodalcell}). 

Finally, substituting~\eref{eq:disc_mtrxhg} and~\eref{eq:disc_mtrxer} into~\eref{eq:projmat1}
yields the following discrete version of the projection operator:
\begin{equation}
\Pi\vm{v}^h=\mat{N}\mat{H}_E\mat{W}_E^\transpose\mat{q}+\mat{N}\mat{G}_E\mat{R}_E^\transpose\mat{q}=
\mat{N}(\mat{H}_E\mat{W}_E^\transpose+\mat{G}_E\mat{R}_E^\transpose)\mat{q},
\label{eq:disc_pi}
\end{equation}
which defines the projection matrix (the matrix form of $\Pi$) as
\begin{equation}
\mat{P}_E=\mat{H}_E\mat{W}_E^\transpose+\mat{G}_E\mat{R}_E^\transpose.
\label{eq:projection_matrix}
\end{equation}

\subsection{Nodally integrated stiffness matrix}\label{sec:stiffness_matrix}

The nodally integrated local stiffness matrix is obtained by discretizing the virtual
element decomposition~\eref{eq:vem_decomp2}. We start by working on
the consistency term (i.e., the first term on the right-hand side
of~\eref{eq:vem_decomp2}). On using~\eref{eq:pi_operator} to obtain
$\bsym{\nabla}\Pi\vm{v}^h=\bsym{\varepsilon}_E(\vm{v}^h)+\bsym{\omega}_E(\vm{v}^h)$
and considering the constitutive relation $\bsym{\sigma}(\Pi\vm{u}^h)=\mathcal{D}:\bsym{\varepsilon}_E(\vm{u}^h)$, 
we can write
\begin{align}\label{eq:consistent_bilinearform2}
a_E(\Pi\vm{u}^h,\Pi\vm{v}^h) &=|E|\bsym{\sigma}(\Pi\vm{u}^h):\bsym{\nabla}\Pi\vm{v}^h\nonumber\\
&=|E|\,\bsym{\varepsilon}_E(\vm{v}^h):\mathcal{D}:\bsym{\varepsilon}_E(\vm{u}^h).
\end{align}
Using the symmetry of tensors $\bsym{\varepsilon}_E$ and $\mathcal{D}$, 
\eref{eq:consistent_bilinearform2} can be written in Voigt notation as
\begin{equation}
a_E(\Pi\vm{u}^h,\Pi\vm{v}^h)= |E|\,\bsym{\upvarepsilon}_E^\transpose(\vm{v}^h)\,\mat{D}\,\bsym{\upvarepsilon}_E(\vm{u}^h),
\label{eq:consistent_bilinearform3}
\end{equation}
where $\mat{D}$ is the constitutive
matrix for an isotropic linear elastic material given by
\begin{equation}
\mat{D} =
\frac{E_\mathrm{Y}}{(1+\nu)(1-2\nu)}\smat{1-\nu & \nu & 0\\ \nu & 1-\nu & 0\\ 0 & 0 & \frac{1-2\nu}{2}},\quad
\mat{D} =
\frac{E_\mathrm{Y}}{1-\nu^2}\smat{1 & \nu & 0 \\ \nu & 1 & 0 \\ 0 & 0 & \frac{1-\nu}{2}}
\end{equation}
for plane strain and plane stress conditions, respectively, where $E_\mathrm{Y}$ is the 
Young's modulus and $\nu$ is the Poisson's ratio; $\bsym{\upvarepsilon}_E(\vm{v}^h)$ 
is defined in~\eref{eq:mtrxer}.

Now, substituting the discrete version of $\bsym{\upvarepsilon}_E(\vm{v}^h)$ as given
in~\eref{eq:disc_mtrxer} into~\eref{eq:consistent_bilinearform3} leads to the
following discrete local consistency bilinear form:
\begin{equation}
a_E(\Pi\vm{u}^h,\Pi\vm{v}^h)=\mat{q}^\transpose|E|\,\mat{W}_E\,\mat{D}\,\mat{W}_E^\transpose\mat{d},
\label{eq:discrete_consistent_bilinear_form}
\end{equation}
where $\mat{d}$ is a column vector similar to $\mat{q}$ given in~\eref{eq:matrix_q}
that contains the nodal coefficients that are associated with the displacement
trial functions.

The discrete local stability bilinear form is obtained by replacing 
the field discretizations
$\vm{u}^h=\mat{N}\mat{d}$ and $\vm{v}^h=\mat{N}\mat{q}$, where $\mat{N}$ is
defined in~\eref{eq:matrix_N}, along with~\eref{eq:disc_pi} into the
second term on the right-hand side of~\eref{eq:vem_decomp2}. This yields
\begin{align}\label{eq:discrete_stability_bilinear_form}
s_E(\vm{u}^h-\Pi\vm{u}^h,\vm{v}^h-\Pi\vm{v}^h) &= s_E(\mat{N}\mat{d}-\mat{N}\mat{P}_E\mat{d},\mat{N}\mat{q}-\mat{N}\mat{P}_E\mat{q}) \nonumber\\
&= \mat{q}^\transpose(\mat{I}_{2m}-\mat{P}_E)^\transpose\,\mat{S}_E\,(\mat{I}_{2m}-\mat{P}_E)\,\mat{d},
\end{align}
where $\mat{I}_{2m}$ is the identity ($2m\times 2m$) matrix and
$\mat{S}_E=s_E(\mat{N}^\transpose,\mat{N})$. Thus, substituting
\eref{eq:discrete_consistent_bilinear_form} and
\eref{eq:discrete_stability_bilinear_form} into the first and
second terms on the right-hand side of~\eref{eq:vem_decomp2},
respectively, gives
\begin{equation}
a_E^h(\vm{u}^h,\vm{v}^h)=\mat{q}^\transpose\left(|E|\,\mat{W}_E\,\mat{D}\,\mat{W}_E^\transpose +
(\mat{I}_{2m}-\mat{P}_E)^\transpose\,\mat{S}_E\,(\mat{I}_{2m}-\mat{P}_E)\right)\mat{d},
\label{eq:discrete_vem_bilinear_form}
\end{equation}
which defines the nodally integrated local stiffness matrix as the sum
of the consistency stiffness matrix $\mat{K}_E^\cons$ and the stability 
stiffness matrix $\mat{K}_E^\stab$, as follows:
\begin{equation}
\mat{K}_E = \mat{K}_E^\cons + \mat{K}_E^\stab,\quad
\mat{K}_E^\cons=|E|\,\mat{W}_E\,\mat{D}\,\mat{W}_E^\transpose,\quad
\mat{K}_E^\stab = (\mat{I}_{2m}-\mat{P}_E)^\transpose\,\mat{S}_E\,(\mat{I}_{2m}-\mat{P}_E).
\label{eq:nodal_stiffness_matrix}
\end{equation}

Regarding the stability stiffness, we make the following choice for $\mat{S}_E$:
\begin{equation}
\mat{S}_E = \mat{I}_{2m}\odot\left(\mat{1}_{2m}\otimes\bsym{\gamma}\right),
\label{eq:se_matrix}
\end{equation}
where $\mat{1}_{2m}$ is the unit ($2m\times 1$) column vector, $\bsym{\gamma} = \mathrm{diag}(\mat{K}_E^\cons)$ is
a column vector containing the diagonal of the matrix $\mat{K}_E^\cons$ and $\odot$ is
the element-wise product. This means that $\mat{S}_E$ is a diagonal matrix
whose diagonal entries are those of the diagonal of the nodally integrated local consistency
stiffness matrix, which is in the spirit of the ``$D$-recipe'' studied
in References~\cite{beirao:HOVEM:2017,dassi:EHOVEM:2018}.

\subsection{Nodally integrated force vector}\label{sec:force_vector}

For linear fields, the simplest approximation for the 
body force vector is constructed by projecting both the body force vector $\vm{b}$
and the test displacements $\vm{v}^h$ onto constants, as follows~\cite{BeiraoDaVeiga-Brezzi-Cangiani-Manzini-Marini-Russo:2013,BeiraodaVeiga-Brezzi-Marini:2013,artioli:AO2DVEM:2017}:
\begin{equation}\label{eq:discrete_element_body_force_linear_form}
\ell_{b,E}^h(\vm{v}^h)=\int_E P_0\vm{b}^h\cdot P_0\vm{v}^h\,\diffx=
P_0\vm{b}^h\cdot P_0\vm{v}^h\int_E\,\diffx = |E|\,\vm{b}_E\cdot\vm{v}_E
= \mat{q}^\transpose |E|\mat{N}_E^\transpose\vm{b}_E,
\end{equation}
where we have used $P_0$ as defined in~\eref{eq:pipo_def}, and
\begin{equation}\label{eq:matrix_Ne}
\mat{N}_E =
\left[(\mat{N}_E)_1 \quad \cdots \quad (\mat{N}_E)_a \quad \cdots \quad
(\mat{N}_E)_m\right], \quad (\mat{N}_E)_a=\smat{\phi_a^E & 0 \\ 0 &
\phi_a^E},
\end{equation}
where $\phi_a^E$ is defined below~\eref{eq:matrix_R}. Hence, the nodal body 
force vector is given by
\begin{equation}
\mat{f}_{b,E}=|E|\,\mat{N}_E^\transpose\vm{b}_E,
\end{equation}
which coincides with the direct integration of the body force vector at the node 
with coordinates $\vm{x}_E$.

The integral that defines the traction force vector is
similar to the integral that defines the body force vector, but it is
one dimension lower. This means that
we can simply apply a direct nodal integration on the Neumann edges.
Proceeding likewise leads to the following nodal traction force vector:
\begin{equation}
\mat{f}_{t,S}=|S|\,\mat{N}_S^\transpose\bar{\vm{t}}_S,
\end{equation}
where 
\begin{equation}
\mat{N}_S =
\left[(\mat{N}_S)_1 \quad \cdots \quad (\mat{N}_S)_a \quad \cdots \quad
(\mat{N}_S)_m\right], \quad (\mat{N}_S)_a=\smat{\phi_a^S & 0 \\ 0 &
\phi_a^S}
\end{equation}
with $\phi_a^S:=\phi_a(\vm{x}_S)$.

\subsection{Nodally integrated mass matrix}
\label{sec:mass_matrix}

With the aim of testing the NIVED approach for an elastodynamic problem,
we show how to compute the approximation of the mass matrix using the virtual
element decomposition. The mass matrix is constructed along the same lines
as the stiffness matrix, that is, the mass matrix is split into a consistency part 
and a stability part. However, this construction is done using an $L^2$-projection 
on $[\mathcal{P}(E)]^2$ of any function $\vm{v}^h\in [\mathcal{W}(E)]^2$.
For linear and quadratic fields, the $L^2$-projection coincides with the 
elliptic projection $\Pi$~\cite{ahmad:EQPVEM:2013,BeiraodaVeiga-Brezzi-Marini-Russo:2014}.
This means, 
\[
\int_E\phi_a\phi_b\diffx=\int_E\Pi\phi_a\Pi\phi_b\diffx+\int_E(\phi_a-\Pi\phi_a)(\phi_b-\Pi\phi_b).
\]
For testing our nodal integration approach, we use Newmark's constant average acceleration
time integration method, which is implicit and unconditionally stable. For such integrators,
the stability part of the mass matrix is not needed~\cite{vacca:VEMPAR:2015}.
Thus, generally we only require the consistency term $\int_E\Pi\phi_a\Pi\phi_b \, \diffx$.
Hence, the expression for the
nodally integrated mass matrix for a solid
material of density $\rho$, is:
\begin{equation}
\mat{M}_E=\rho|E|\int_E\mat{P}^\transpose\mat{N}^\transpose\mat{N}\mat{P}\,\diffx=
\rho|E|\mat{P}^\transpose\left(\int_E\mat{N}^\transpose\mat{N}\,\diffx\right)\mat{P}=
\rho|E|\mat{P}^\transpose\mat{N}_E^\transpose\mat{N}_E\mat{P},
\label{eq:nodal_mass_matrix}
\end{equation}
where the integral has been nodally integrated at $\vm{x}_E$ leading 
to $\mat{N}_E$ as defined in~\eref{eq:matrix_Ne}. 

For a detailed treatment of time-dependent problems using the virtual element
framework, the interested
reader is referred to References~\cite{vacca:VEMPAR:2015,vacca:VEMHYP:2017,park:VEMELASTDYN:2019}.

\section{PATCH TEST SATISFACTION} 
\label{sec:patch_test}

In the patch test, the exact linear field $\vm{u}(\vm{x})=\vm{p}(\vm{x}) \in [\mathcal{P}(E)]^2$ 
is imposed on the entire boundary of the domain. Hence, the numerical solution for the patch test 
must be $\vm{u}^h(\vm{x})=\vm{p}(\vm{x})$ on the whole domain. Using this in~\eref{eq:vem_decomp2},
we get
\begin{equation}\label{eq:pt1}
a_E(\vm{p},\vm{v}^h)=a_E(\Pi\vm{p},\Pi\vm{v}^h)+a_E(\vm{p}-\Pi\vm{p},\vm{v}^h-\Pi\vm{v}^h).
\end{equation}
The second term on the right-handside of~\eref{eq:pt1} is zero 
because $\vm{p}-\Pi\vm{p}=\vm{p}-\vm{p}=0$. Then,
\begin{equation}\label{eq:pt2}
a_E(\vm{p},\vm{v}^h)=\int_E\bsym{\sigma}(\vm{p}):\bsym{\varepsilon}(\Pi\vm{v}^h)\,\diffx.
\end{equation}
Note that $\bsym{\sigma}(\vm{p})$ is a constant field. Let $\bsym{\sigma}^c$ denote this 
constant stress field and observe that $\bsym{\varepsilon}(\Pi\vm{v}^h)=\bsym{\varepsilon}_E(\vm{v}^h)$
is the nodal strain (cell average) defined in~\eref{eq:average_symstrain} and thus 
a constant field as well. Hence, \eref{eq:pt2} becomes
\begin{equation}\label{eq:pt3}
a_E(\vm{p},\vm{v}^h)=|E|\bsym{\sigma}^c:\bsym{\varepsilon}_E(\vm{v}^h)=
\bsym{\sigma}^c:\frac{1}{2}\int_{\partial E}\left(\vm{v}^h\otimes\vm{n}+\vm{n}\otimes\vm{v}^h\right)\,\diffs.
\end{equation}

We now substitute $\vm{v}^h(\vm{x}) =\phi_a(\vm{x})\vm{v}_a$ into~\eref{eq:pt3} and sum 
over all the nodal cells to obtain
\begin{equation}\label{eq:pt4}
\sum_{E\in\mathcal{T}^h} a_E(\vm{p},\vm{v}^h)=
\sum_{E\in\mathcal{T}^h}\left(\int_{\partial E}\mat{B}_a\bsym{\sigma}^c\,\diffs\right)\vm{v}_a=\vm{f}_a\cdot\vm{v}_a,
\end{equation}
where
\begin{equation}\label{eq:pt5}
\mat{B}_a=
\smat{\phi_a n_1 & 0 \\ 0 & \phi_a n_2 \\ \phi_a n_2 & \phi_a n_1},
\end{equation}
and $\vm{f}_a$ is an interior nodal force. Since the patch test produces a 
state of constant strains and stresses, all the interior nodal forces must 
be identically equal to zero. Therefore, for the patch test to be satisfied 
it suffices to show that
\begin{equation}\label{eq:pt6}
\vm{f}_a=\sum_{E\in\mathcal{T}^h}\int_{\partial E}\mat{B}_a\bsym{\sigma}^c\,\diffs=\vm{0},
\end{equation}
where the assembly is over all the nodal cells that have a non-zero intersection with 
the support of $\phi_a$. To compute~\eref{eq:pt6}, we choose a 1-point Gauss rule over 
the edges of the integration cells, and since the evaluation of $\mat{B}_a$ 
at a given interior edge will arise from two adjacent cells in the assembly, the two 
contributions cancel each other. Thus, the net contribution to $\vm{f}_a$ from all 
the interior edges vanishes, and hence~\eref{eq:pt6} is satisfied.


\section{EXTENSION TO NONLINEAR ANALYSIS: THE VISCOELASTIC CASE} 
\label{sec:nonlinear_nived}

In this section, the proposed NIVED scheme is extended to solid mechanics
problems in which the constitutive law for the solid material depends on the history
of deformation. In particular, we consider the case of Generalized Maxwell viscoelastic 
model describing a standard linear isotropic solid~\cite{zienkiewicz:2005:FEM2}, wherein
the Cauchy stress in vector form and as a function of time $t$ is split as
\begin{equation}\label{eq:stress_split}
\bsym{\upvarsigma}(t)=\mat{s}(t)+\mat{m}K\mat{m}^\transpose\bsym{\upvarepsilon},
\end{equation}
where $K$ is the material's Bulk modulus, $\bsym{\upvarepsilon}$ is the small strain tensor 
in vector form, $\mat{m}=[1\quad 1 \quad 0]^\transpose$ 
and $\mat{s}(t)$ is the deviatoric stress given by
\begin{equation}\label{eq:devstressintform}
\mat{s}(t)= \int_{-\infty}^{t}2G(t-t')\frac{\partial\mat{e}}{\partial t'}\,dt',
\end{equation}
where $\mat{e}=\bsym{\upvarepsilon}-1/3\,\mat{m}\mat{m}^\transpose\bsym{\upvarepsilon}$ is the
deviatoric strain and $G(t)$ is the shear relaxation modulus given by a one-term
Prony series that is of the form
\begin{equation}\label{eq:relaxmodulus}
G(t)=G[\mu_0+\mu_1\exp(-t/\lambda_1)],
\end{equation}
where $G=E/(2(1+\nu))$, $\mu_0>0$ and $\mu_1>0$ are parameters that satisfy $\mu_0+\mu_1=1$, and
$\lambda_1$ is a relaxation time. Substituting~\eref{eq:relaxmodulus} into~\eref{eq:devstressintform}
leads to
\begin{equation}\label{eq:devstress}
\mat{s}(t)= 2G\left[\mu_0\mat{e}(t)+\mu_1\mat{q}^{(1)}(t)\right],
\end{equation}
where $\mat{q}^{(1)}(t)$ is the dimensionless partial deviatoric strain given by
\begin{equation}\label{eq:pardevstrain}
\mat{q}^{(1)}(t)= \int_{-\infty}^{t}\exp(-(t-t')/\lambda_1)\frac{\partial\mat{e}}{\partial t'}\,dt'.
\end{equation}
The term $\mat{q}^{(1)}(t)$ is solved using the recursion formula~\cite{zienkiewicz:2005:FEM2}
\begin{equation}\label{eq:pardevstrainrecform}
\mat{q}_{n+1}^{(1)}= \exp(-\Delta t/\lambda_1)\mat{q}_n^{(1)}+\Delta\mat{q}_{n+1}^{(1)},
\end{equation}
where 
\begin{equation}\label{eq:deltapardevstrain}
\Delta\mat{q}_{n+1}^{(1)}= \frac{\lambda_1}{\Delta t}(1-\exp(-\Delta t/\lambda_1))(\mat{e}_{n+1}-\mat{e}_n).
\end{equation}
Thus, using the recursion formula, the stresses can be written as
\begin{equation}\label{eq:devstressrecform}
\mat{s}_{n+1}= 2G\left[\mu_0\mat{e}_{n+1}+\mu_1\mat{q}_{n+1}^{(1)}\right],
\end{equation}
and
\begin{equation}\label{eq:stress_split_recform}
\bsym{\upvarsigma}_{n+1}=\mat{s}_{n+1}+\mat{m}K\mat{m}^\transpose\bsym{\upvarepsilon}_{n+1}.
\end{equation}

On considering~\eref{eq:stress_split_recform} along with~\eref{eq:devstressrecform} and
\eref{eq:pardevstrainrecform}, it is realized that $\bsym{\upvarsigma}_{n+1}$ is
determined by the knowledge of the time increment $\Delta t$ and the strain 
$\bsym{\upvarepsilon}$ at times $t_n$ and $t_{n+1}$.

The numerical solution of the equilibrium equation for the standard linear isotropic
solid described by the Generalized Maxwell viscoelastic model is sought
incrementally using a Newton-Raphson scheme. Under the NIVED framework,
the constitutive equation~\eref{eq:stress_split_recform} is computed using the vector 
form of the nodal strain~\eref{eq:average_symstrain} and the base equilibrium equation
for the Newton-Raphson scheme is obtained by linearizing the following VED 
representation of the virtual work:
\begin{align}\label{eq:nonlinear_ved_vwork}
&|E|\,\bsym{\upvarepsilon}_E^\transpose(\vm{v}^h)
\bsym{\upvarsigma}_{n+1}(\bsym{\upvarepsilon}_E(\vm{u}_{n+1}^h))
+\int_E \left[\bsym{\upvarepsilon}(\vm{v}^h)-\bsym{\upvarepsilon}_E(\vm{v}^h)\right]^\transpose
\bsym{\upvarsigma}_{n+1}(\bsym{\upvarepsilon}(\vm{u}_{n+1}^h)-\bsym{\upvarepsilon}_E(\vm{u}_{n+1}^h))\,\diffx\nonumber\\
&-|E|\,(P_0\vm{v}^h)^\transpose\vm{b}_E - |S|\,(P_0\vm{v}^h)^\transpose\bar{\vm{t}}_S = 0.
\end{align}
The linearization of~\eref{eq:nonlinear_ved_vwork} yields
\begin{align}\label{eq:linearized_ved_vwork}
&|E|\,\bsym{\upvarepsilon}_E^\transpose(\vm{v}^h)\mat{D}_\tangent|_{n+1}\bsym{\upvarepsilon}_E(\Delta\vm{u}^h)
+\int_E \left[\bsym{\upvarepsilon}(\vm{v}^h)-\bsym{\upvarepsilon}_E(\vm{v}^h)\right]^\transpose
\mat{D}_\tangent|_{n+1}\left[\bsym{\upvarepsilon}(\Delta\vm{u}^h)-\bsym{\upvarepsilon}_E(\Delta\vm{u}^h)\right]\,\diffx\nonumber\\
&=-\Big\{|E|\,\bsym{\upvarepsilon}_E^\transpose(\vm{v}^h)
\bsym{\upvarsigma}_{n+1}(\bsym{\upvarepsilon}_E(\vm{u}_{n+1}^h))
+\int_E \left[\bsym{\upvarepsilon}(\vm{v}^h)-\bsym{\upvarepsilon}_E(\vm{v}^h)\right]^\transpose
\bsym{\upvarsigma}_{n+1}(\bsym{\upvarepsilon}(\vm{u}_{n+1}^h)-\bsym{\upvarepsilon}_E(\vm{u}_{n+1}^h))\,\diffx\Big.\nonumber\\
&\Big.\quad\quad\quad-|E|\,(P_0\vm{v}^h)^\transpose\vm{b}_E - |S|\,(P_0\vm{v}^h)^\transpose\bar{\vm{t}}_S\Big\},
\end{align}
where $\mat{D}_\tangent$ is the tangent moduli that is obtained by deriving~\eref{eq:stress_split_recform},
as follows:
\begin{equation}\label{eq:tangent_moduli}
\mat{D}_\tangent|_{n+1}=\frac{\partial\bsym{\upvarsigma}_{n+1}}{\partial\bsym{\upvarepsilon}_{n+1}}=
2G\Big[\mu_0+\mu_1\frac{\lambda_1}{\Delta t}(1-\exp(-\Delta t/\lambda_1))\Big]
\Big[\mat{I}-\frac{1}{3}\mat{m}\mat{m}^\transpose\Big]+K\mat{m}\mat{m}^\transpose.
\end{equation}
Adopting the same argument used
to obtain the VED decomposition~\eref{eq:vem_decomp2}, we approximate the second term
on both sides of~\eref{eq:linearized_ved_vwork} by computable quantities, thus giving
\begin{align}\label{eq:linearized_ved_vwork_final}
&|E|\,\bsym{\upvarepsilon}_E^\transpose(\vm{v}^h)\mat{D}_\tangent|_{n+1}\bsym{\upvarepsilon}_E(\Delta\vm{u}^h)
+s_{E,\tangent}(\vm{v}^h-\Pi\vm{v}^h,\Delta\vm{u}^h-\Pi\Delta\vm{u}^h)\nonumber\\
&=-\Big\{|E|\,\bsym{\upvarepsilon}_E^\transpose(\vm{v}^h)
\bsym{\upvarsigma}_{n+1}+s_E(\vm{v}^h-\Pi\vm{v}^h,\vm{u}^h-\Pi\vm{u}^h)
-|E|\,(P_0\vm{v}^h)^\transpose\vm{b}_E - |S|\,(P_0\vm{v}^h)^\transpose\bar{\vm{t}}_S\Big\}.
\end{align}

We remark that the linearized virtual work~\eref{eq:linearized_ved_vwork_final} is 
supported by the virtual element framework for inelastic analysis presented in 
Reference~\cite{artioli:AOVEMIN:2017}. After discretizing~\eref{eq:linearized_ved_vwork_final}, 
relying on the arbitrariness of nodal variations and summing through all the cells, 
we finally obtain the following Newton-Raphson iterations to solve 
the equilibrium state $\mat{d}_{n+1}^{(k)}=\mat{d}_{n+1}^{(k-1)}+\Delta\mat{d}^{(k)}$ 
at time $t_{n+1}$ with a time increment $\Delta t = t_{n+1}-t_n$:
\begin{subequations}\label{eq:nr_scheme}
\begin{align}
\sum_{E\in\mathcal{T}^h}\Big(\mat{K}_{E,\tangent}^\cons+\mat{K}_{E,\tangent}^\stab\Big)_{n+1}^{(k-1)}\Delta\mat{d}^{(k)} 
= -\sum_{E\in\mathcal{T}^h, S\in\mathcal{E}^h}\Big(\mat{f}_E^\cons+&\mat{f}_E^\stab-\mat{f}_{b,E}-\mat{f}_{t,S}\Big)_{n+1}^{(k-1)},\label{eq:nr_scheme_a}\\ 
\intertext{with}
\mat{K}_{E,\tangent}^\cons=|E|\mat{W}_E\mat{D}_\tangent|_{n+1}\mat{W}_E^\transpose,
\quad\mat{K}_{E,\tangent}^\stab =(\mat{I}_{2m}-\mat{P}_E)^\transpose \,\mat{I}_{2m}&\odot\left(\mat{1}_{2m}\otimes\bsym{\gamma}\right)\,(\mat{I}_{2m}-\mat{P}_E),\label{eq:nr_scheme_b}\\ 
\mat{f}_E^\cons=|E|\mat{W}_E\bsym{\upvarsigma}_{n+1}, \quad \mat{f}_E^\stab=(\mat{I}_{2m}-\mat{P}_E)^\transpose\,\mat{I}_{2m}
\odot (\mat{1}_{2m}\otimes &\bsym{\gamma})\,(\mat{I}_{2m}-\mat{P}_E)\,\mat{d}_{n+1}^{(k-1)},\\
\mat{f}_{b,E}=|E|\mat{N}_E^\transpose\vm{b}_E, \quad \mat{f}_{t,S} =|S|\mat{N}_S^\transpose &\bar{\vm{t}}_S,\label{eq:nr_scheme_c}
\end{align}
\end{subequations}
where $\bsym{\gamma} = \mathrm{diag}(\mat{K}_{E,\tangent}^\cons)$, and $\mat{W}_E$, $\mat{P}_E$, $\mat{I}_{2m}$, $\mat{1}_{2m}$, $\mat{N}_E$, $\mat{N}_S$, $\vm{b}_E$ and $\bar{\vm{t}}_S$ are the same as those used in the linear NIVED formulation (see Section~\ref{sec:nived}).


\section{NUMERICAL EXPERIMENTS} \label{sec:num_experiments}

In this section, numerical experiments are presented to demonstrate the accuracy and
convergence of the NIVED approach. Results are compared with a meshfree Galerkin
method that uses maximum-entropy basis functions and standard Gauss quadrature
with integration points defined in the interior of the cells of a background mesh
of 3-node triangles. We denote this method by the acronym MEM. 
When available, we also provide comparisons of our results 
to reference solutions obtained with the finite element method (FEM). The construction 
of the nodal representative polygonal cells for the NIVED method 
is carried out using the same background mesh of 3-node triangles used 
for the MEM method, that is, the construction method depicted 
in~\fref{fig:nodalpartition}(b) is considered. 
This allows a direct comparison of the NIVED and MEM methods since the number 
of degrees of freedom remains the same in both methods. 

For numerical integration in the NIVED approach, we use a 1-point Gauss rule 
per edge of the polygonal cell and at these points only the evaluation of 
basis functions is required. In the MEM approach, we consider 1-point, 3-point, 
6-point and 12-point Gauss rules in the interior of the 3-node triangular cell.
In contrast to the NIVED method, at these interior integration points the 
evaluation of both the basis functions and their derivatives is required. 
The maximum-entropy basis functions in the 
NIVED and MEM methods are performed using $\gamma=2.0$ in 
all the numerical experiments that follow.

\subsection{Patch test}\label{sec:numexamples_patchtest}

\begin{figure}[!tbhp]
\centering
\epsfig{file = 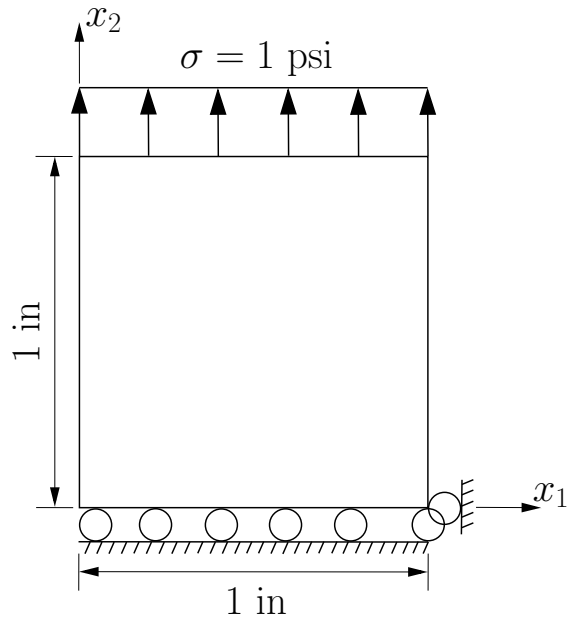, width = 0.5\textwidth}
\caption{Patch test problem.}
\label{fig:patchtestproblem}
\end{figure}

\begin{figure}[!tbhp]
\centering
\mbox{
\subfigure[]{\label{fig:patchtestmeshes_a}
\epsfig{file = 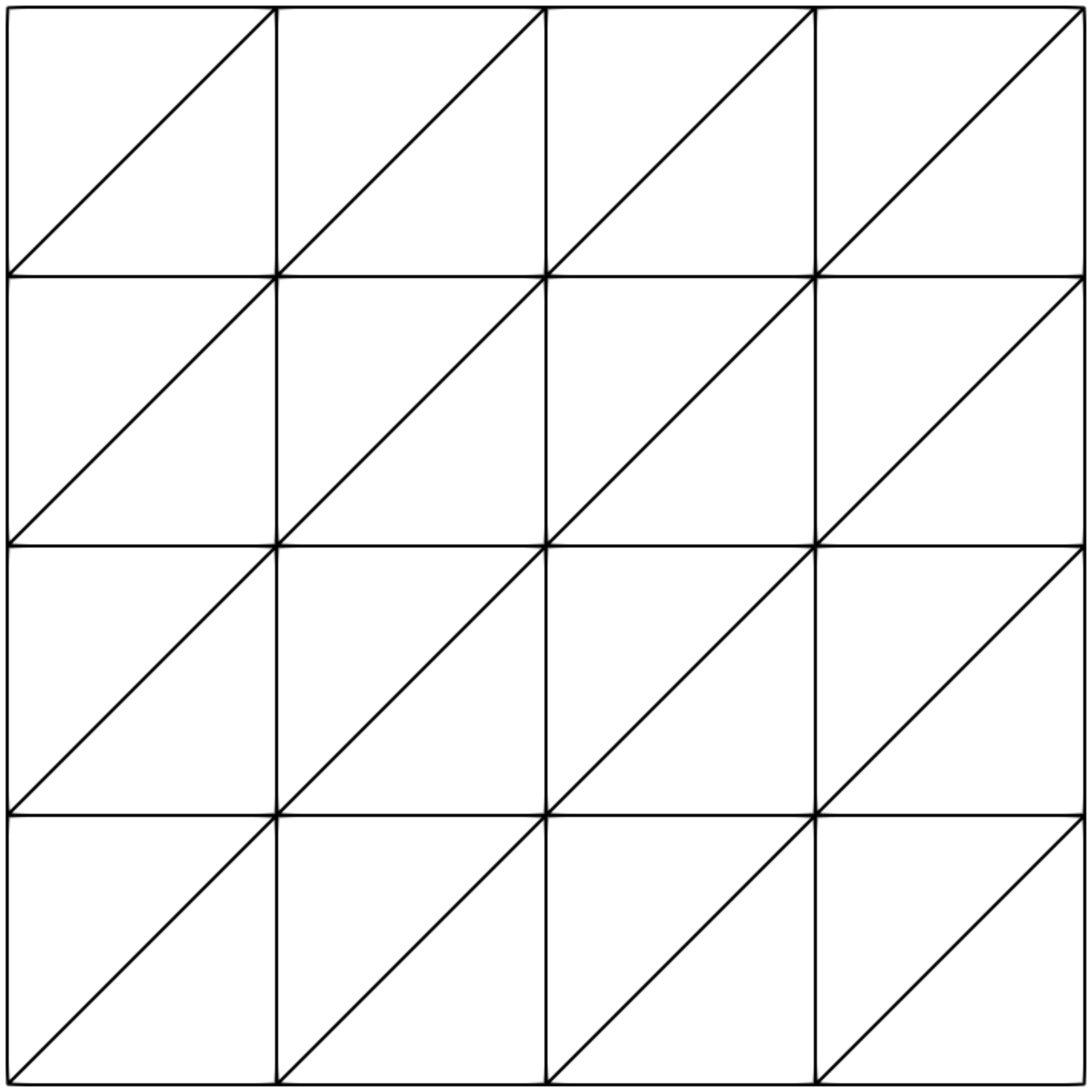, width = 0.28\textwidth}}

\subfigure[]{\label{fig:patchtestmeshes_b}
\epsfig{file = 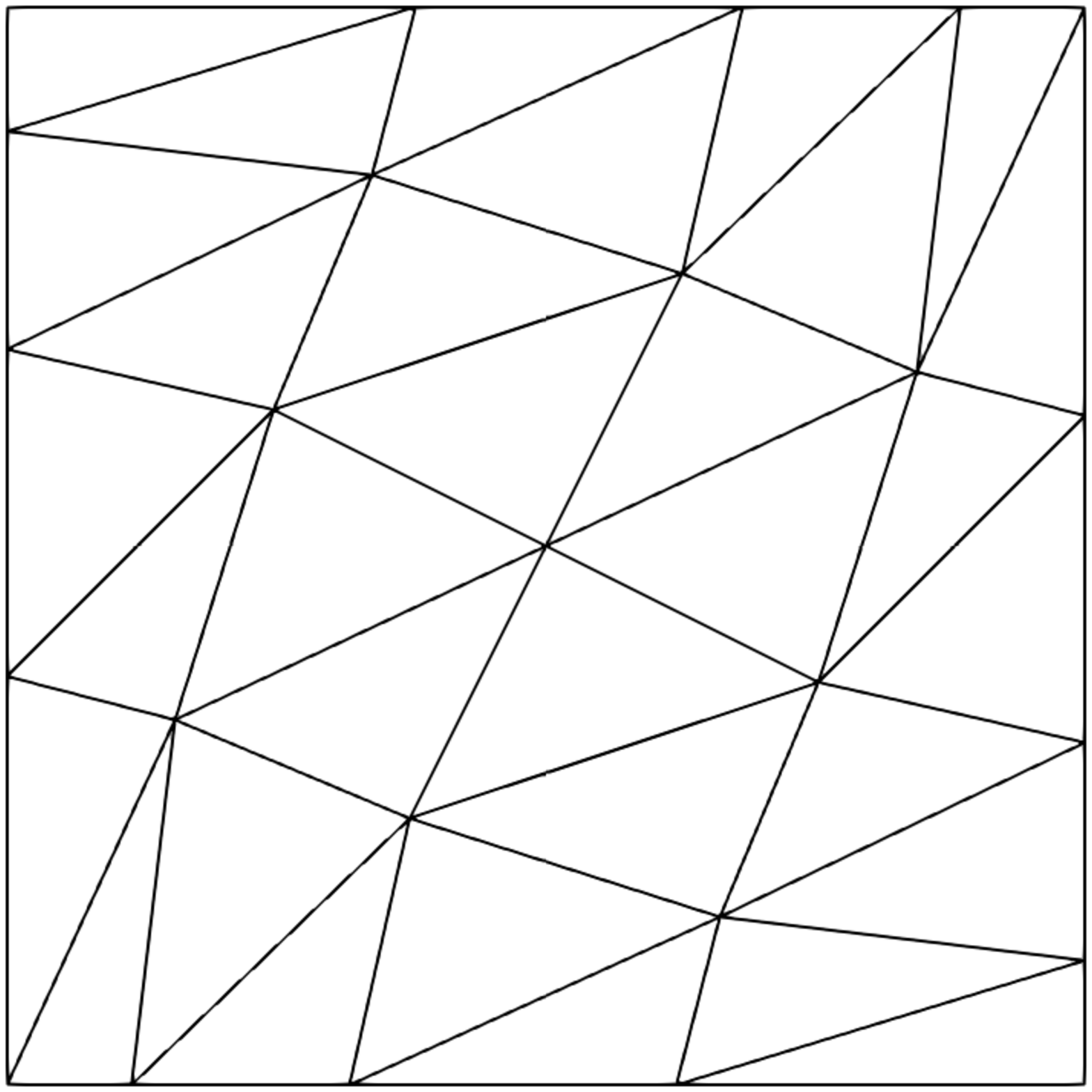, width = 0.28\textwidth}}

\subfigure[]{\label{fig:patchtestmeshes_c}
\epsfig{file = 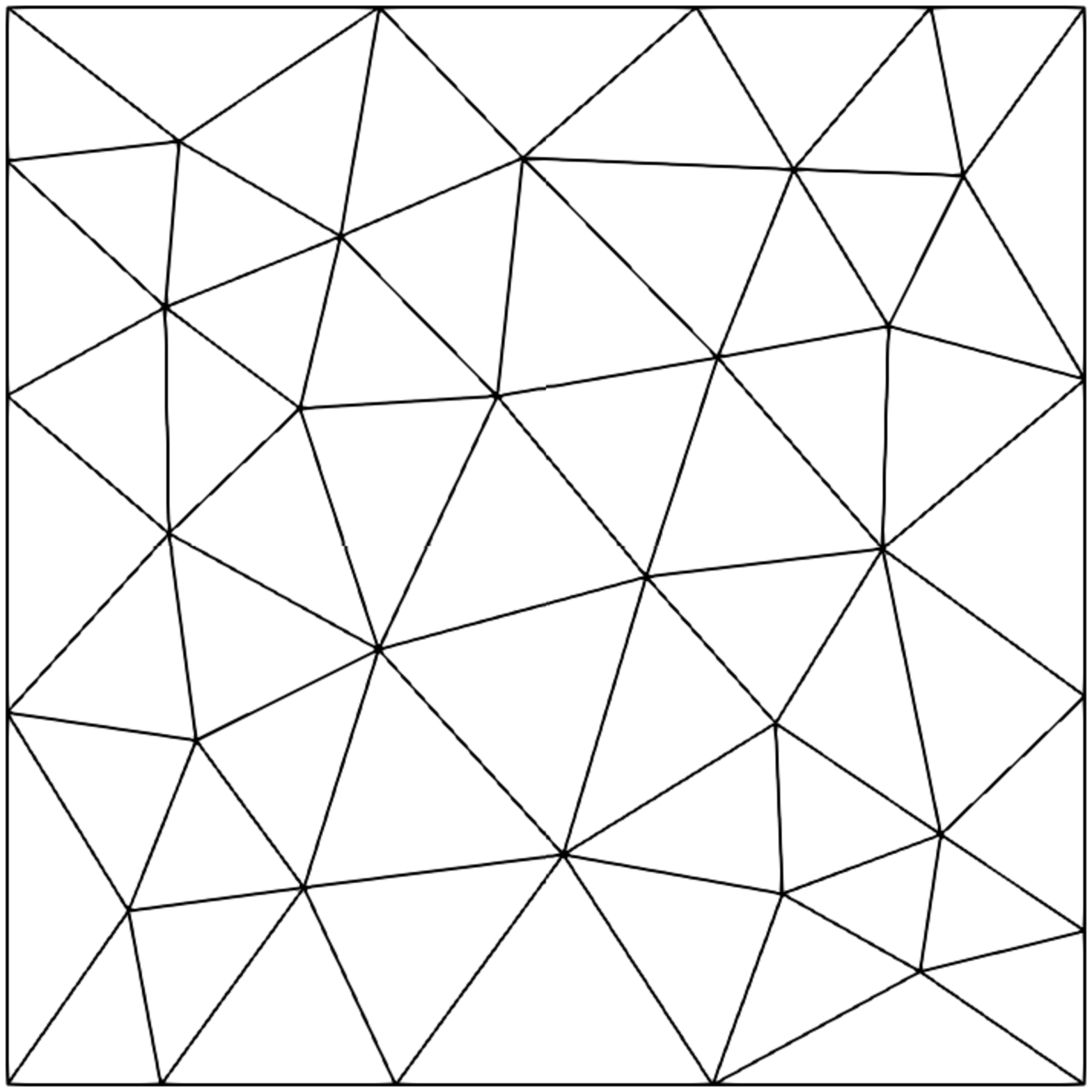, width = 0.28\textwidth}}
}
\mbox{
\subfigure[]{\label{fig:patchtestmeshes_d}
\epsfig{file = 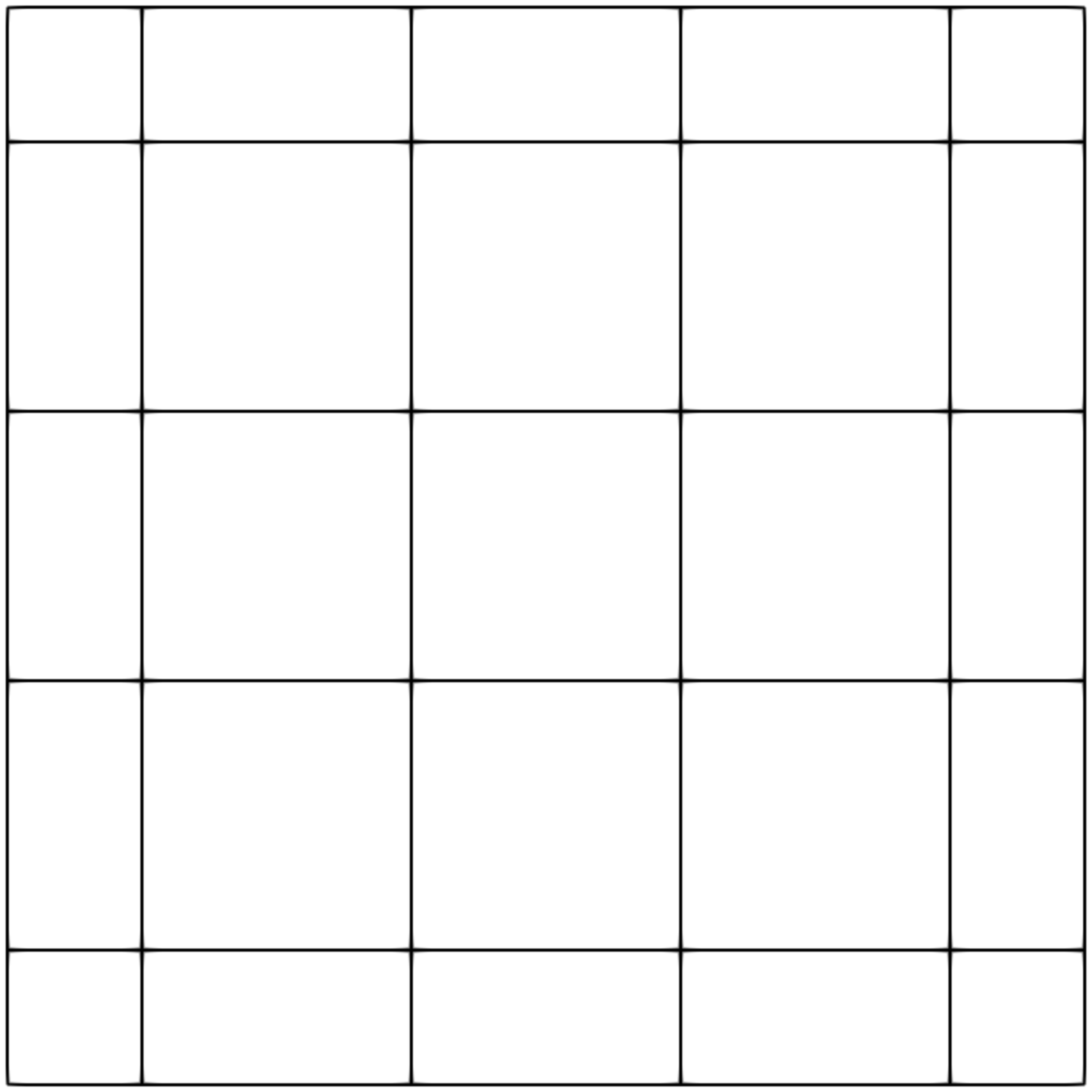, width = 0.28\textwidth}}

\subfigure[]{\label{fig:patchtestmeshes_e}
\epsfig{file = 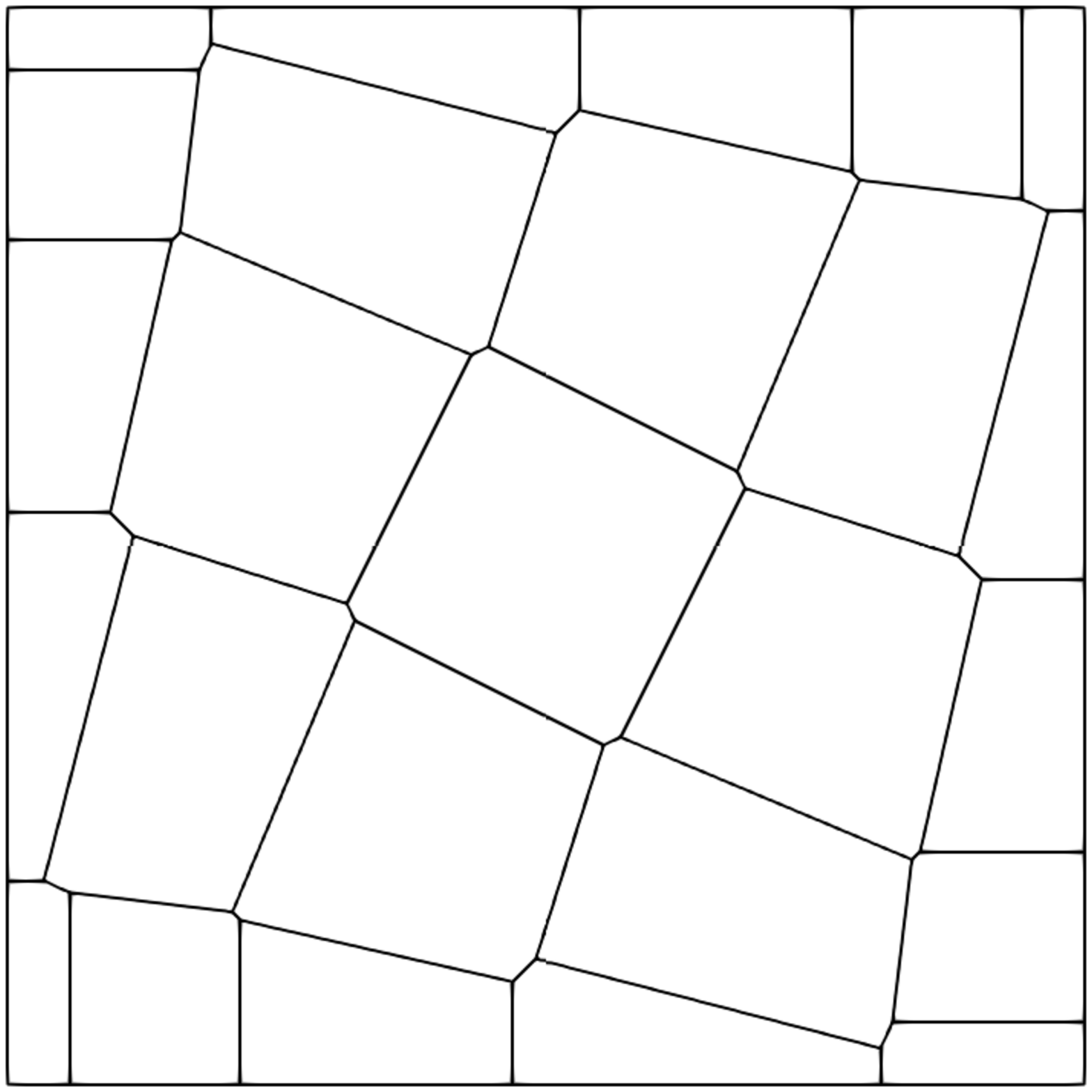, width = 0.28\textwidth}}

\subfigure[]{\label{fig:patchtestmeshes_f}
\epsfig{file = 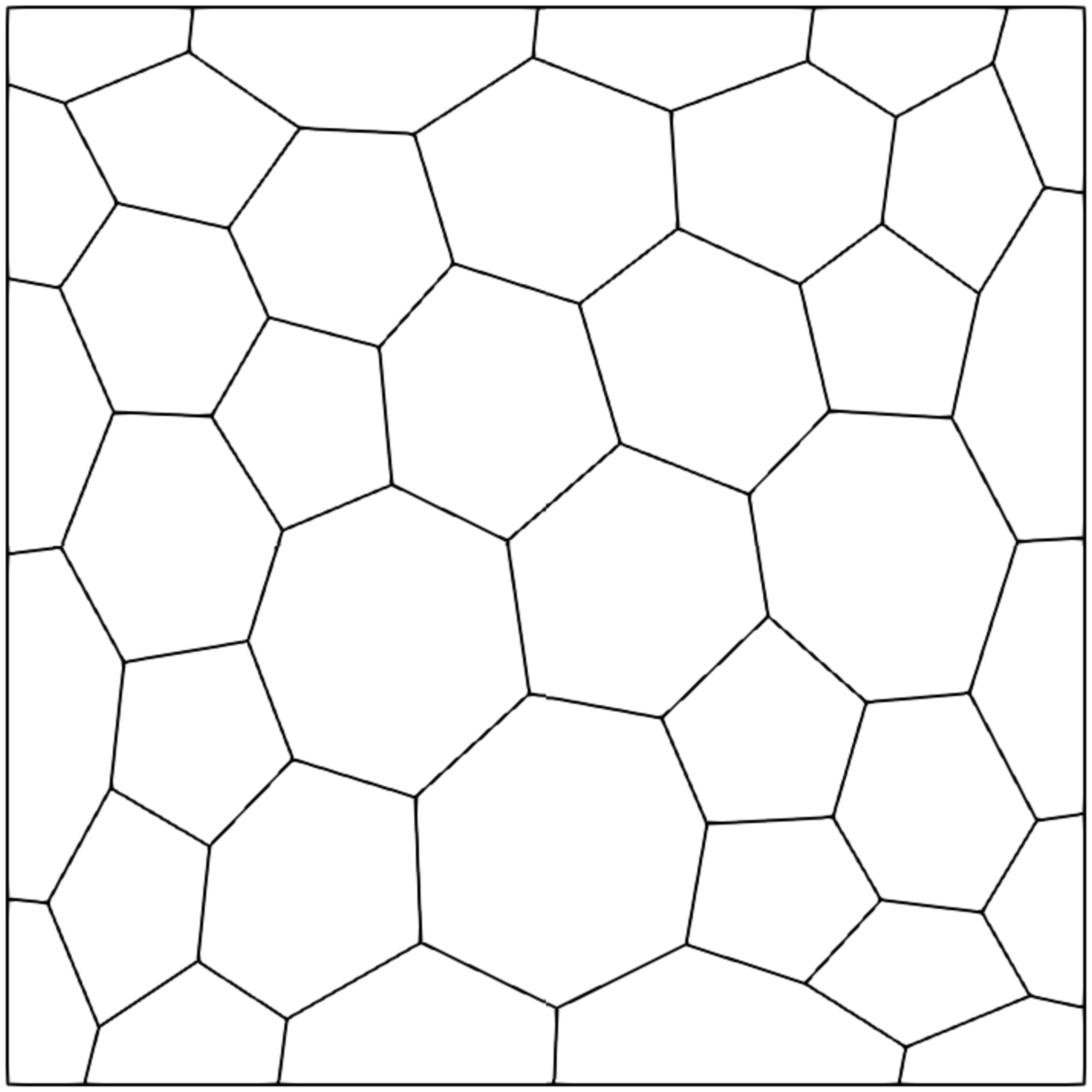, width = 0.28\textwidth}}
}
\caption{Background meshes used for the patch test.
(a) Regular mesh, (b) distorted mesh and (c) unstructured mesh for the MEM
approach; (d) regular mesh, (d) distorted mesh and (e) unstructured
mesh for the NIVED approach.}
\label{fig:patchtestmeshes}
\end{figure}

We solve the boundary-value problem~\eref{eq:strongform} for the patch test
that is schematically depicted in \fref{fig:patchtestproblem}. Plane stress 
condition is assumed with the following material
parameters: $E_\mathrm{Y}=3\times 10^7$ psi and $\nu = 0.3$. The exact solution 
for this problem is the linear field 
$\vm{u}=\smat{\nu(1-x_1)/E_\mathrm{Y} & x_2/E_\mathrm{Y}}^\transpose$.
The Dirichlet boundary conditions are imposed as follows: on the bottom boundary, 
the vertical displacement is fixed, and at its right corner, 
the displacement is additionally fixed in the horizontal direction. 
The Neumann boundary conditions are 
$\bar{\vm{t}}=\smat{0 & 0}^\transpose$
on the left and right boundaries, and
$\bar{\vm{t}}=\smat{0 & \sigma}^\transpose$ on the top boundary.
The Gaussian prior weight function is used for the maximum-entropy basis functions.
The background meshes used in this study are depicted in \fref{fig:patchtestmeshes}.
Numerical results for the relative error in the $L^2$ norm and the $H^1$ seminorm 
are presented in Tables~\ref{table:patchtest_L2norm} and 
\ref{table:patchtest_H1norm}, respectively, for the MEM and NIVED approaches. 
Numerical results confirm that the patch test is met to machine
precision only for the NIVED method.

\begin{table}[!tbhp]
\setlength{\arrayrulewidth}{.15em}
\newcolumntype{a}{>{\columncolor{white}}c}
\begin{center}
\small\addtolength{\tabcolsep}{-1pt}
\caption{Relative error in the $L^{2}$ norm for the patch test.}
\vspace*{0pt}
\renewcommand{\arraystretch}{1.5}
\begin{tabular}{aaaaa}
\hline
\rowcolor{white}
Method & Gauss rule & Regular & Distorted  & Unstructured \\
\hline
\rowcolor{white}
MEM & 1-pt (interior of the cell) & $1.0 \times 10^{-2}$ & $2.0 \times 10^{-2}$ & $1.7 \times 10^{-2}$ \\
\rowcolor{white}
MEM & 3-pt (interior of the cell) & $2.3 \times 10^{-3}$ & $1.6 \times 10^{-3}$ & $1.6 \times 10^{-3}$ \\
\rowcolor{white}
MEM & 6-pt (interior of the cell) & $5.0 \times 10^{-5}$ & $8.0 \times 10^{-4}$ & $1.2 \times 10^{-3}$ \\
\rowcolor{white}
MEM & 12-pt (interior of the cell) & $2.2 \times 10^{-7}$ & $3.0 \times 10^{-4}$ & $5.0 \times 10^{-4}$ \\
\rowcolor{white}
NIVED & 1-pt (per edge of the cell) & $3.6 \times 10^{-15}$ & $4.1 \times 10^{-15}$ & $2.5 \times 10^{-15}$ \\
\hline
\end{tabular}
\label{table:patchtest_L2norm}
\end{center}
\end{table}
\begin{table}[!tbhp]
\setlength{\arrayrulewidth}{.15em}
\newcolumntype{a}{>{\columncolor{white}}c}
\begin{center}
\small\addtolength{\tabcolsep}{-1pt}
\caption{Relative error in the $H^{1}$ seminorm for the patch test.}
\vspace*{0pt}
\renewcommand{\arraystretch}{1.5}
\begin{tabular}{aaaaa}
\hline
\rowcolor{white}
Method & Gauss rule & Regular & Distorted & Unstructured \\
\hline
\rowcolor{white}
MEM & 1-pt (interior of the cell) & $1.4 \times 10^{-2}$ & $5.4 \times 10^{-2}$ & $2.5 \times 10^{-2}$ \\
\rowcolor{white}
MEM & 3-pt (interior of the cell) & $2.6 \times 10^{-3}$ & $5.3 \times 10^{-3}$ & $4.8 \times 10^{-3}$ \\
\rowcolor{white}
MEM & 6-pt (interior of the cell) & $5.4 \times 10^{-5}$ & $1.9 \times 10^{-3}$ & $1.3 \times 10^{-3}$ \\
\rowcolor{white}
MEM & 12-pt (interior of the cell) & $2.3 \times 10^{-7}$ & $7.7 \times 10^{-4}$ & $4.5 \times 10^{-4}$ \\
\rowcolor{white}
NIVED & 1-pt (per edge of the cell)  & $3.6 \times 10^{-15}$ & $5.2 \times 10^{-15}$ & $7.8 \times 10^{-15}$ \\
\hline
\end{tabular}
\label{table:patchtest_H1norm}
\end{center}
\end{table}

\subsection{Numerical stability}\label{sec:numexamples_stability}
\begin{figure}[!tbhp]
\centering
\mbox{
\subfigure[]{\label{fig:numstability_quartic_prior_a}
\epsfig{file = 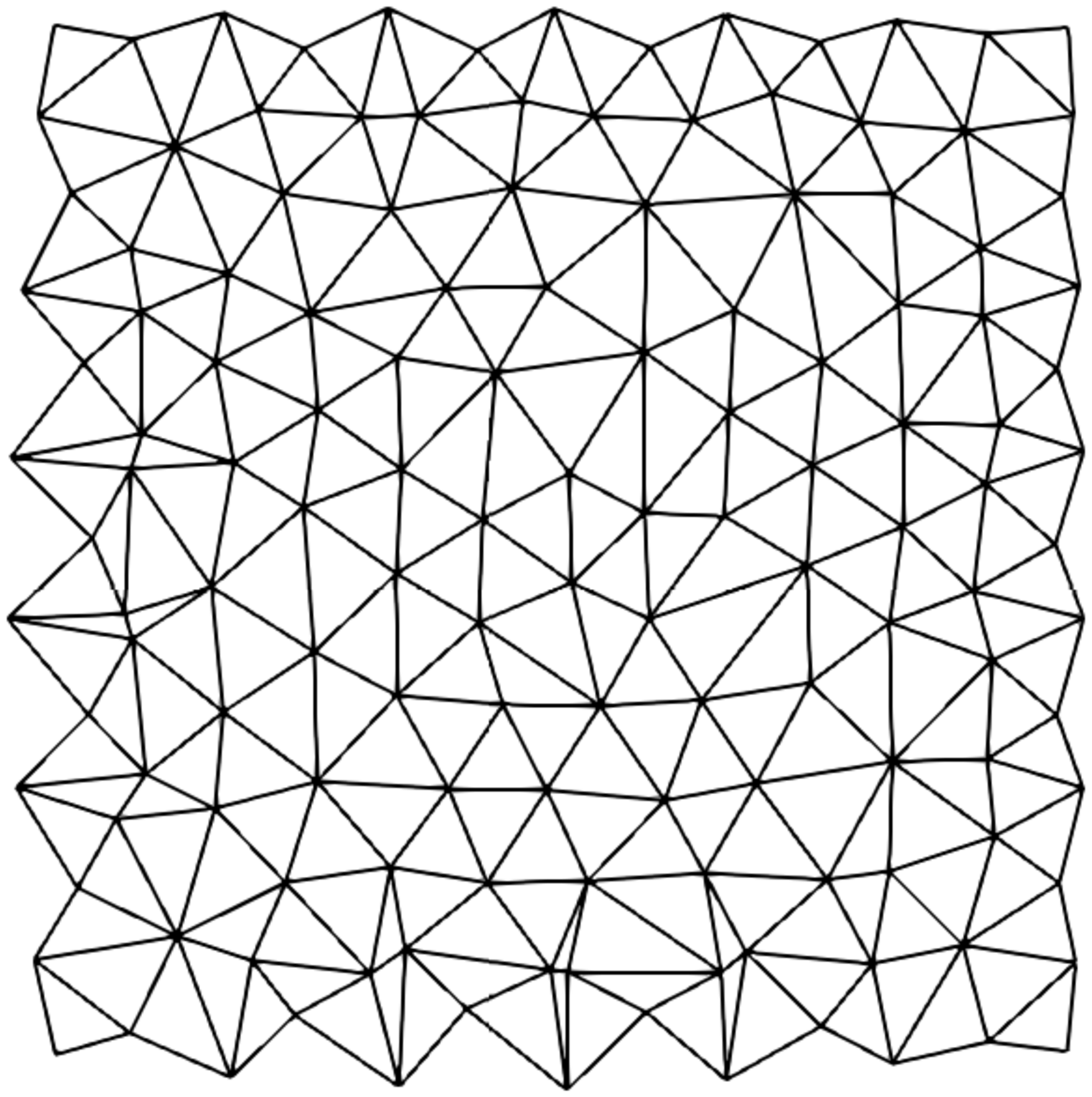, width = 0.26\textwidth}}
\subfigure[]{\label{fig:numstability_quartic_prior_b}
\epsfig{file = 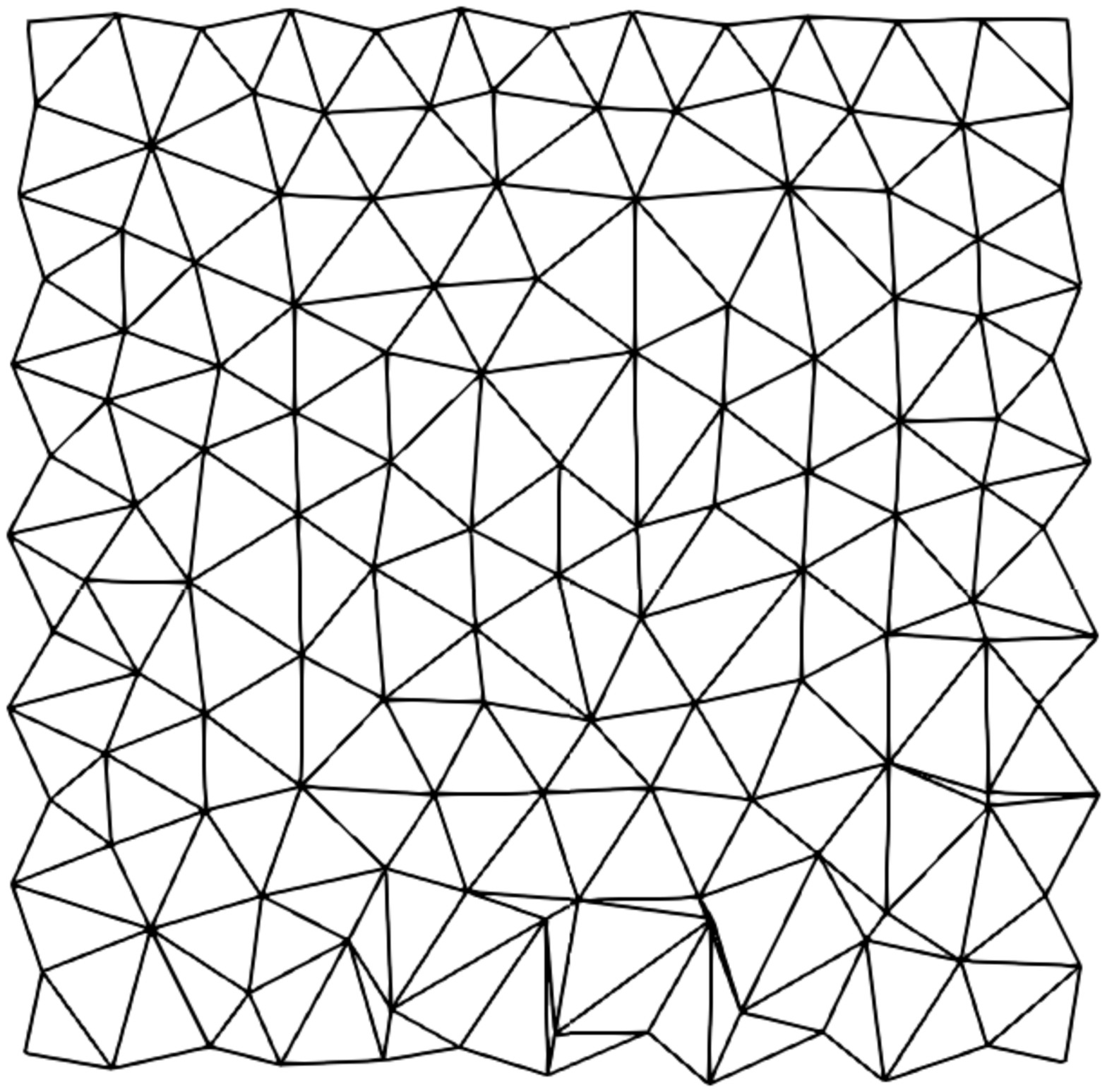, width = 0.26\textwidth}}
\subfigure[]{\label{fig:numstability_quartic_prior_c}
\epsfig{file = 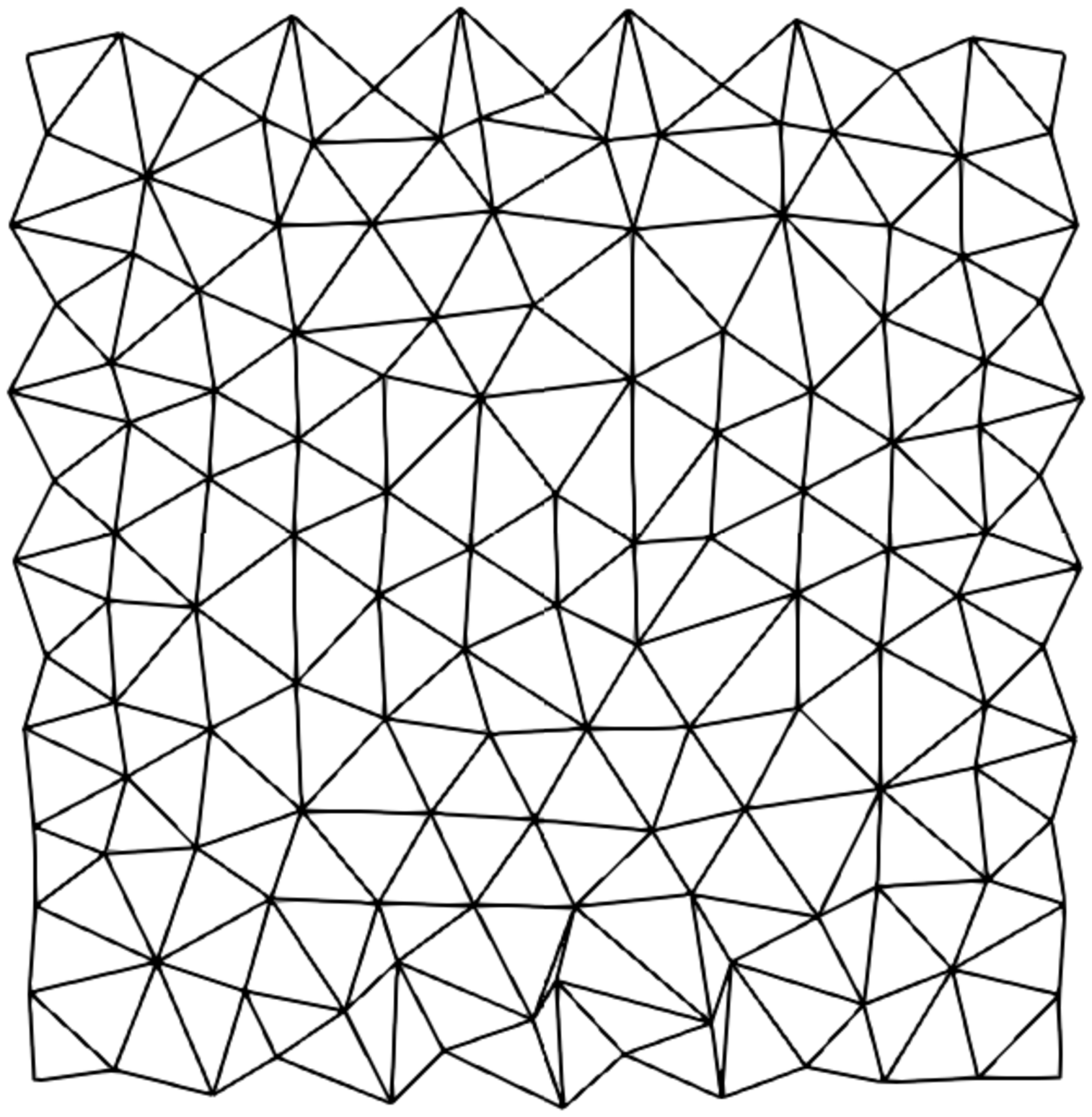, width = 0.26\textwidth}}
}
\mbox{
\subfigure[]{\label{fig:numstability_quartic_prior_d}
\epsfig{file = 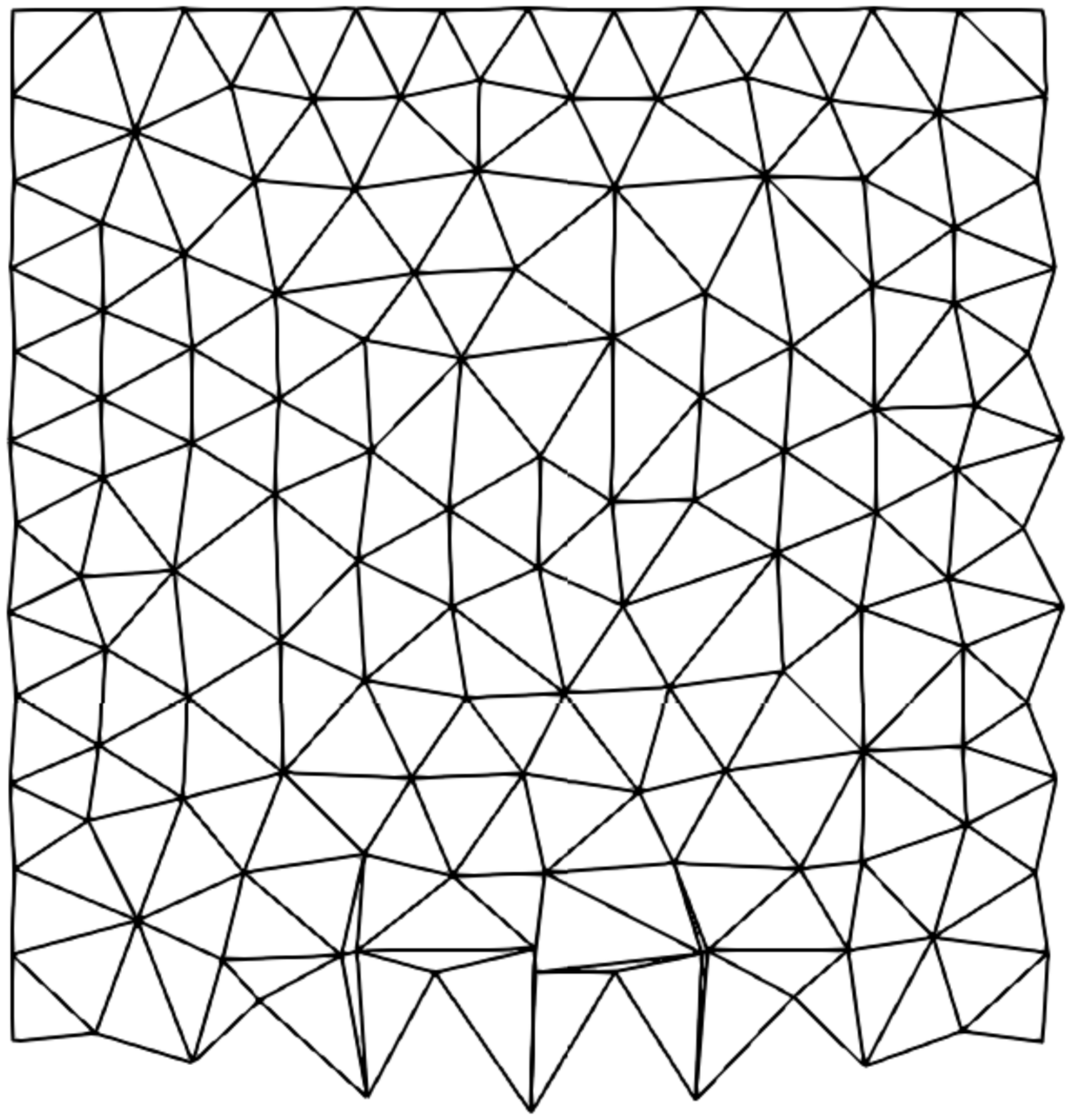, width = 0.26\textwidth}}
\subfigure[]{\label{fig:numstability_quartic_prior_e}
\epsfig{file = 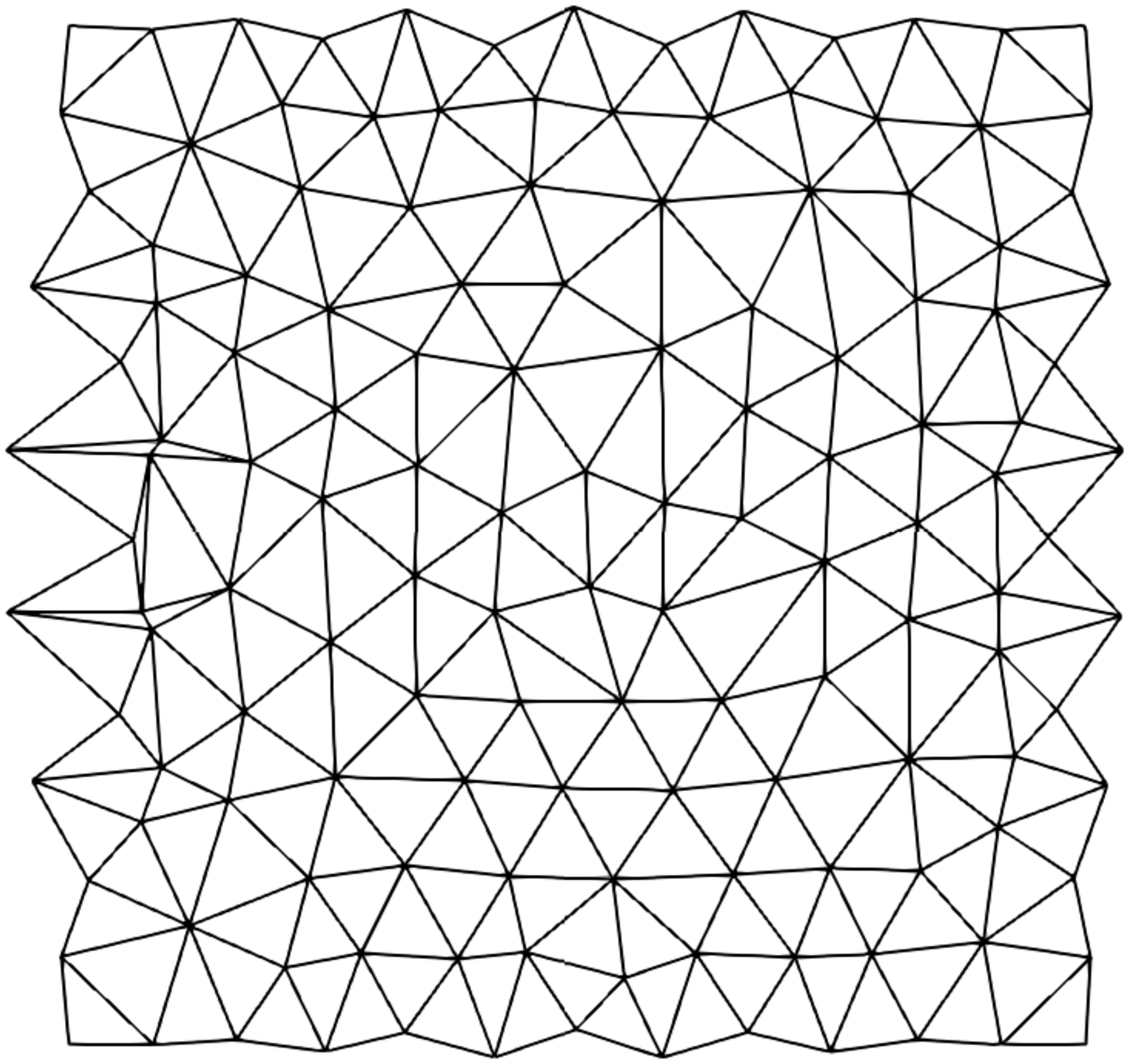, width = 0.26\textwidth}}
\subfigure[]{\label{fig:numstability_quartic_prior_f}
\epsfig{file = 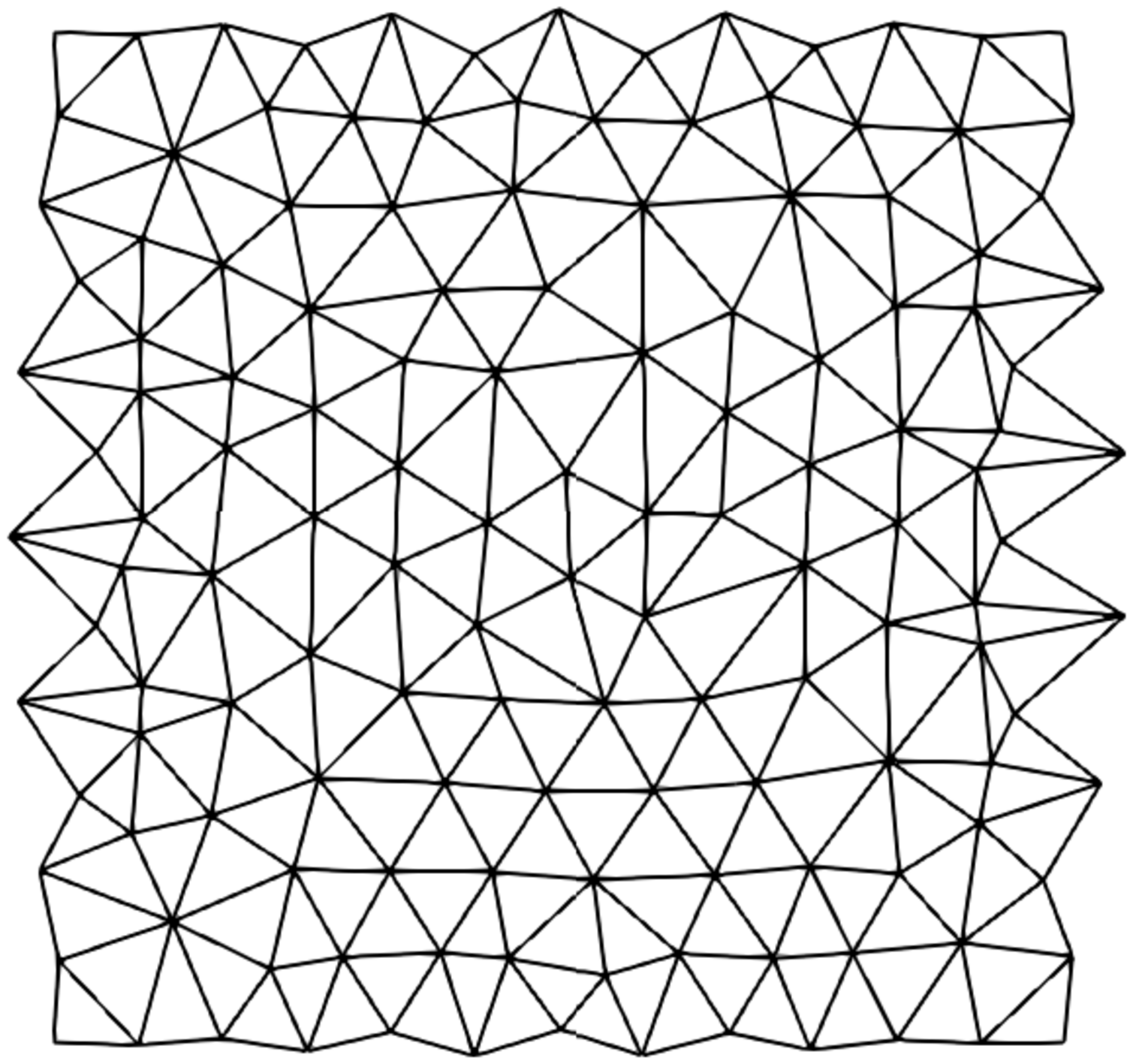, width = 0.26\textwidth}}
}
\mbox{
\subfigure[]{\label{fig:numstability_quartic_prior_g}
\epsfig{file = 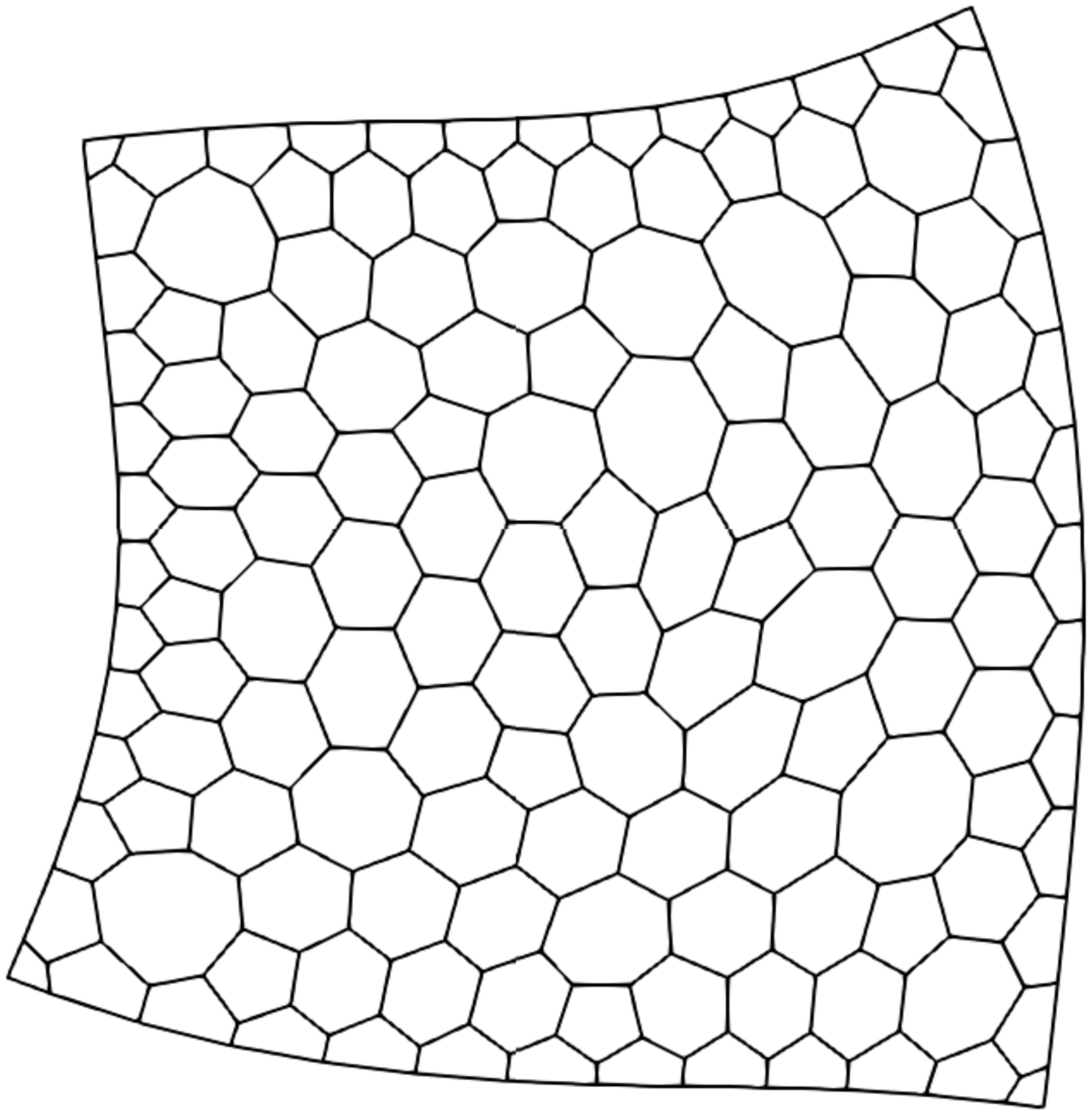, width = 0.26\textwidth}}
\subfigure[]{\label{fig:numstability_quartic_prior_h}
\epsfig{file = 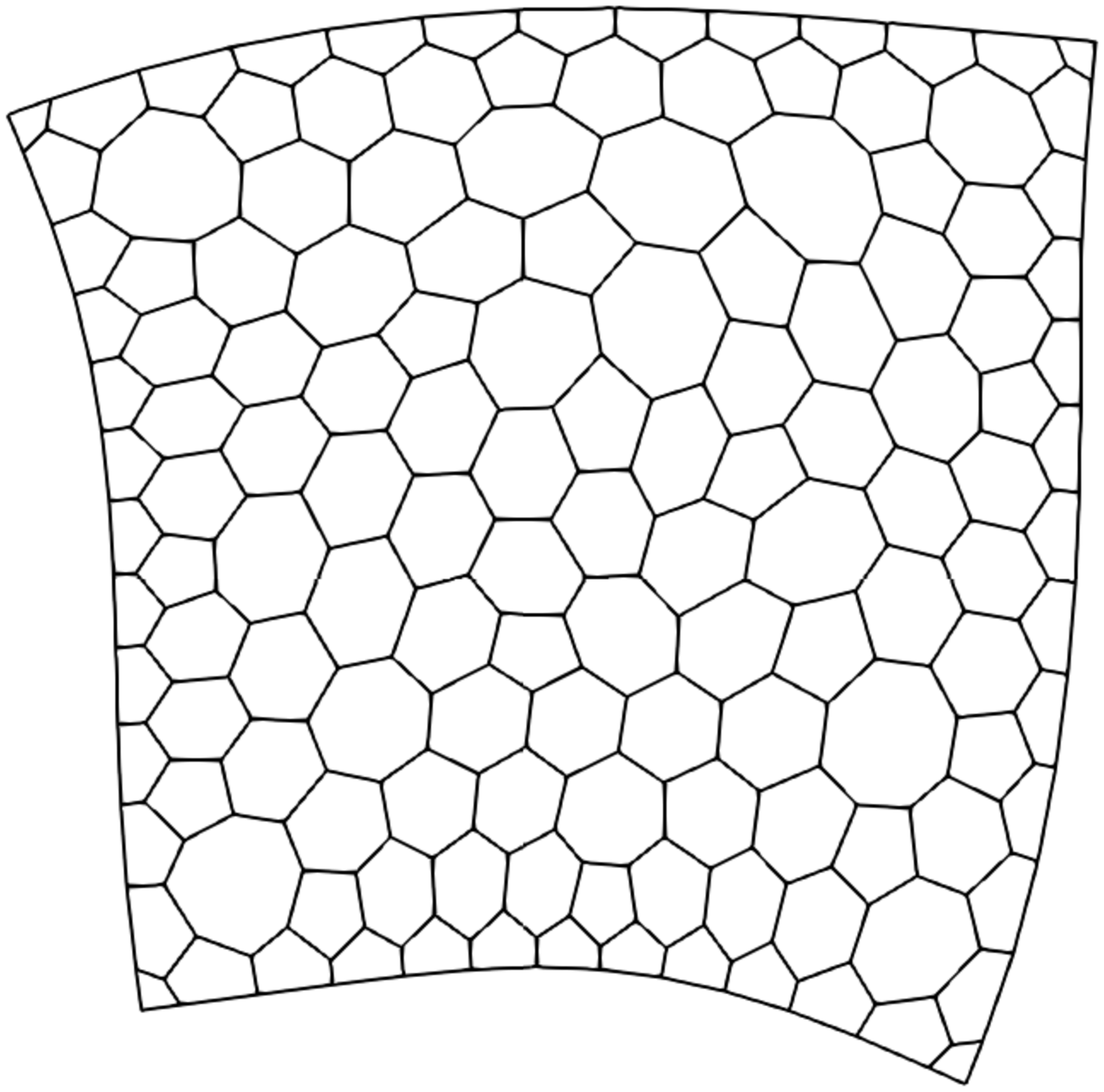, width = 0.26\textwidth}}
\subfigure[]{\label{fig:numstability_quartic_prior_i}
\epsfig{file = 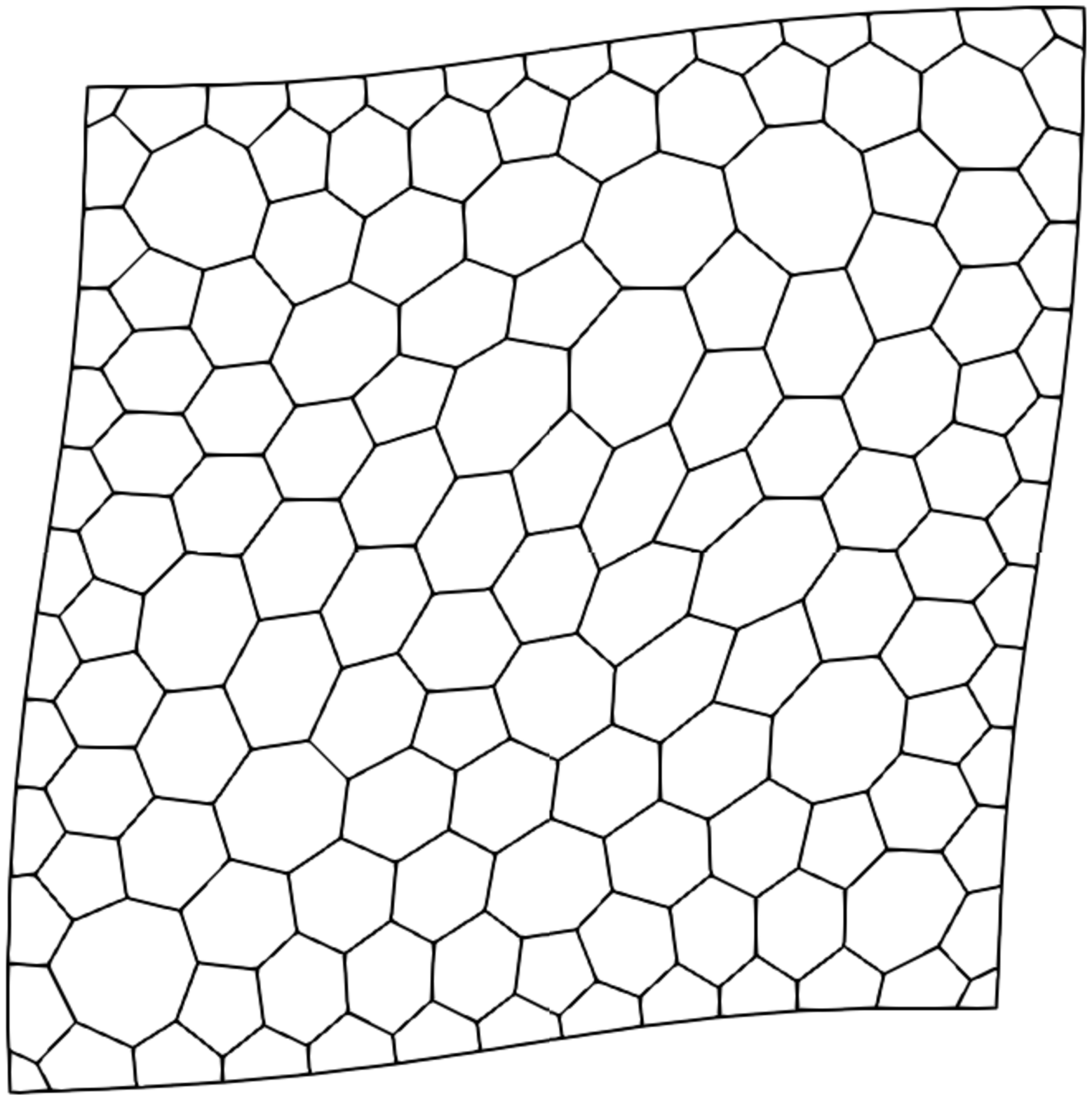, width = 0.26\textwidth}}
}
\caption{Eigenvalue analyses for the MEM and NIVED methods with
the maximum-entropy basis functions computed using the quartic prior. 
Depiction of the three mode shapes that follow
the three rigid-body modes. (a)-(c) MEM (1-pt), (d)-(f) MEM (12-pt), (g)-(i) NIVED.
Instabilities in the MEM method are evidenced by the presence of nonsmooth mode shapes.}
\label{fig:numstability_quartic_prior}
\end{figure}

To assess the stability of the NIVED method,
eigenvalue analyses are performed on a unit square domain. 
As a reference for comparison, we include the MEM approach
in the analyses. The quartic polynomials and the Gaussian radial basis function
are considered as the prior weight function
in the evaluation of the maximum-entropy basis functions in these analyses.
Three zero eigenvalues, which correspond to the three zero-energy rigid-body modes, 
are obtained for both methods. 
The three mode shapes that follow the three rigid-body 
modes are depicted in \fref{fig:numstability_quartic_prior} for the MEM and NIVED methods
using the quartic polynomials as the prior weight function. For the quartic prior,
instabilities are observed for the MEM with a 1-point Gauss rule since it exhibits nonsmooth
mode shapes (Figures~\ref{fig:numstability_quartic_prior_a}--\subref{fig:numstability_quartic_prior_c}) 
and even a 12-point Gauss rule is not sufficient for removing the instabilities 
(Figures~\ref{fig:numstability_quartic_prior_d}--\subref{fig:numstability_quartic_prior_f}).
In stark contrast, Figures~\ref{fig:numstability_quartic_prior_g}--\subref{fig:numstability_quartic_prior_i} 
show the smooth mode shapes that are obtained in the NIVED approach when using the quartic prior.

\begin{figure}[!tbhp]
\centering
\mbox{
\subfigure[]{\label{fig:numstability_gaussian_prior_a}
\epsfig{file = 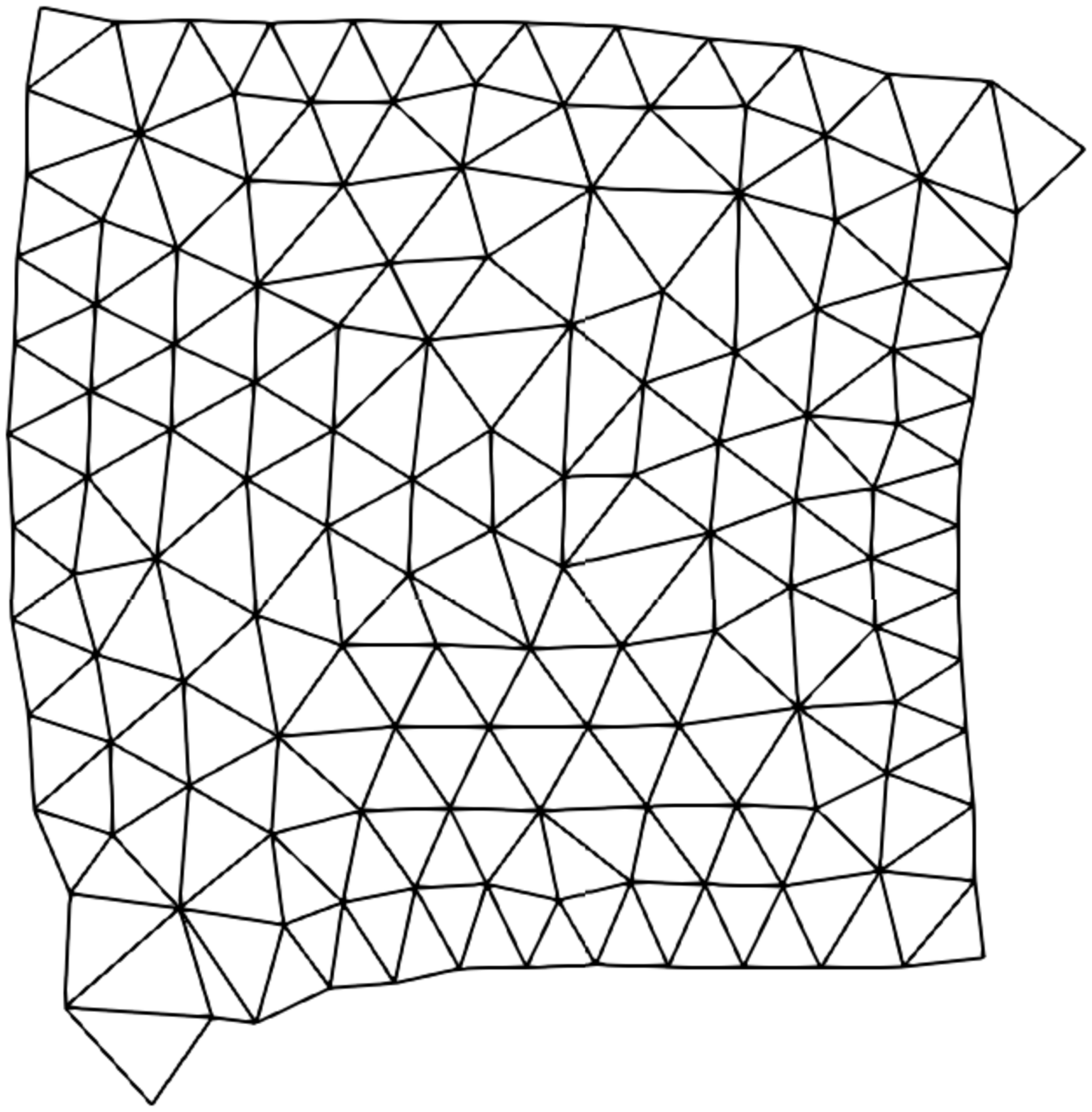, width = 0.26\textwidth}}
\subfigure[]{\label{fig:numstability_gaussian_prior_b}
\epsfig{file = 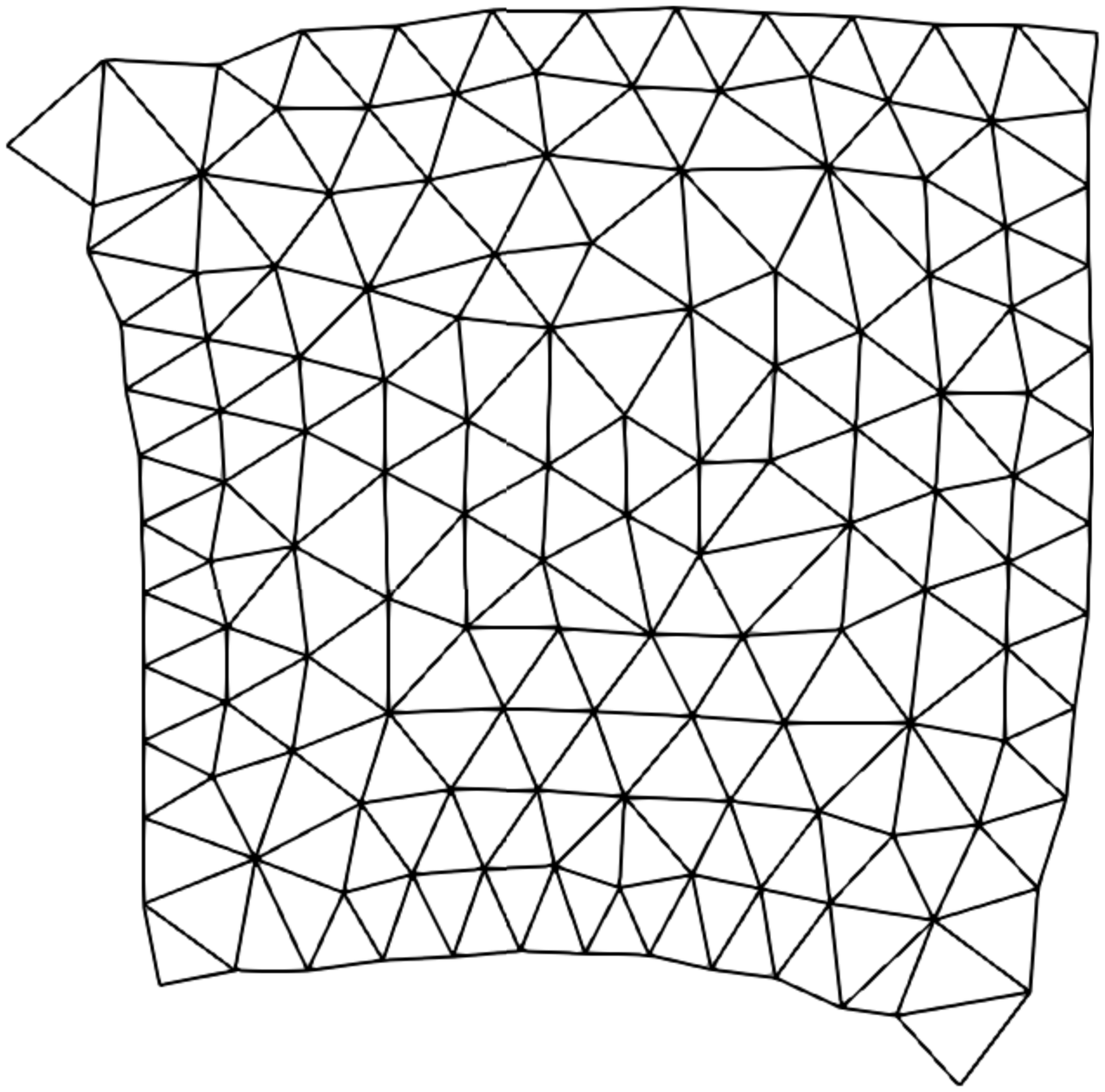, width = 0.26\textwidth}}
\subfigure[]{\label{fig:numstability_gaussian_prior_c}
\epsfig{file = 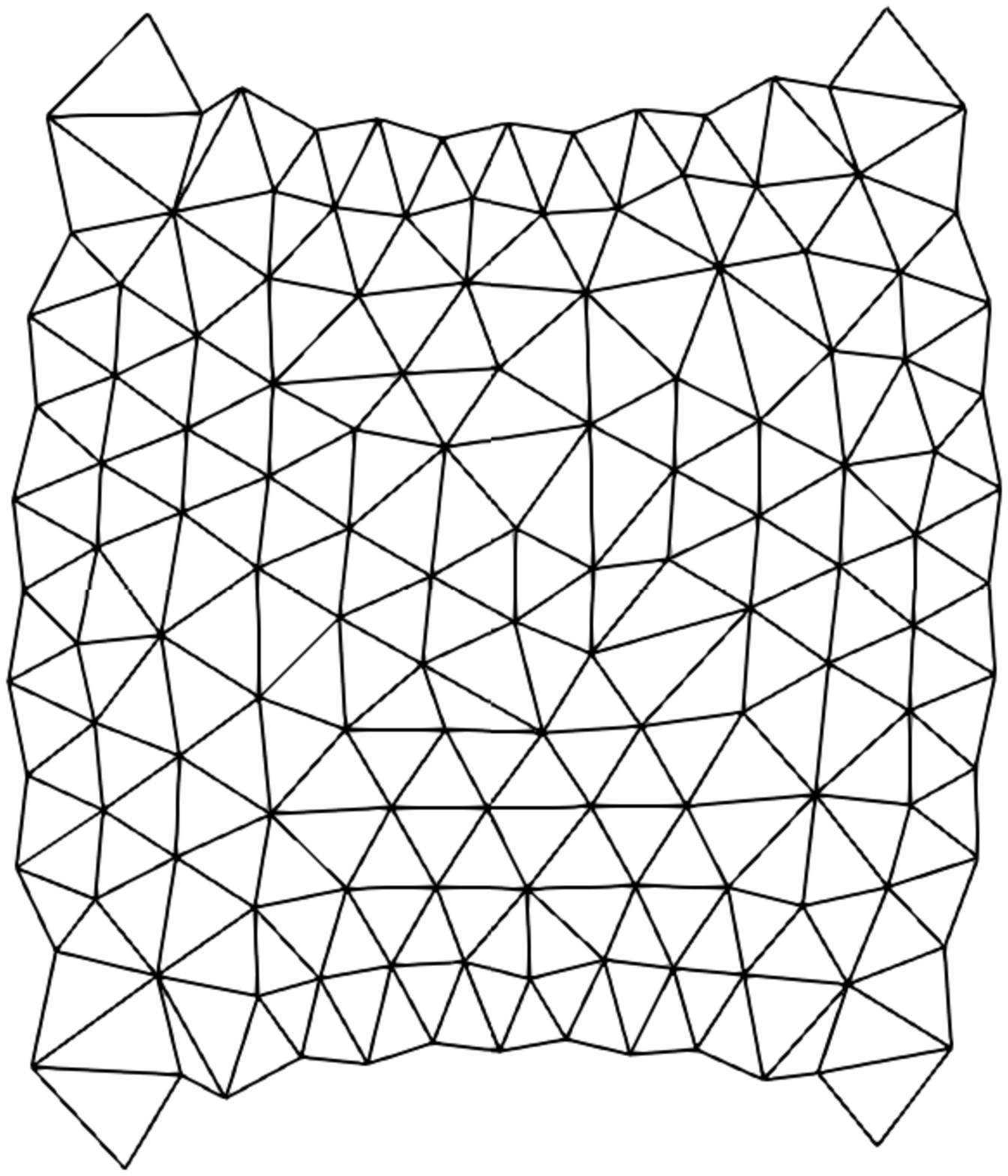, width = 0.26\textwidth}}
}
\mbox{
\subfigure[]{\label{fig:numstability_gaussian_prior_d}
\epsfig{file = 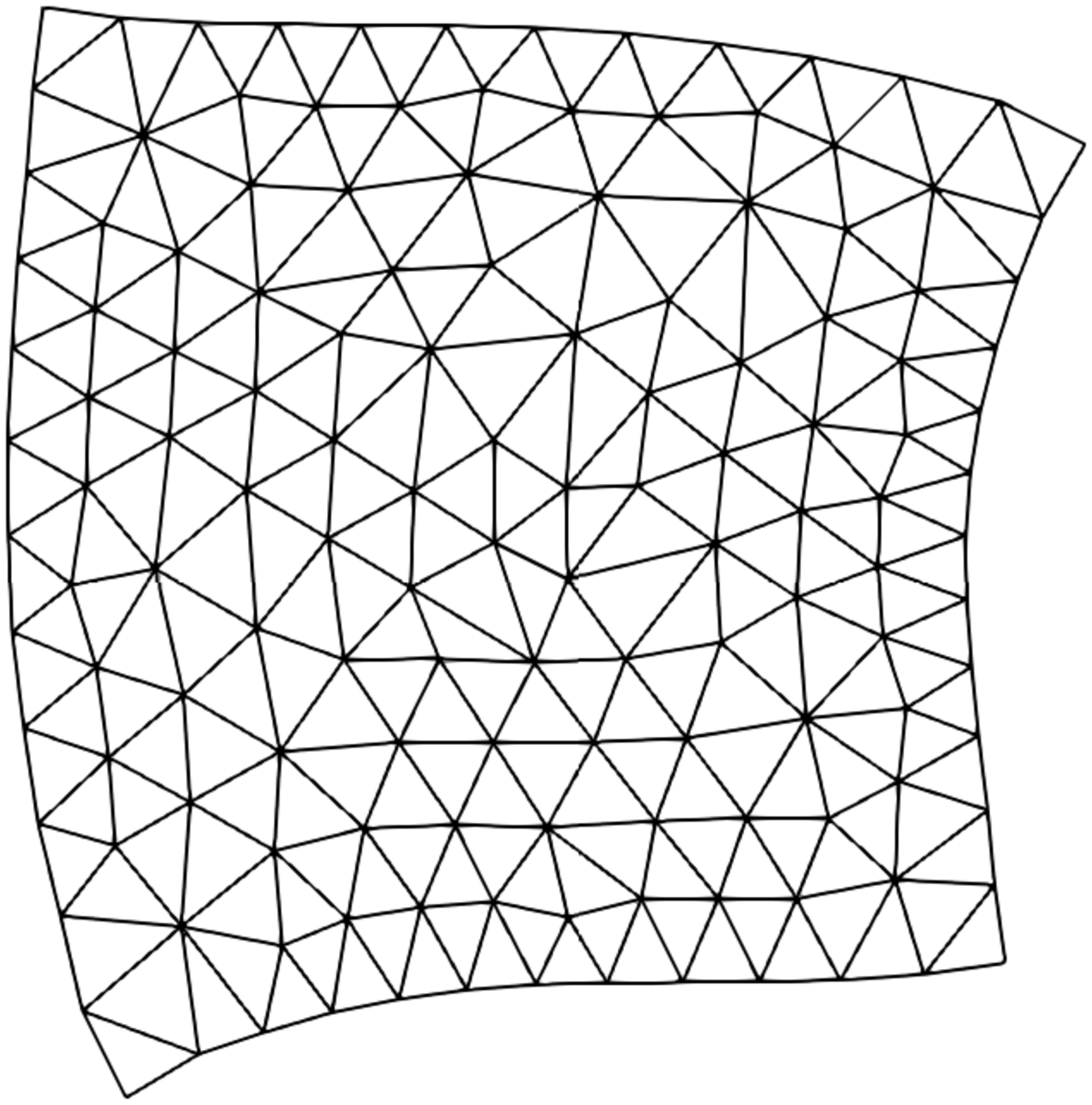, width = 0.26\textwidth}}
\subfigure[]{\label{fig:numstability_gaussian_prior_e}
\epsfig{file = 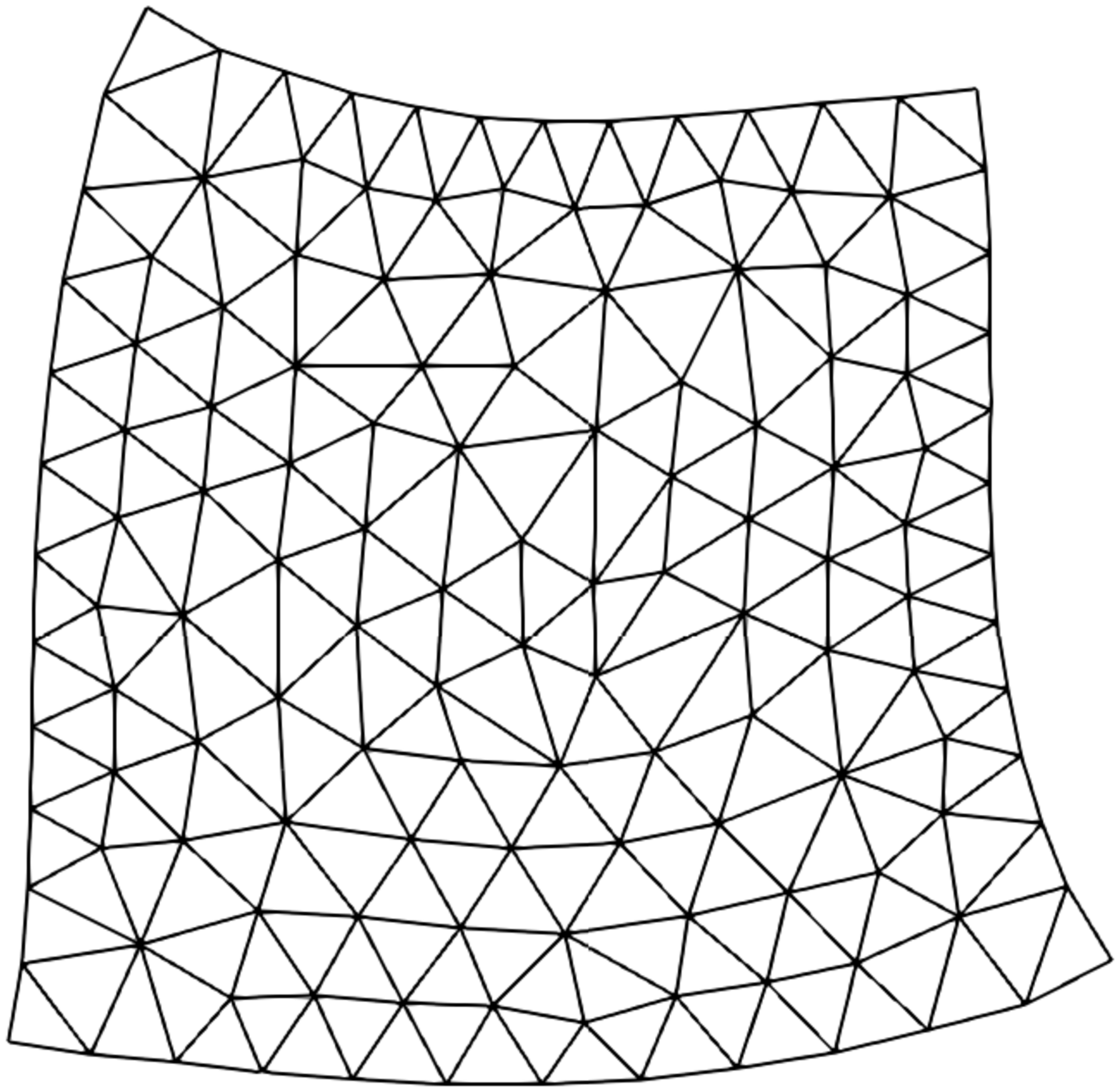, width = 0.26\textwidth}}
\subfigure[]{\label{fig:numstability_gaussian_prior_f}
\epsfig{file = 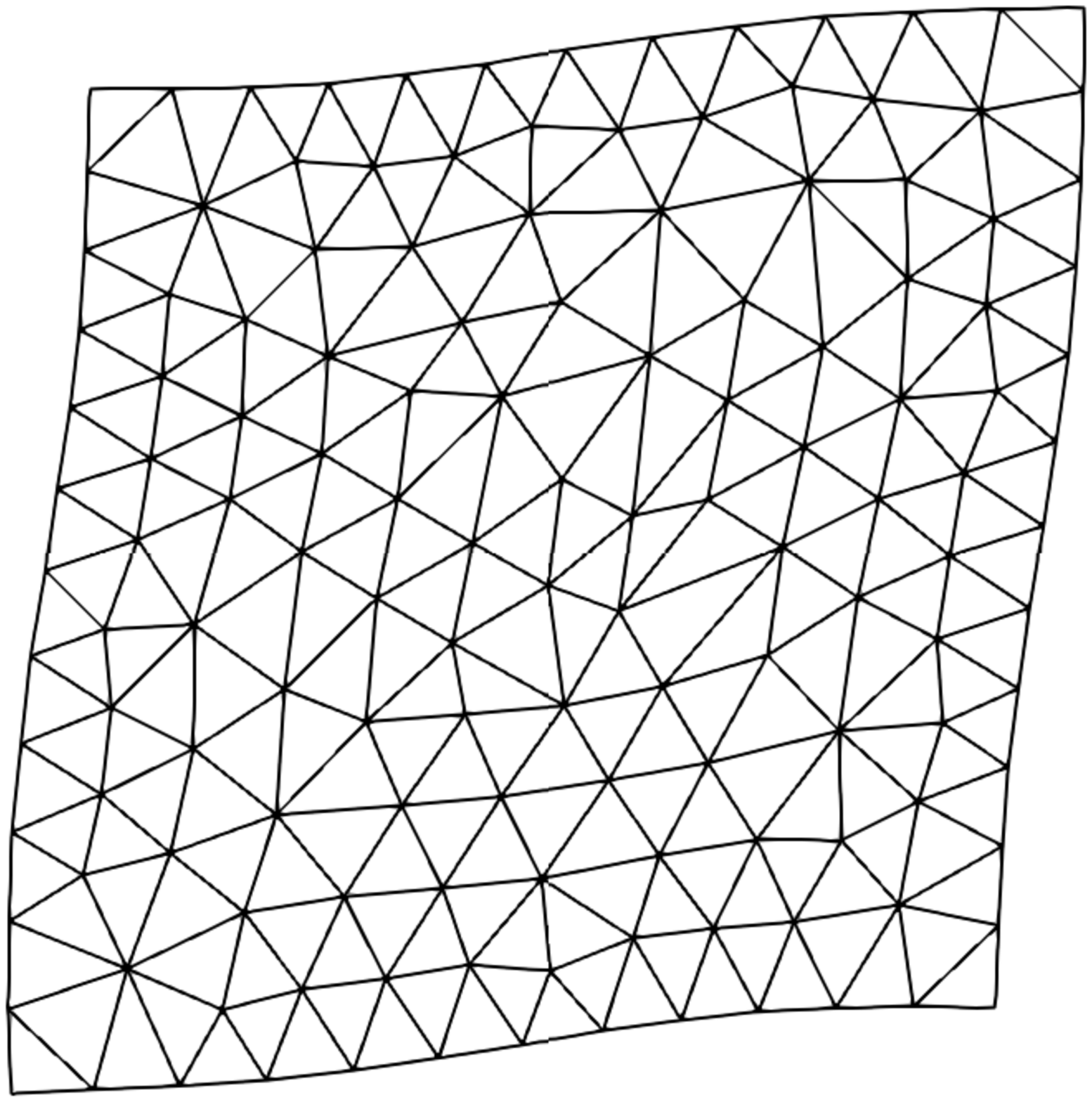, width = 0.26\textwidth}}
}
\mbox{
\subfigure[]{\label{fig:numstability_gaussian_prior_g}
\epsfig{file = 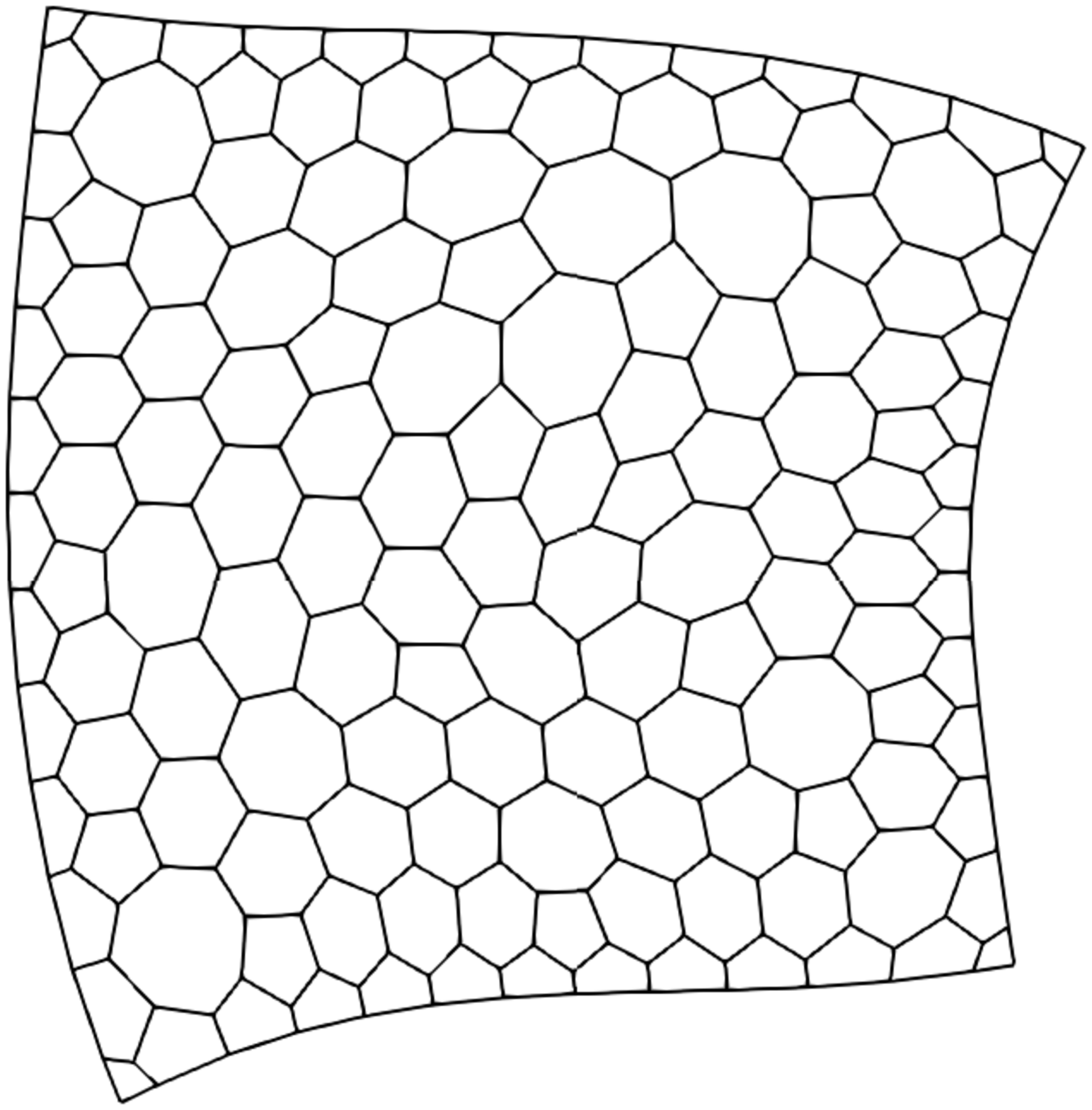, width = 0.26\textwidth}}
\subfigure[]{\label{fig:numstability_gaussian_prior_h}
\epsfig{file = 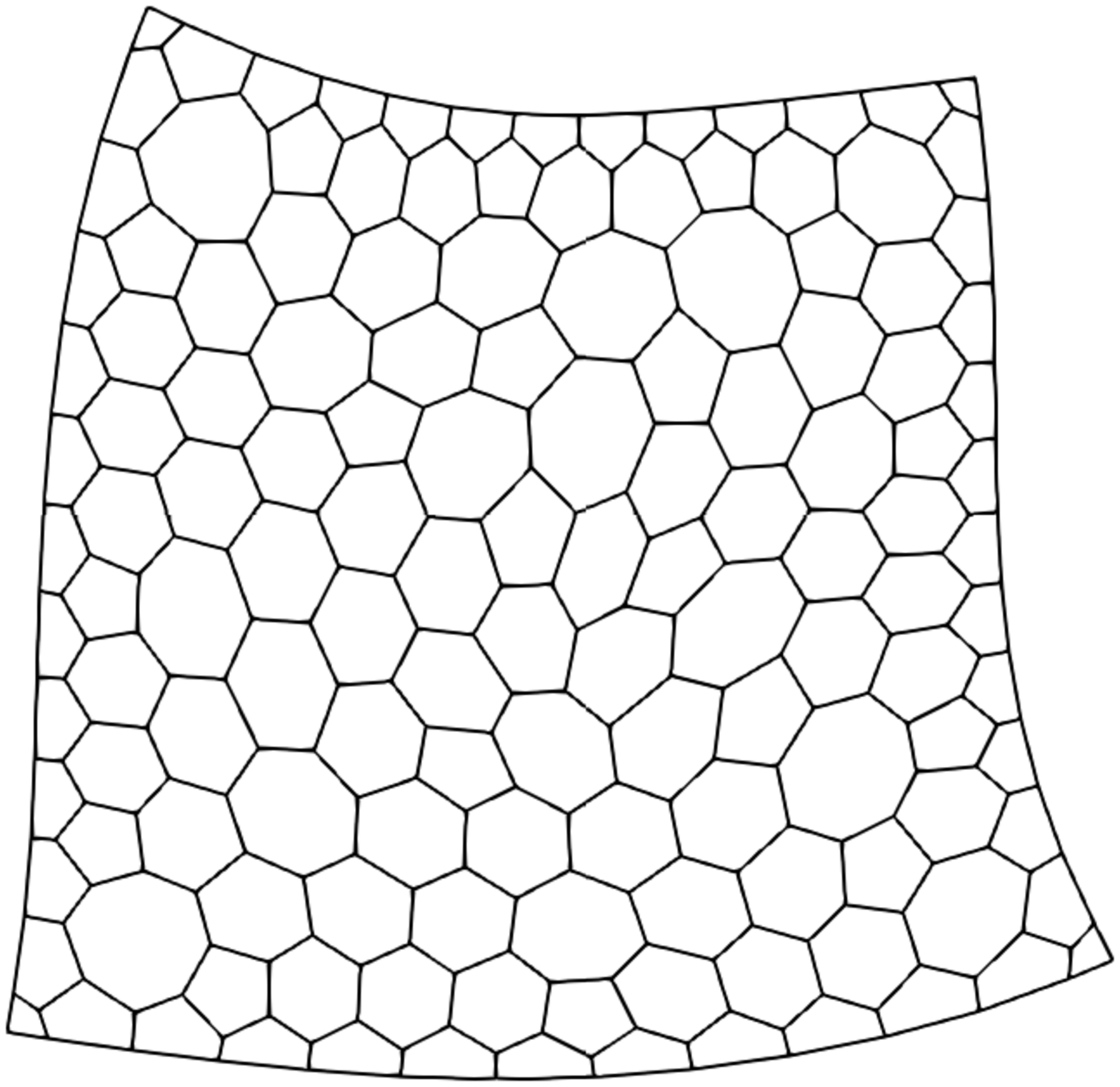, width = 0.26\textwidth}}
\subfigure[]{\label{fig:numstability_gaussian_prior_i}
\epsfig{file = 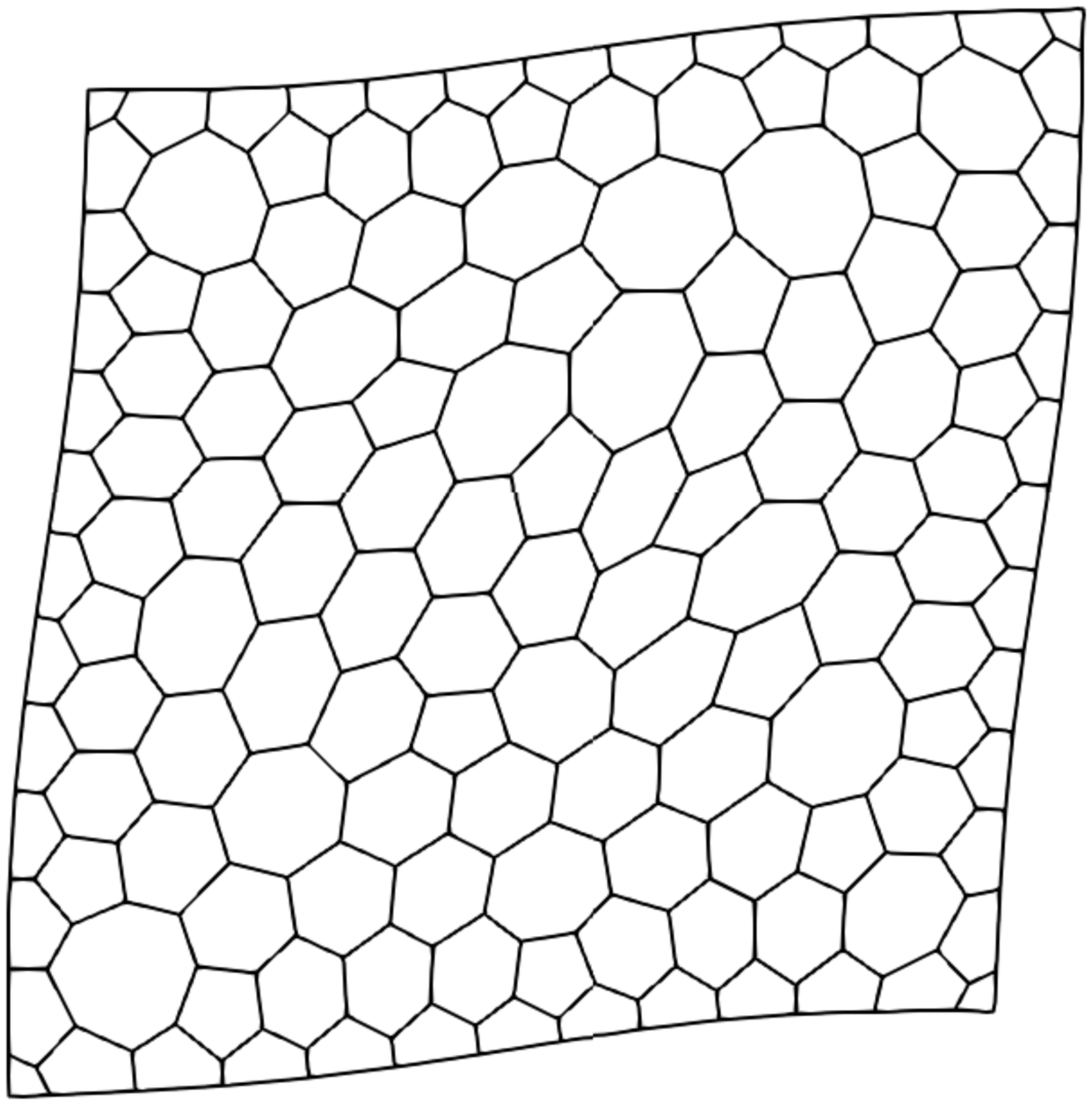, width = 0.26\textwidth}}
}
\caption{Eigenvalue analyses for the MEM and NIVED methods with
the maximum-entropy basis functions computed using the Gaussian prior. 
Depiction of the three mode shapes that follow
the three rigid-body modes. (a)-(c) MEM (1-pt), (d)-(f) MEM (6-pt), (g)-(i) NIVED.
Instabilities in the MEM method are evidenced by the presence of nonsmooth mode shapes.}
\label{fig:numstability_gaussian_prior}
\end{figure}

Similarly, the three mode shapes that follow the three rigid-body 
modes are depicted in \fref{fig:numstability_gaussian_prior} for the MEM and NIVED methods
using the Gaussian radial basis function as the prior weight function. For the Gaussian prior,
instabilities are observed for the MEM with a 1-point Gauss rule since it exhibits nonsmooth
mode shapes (Figures~\ref{fig:numstability_gaussian_prior_a}--\subref{fig:numstability_gaussian_prior_c}), 
but a 6-point Gauss rule can effectively remove the instabilities 
(Figures~\ref{fig:numstability_gaussian_prior_d}--\subref{fig:numstability_gaussian_prior_f}).
The NIVED approach using the Gaussian prior 
is also free of instabilities as revealed by the smooth mode shapes depicted in 
Figures~\ref{fig:numstability_gaussian_prior_g}--\subref{fig:numstability_gaussian_prior_i}.

\subsection{Cantilever beam}\label{sec:numexamples_beam}
\begin{figure}[!tbhp]
\centering
\epsfig{file = 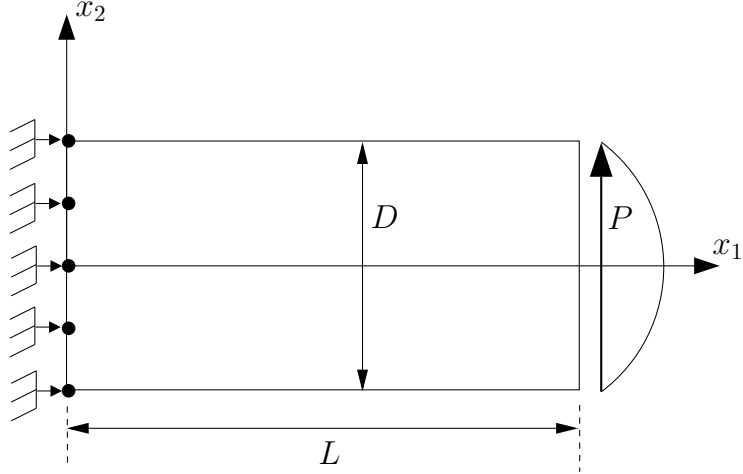, width = 0.65\textwidth}
\caption{Cantilever beam problem.}
\label{fig:beamproblem}
\end{figure}

\begin{figure}[!tbhp]
\centering
\mbox{
\subfigure[]{\label{fig:beam_mesh_a} \epsfig{file = 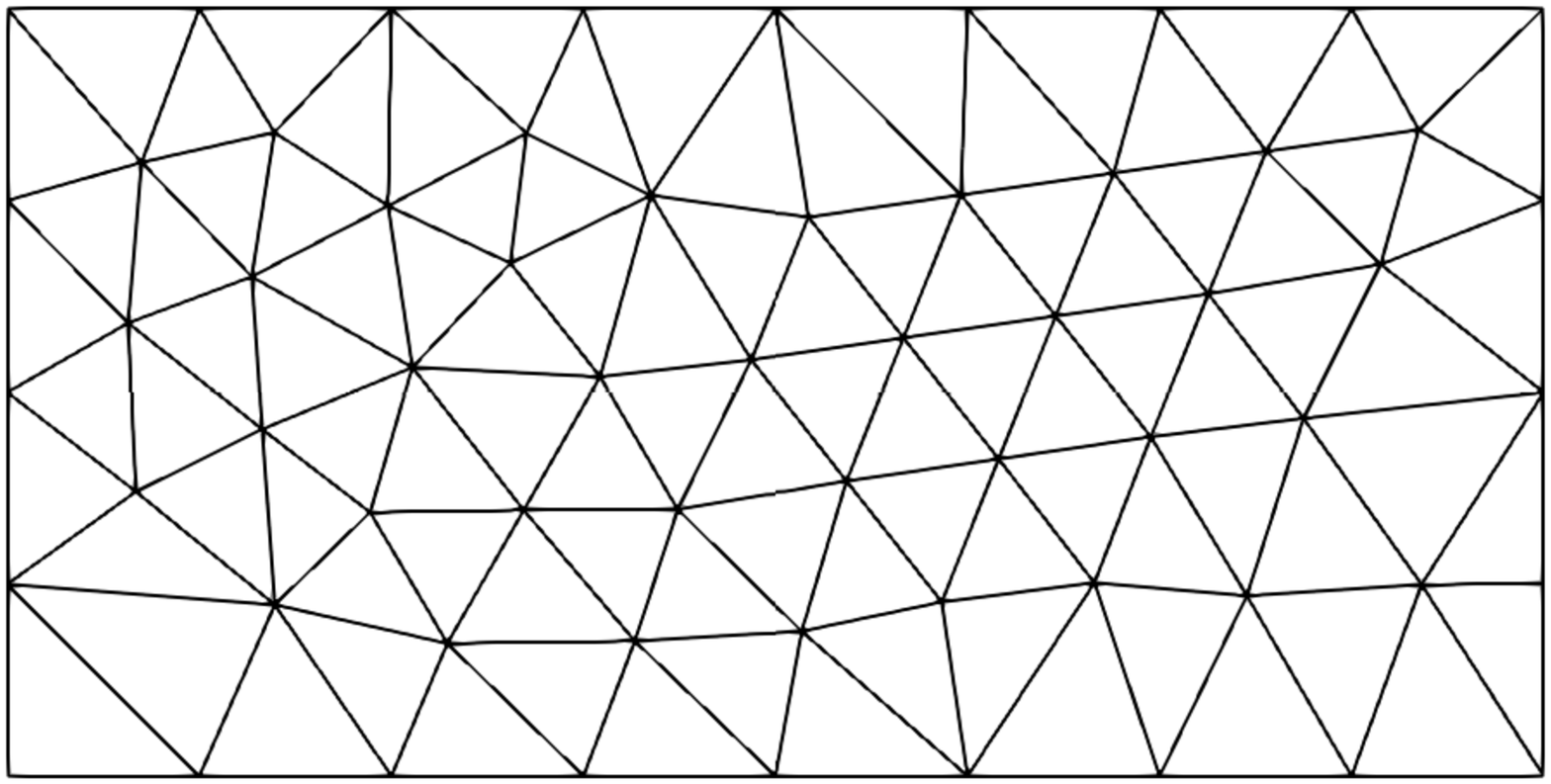, width = 0.23\textwidth}}
\subfigure[]{\label{fig:beam_mesh_b} \epsfig{file = 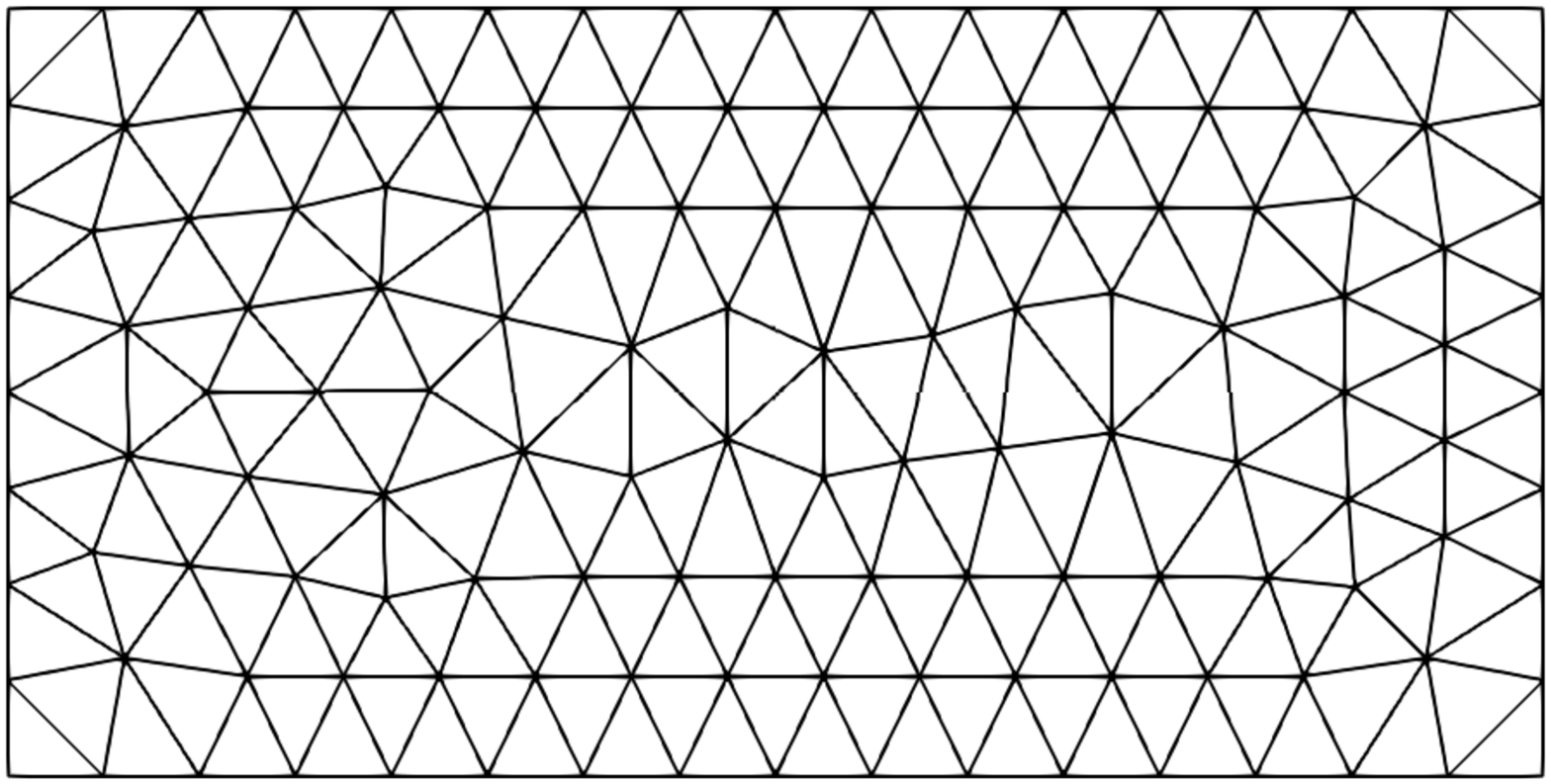, width = 0.23\textwidth}}
\subfigure[]{\label{fig:beam_mesh_c} \epsfig{file = 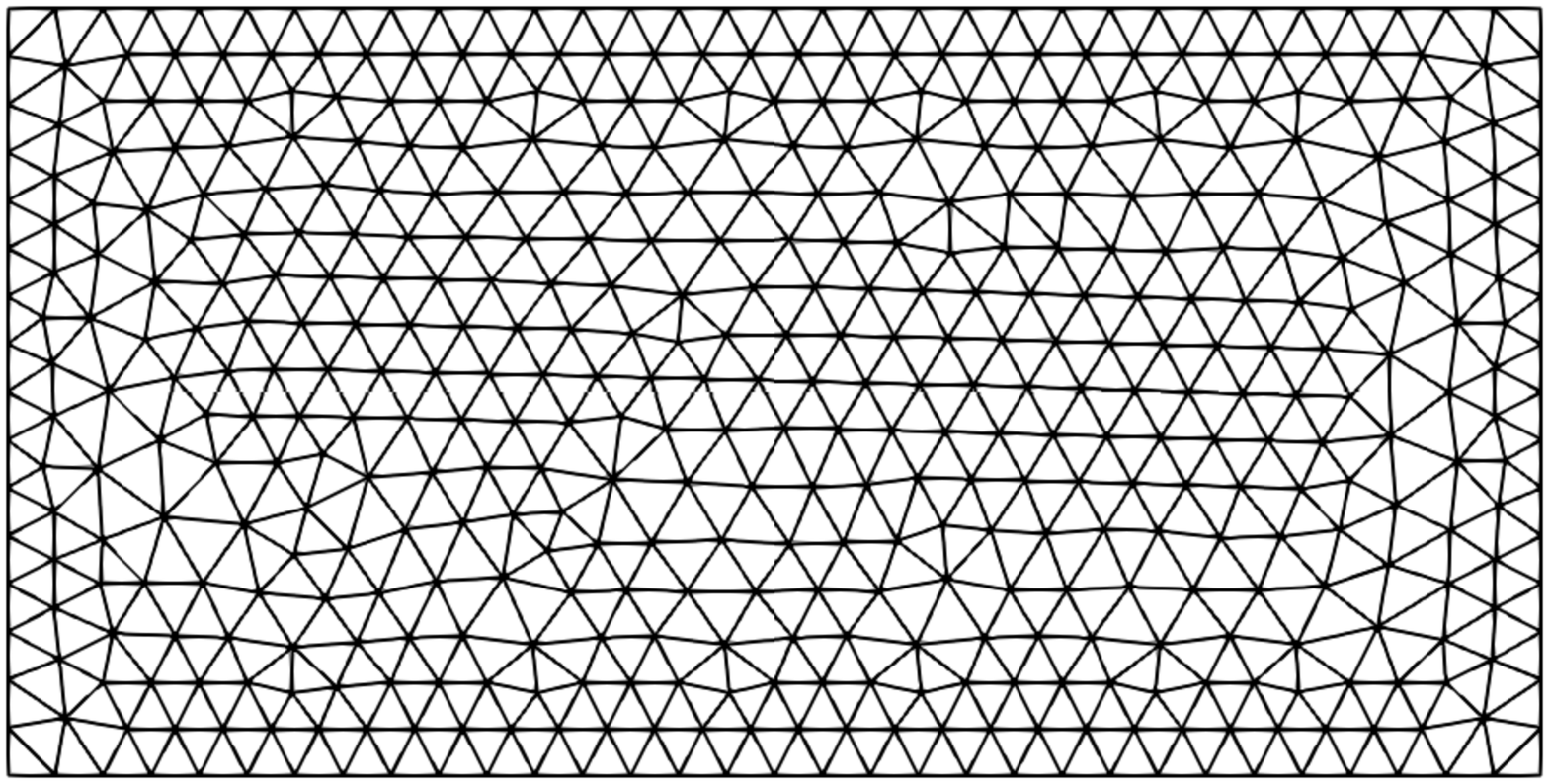, width = 0.23\textwidth}}
\subfigure[]{\label{fig:beam_mesh_d} \epsfig{file = 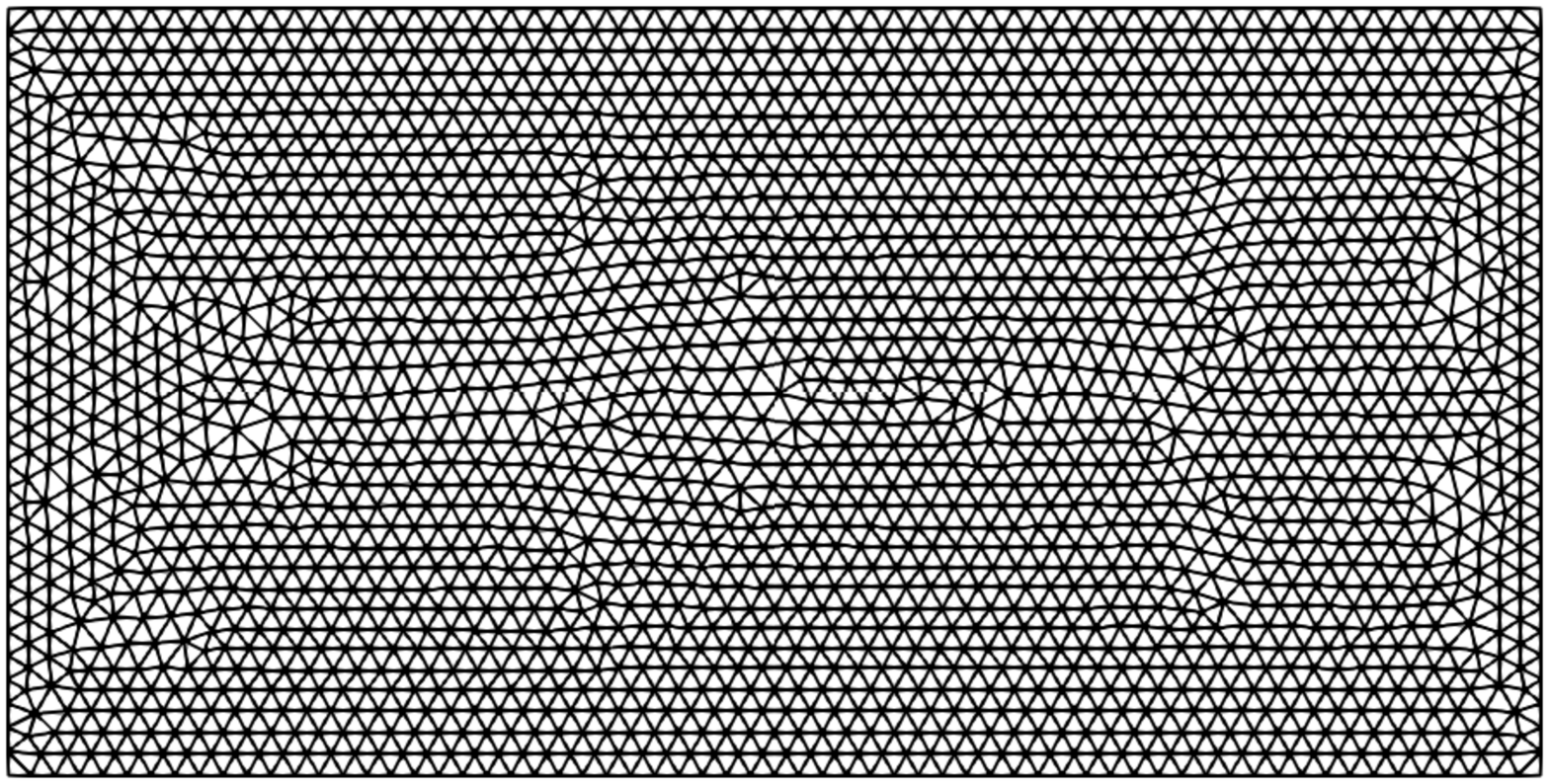, width = 0.23\textwidth}}
}
\mbox{
\subfigure[]{\label{fig:beam_mesh_e} \epsfig{file = 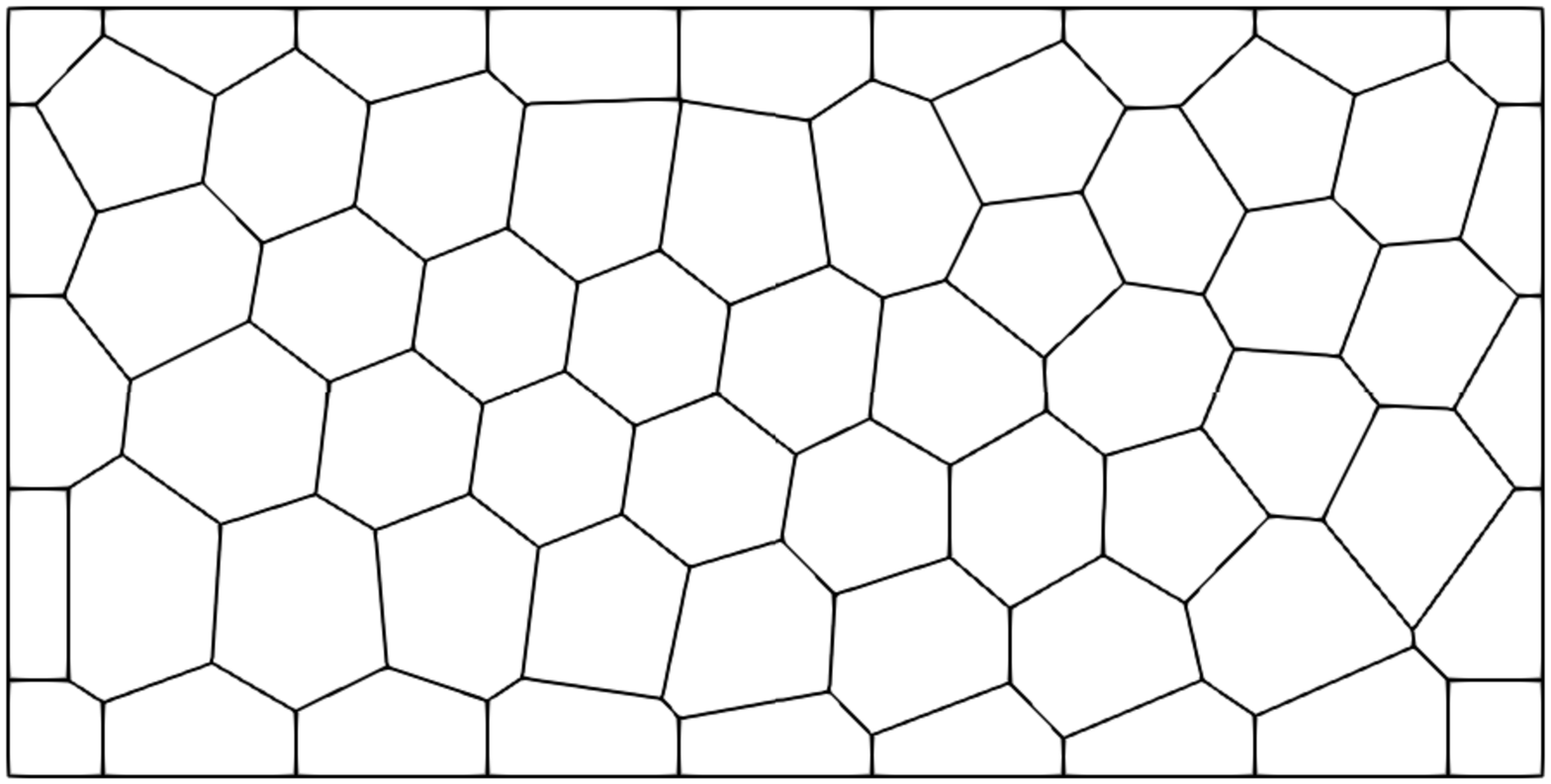, width = 0.23\textwidth}}
\subfigure[]{\label{fig:beam_mesh_f} \epsfig{file = 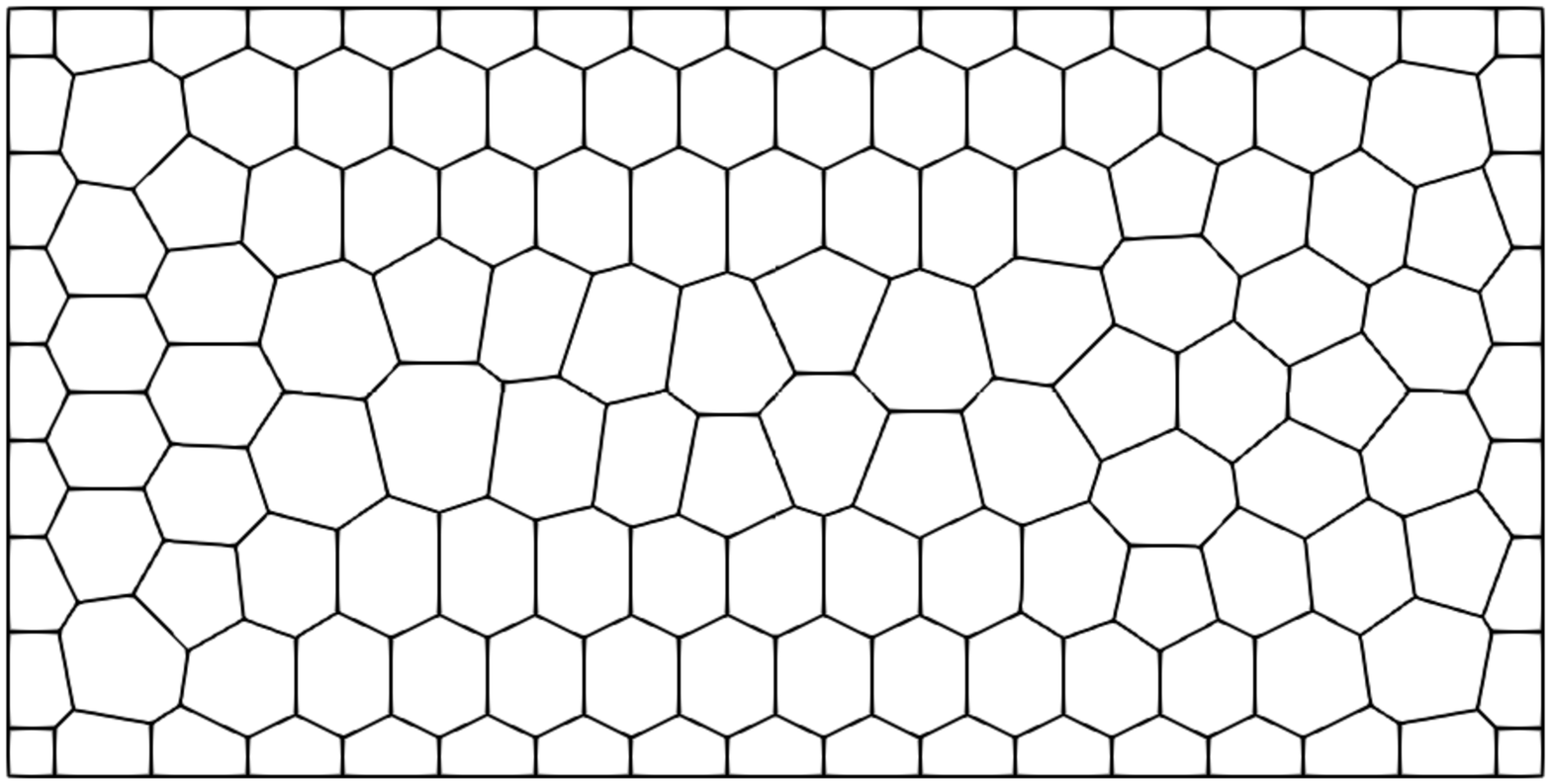, width = 0.23\textwidth}}
\subfigure[]{\label{fig:beam_mesh_g} \epsfig{file = 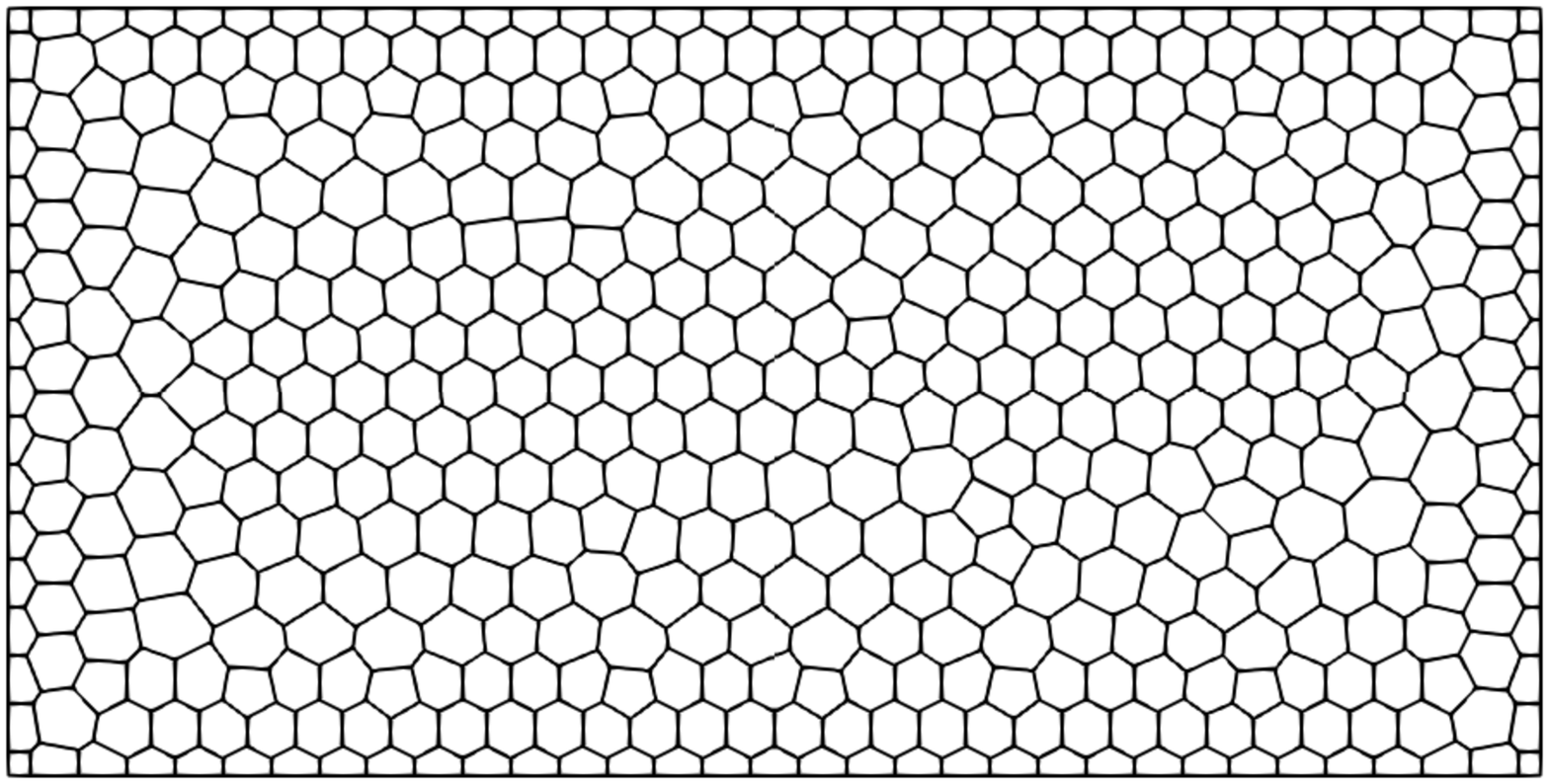, width = 0.23\textwidth}}
\subfigure[]{\label{fig:beam_mesh_h} \epsfig{file = 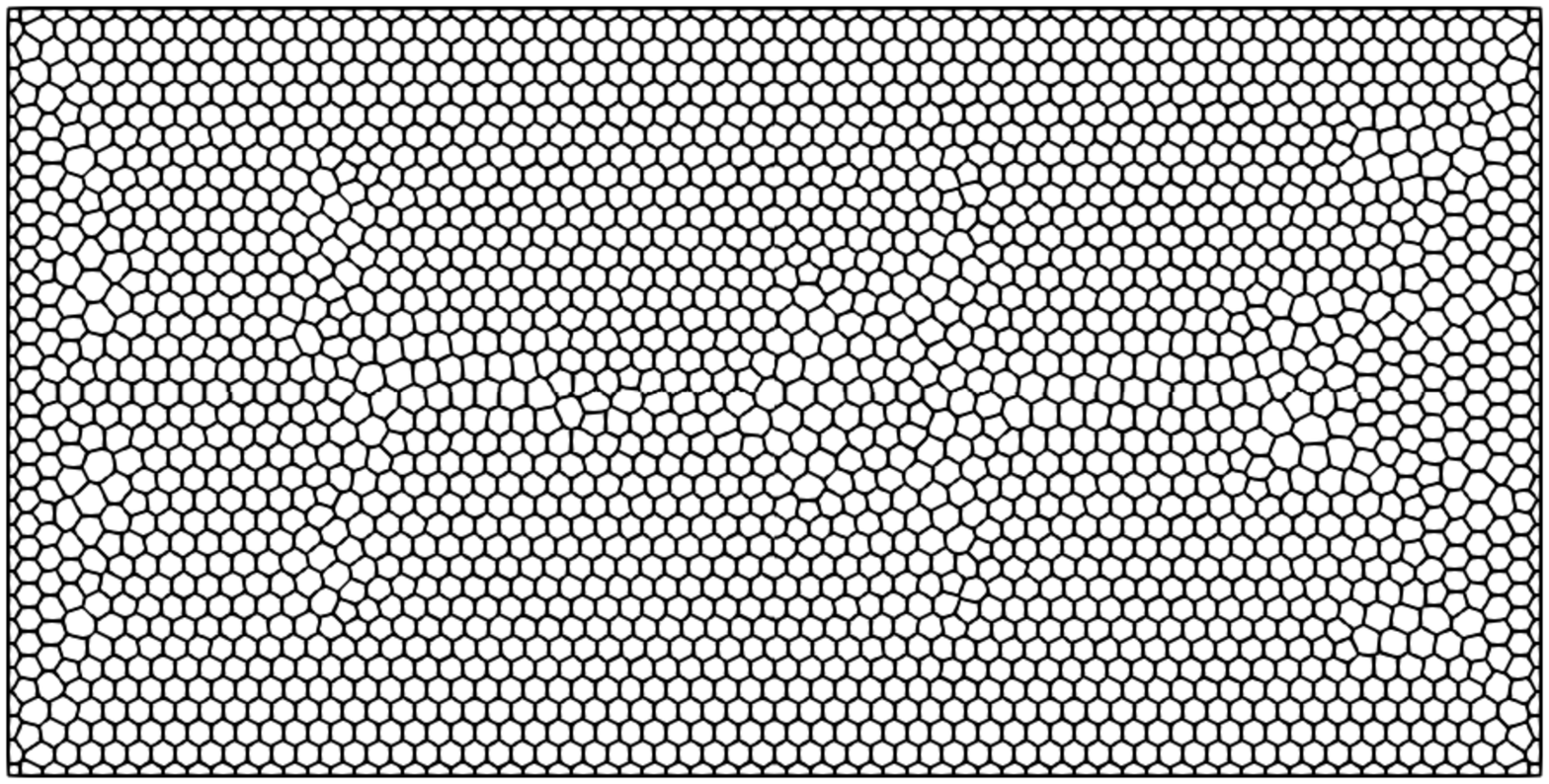, width = 0.23\textwidth}}
}
\caption{Sequence of background integration meshes used for the cantilever beam problem.
         (a)-(d) Meshes for the MEM method and (e)-(h) meshes for the NIVED method.}
\label{fig:beam_mesh}
\end{figure}

We conduct a convergence study 
for the problem of a cantilever beam
of unit thickness subjected to a parabolic end load $P=-1000$ lbf.
A schematic representation of the problem is shown in \fref{fig:beamproblem}. 
Plane strain condition is assumed with material parameters given by 
$E_\mathrm{Y}=10^7$ psi and $\nu = 0.3$. The essential boundary 
conditions on the clamped edge are applied according to the analytical 
solution given by Timoshenko and Goodier~\cite{timoshenko:1970:TOE}:
\begin{equation*}\label{beam_exact_sol}
\vm{u}=
\left[
\begin{array}{c}
-\frac{Px_2}{6\bar{E}_\mathrm{Y}I}\left((6L-3x_1)x_1 + (2+\bar{\nu})x_2^2 - \frac{3D^2}{2}(1+\bar{\nu})\right)\\
\frac{P}{6\bar{E}_\mathrm{Y}I}\left(3\bar{\nu}x_2^{2}(L-x_1)+(3L-x_1)x_1^{2}\right)
\end{array}
\right],
\end{equation*}
where $\bar{E}_\mathrm{Y}=E_\mathrm{Y}/\left(1-\nu^{2}\right)$ and $\bar{\nu}=\nu/\left(1-\nu\right)$;
$L=8$ inch is the length of the beam, $D=4$ inch is the height of 
the beam, and $I=D^3/12$ is the second-area moment of the beam section.
The exact stress field is:
\begin{equation*}
\left[\begin{array}{c}
\sigma_{11}\\
\sigma_{22}\\
\sigma_{12}
\end{array}
\right]
=\left[\begin{array}{c}
\frac{-P(L-x_1)x_2}{I}\\
0\\
\frac{P}{2I}\left(\frac{D^2}{4}-x_2^2\right)
\end{array}
\right].
\end{equation*}
The Neumann boundary conditions are applied using the exact 
stress field, which gives $\bar{\vm{t}}=\smat{0 & 0}^\transpose$
on the top and bottom boundaries, and
$\bar{\vm{t}}=\smat{0 & \sigma_{12}}^\transpose$
on the right boundary. For the evaluation of the maximum-entropy
basis functions, the Gaussian prior is used.
The sequence of background integration meshes used in the 
study are shown in~\fref{fig:beam_mesh}. 

\begin{figure}[!tbhp]
\centering
\subfigure[]{\label{fig:beam_norms_a} \epsfig{file = 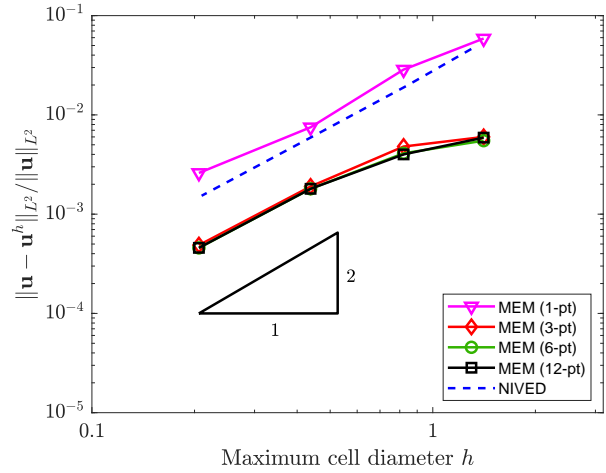, width = 0.58\textwidth}}
\subfigure[]{\label{fig:beam_norms_b} \epsfig{file = 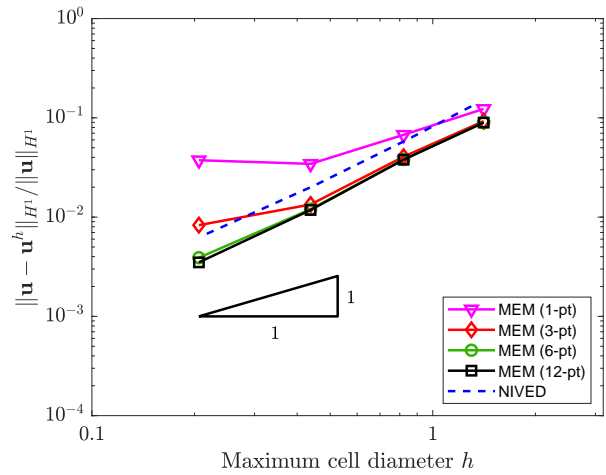, width = 0.58\textwidth}}
\caption{Convergence rates for the cantilever beam problem. Optimal rate of $2$ in the
         $L^2$ norm is delivered by the MEM method using a 3-point Gauss rule, but
         its convergence in the $H^1$ seminorm behaves erratically for 1-point and 3-point Gauss
         rules due to integration errors; a 6-point Gauss rule is needed to recover the optimal
         rate of $1$. Optimal rates of $2$ and $1$ in the $L^2$ norm and the $H^1$ seminorm,
         respectively, are delivered by the NIVED approach.}
\label{fig:beam_norms}
\end{figure}

The convergence of the MEM and NIVED methods in the $L^2$ norm and the $H^1$ 
seminorm of the error are shown in \fref{fig:beam_norms}. The MEM approach 
needs a 3-point Gauss rule to deliver the optimal rate in the $L^2$ norm. 
Regarding the $H^1$ seminorm, the convergence of the MEM method behaves 
erratically for 1-point and 3-point Gauss rules due to integration errors,
and a 6-point Gauss rule is needed to recover the optimal rate. 
In contrast, the proposed NIVED method delivers the optimal rate 
of convergence in both the $L^2$ norm and the $H^1$ seminorm of the error.

\subsection{Infinite plate with a circular hole}\label{sec:numexamples_infiniteplate}

\begin{figure}[!tbhp]
\centering
\mbox{
\subfigure[]{\label{fig:plate_hole_a} \epsfig{file = 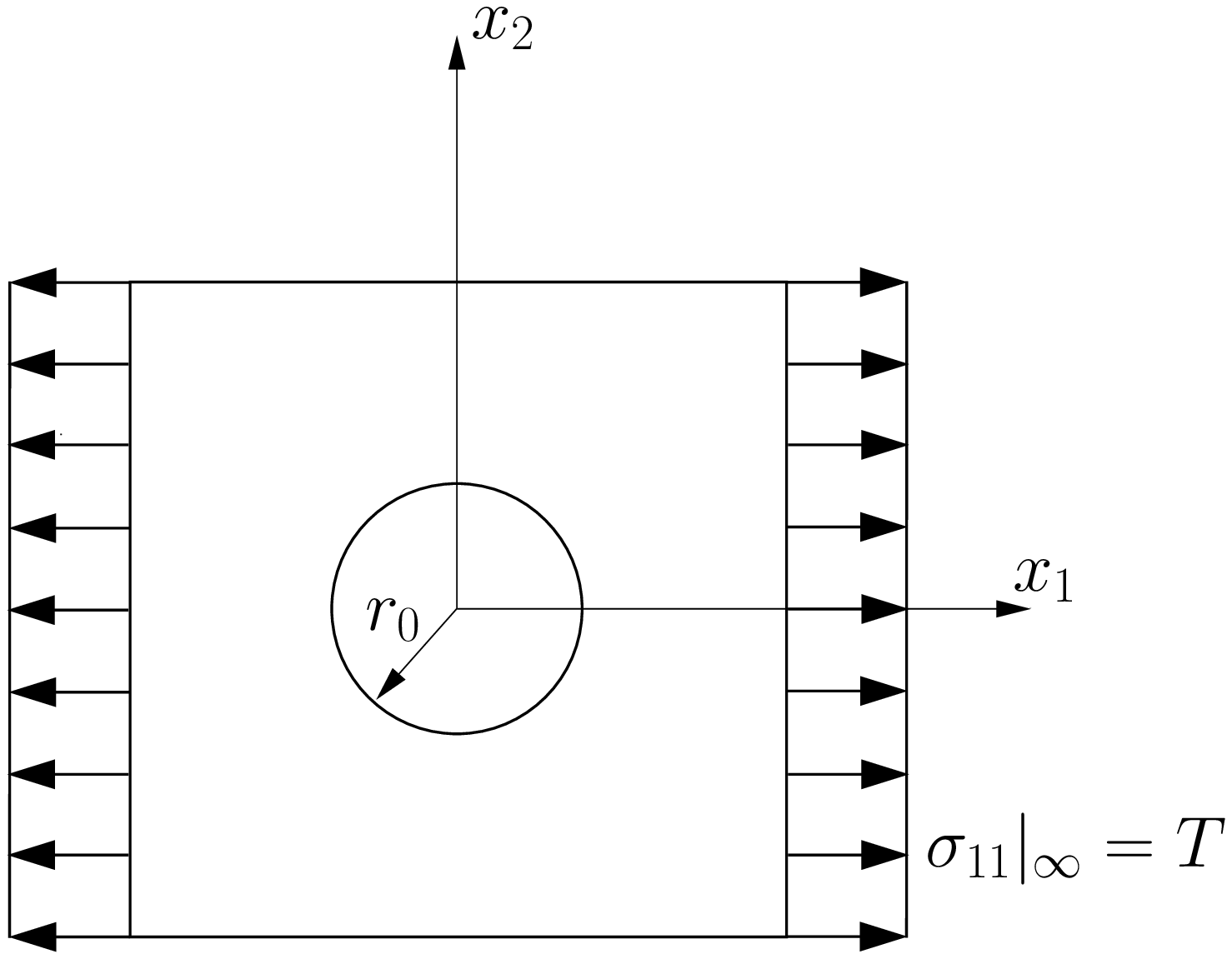, width = 0.5\textwidth}}
\subfigure[]{\label{fig:plate_hole_b} \epsfig{file = 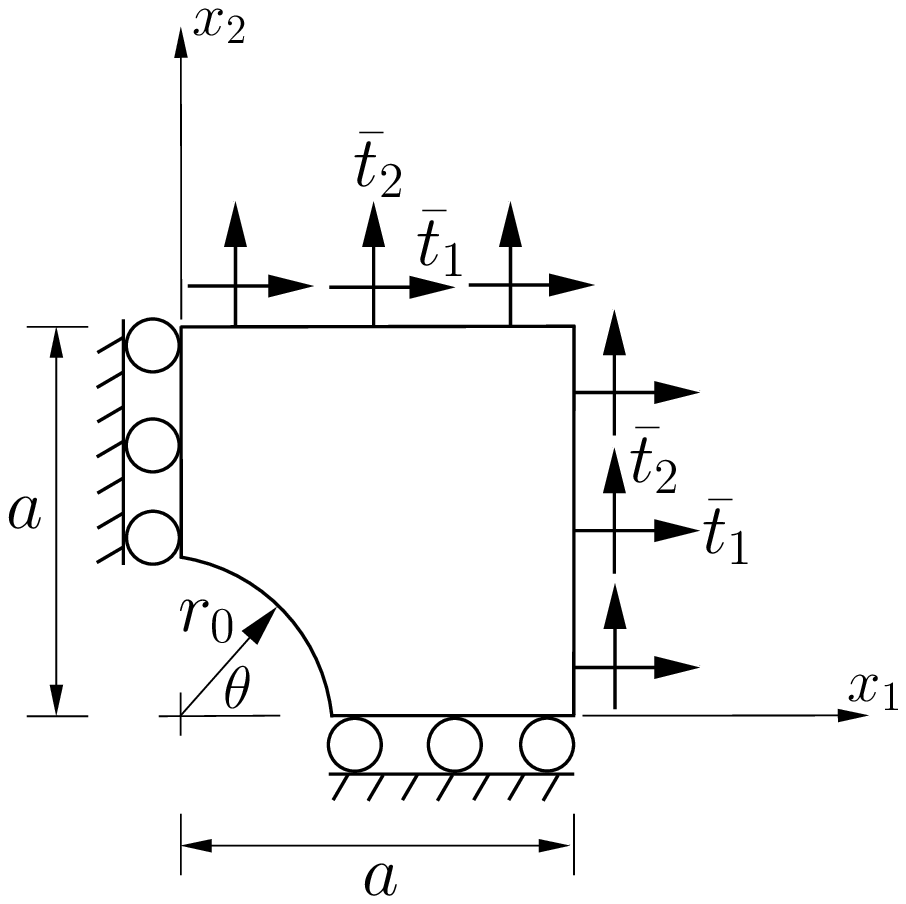, width = 0.4\textwidth}}
}
\caption{Infinite plate with a circular hole. (a) Infinite plate and (b) domain of analysis.}
\label{fig:plate_hole}
\end{figure}

\begin{figure}[!tbhp]
\centering
\mbox{
\subfigure[]{\label{fig:plate_hole_mesh_a} \epsfig{file = 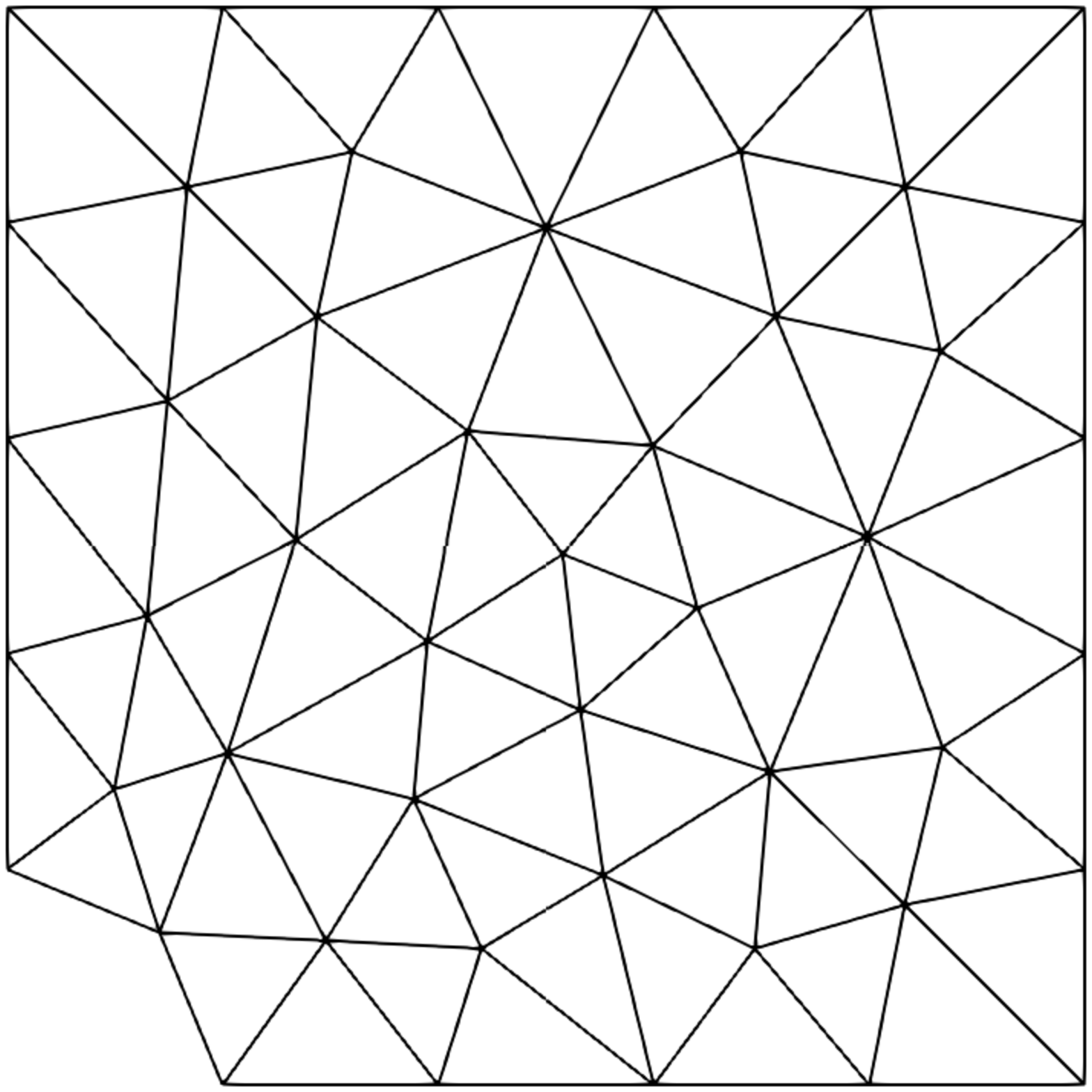, width = 0.15\textwidth}}
\subfigure[]{\label{fig:plate_hole_mesh_b} \epsfig{file = 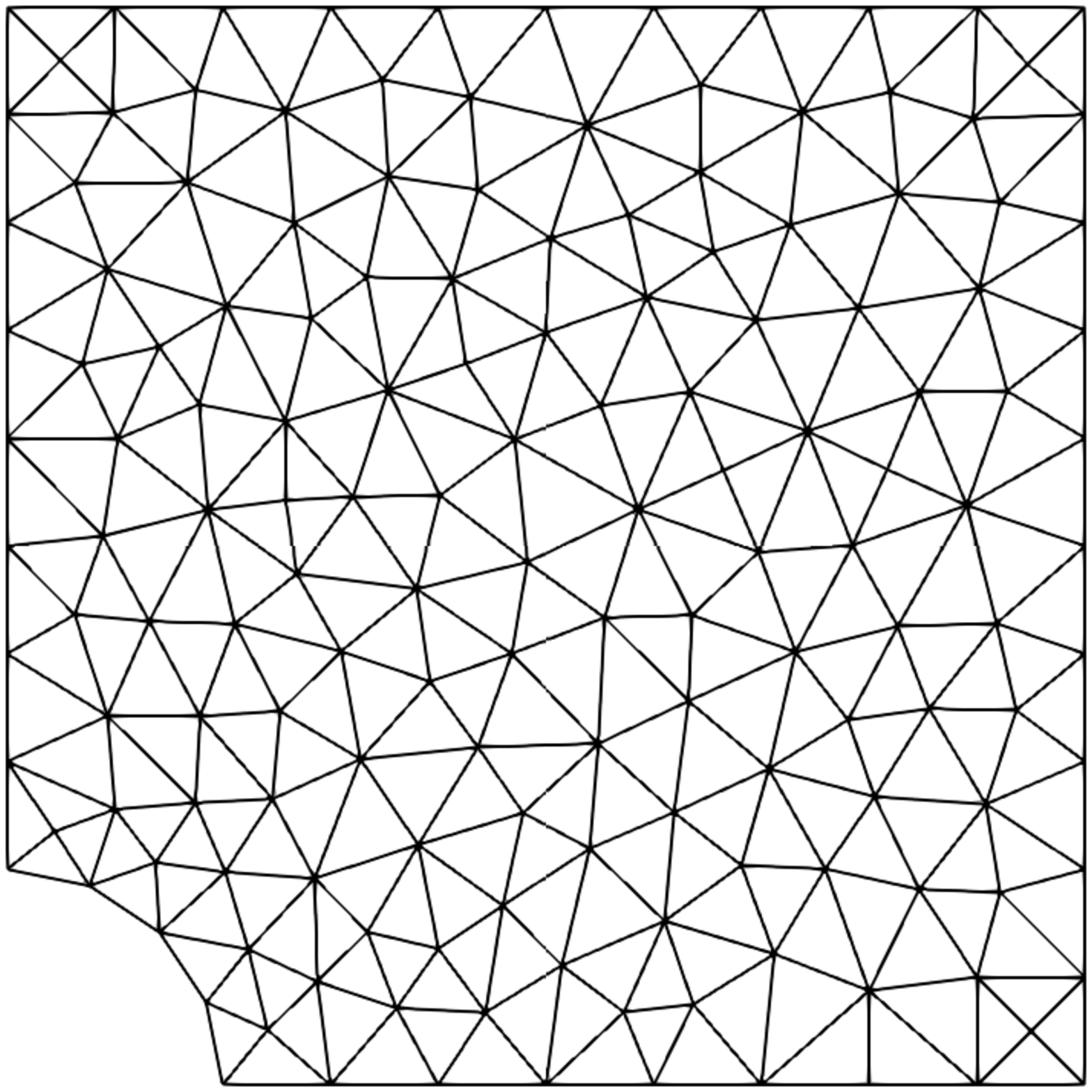, width = 0.15\textwidth}}
\subfigure[]{\label{fig:plate_hole_mesh_c} \epsfig{file = 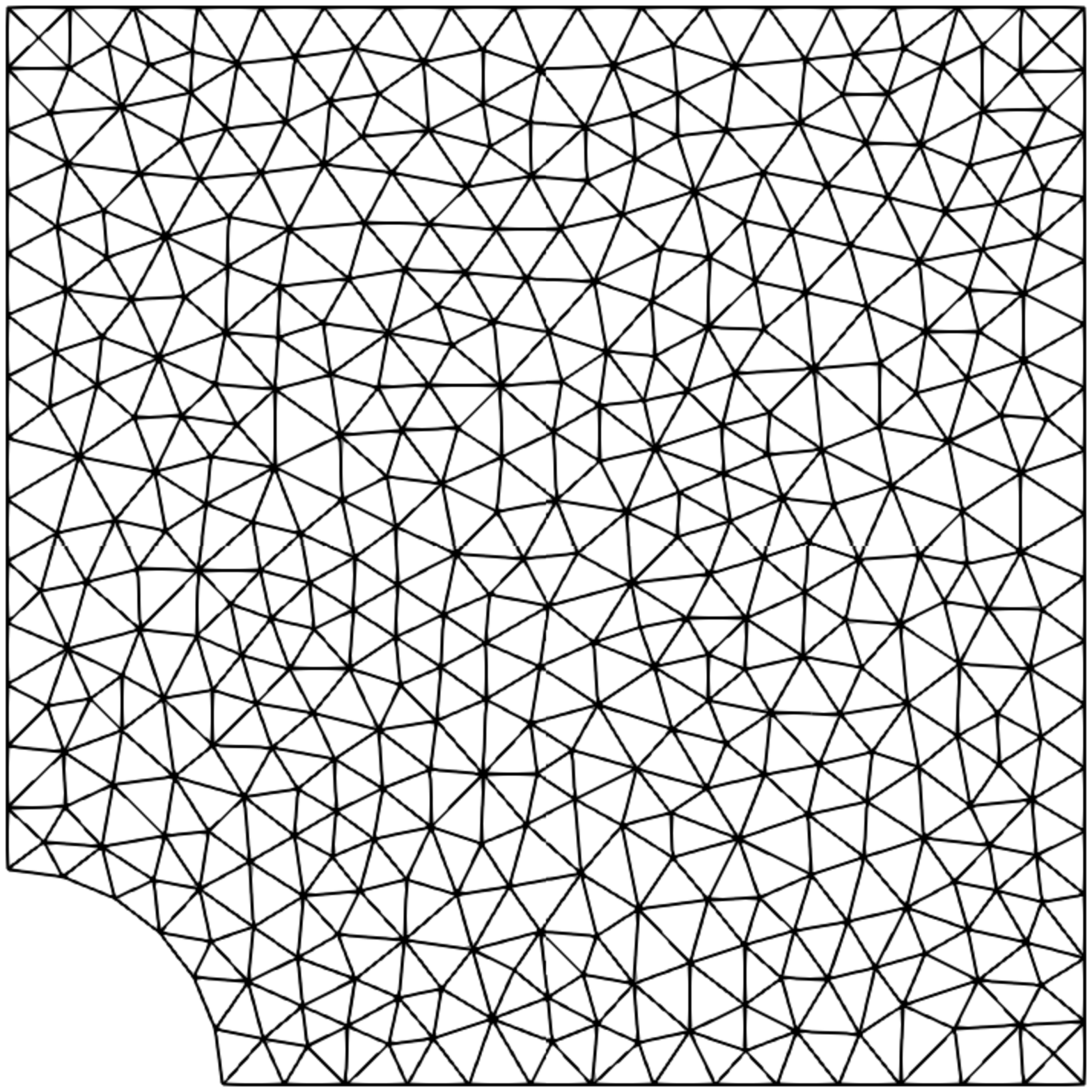, width = 0.15\textwidth}}
\subfigure[]{\label{fig:plate_hole_mesh_d} \epsfig{file = 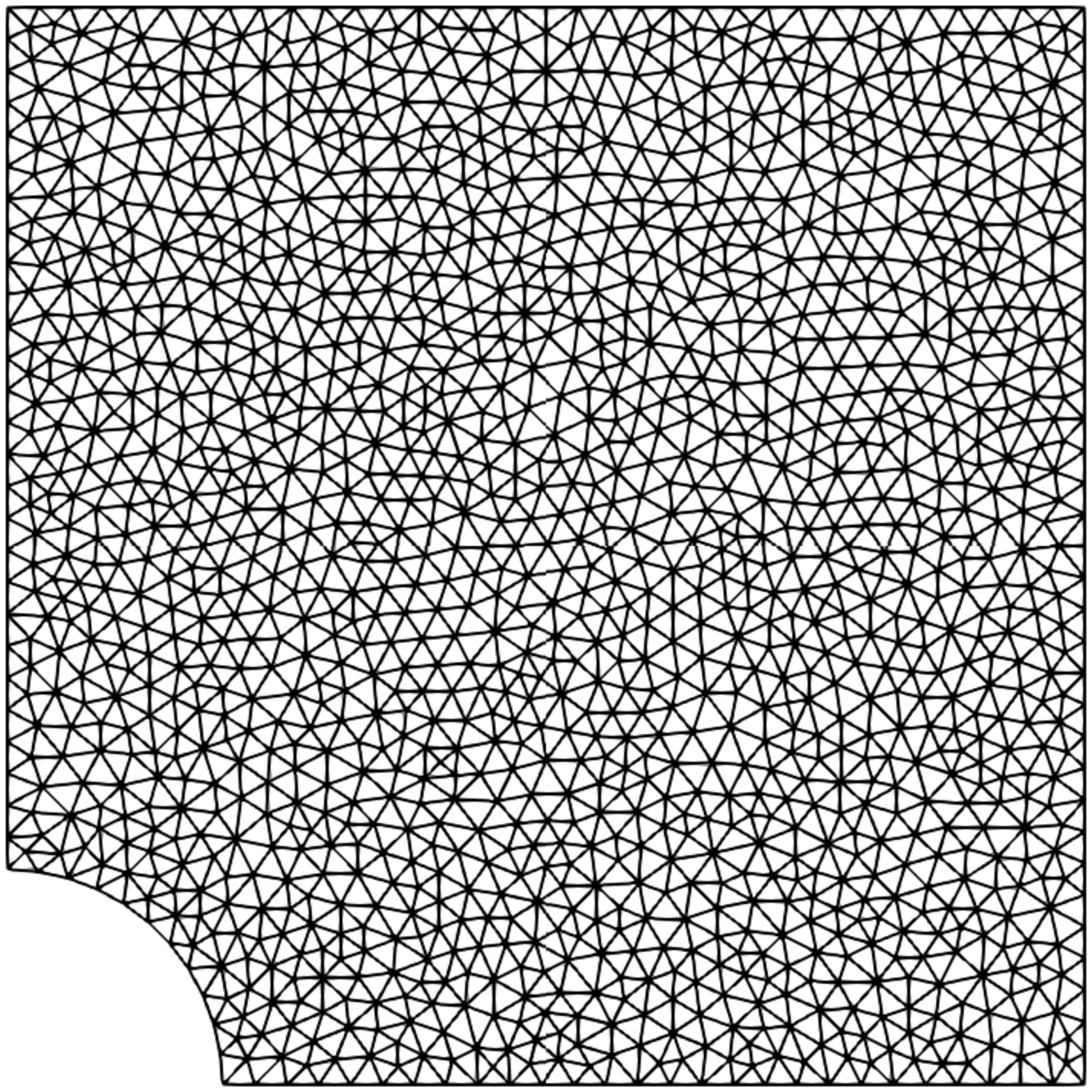, width = 0.15\textwidth}}
\subfigure[]{\label{fig:plate_hole_mesh_e} \epsfig{file = 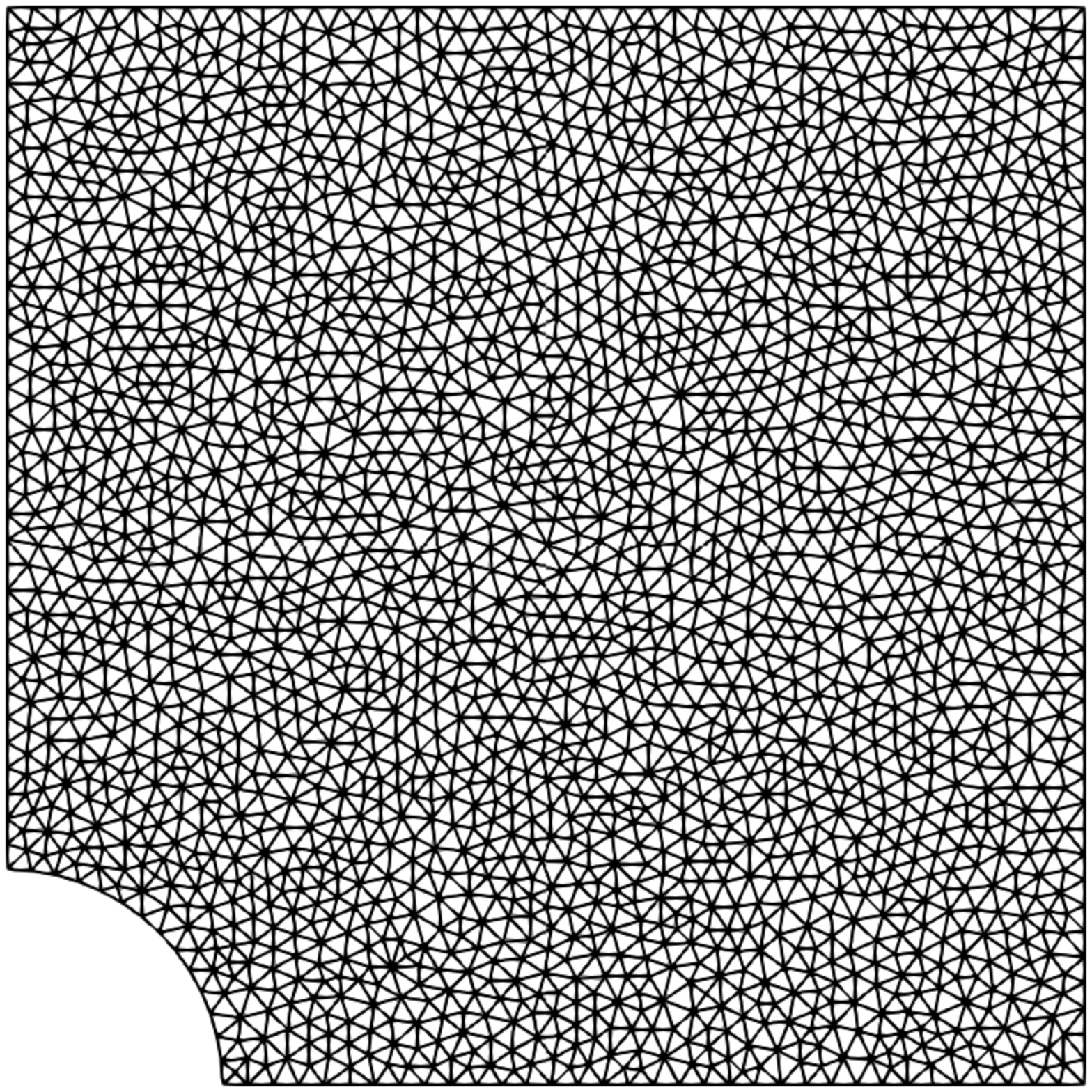, width = 0.15\textwidth}}
\subfigure[]{\label{fig:plate_hole_mesh_f} \epsfig{file = 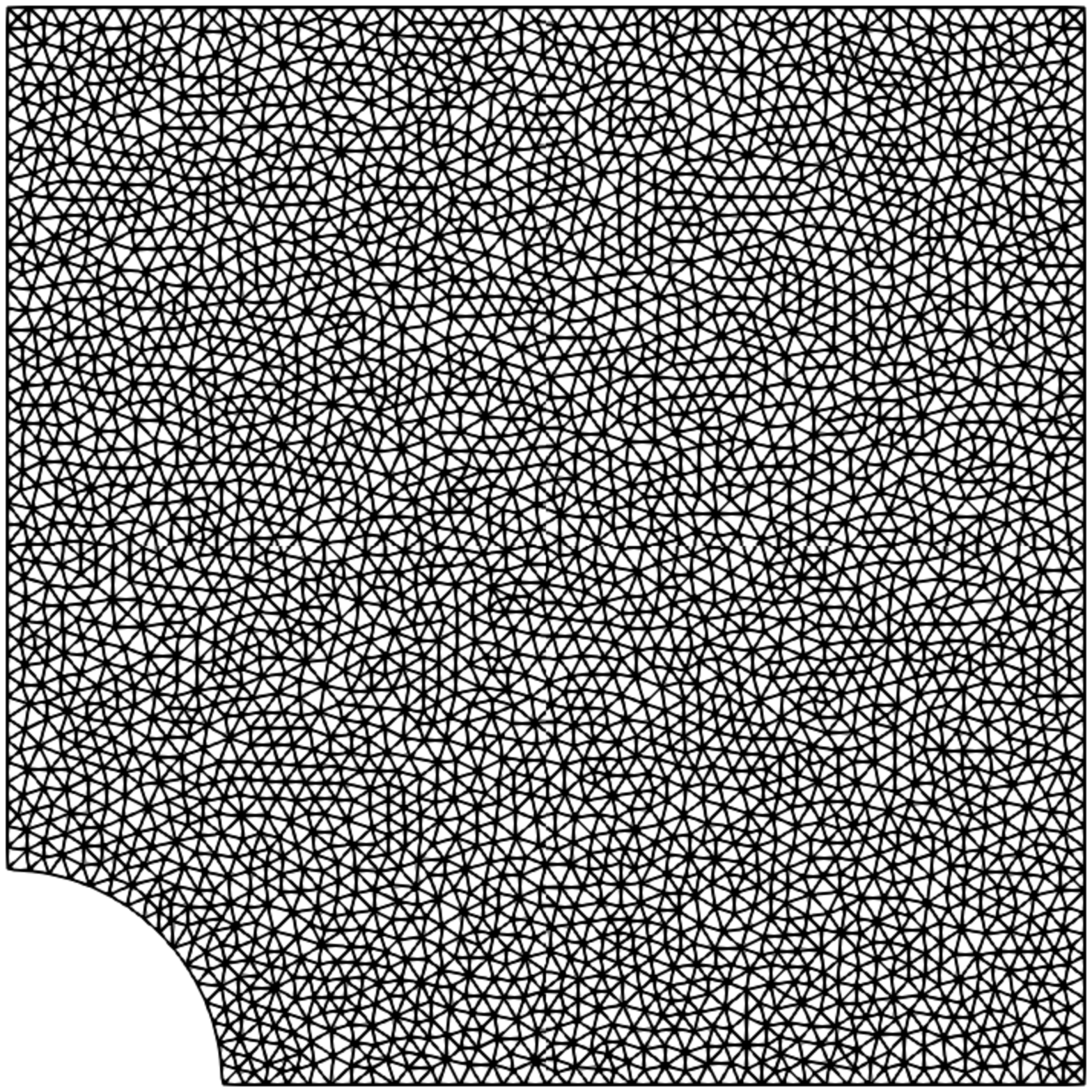, width = 0.15\textwidth}}
}
\mbox{
\subfigure[]{\label{fig:plate_hole_mesh_g} \epsfig{file = 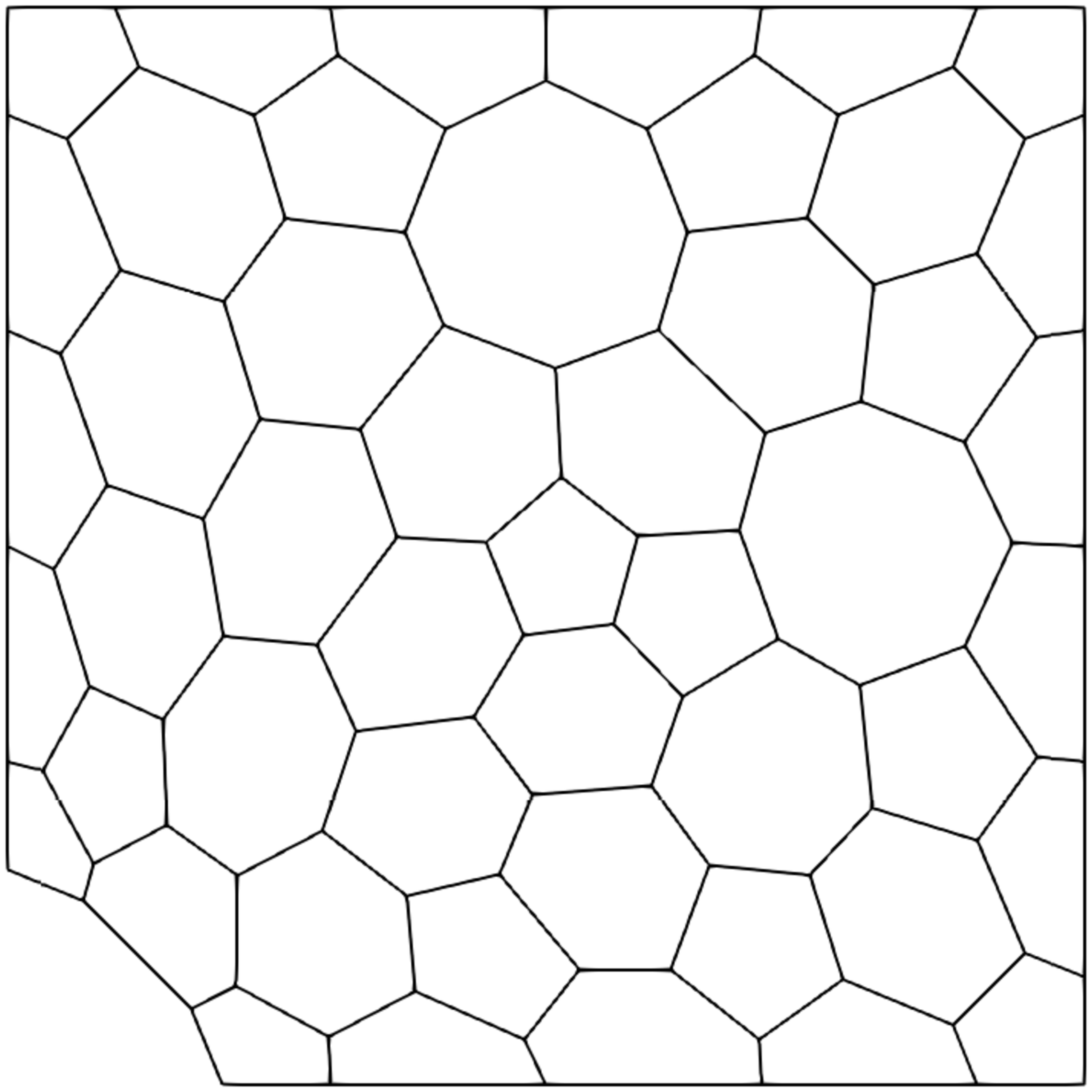, width = 0.15\textwidth}}
\subfigure[]{\label{fig:plate_hole_mesh_h} \epsfig{file = 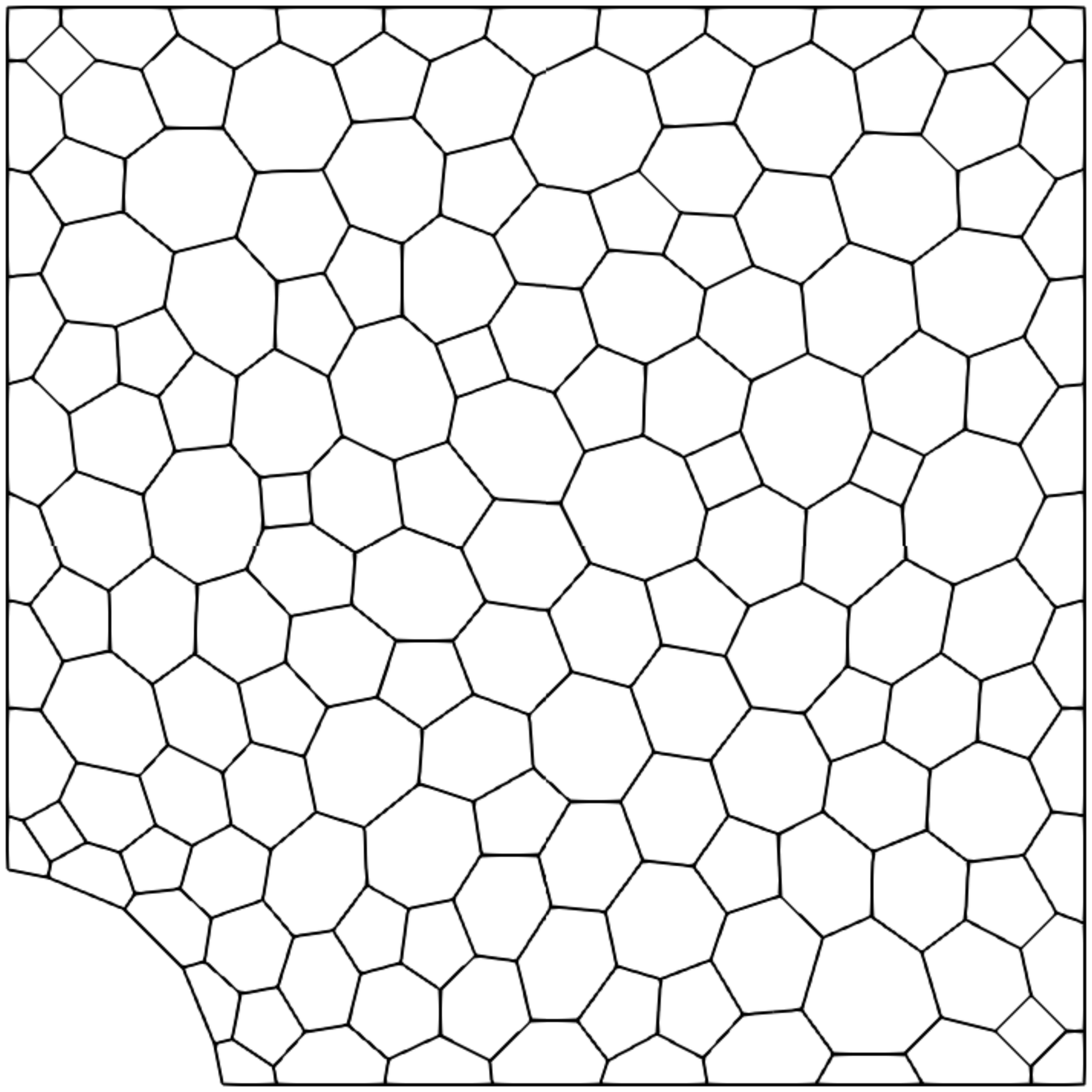, width = 0.15\textwidth}}
\subfigure[]{\label{fig:plate_hole_mesh_i} \epsfig{file = 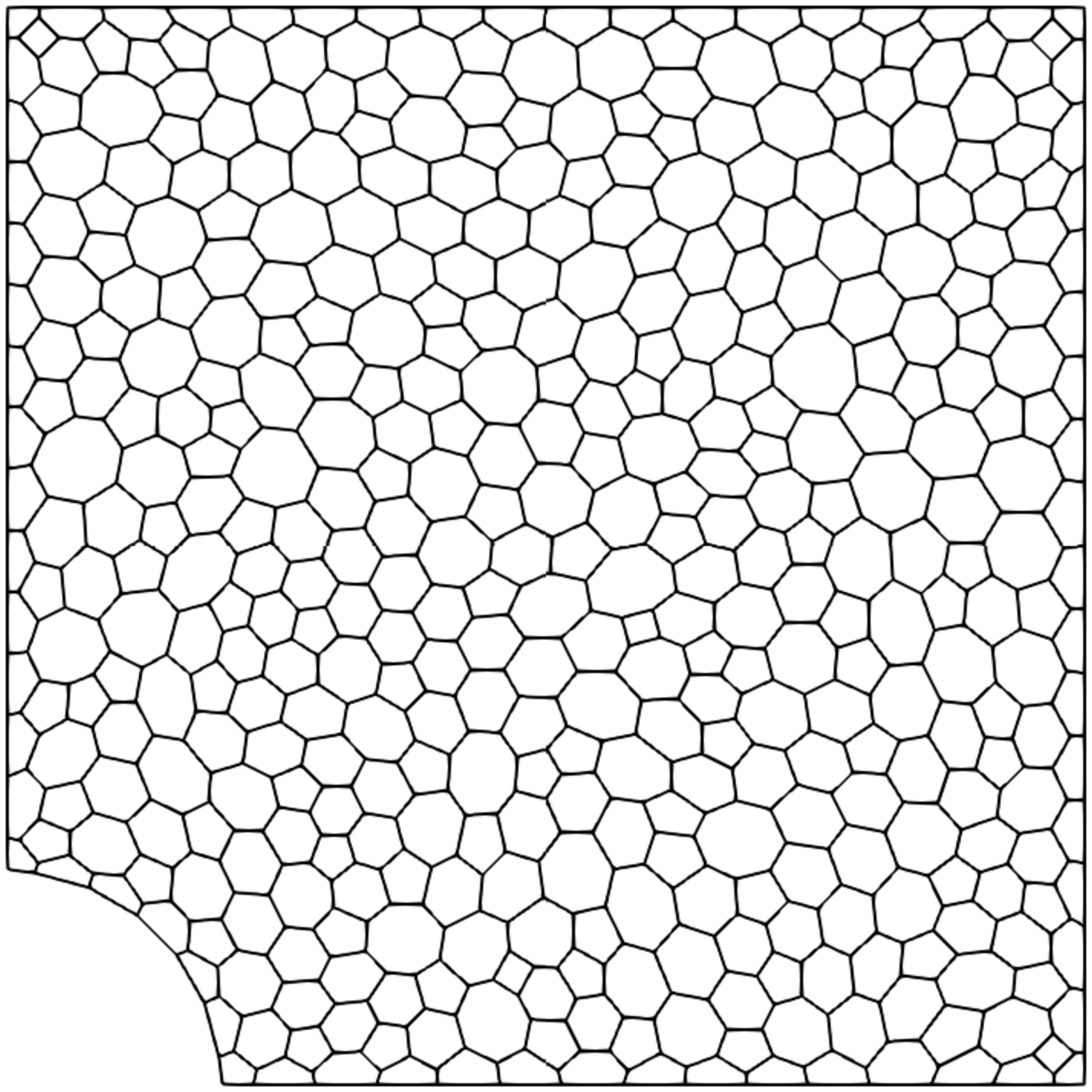, width = 0.15\textwidth}}
\subfigure[]{\label{fig:plate_hole_mesh_j} \epsfig{file = 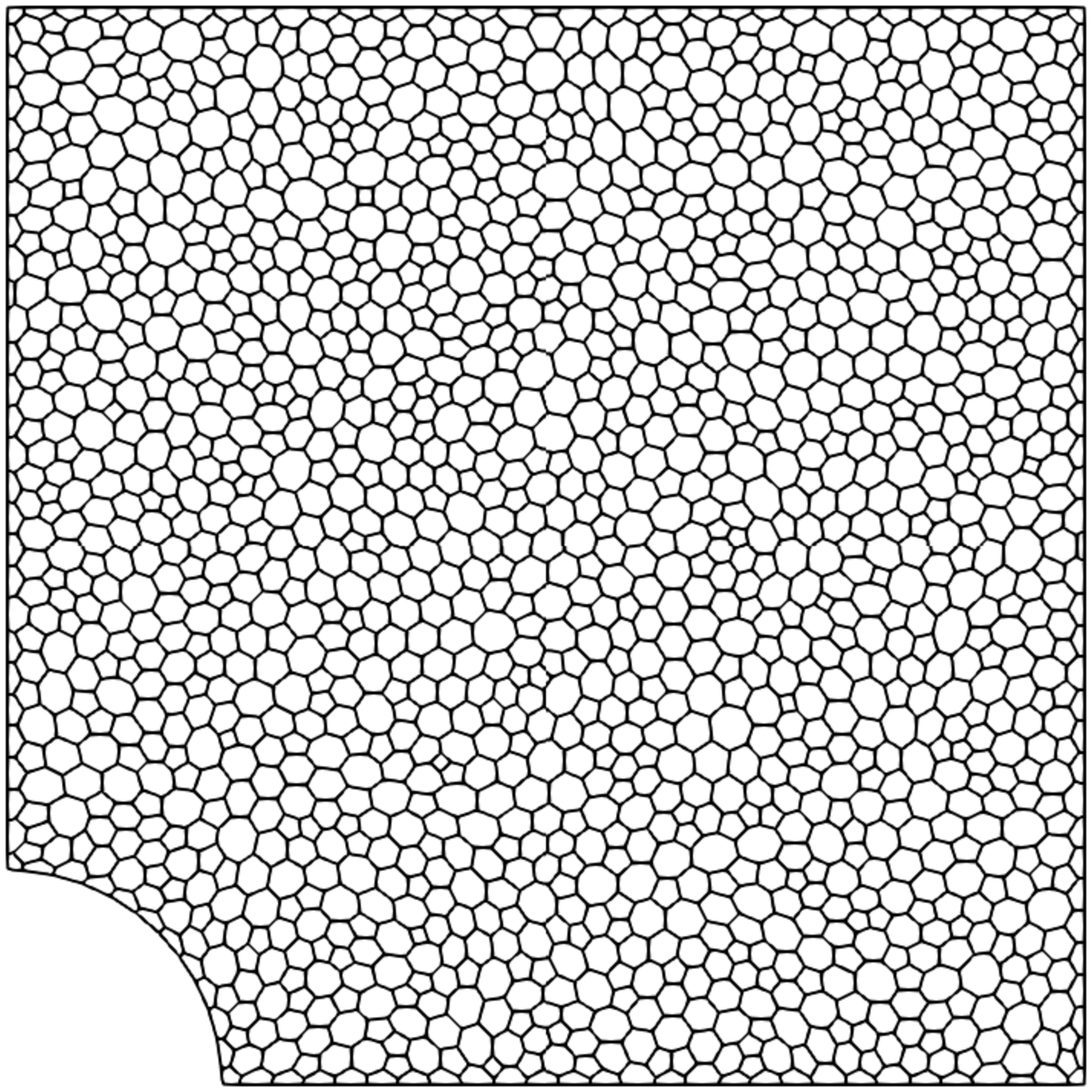, width = 0.15\textwidth}}
\subfigure[]{\label{fig:plate_hole_mesh_k} \epsfig{file = 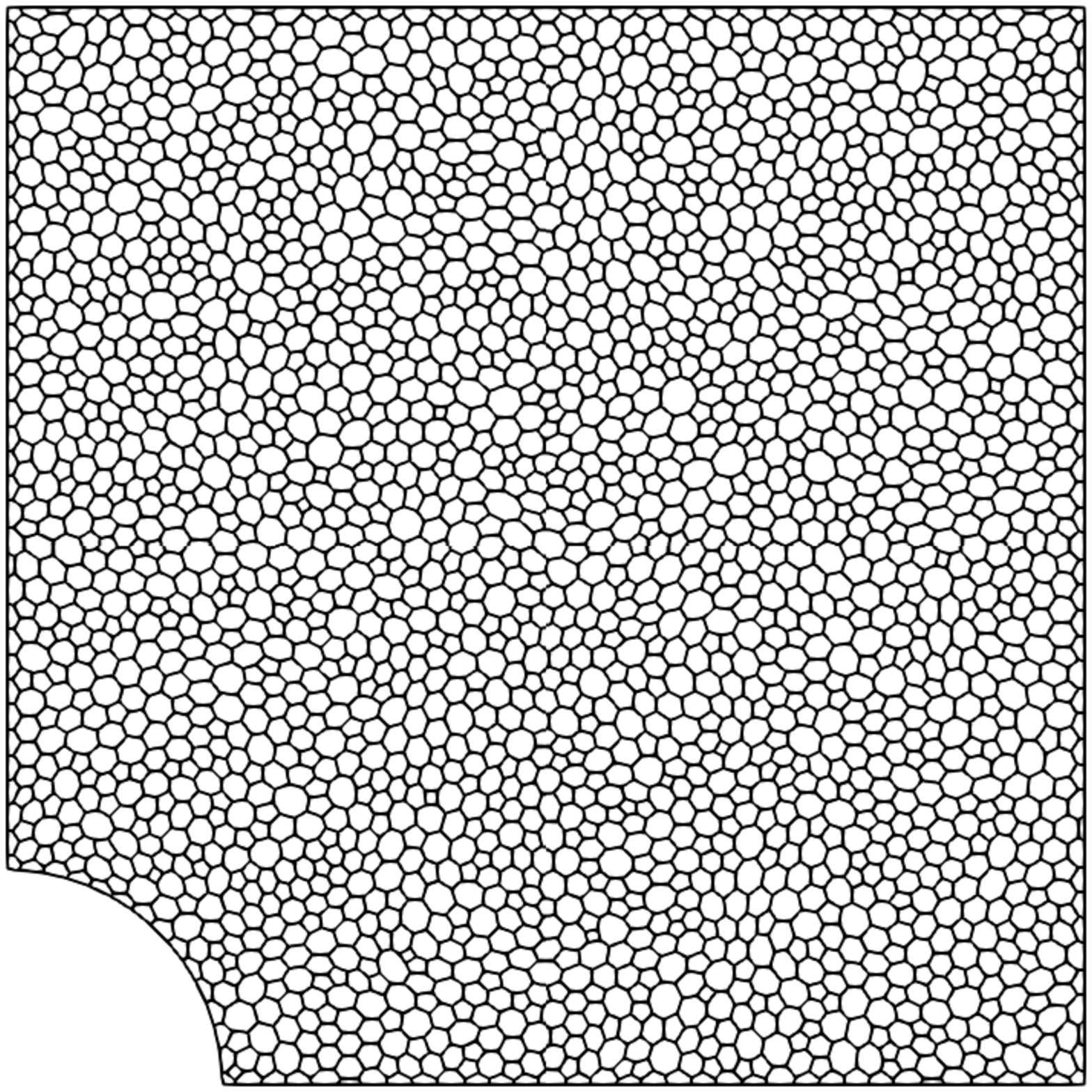, width = 0.15\textwidth}}
\subfigure[]{\label{fig:plate_hole_mesh_l} \epsfig{file = 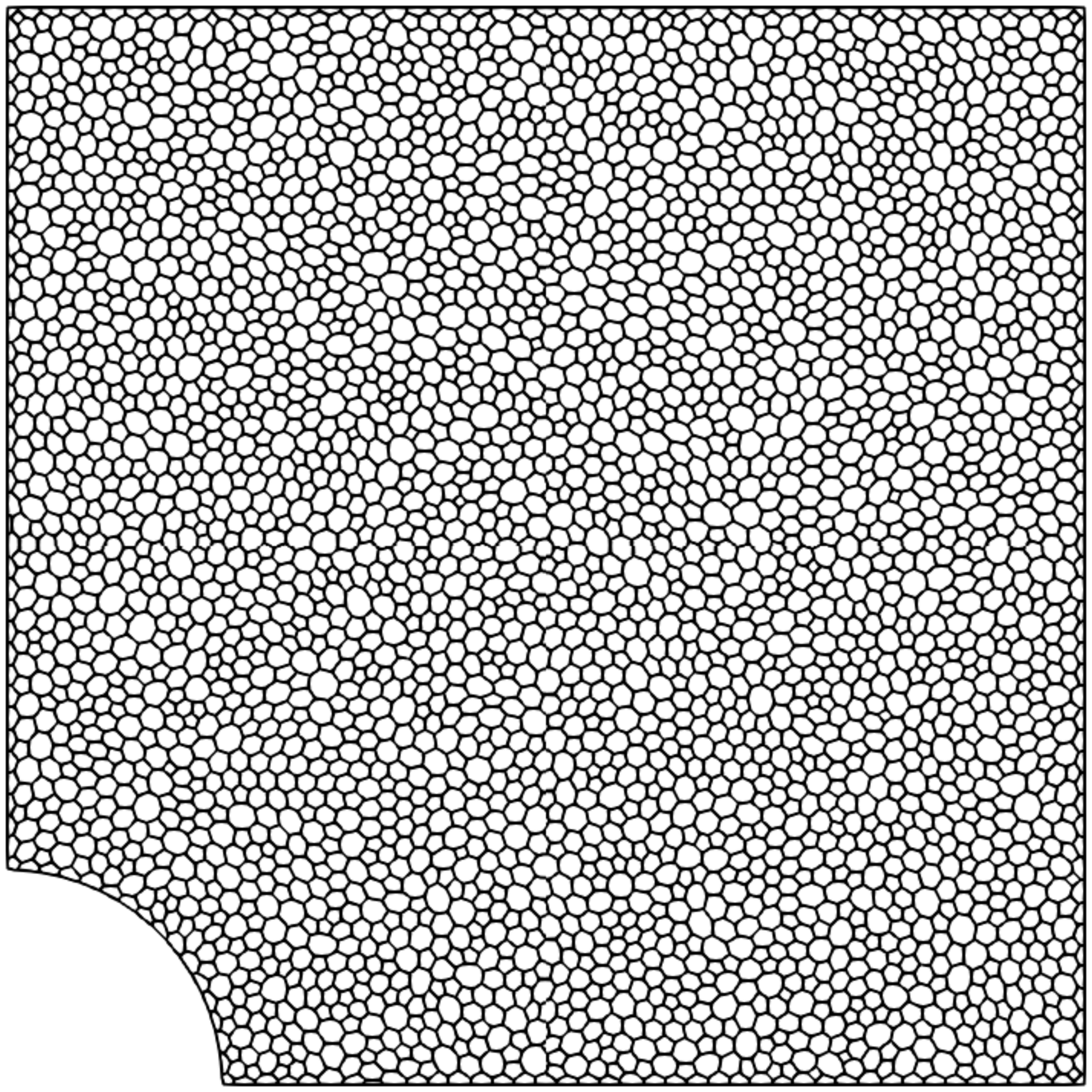, width = 0.15\textwidth}}
}
\caption{Sequence of background integration meshes used for the infinite plate with a circular hole problem.
         (a)-(f) Meshes for the MEM method and (g)-(l) meshes for the NIVED method.}
\label{fig:plate_hole_mesh}
\end{figure}

The rates of convergence are studied for the problem of an infinite plate with 
a circular hole that is loaded at infinity according to the tractions
$\sigma_{11}=T$ and $\sigma_{22}=\sigma_{12}=0$ (\fref{fig:plate_hole_a}). 
Due to the symmetry of the geometry and boundary conditions, we consider
the domain of analysis shown in \fref{fig:plate_hole_b}. Plane stress condition is
assumed with $E_\mathrm{Y}=10^3$ psi and $\nu = 0.3$.
The exact solution on the domain of analysis is given by~\cite{timoshenko:1970:TOE}
\begin{equation*}
\vm{u}=\left[
\begin{array}{c}
\frac{T}{4G} \bigg( \frac{\kappa+1}{2}r\cos{\theta}+
 \frac{r_0^2}{r}\Big( (\kappa+1)\cos{\theta}+\cos{3\theta}\Big)-\frac{r_0^4}{r^3}\cos{3\theta} \bigg)\\
 \frac{T}{4G} \bigg( \frac{\kappa-3}{2}r\sin{\theta}+
 \frac{r_0^2}{r}\Big( (\kappa-1)\sin{\theta}+\sin{3\theta}\Big)-\frac{r_0^4}{r^3}\sin{3\theta} \bigg)
\end{array}
\right],
\end{equation*}
where $G=E_\mathrm{Y}/(2(1+\nu))$ and $\kappa = (3-\nu)/(1+\nu)$. The exact stress field is:
\begin{equation*}
\left[\begin{array}{c}
\sigma_{11}\\
\sigma_{22}\\
\sigma_{12}
\end{array}
\right]
=\left[\begin{array}{c}
T\bigg(1-\frac{r_0^2}{r^2}\left( \frac{3}{2}\cos{2\theta}
+\cos{4\theta}\right)+\frac{3r_0^4}{2r^4}\cos{4\theta} \bigg)\\
-T\bigg(\frac{r_0^2}{r^2}\left( \frac{1}{2}\cos{2\theta}-\cos{4\theta}\right)
+\frac{3r_0^4}{2r^4}\cos{4\theta} \bigg)\\
-T\bigg(\frac{r_0^2}{r^2}\left( \frac{1}{2}\sin{2\theta}
+\sin{4\theta}\right)-\frac{3r_0^4}{2r^4}\sin{4\theta}\bigg)
\end{array}
\right],
\end{equation*}
where $r$ is the radial distance from the center $(x_1=0,x_2=0)$ to a
point $(x_1,x_2)$ in the domain of analysis. In the computations, the
following data are used: $T=100$ psi, $r_0=1$ inch and $a=5$ inch. The Dirichlet 
boundary conditions on the domain of analysis are imposed as follows: $\bar{u}_1 = 0$ 
on the left side and $\bar{u}_2 = 0$ on the bottom side. 
The Neumann boundary conditions are prescribed using the exact 
stresses, as follows:
$\bar{\vm{t}}=\smat{\sigma_{12} & \sigma_{22}}^\transpose$ on the top side and
$\bar{\vm{t}}=\smat{\sigma_{11} & \sigma_{12}}^\transpose$ on the right side.
The sequence of background integration meshes used in 
the study are shown in~\fref{fig:plate_hole_mesh}. The Gaussian prior is
used for the evaluation of the maximum-entropy basis functions. 

\begin{figure}[!tbhp]
\centering
\subfigure[]{\label{fig:ex_plate_hole_norms_a} \epsfig{file = 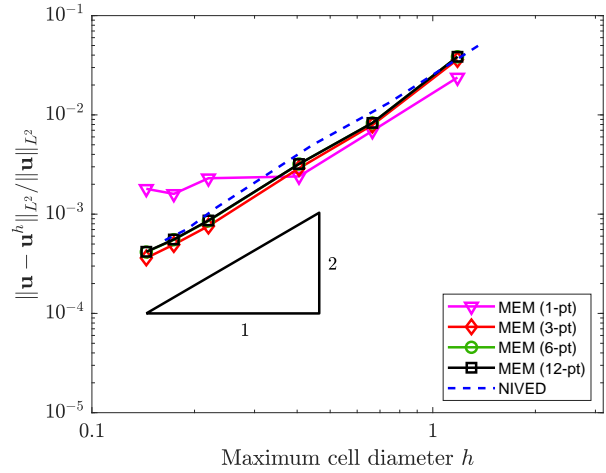, width = 0.58\textwidth}}
\subfigure[]{\label{fig:ex_plate_hole_norms_b} \epsfig{file = 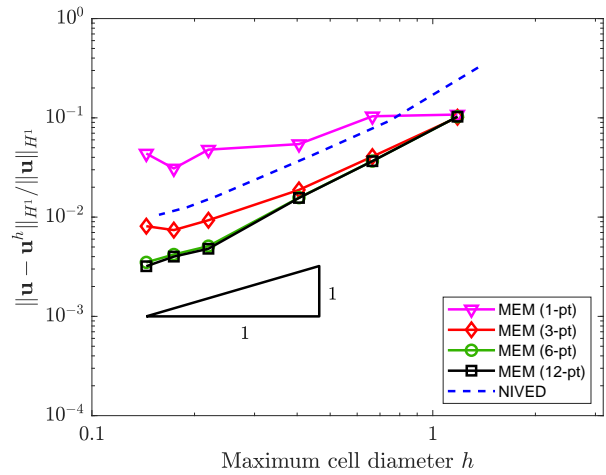, width = 0.58\textwidth}}
\caption{Rates of convergence for the infinite plate with a circular hole
         problem. For the MEM approach, an erratic convergence in the $L^2$ norm is 
         obtained using a 1-point Gauss rule and the optimal rate of 2 is recovered 
         using a 3-point Gauss rule. Its convergence in the $H^1$ seminorm is optimal 
         when using at least a 6-point Gauss rule. The optimal rates of 2 and 1 are 
         obtained in the $L^2$ norm and the $H^1$ seminorm,
         respectively, for the NIVED method.}
\label{fig:ex_plate_hole_norms}
\end{figure}

The convergence rates that are delivered by the
MEM and NIVED approaches in the $L^2$ norm and the $H^1$ seminorm of
the error are compared in \fref{fig:ex_plate_hole_norms}.
We observe that the convergence of the MEM using a 1-point Gauss rule 
behaves erratically in the $L^2$ norm and that the optimal rate of 2 
is recovered using a 3-point Gauss rule. 
Regarding its convergence in the $H^1$ seminorm, we observe that 
it is optimal when using at least a 6-point Gauss rule. 
On the other hand, the NIVED approach delivers the optimal rates 
of 2 and 1 in the $L^2$ norm and the $H^1$ seminorm, respectively.

\subsection{$L$-shaped domain under traction load}\label{sec:numexamples_lshaped}

\begin{figure}[!tbhp]
\centering
\epsfig{file = 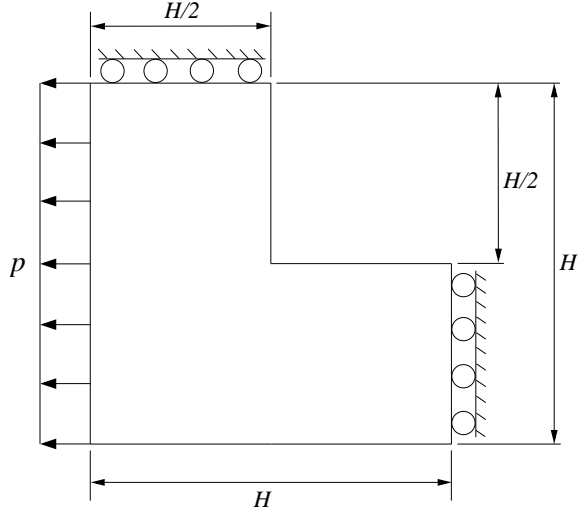, width = 0.5\textwidth}
\caption{$L$-shaped domain under traction load.}
\label{fig:lshapeddomain}
\end{figure}

\begin{figure}[!tbhp]
\centering
\mbox{
\subfigure[]{\label{fig:lshaped_mesh_a} \epsfig{file = 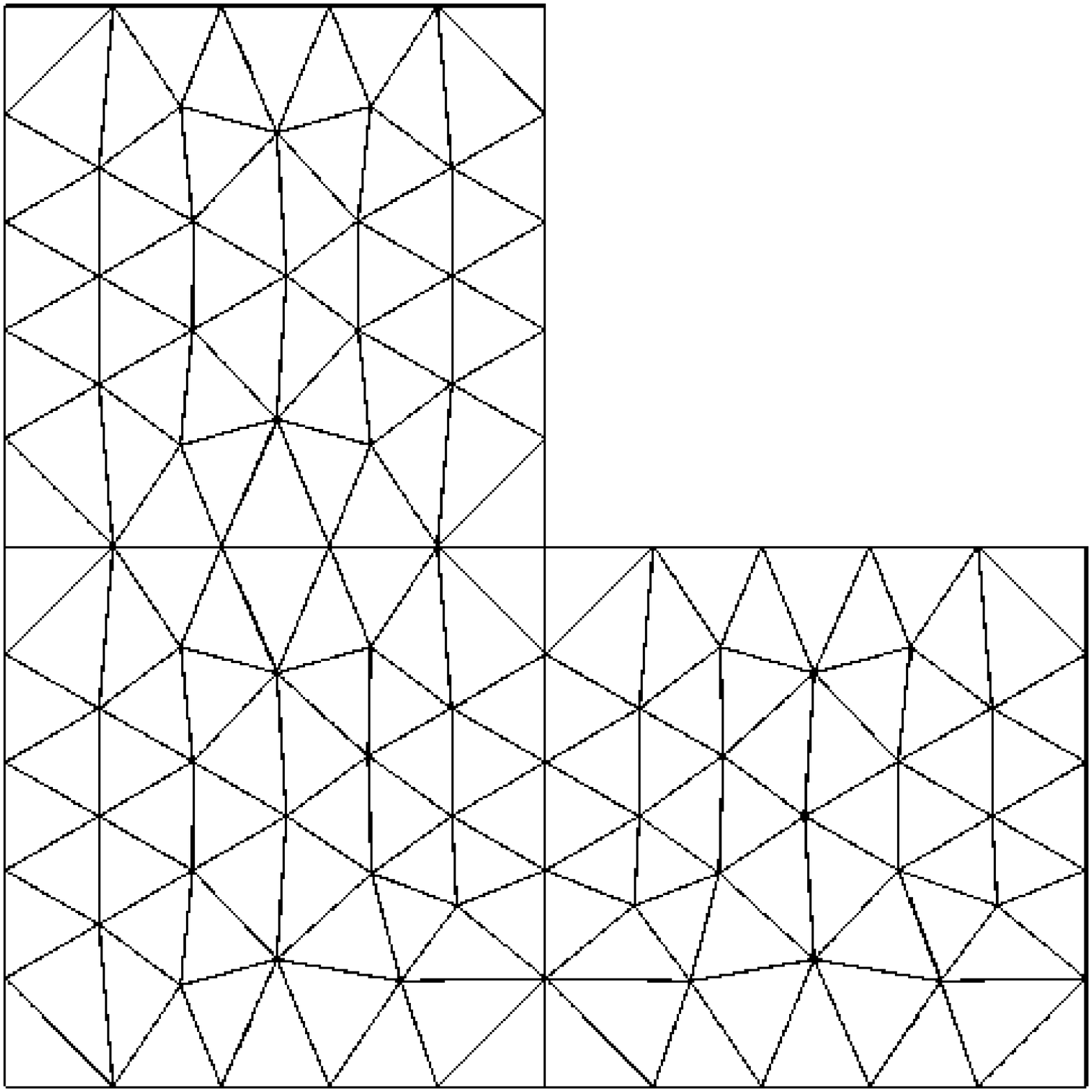, width = 0.185\textwidth}}
\subfigure[]{\label{fig:lshaped_mesh_b} \epsfig{file = 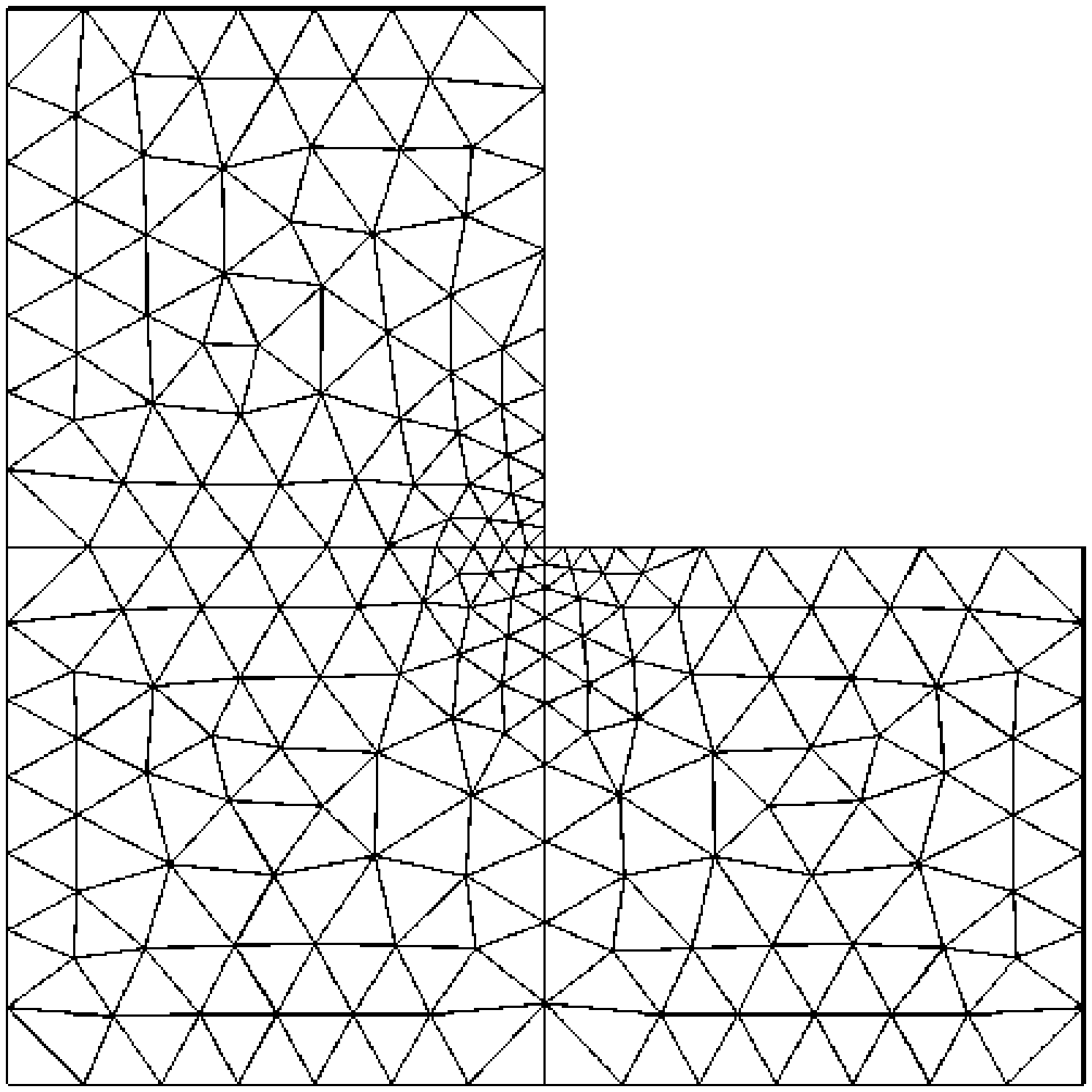, width = 0.185\textwidth}}
\subfigure[]{\label{fig:lshaped_mesh_c} \epsfig{file = 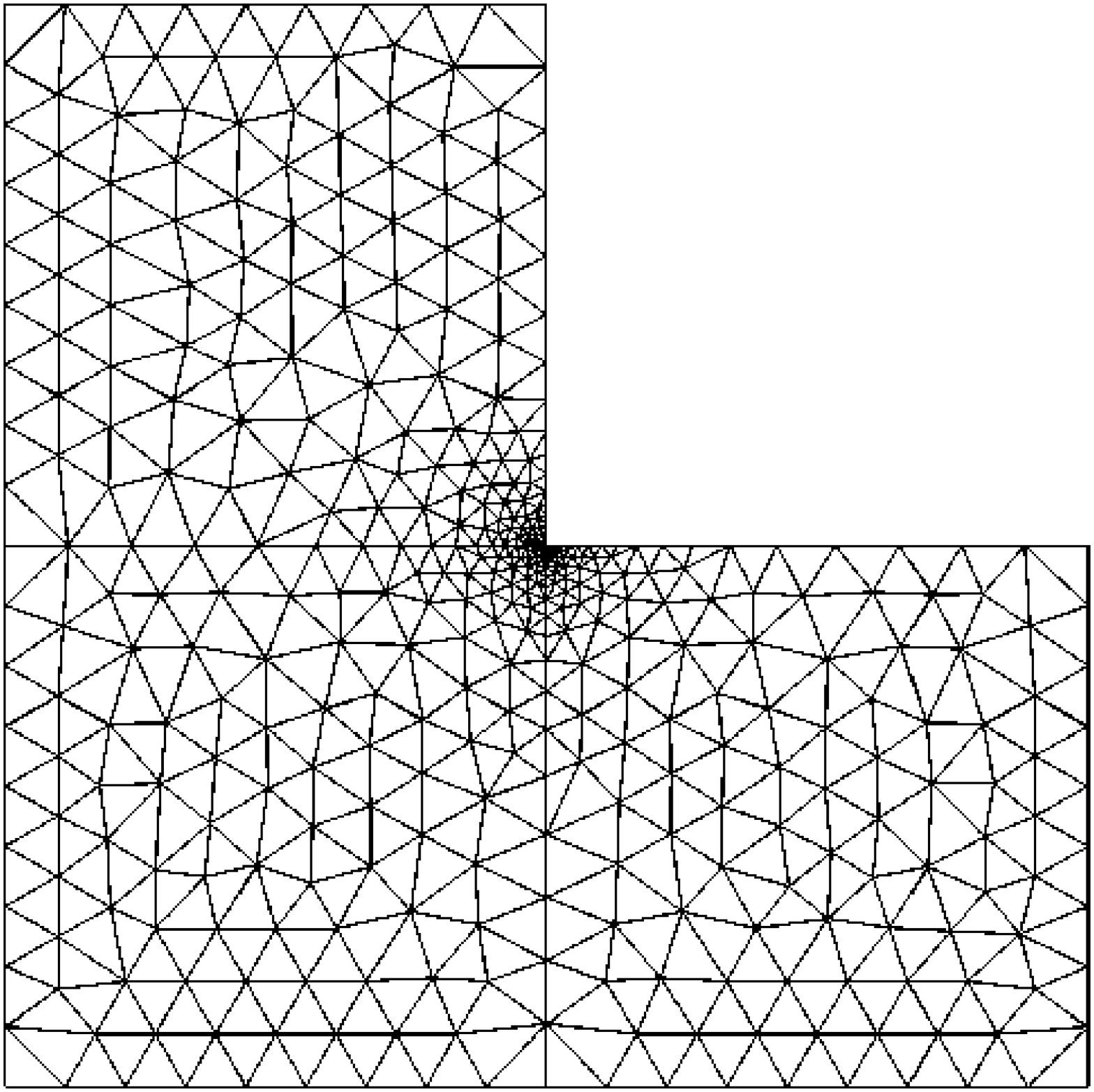, width = 0.185\textwidth}}
\subfigure[]{\label{fig:lshaped_mesh_d} \epsfig{file = 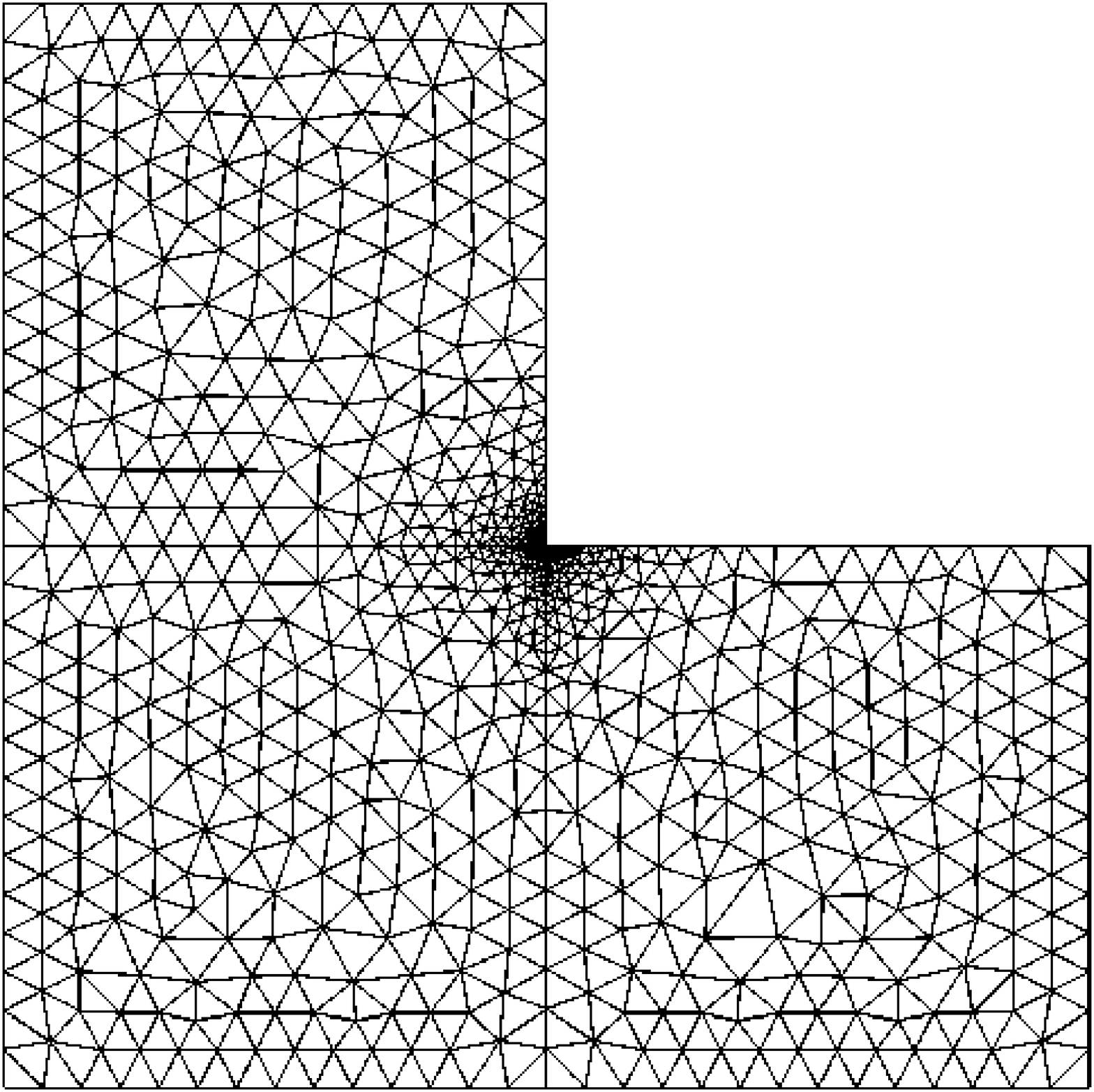, width = 0.185\textwidth}}
\subfigure[]{\label{fig:lshaped_mesh_e} \epsfig{file = 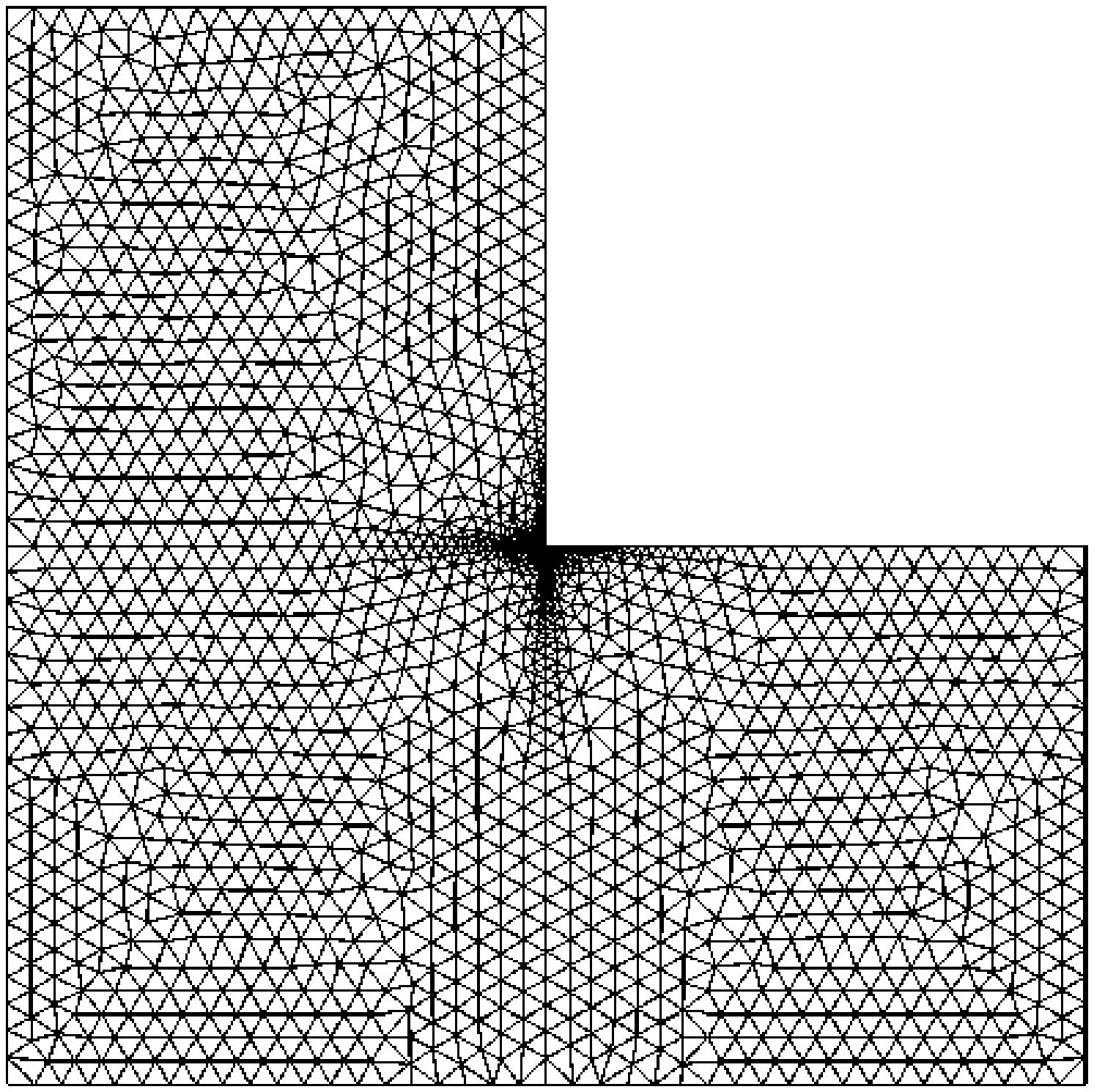, width = 0.185\textwidth}}
}
\mbox{
\subfigure[]{\label{fig:lshaped_mesh_f} \epsfig{file = 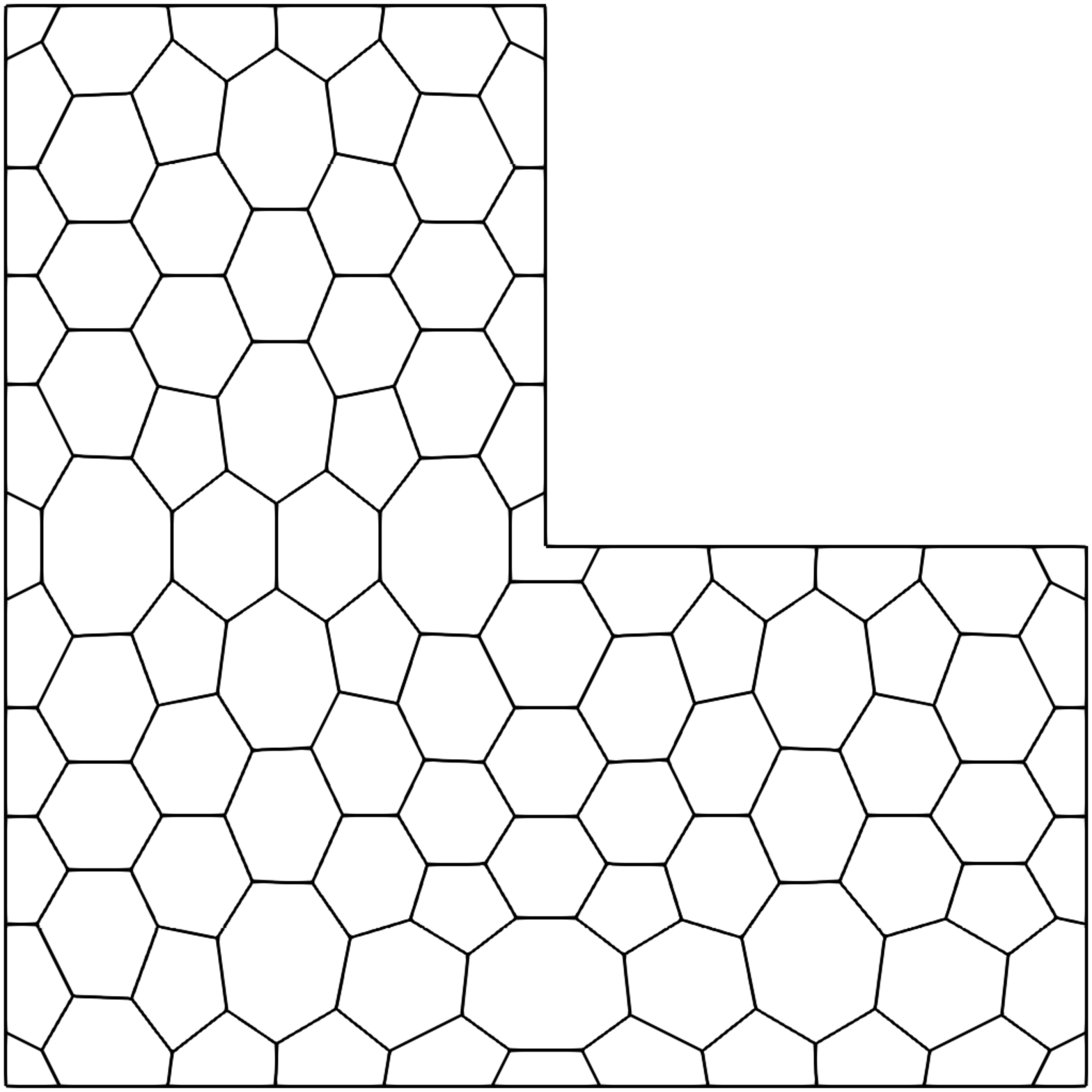, width = 0.185\textwidth}}
\subfigure[]{\label{fig:lshaped_mesh_g} \epsfig{file = 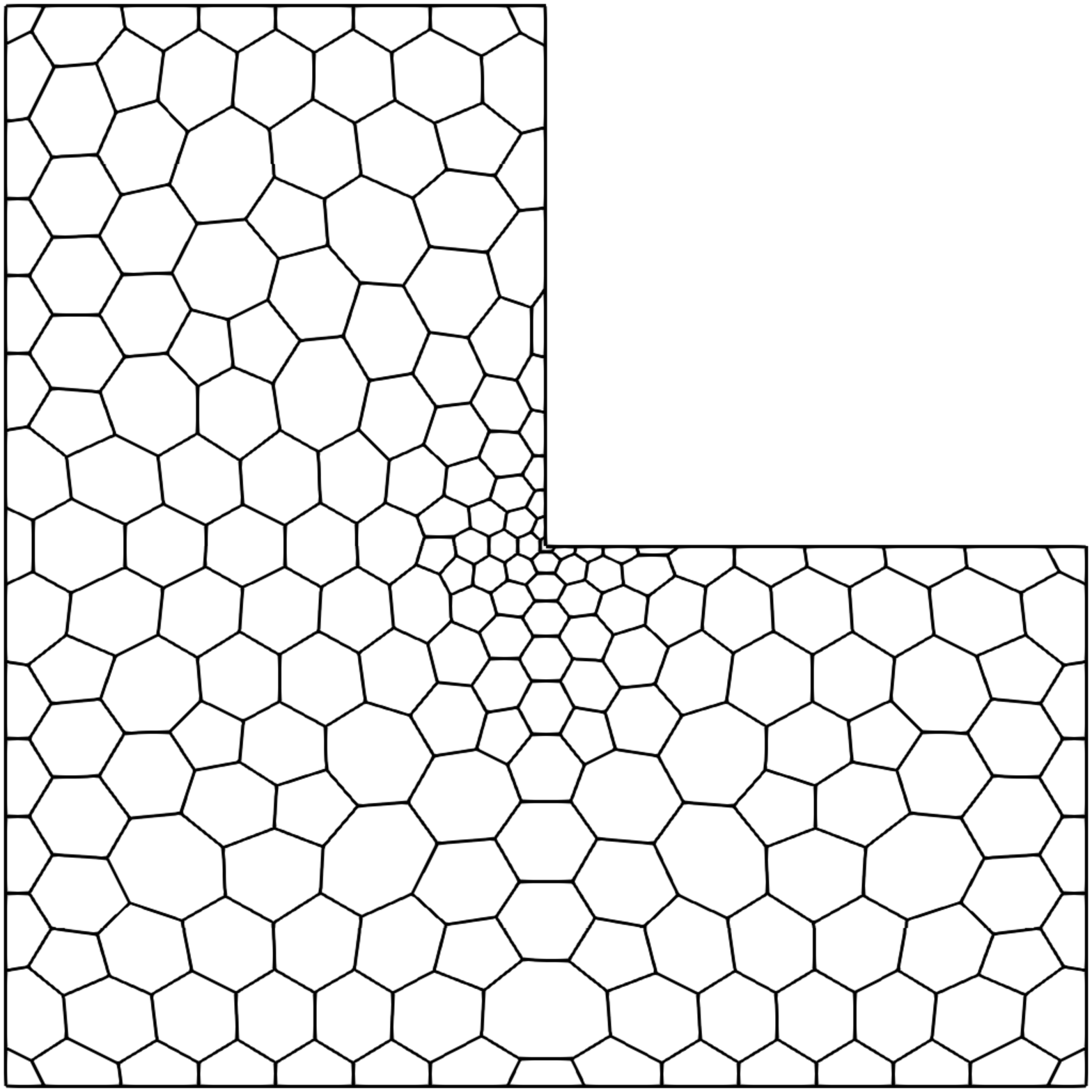, width = 0.185\textwidth}}
\subfigure[]{\label{fig:lshaped_mesh_h} \epsfig{file = 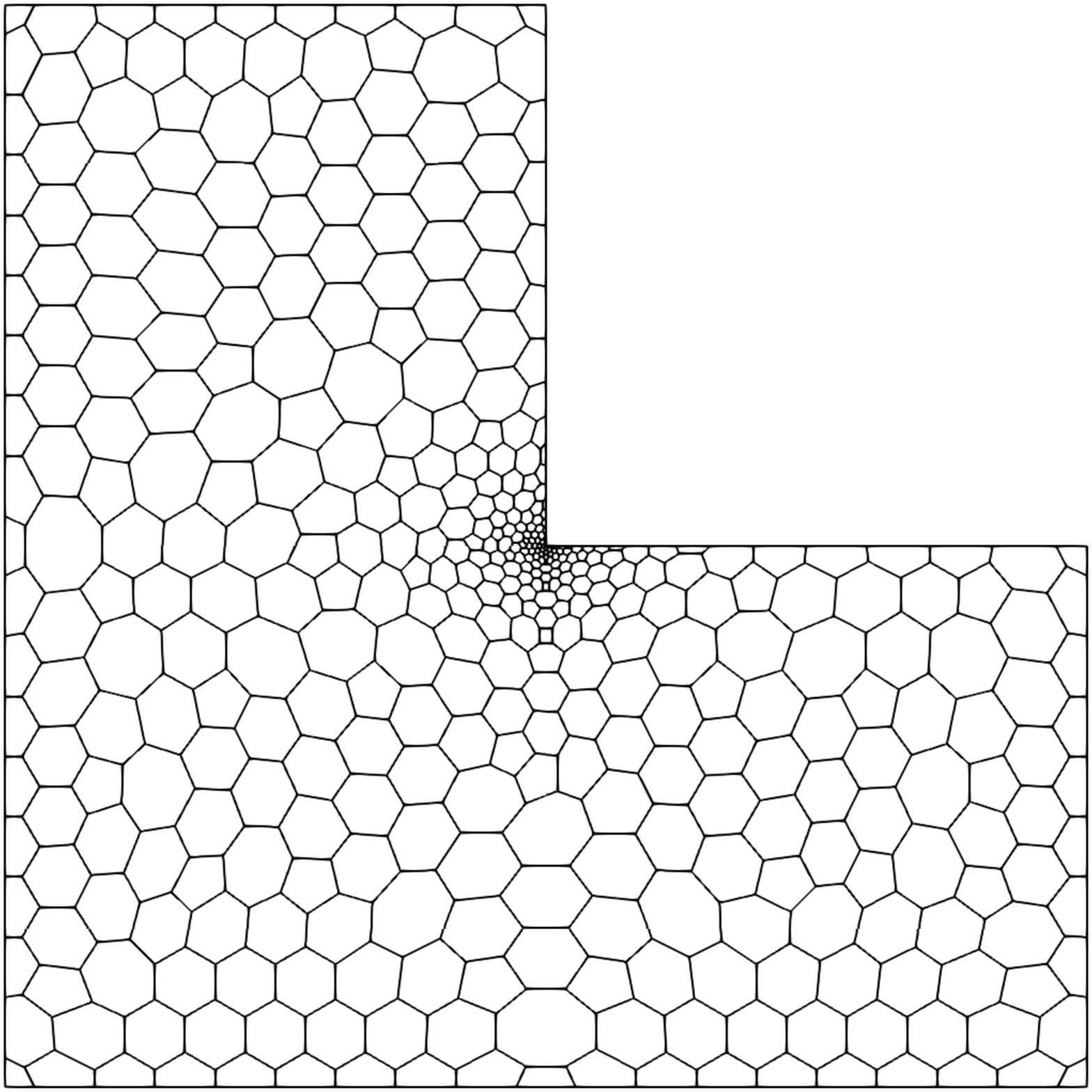, width = 0.185\textwidth}}
\subfigure[]{\label{fig:lshaped_mesh_i} \epsfig{file = 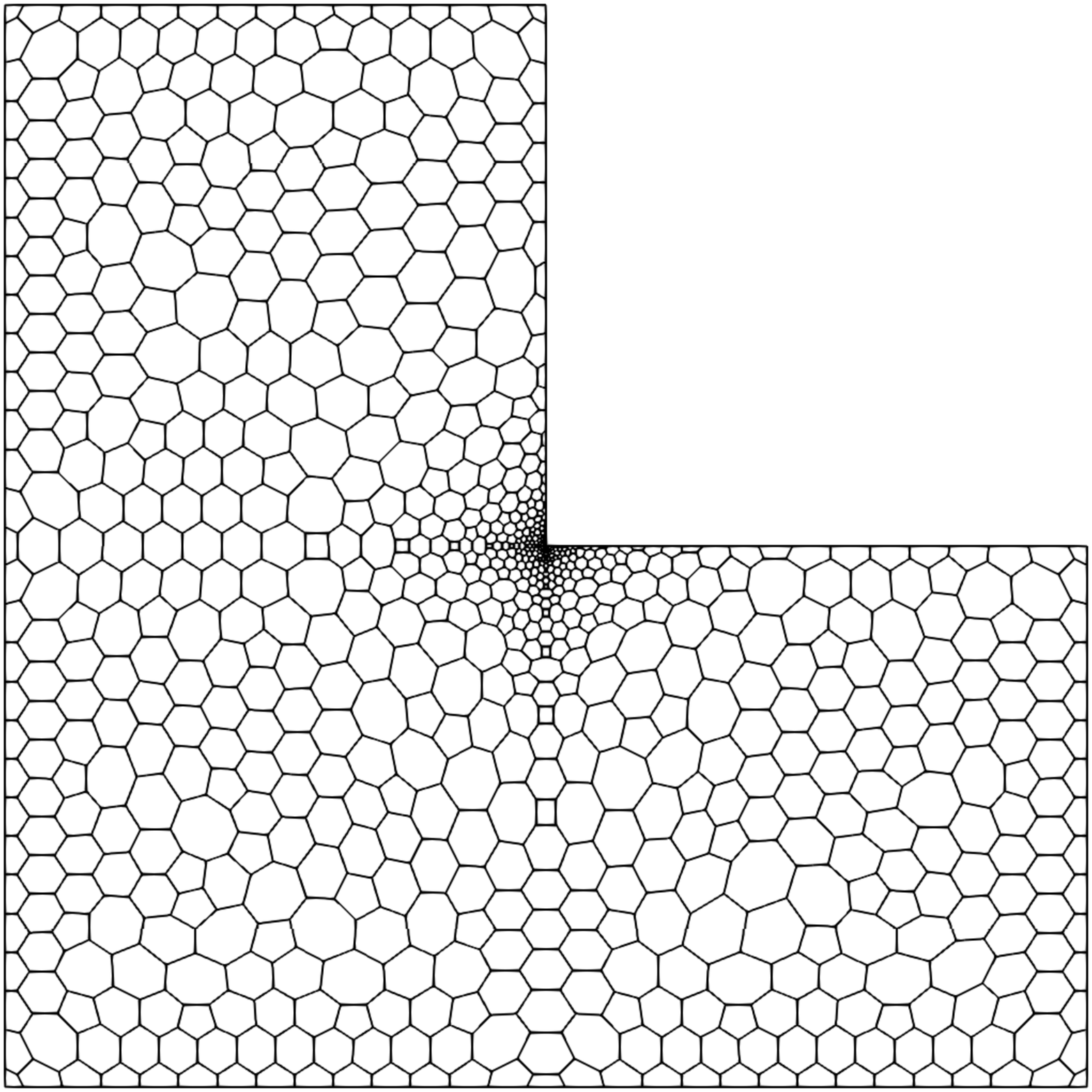, width = 0.185\textwidth}}
\subfigure[]{\label{fig:lshaped_mesh_j} \epsfig{file = 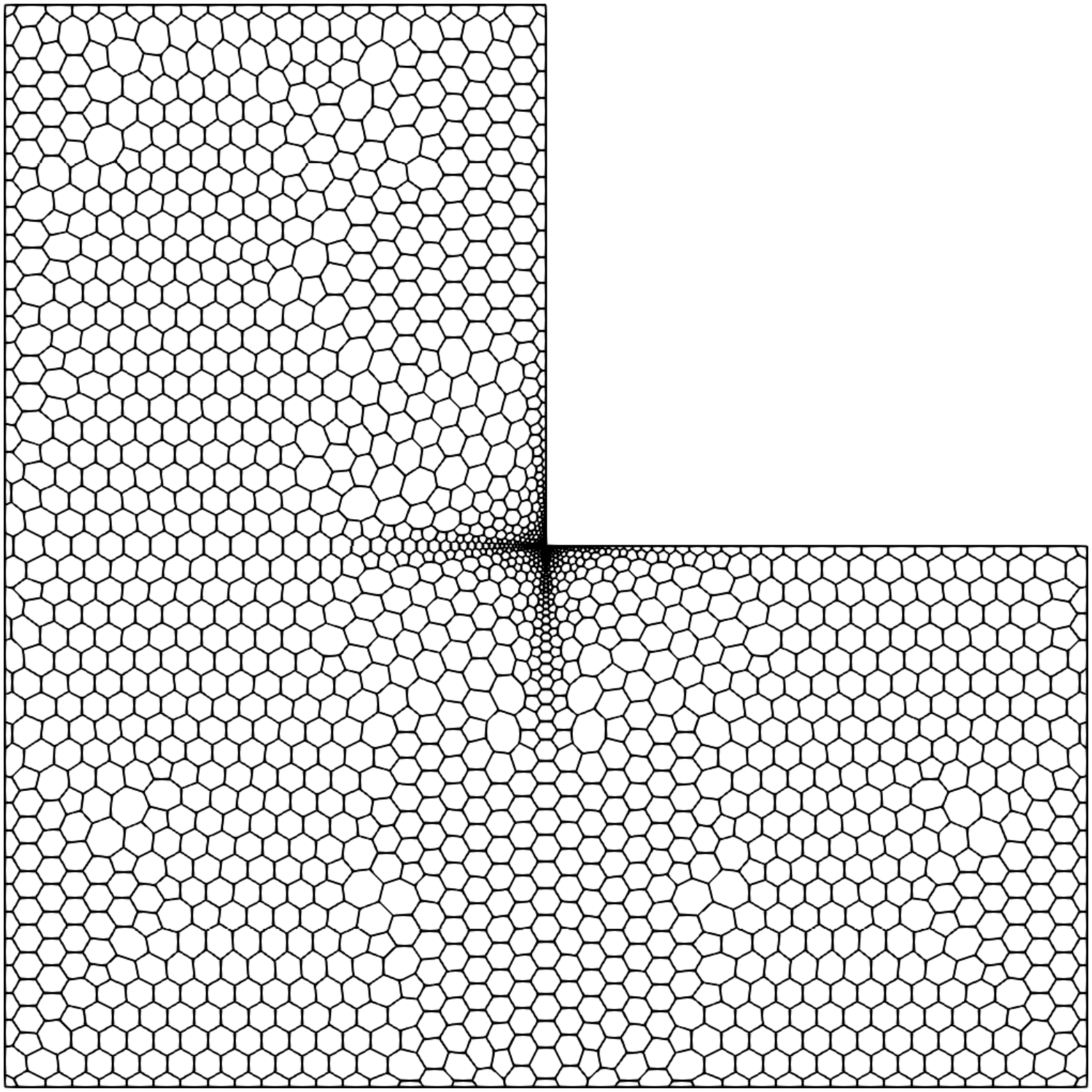, width = 0.185\textwidth}}
}
\caption{Sequence of background integration meshes used for the $L$-shaped domain problem.
         (a)-(e) Meshes for the MEM method and (f)-(j) meshes for the NIVED method.}
\label{fig:lshaped_mesh}
\end{figure}

\begin{figure}[!tbhp]
\centering
\epsfig{file = 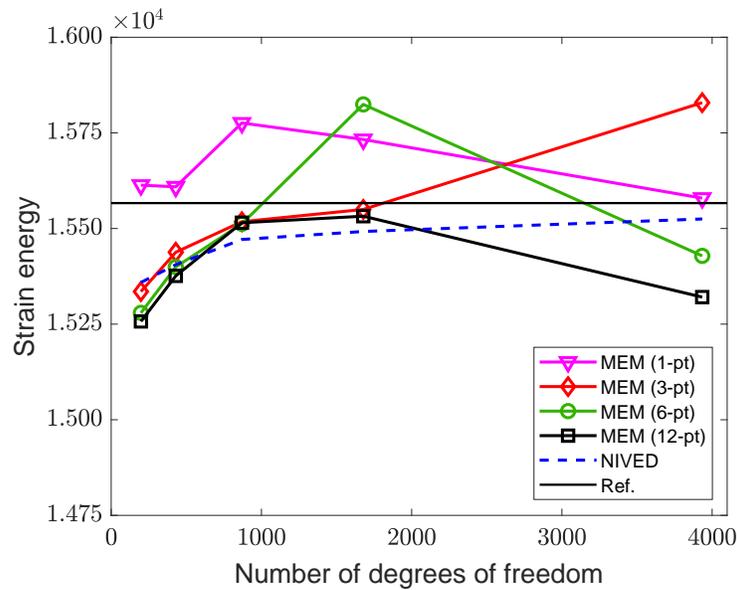, width = 0.7\textwidth}
\caption{Convergence of the strain energy for the $L$-shaped domain problem. 
         The MEM method does not uniformly approach the reference value 
         of the strain energy, whereas the proposed NIVED scheme does.}
\label{fig:ex_lshaped_strain_energy}
\end{figure}

An $L$-shaped domain under traction load is considered. The geometry and boundary
conditions are shown in~\fref{fig:lshapeddomain}, where $H=100$ inch and $p=1$ psi.
The thickness is 1 inch. A stress singularity occurs at the re-entrant corner 
and the exact solution is not known. However, an estimation of the exact
strain energy using Richardson's extrapolation 
is known and has the value 15566.46~\cite{beckers:NCPE:1993}. We use this value 
to assess the convergence of the strain energy. Plane stress condition is assumed 
with $E_\mathrm{Y}=1$ psi and $\nu = 0.3$. The sequence of background meshes used 
in this study are depicted in \fref{fig:lshaped_mesh}. The convergence of
the strain energy for the MEM with several Gauss integration rules and for
the NIVED approach is shown in \fref{fig:ex_lshaped_strain_energy}. We observe 
that even with a 12-point Gauss rule, the MEM does not uniformly approach 
the reference value of the strain energy. In contrast, the proposed 
NIVED scheme uniformly approaches this reference value.

\subsection{Manufactured elastostatic problem}\label{sec:numexamples_manufactured_elastostatic}

\begin{figure}[!tbhp]
\centering
\mbox{
\subfigure[]{\label{fig:manufactured_mesh_a} \epsfig{file = 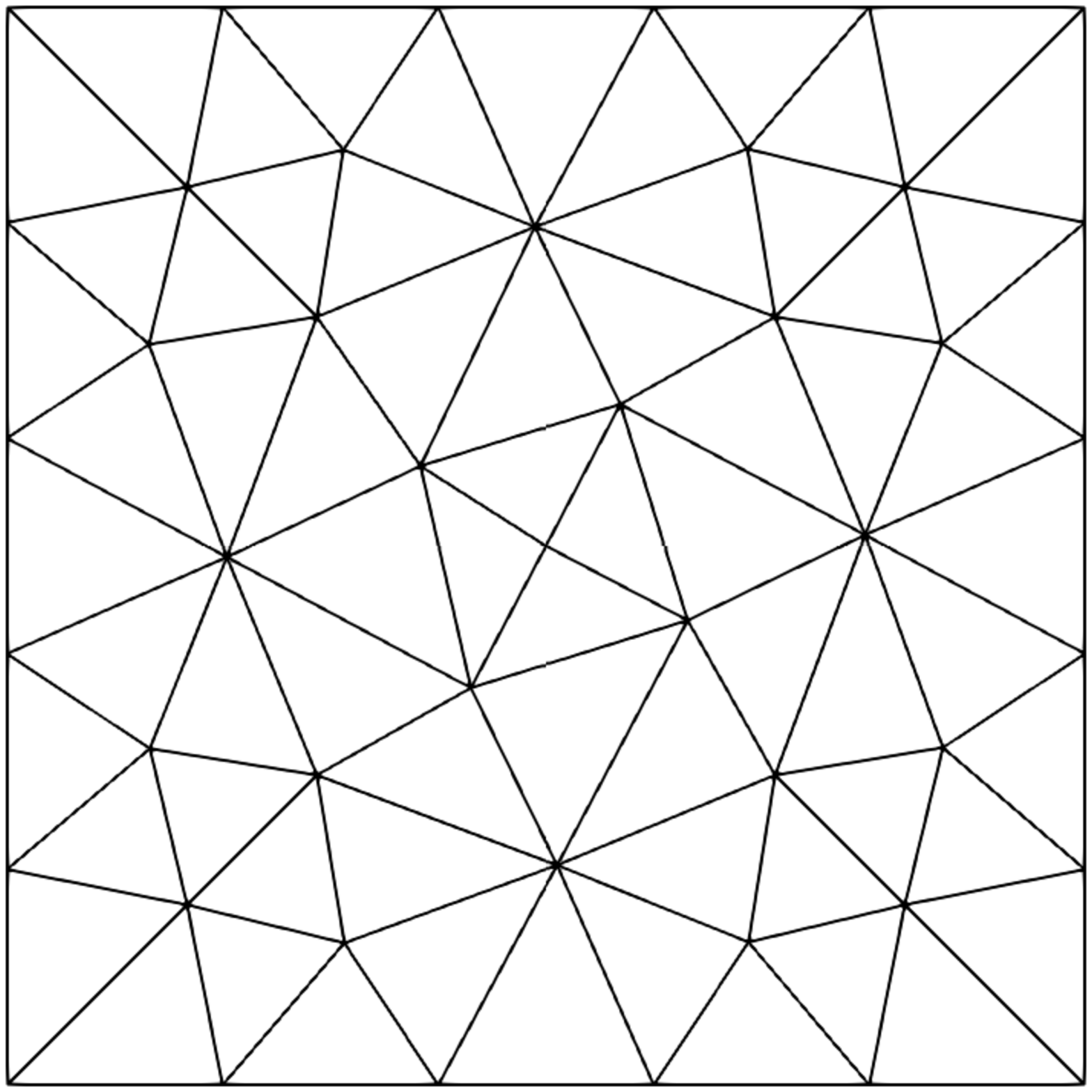, width = 0.15\textwidth}}
\subfigure[]{\label{fig:manufactured_mesh_b} \epsfig{file = 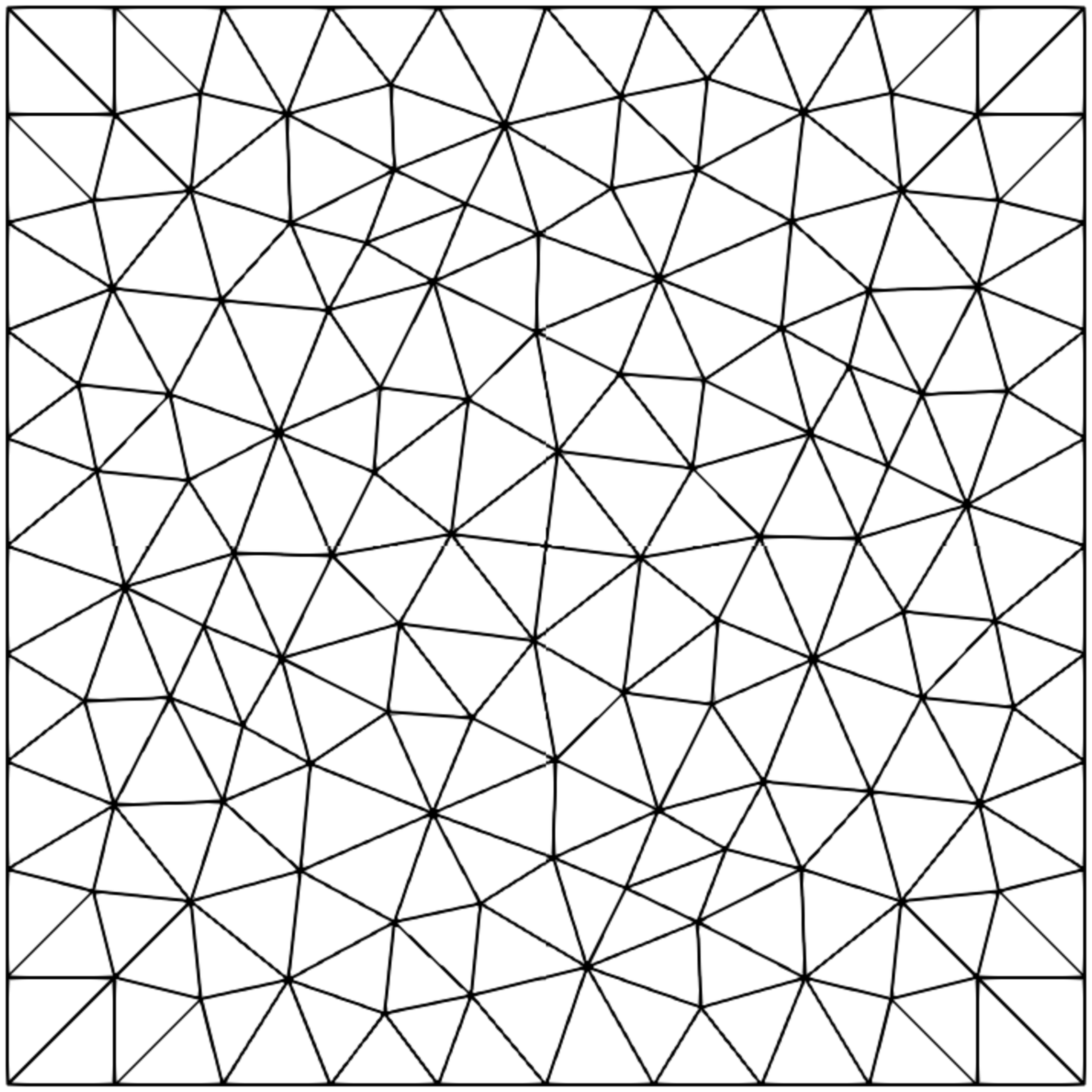, width = 0.15\textwidth}}
\subfigure[]{\label{fig:manufactured_mesh_c} \epsfig{file = 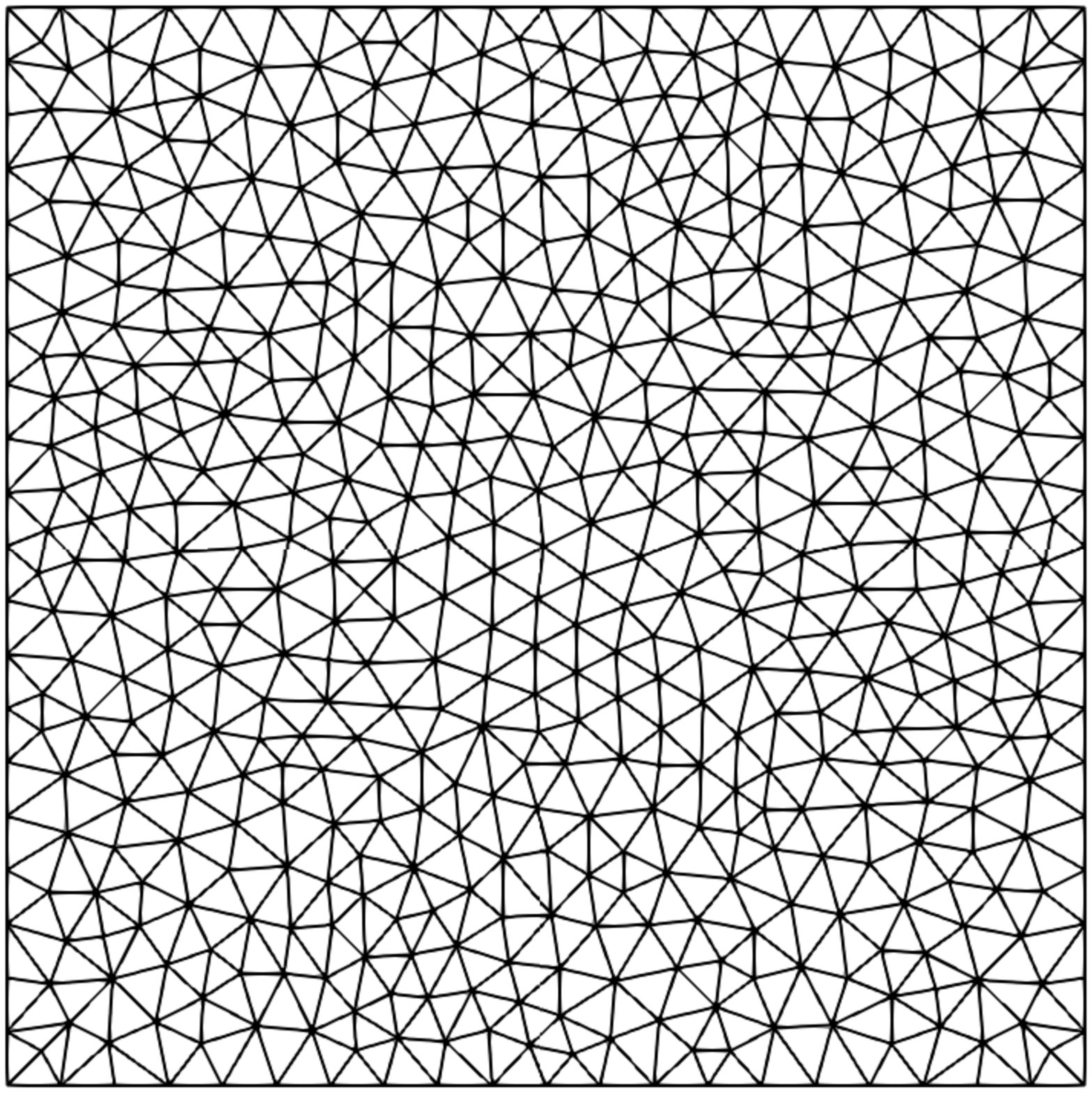, width = 0.15\textwidth}}
\subfigure[]{\label{fig:manufactured_mesh_d} \epsfig{file = 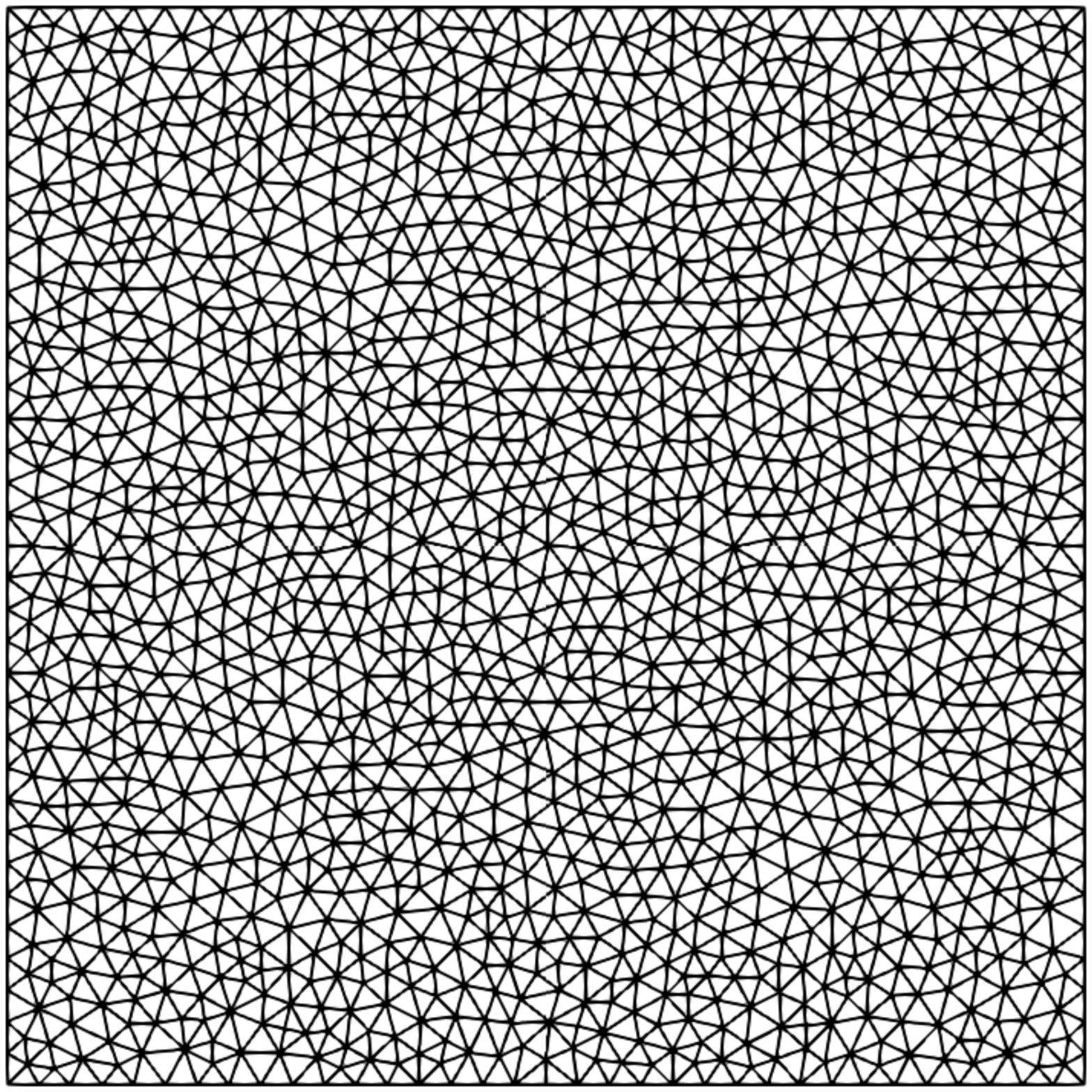, width = 0.15\textwidth}}
\subfigure[]{\label{fig:manufactured_mesh_e} \epsfig{file = 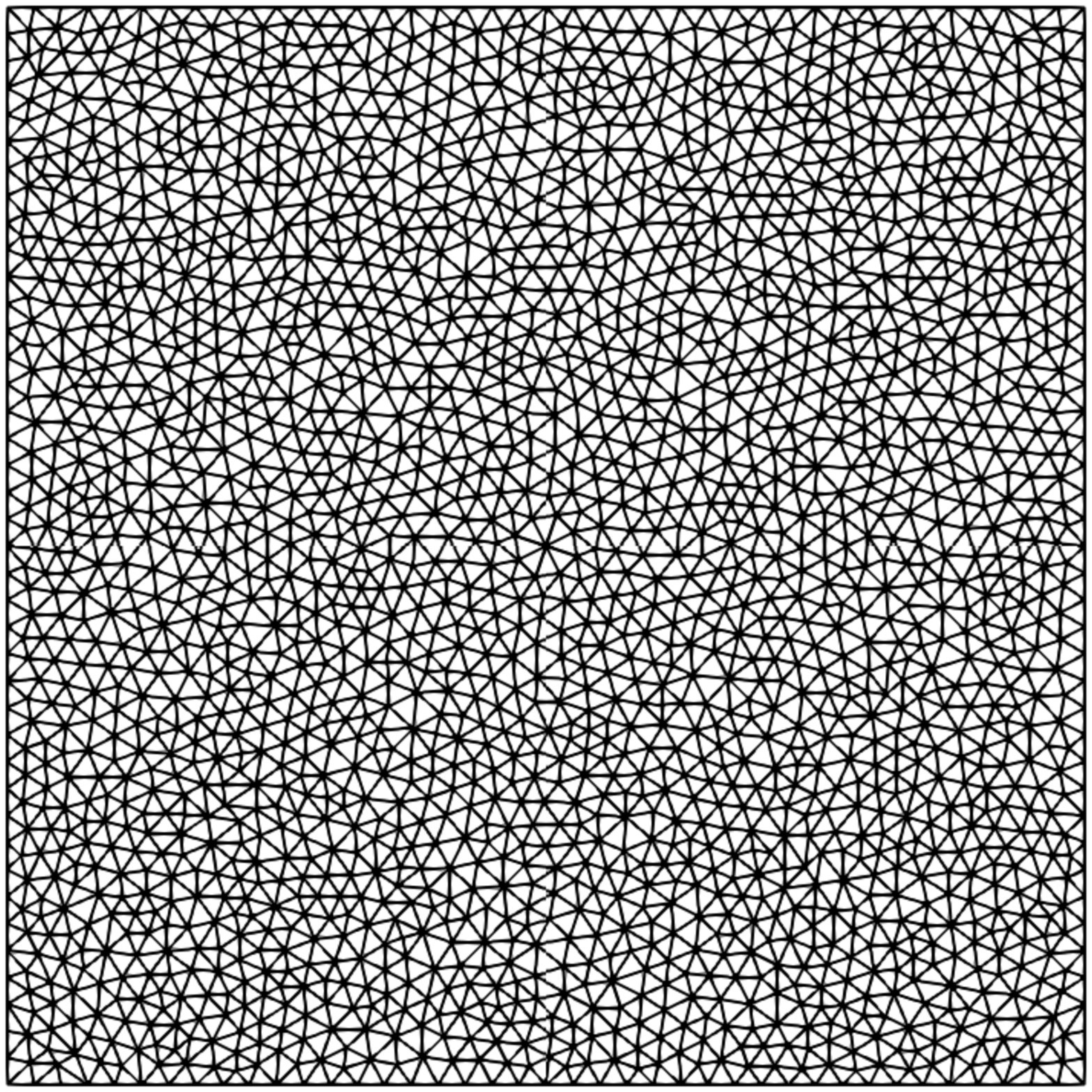, width = 0.15\textwidth}}
\subfigure[]{\label{fig:manufactured_mesh_f} \epsfig{file = 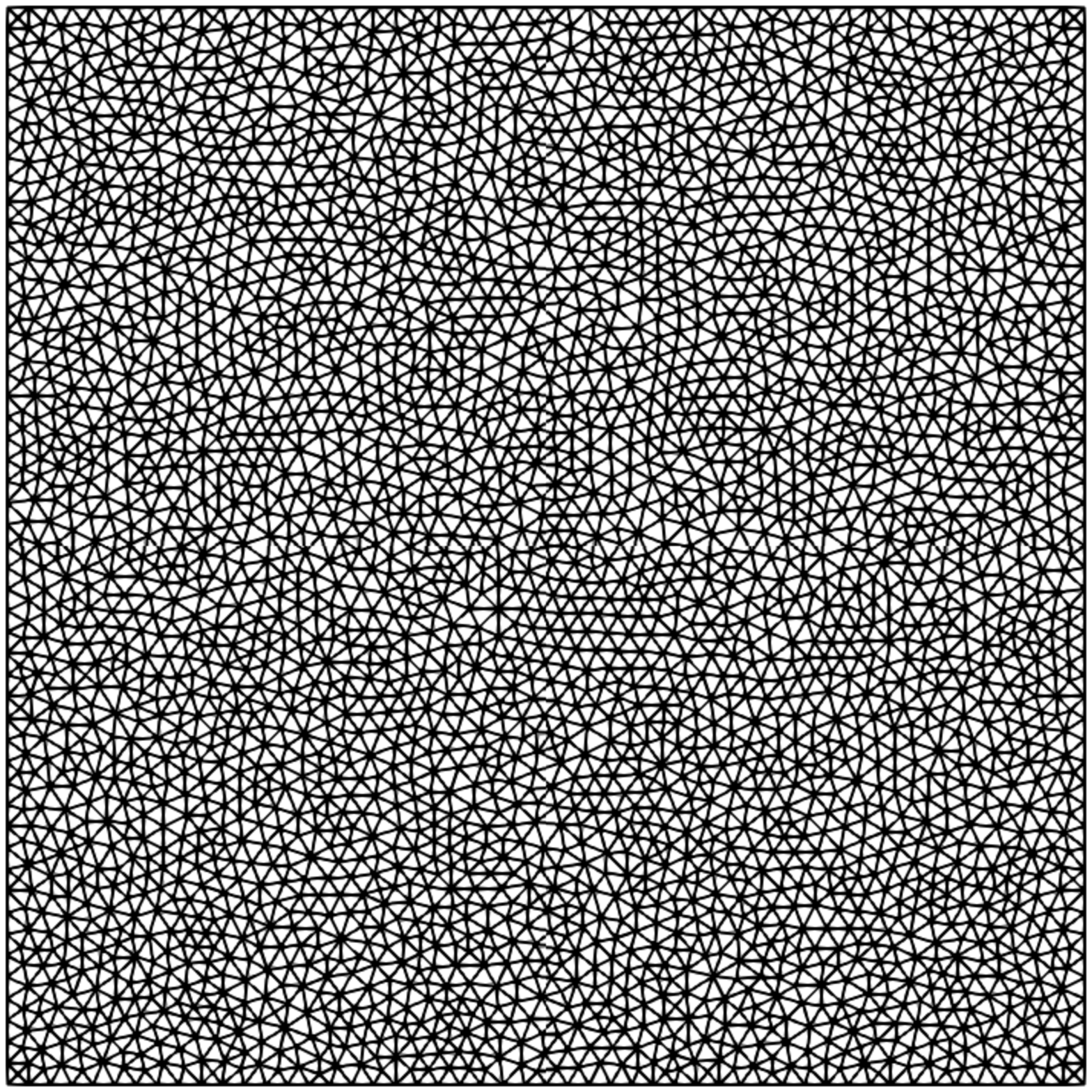, width = 0.15\textwidth}}
}
\mbox{
\subfigure[]{\label{fig:manufactured_mesh_g} \epsfig{file = 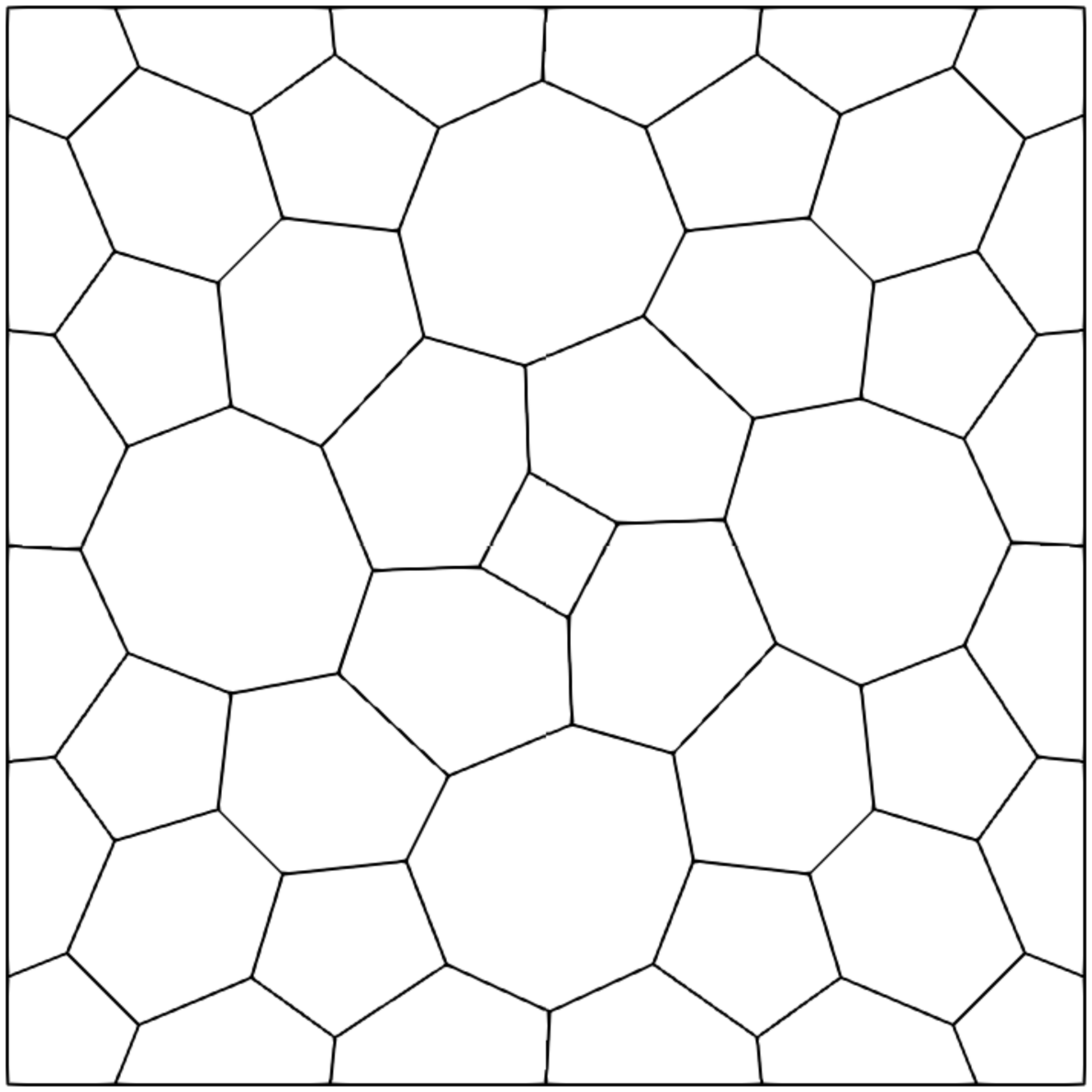, width = 0.15\textwidth}}
\subfigure[]{\label{fig:manufactured_mesh_h} \epsfig{file = 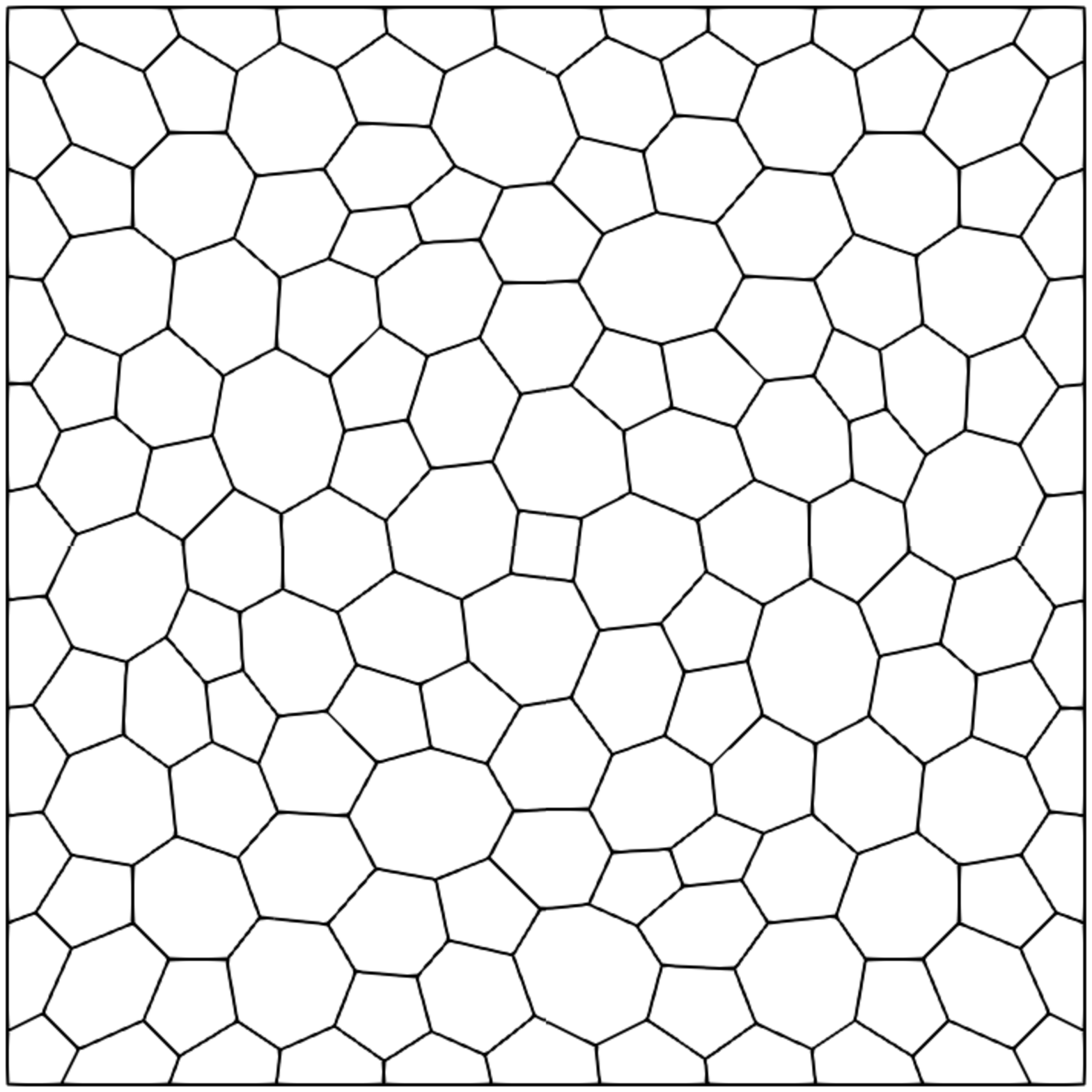, width = 0.15\textwidth}}
\subfigure[]{\label{fig:manufactured_mesh_i} \epsfig{file = 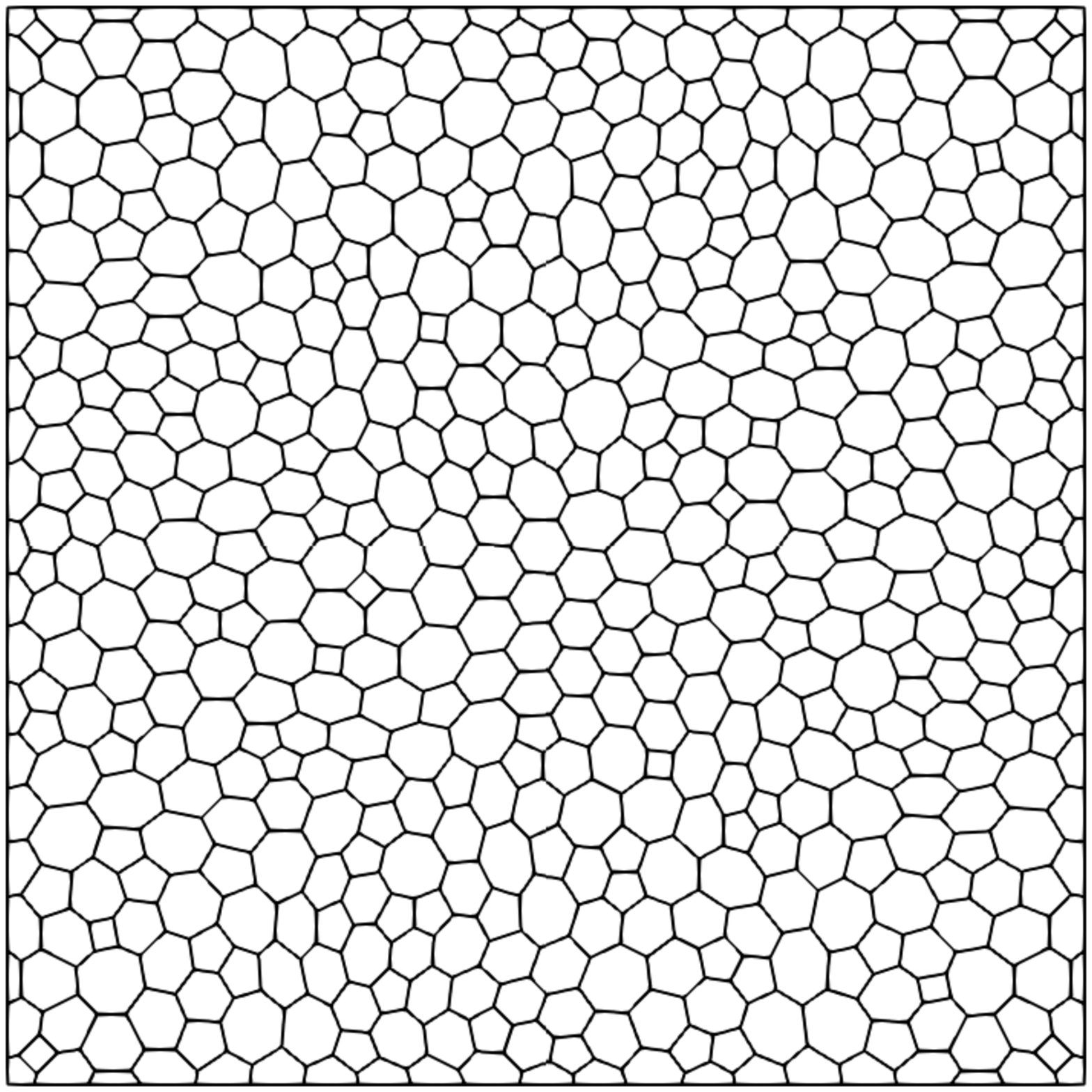, width = 0.15\textwidth}}
\subfigure[]{\label{fig:manufactured_mesh_j} \epsfig{file = 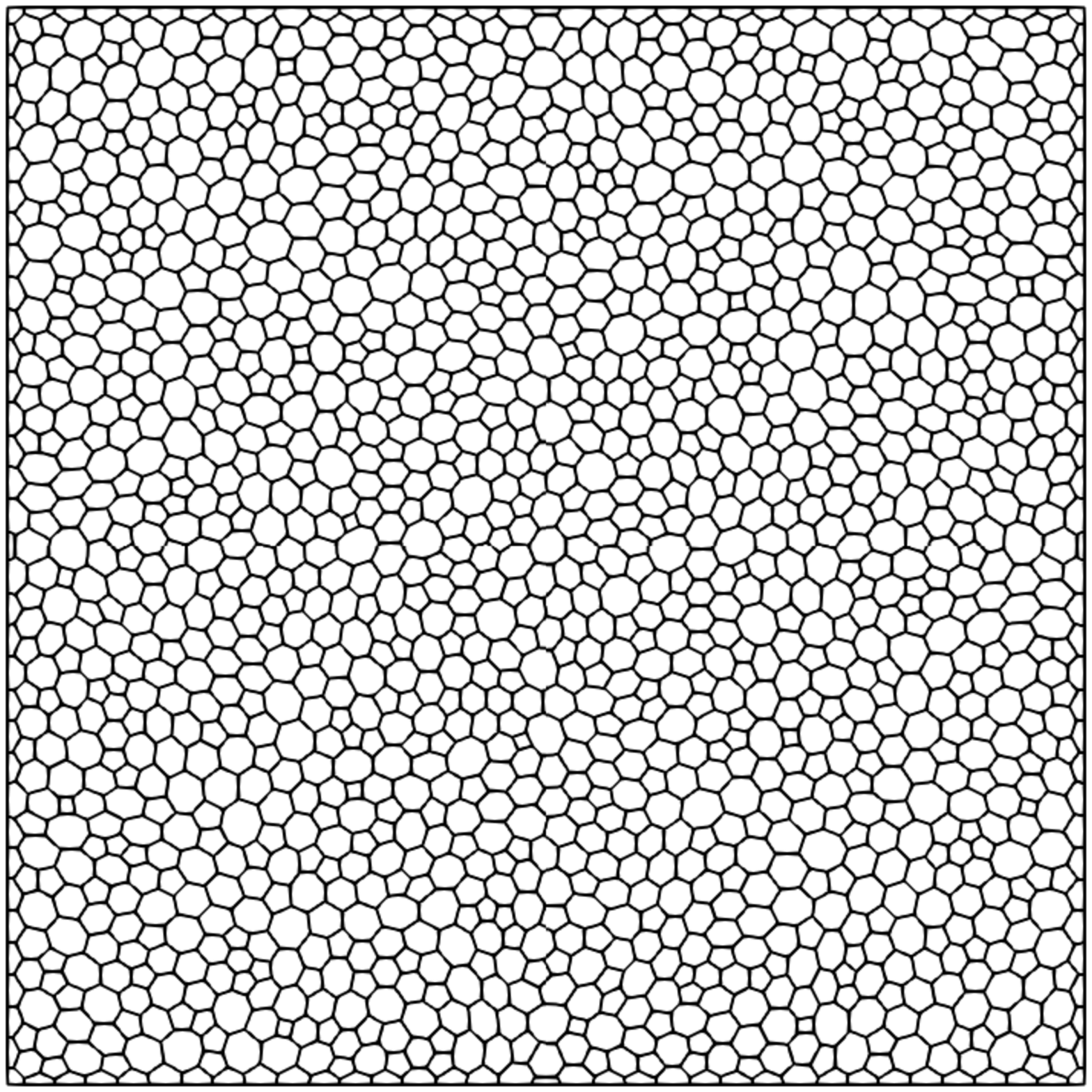, width = 0.15\textwidth}}
\subfigure[]{\label{fig:manufactured_mesh_k} \epsfig{file = 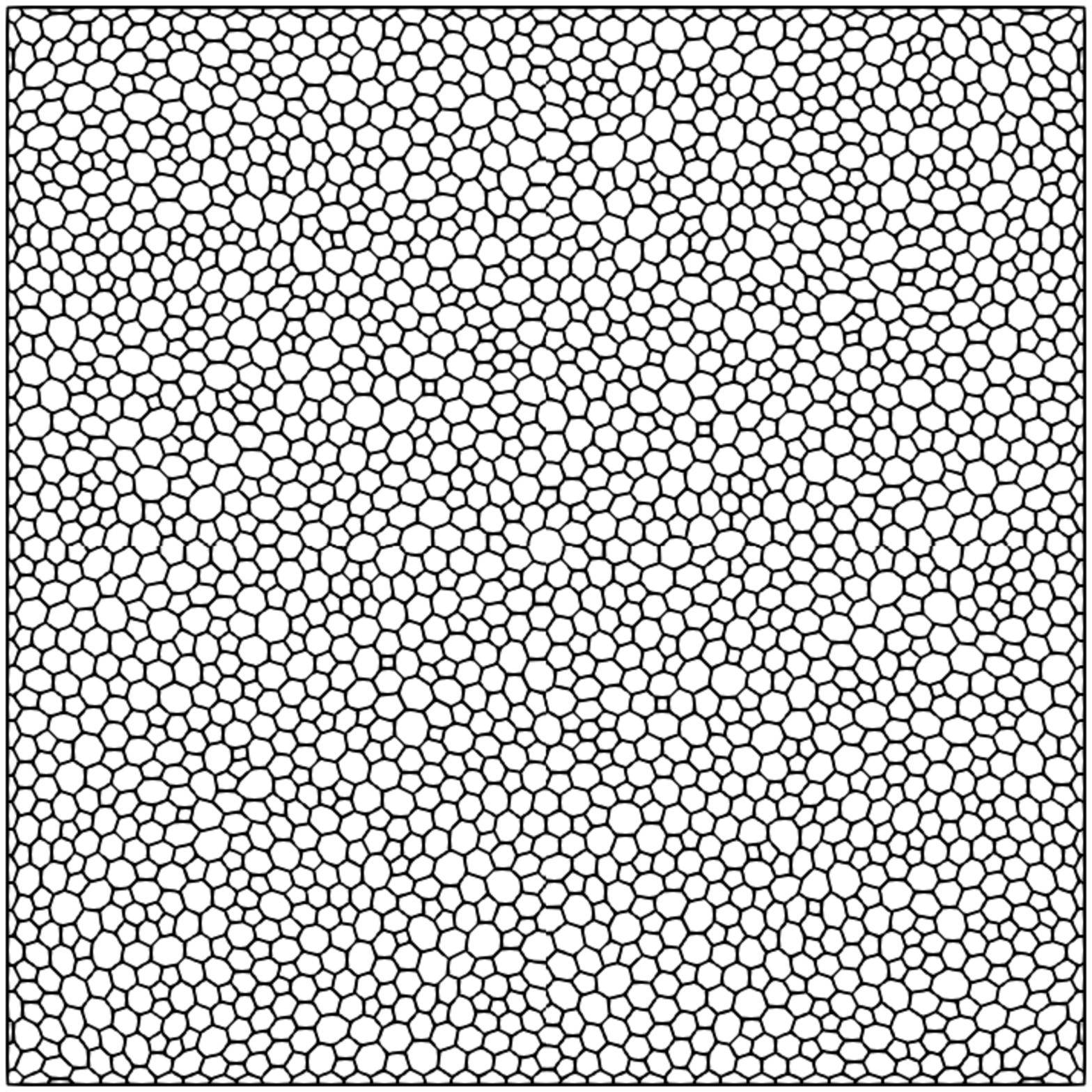, width = 0.15\textwidth}}
\subfigure[]{\label{fig:manufactured_mesh_l} \epsfig{file = 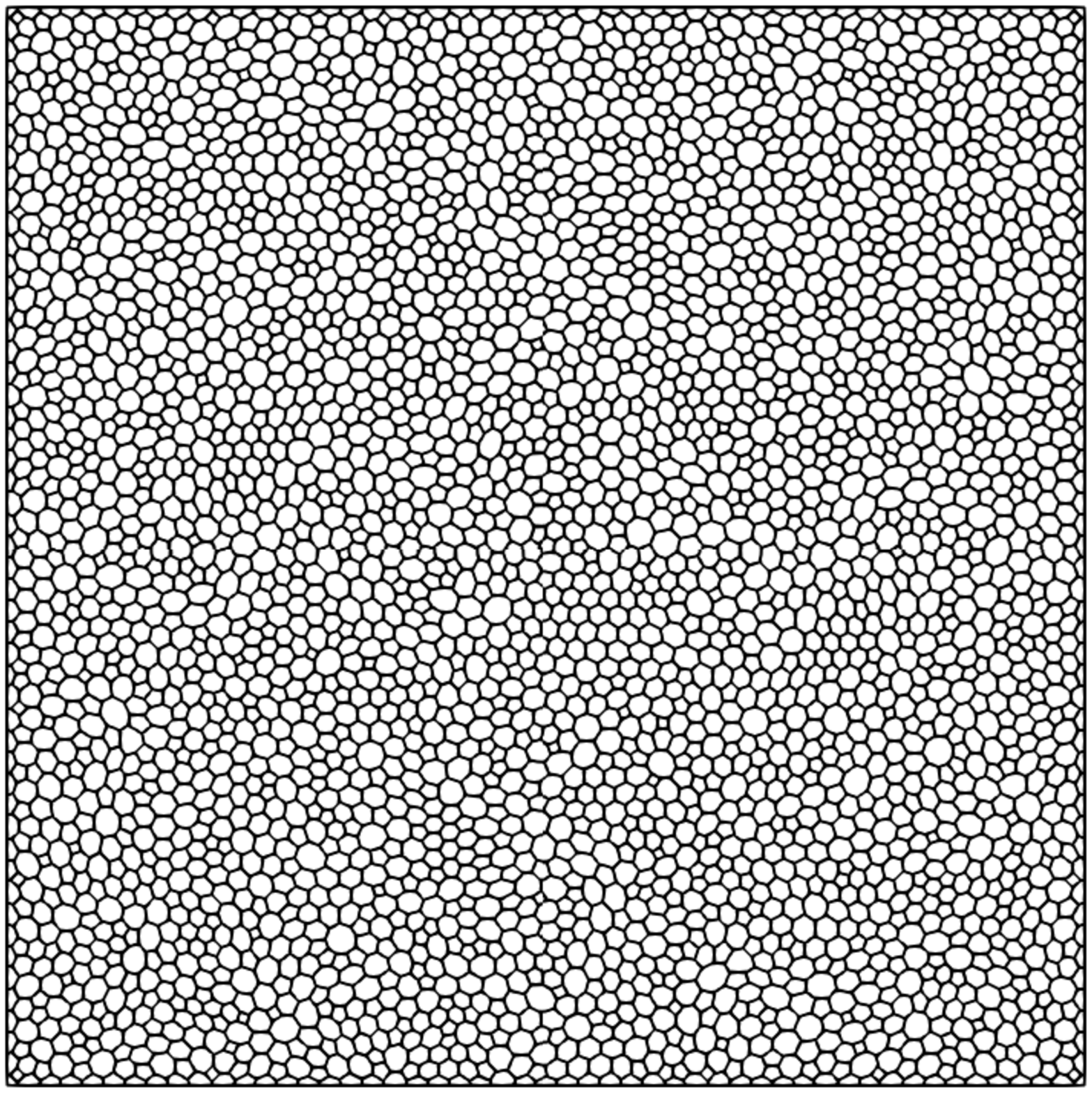, width = 0.15\textwidth}}
}
\caption{Sequence of background integration meshes used for the manufactured elastostatic problem.
         (a)-(f) Meshes for the MEM method and (g)-(l) meshes for the NIVED method.}
\label{fig:manufactured_mesh}
\end{figure}

In this example, the manufactured elastostatic problem 
found in Reference~\cite{duan:2014:CEF} is used to assess
the performance of the NIVED method. The analysis is conducted
on a $2\times 2$ square domain. Plane stress condition is considered
with $E_\mathrm{Y}=10^5$ psi and $\nu = 0.3$. The entire domain boundary 
is prescribed with the following Dirichlet boundary conditions:
\begin{equation*}
u_1(\vx) = \sin{x_1}\cos{x_2},\quad u_2(\vx) = e^{x_1}e^{x_2},
\end{equation*}
which correspond to the exact solutions of the linear elastostatic 
problem~\eref{eq:strongform} manufactured with a body force
given by
\begin{equation*}
\vm{b} = \left[
\begin{array}{c}
(\sin{x_1}\cos{x_2})\mat{D}(1,1) - (e^{x_1}e^{x_2})\mat{D}(1,2) - (e^{x_1}e^{x_2}-\sin{x_1}\cos{x_2})\mat{D}(3,3)\\
(\cos{x_1}\sin{x_2})\mat{D}(2,1) - (e^{x_1}e^{x_2})\mat{D}(2,2) - (e^{x_1}e^{x_2}-\cos{x_1}\sin{x_2})\mat{D}(3,3)
\end{array}
\right].
\end{equation*}
The exact stress field is
\begin{equation*}
\left[\begin{array}{c}
\sigma_{11}\\
\sigma_{22}\\
\sigma_{12}
\end{array}\right] = \left[
\begin{array}{c}
\mat{D}(1,1)\cos{x_1}\cos{x_2} + \mat{D}(1,2)e^{x_1}e^{x_2}\\
\mat{D}(2,1)\cos{x_1}\cos{x_2} + \mat{D}(2,2)e^{x_1}e^{x_2}\\
\mat{D}(3,3)(e^{x_1}e^{x_2}-\sin{x_1}\sin{x_2}).
\end{array}
\right].
\end{equation*}
\fref{fig:manufactured_mesh} depicts the background integration meshes used in 
the study. The Gaussian prior is used for the evaluation of the maximum-entropy 
basis functions.

\begin{figure}[!tbhp]
\centering
\subfigure[]{\label{fig:ex_manufactured_static_norms_a} \epsfig{file = 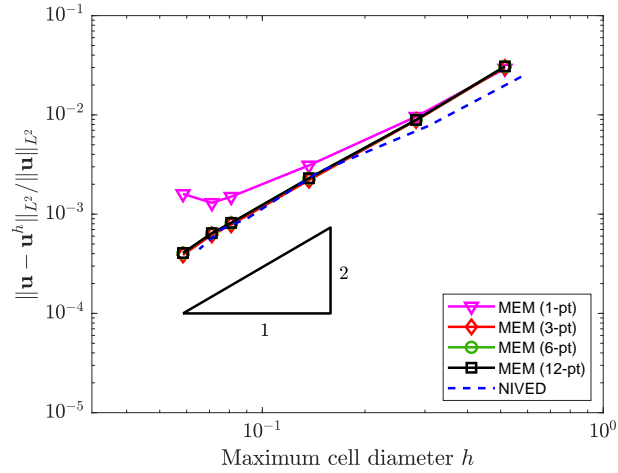, width = 0.58\textwidth}}
\subfigure[]{\label{fig:ex_manufactured_static_norms_b} \epsfig{file = 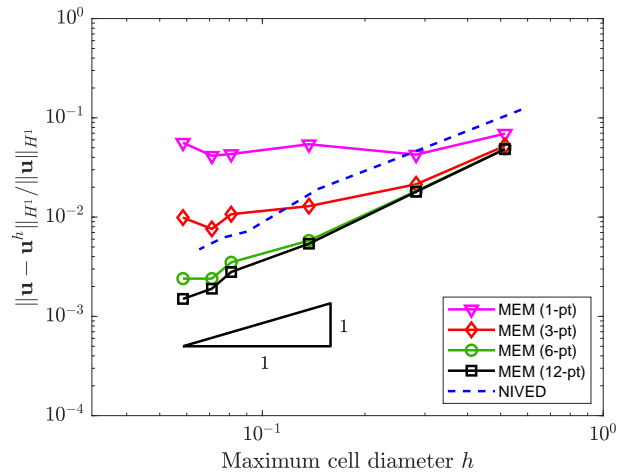, width = 0.58\textwidth}}
\caption{Rates of convergence for the manufactured elastostatic problem. The MEM approach 
         does not convergence in the $L^2$ norm with a 1-point Gauss rule, but 
         the optimal rate of 2 is recovered using a 3-point Gauss rule. 
         Its convergence in the $H^1$ seminorm is optimal 
         when using at least a 12-point Gauss rule. The optimal rates of 2 and 1 are 
         obtained in the $L^2$ norm and the $H^1$ seminorm,
         respectively, for the NIVED method.}
\label{fig:ex_manufactured_static_norms}
\end{figure}

The convergence in the $L^2$ norm and the $H^1$ seminorm for the
MEM and NIVED methods are compared in \fref{fig:ex_manufactured_static_norms}.
The $L^2$ norm is not convergent for the MEM approach using a 1-point Gauss rule, but 
the optimal rate of 2 is recovered using a 3-point Gauss rule. 
The convergence of the MEM method in the $H^1$ seminorm is not good for
1-, 3- and 6-point Gauss rules, but the optimal rate of 1 is recovered
when using a 12-point Gauss rule. The plots show that the NIVED
approach delivers the optimal rates of 2 and 1 in 
the $L^2$ norm and the $H^1$ seminorm,
respectively.

\begin{figure}[!tbhp]
\centering
\mbox{
\subfigure[]{\label{fig:ex_manufactured_static_CPU_cost_a}
\epsfig{file = 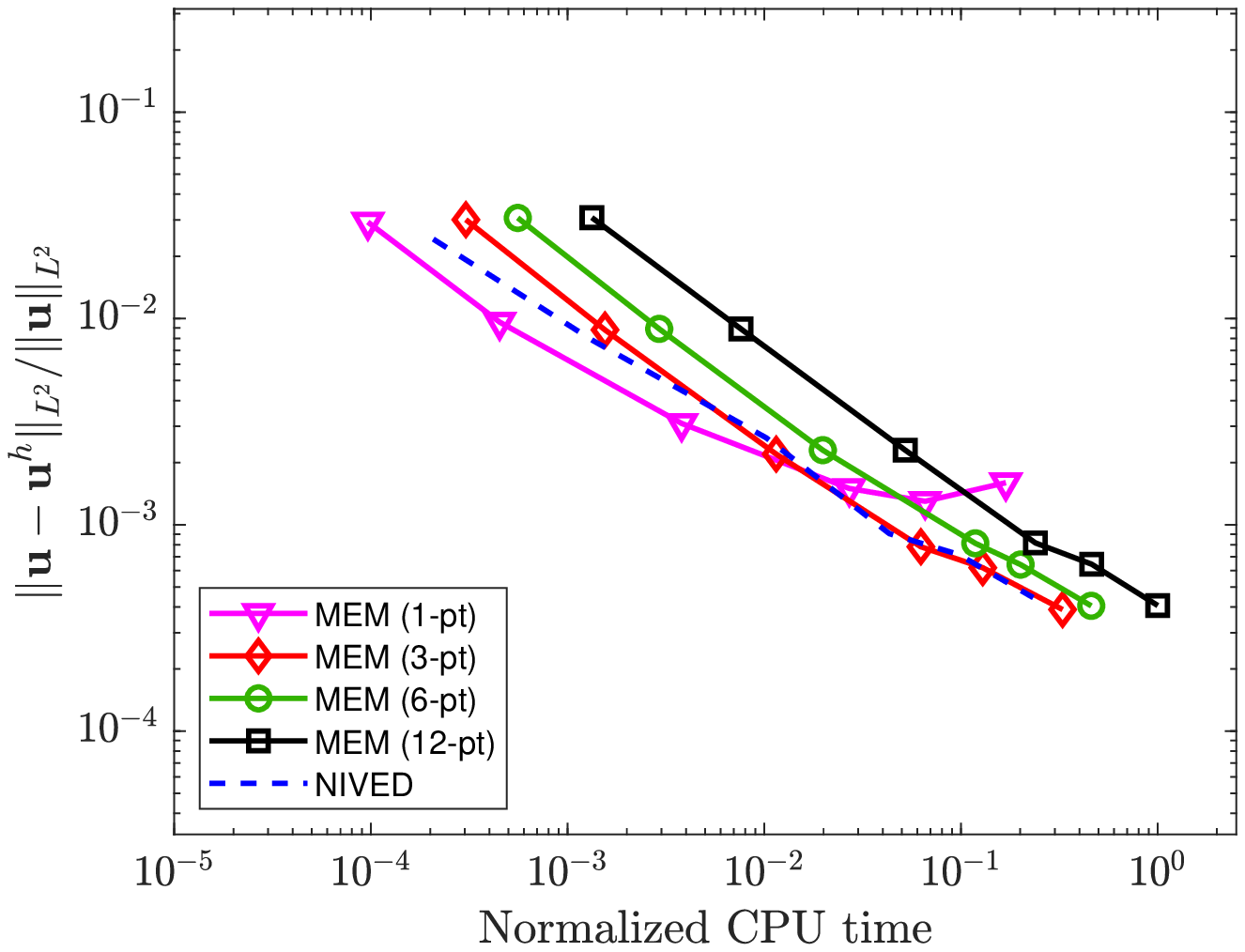, width = 0.5\textwidth}}
\subfigure[]{\label{fig:ex_manufactured_static_CPU_cost_b}
\epsfig{file = 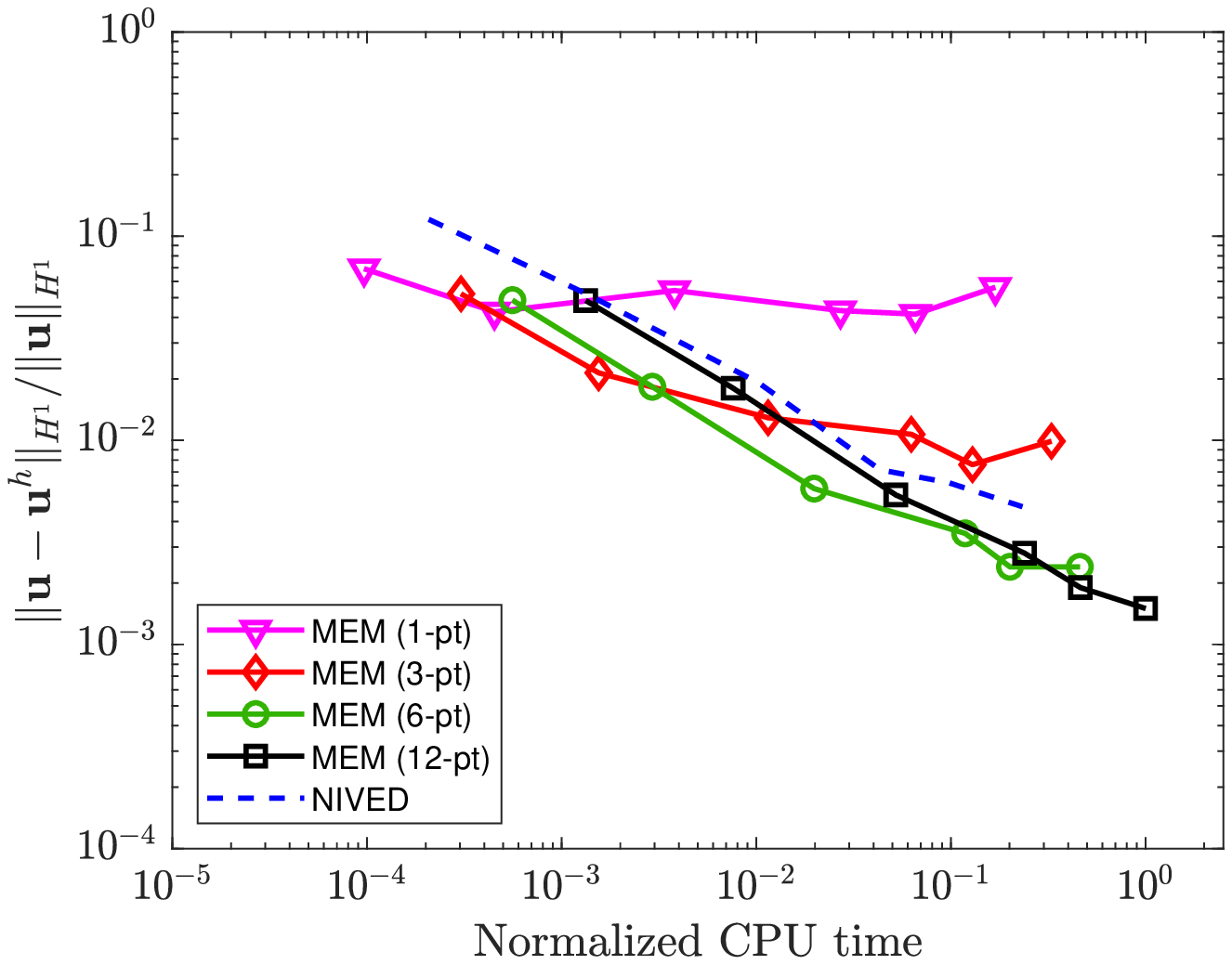, width = 0.5\textwidth}}
}
\subfigure[]{\label{fig:ex_manufactured_static_CPU_cost_c}
\epsfig{file = 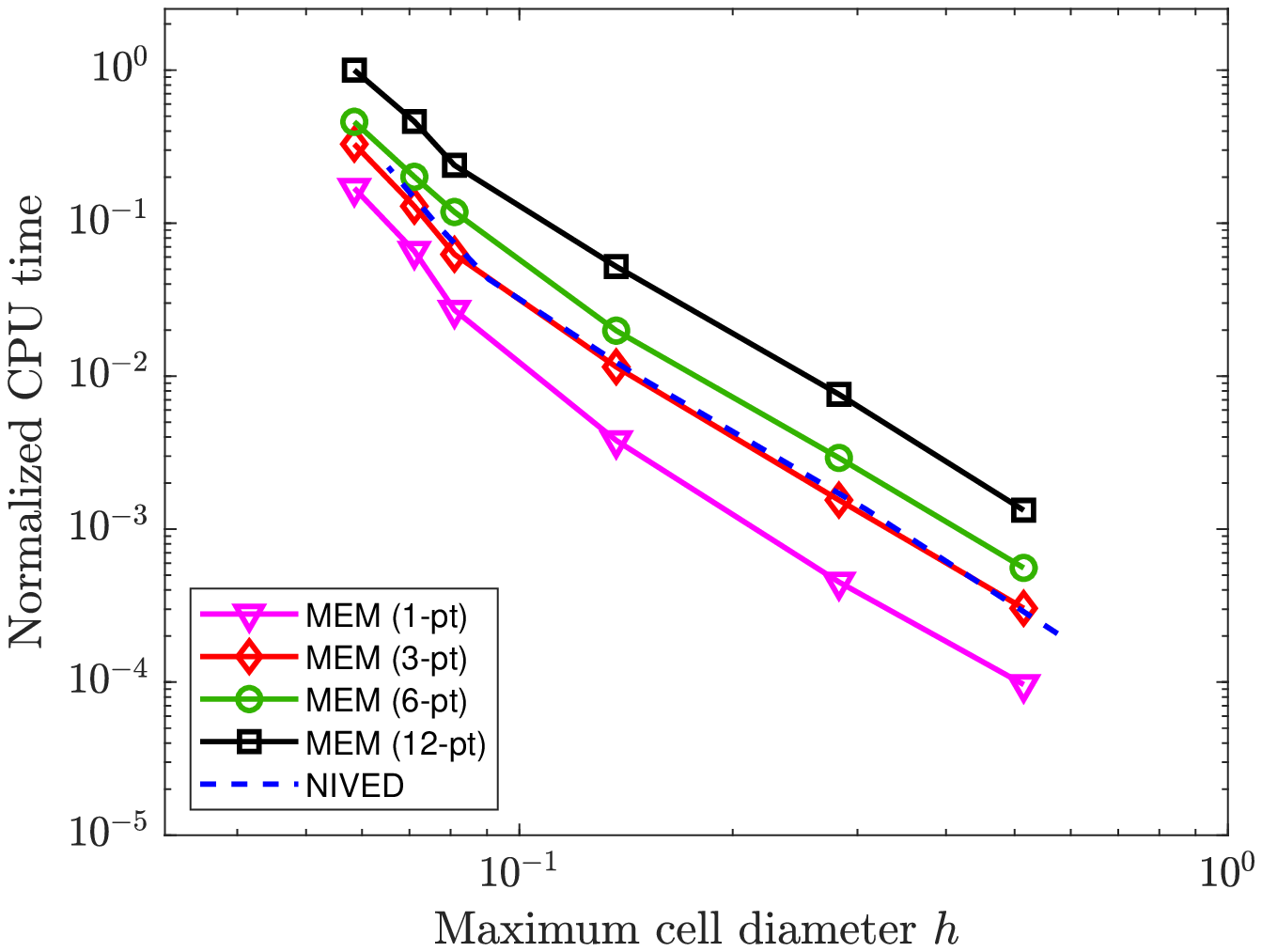, width = 0.5\textwidth}}
\caption{Computational cost of the MEM and NIVED methods for the manufactured elastostatic problem. 
         (a)-(b) Computational cost of the accuracy and (c) computational cost
         of the cell refinement. The computational cost of the NIVED method 
         is similar to the computational cost of the MEM method with a 3-point 
         Gauss rule.}
\label{fig:ex_manufactured_static_CPU_cost}
\end{figure}

The computational cost of the MEM and NIVED methods are compared in
\fref{fig:ex_manufactured_static_CPU_cost}. The methods are assessed
in terms of accuracy and cell refinement using the normalized CPU time, 
which is defined as the ratio of the CPU time of a particular model 
analyzed to the maximum CPU time of any of the models analyzed.
In terms of accuracy, we observe that the computational cost
of the NIVED scheme is similar to the computational cost of the MEM
method with 12-point Gauss rule (\fref{fig:ex_manufactured_static_CPU_cost_b}).
This occurs because for a given level of accuracy, the MEM requires a coarser 
mesh than the NIVED scheme, as shown in~\fref{fig:ex_manufactured_static_norms_b}.
In terms of cell refinement, we observe that 
the computational cost of the NIVED method 
is similar to the computational cost of the MEM method 
with a 3-point Gauss rule (\fref{fig:ex_manufactured_static_CPU_cost_c}). 
However, the MEM with a 3-point
Gauss rule does not converge in the $H^1$ seminorm, which means that 
the computational cost of a convergent MEM approach is greater 
than the computational cost of the NIVED method.
In fact, from \fref{fig:ex_manufactured_static_norms} we know that
a 12-point Gauss rule is needed for optimal $H^1$-convergence of
the MEM method, but this rule is the most expensive as 
shown in \fref{fig:ex_manufactured_static_CPU_cost_c}.
Combining both the cost in terms of accuracy and cell refinement,
we conclude that in general the NIVED scheme can be expected to be faster but less accurate than
a convergent MEM approach with the same number of degrees of freedom.

\subsection{Manufactured elastodynamic problem}\label{sec:numexamples_manufactured_elastodynamic}

A manufactured elastodynamic problem that is found in Reference~\cite{Duan:QCNI:2014} is
considered. The domain of analysis is a $2\times 2$ 
square domain and the background meshes considered are the same used for the manufactured 
elastostatic problem (\fref{fig:manufactured_mesh}). The material parameters
are $E_\mathrm{Y} = 1 \times 10^{5}$ psi, $\nu = 0.3$ and $\rho = 800$ lb/$\mathrm{in}^3$, 
and plane stress condition is assumed. The Gaussian prior is used for the evaluation 
of the maximum-entropy basis functions.

The problem is manufactured with the following body force:
\begin{equation*}
\vm{b} = \left[
\begin{array}{c}
x_1\Big[-2E_\mathrm{Y}g(t) +\rho\ddot{g}(t)\left[-(1+\nu)L^2+ \nu x_2^2 + \frac{1}{3}x_1^2 \right]\Big]\\
x_2\Big[2E_\mathrm{Y}g(t) +\rho\ddot{g}(t)\left[(1+\nu)L^2- \nu x_1^2 - \frac{1}{3}x_2^2 \right] \Big]
\end{array}
\right],
\end{equation*}
where $L$ is the side length of the domain of analysis and
\begin{equation*}
g(t)=\alpha\Big[1-e^{-\frac{\beta t^2}{2}}\Big].
\end{equation*}
The exact displacement field solution is
\begin{equation*}
\vm{u}=\left[
\begin{array}{c}
x_1g(t)\left[-(1+\nu)L^2+\nu x_2^2+\frac{1}{3}x_1^2\right]\\
x_2g(t)\left[(1+\nu)L^2-\nu x_1^2-\frac{1}{3}x_2^2\right]
\end{array}
\right],
\end{equation*}
and the exact stress field is
\begin{equation*}
\left[\begin{array}{c}
\sigma_{11}\\
\sigma_{22}\\
\sigma_{12}
\end{array}\right]=\left[
\begin{array}{c}
E_\mathrm{Y}g(t)\left[x_1^2 - L^2 \right]\\
E_\mathrm{Y}g(t)\left[L^2 - x_2^2 \right]\\
0
\end{array}
\right].
\end{equation*}

At the initial condition $t=0$ seconds, the body is at rest. The exact displacement 
field solution is used to impose the Dirichlet boundary conditions along the entire
boundary of the domain. The following values for the parameters are 
used in the computations: $\alpha = 0.001$, $\beta = 0.001$. 
The Newmark's constant average acceleration method with a time 
increment $\Delta t = 0.01$ seconds is used as the time integration algorithm.

\begin{figure}[!tbhp]
\centering
\subfigure[]{\label{fig:ex_manufactured_dynamic_norms_a} \epsfig{file = 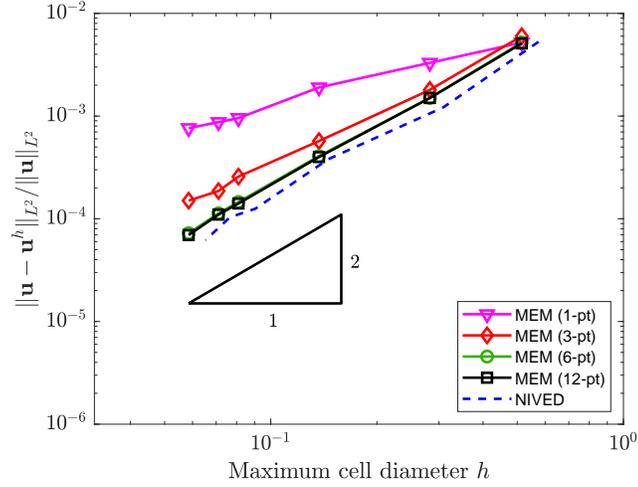, width = 0.6\textwidth}}
\subfigure[]{\label{fig:ex_manufactured_dynamic_norms_b} \epsfig{file = 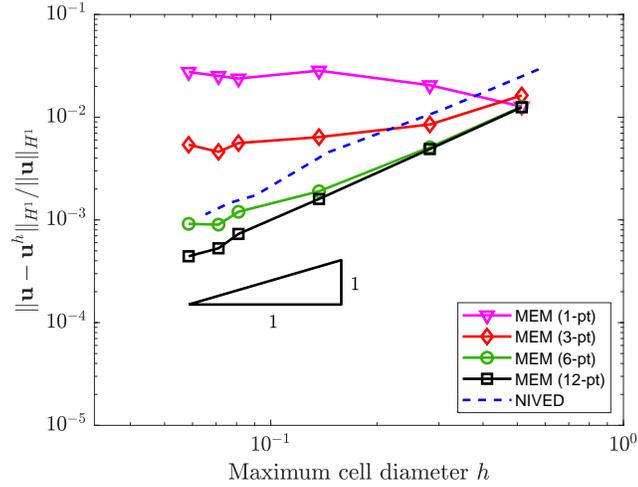, width = 0.6\textwidth}}
\caption{Rates of convergence for the manufactured elastodynamic problem. The MEM approach 
         needs a 3-point Gauss rule to converge with an optimal rate of 2 
         in the $L^2$ norm, and at least a 12-point Gauss rule to converge 
         with an optimal rate of 1 in the $H^1$ seminorm. The optimal rates 
         of 2 and 1 are obtained in the $L^2$ norm and the $H^1$ seminorm,
         respectively, for the NIVED method.}
\label{fig:ex_manufactured_dynamic_norms}
\end{figure}

The convergence of the MEM and NIVED methods in the $L^2$ norm and 
the $H^1$ seminorm after 100 time steps (i.e., $t=1$ s) are compared
in \fref{fig:ex_manufactured_dynamic_norms}.
The MEM approach exhibits suboptimal convergence in the $L^2$ norm
with a 1-point Gauss rule, but 
the optimal rate of 2 is recovered using a 3-point Gauss rule. 
Its convergence in the $H^1$ seminorm is optimal 
when using at least a 12-point Gauss rule. The optimal rates of 2 and 1 are 
obtained in the $L^2$ norm and the $H^1$ seminorm,
respectively, for the NIVED method.

\begin{figure}[!tbhp]
\centering
\mbox{
\subfigure[]{\label{fig:ex_manufactured_dynamic_CPU_cost_a}
\epsfig{file = 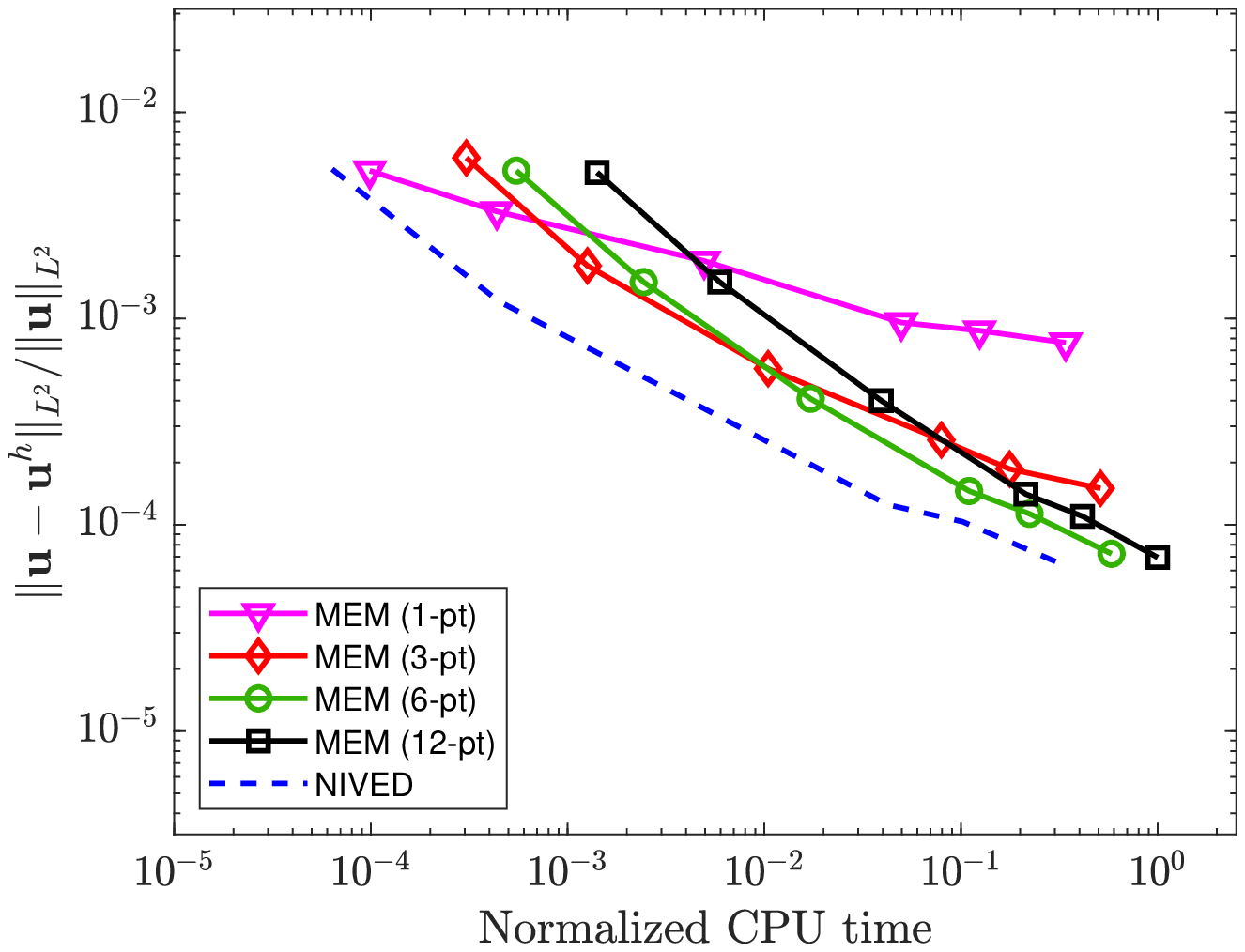, width = 0.5\textwidth}}
\subfigure[]{\label{fig:ex_manufactured_dynamic_CPU_cost_b}
\epsfig{file = 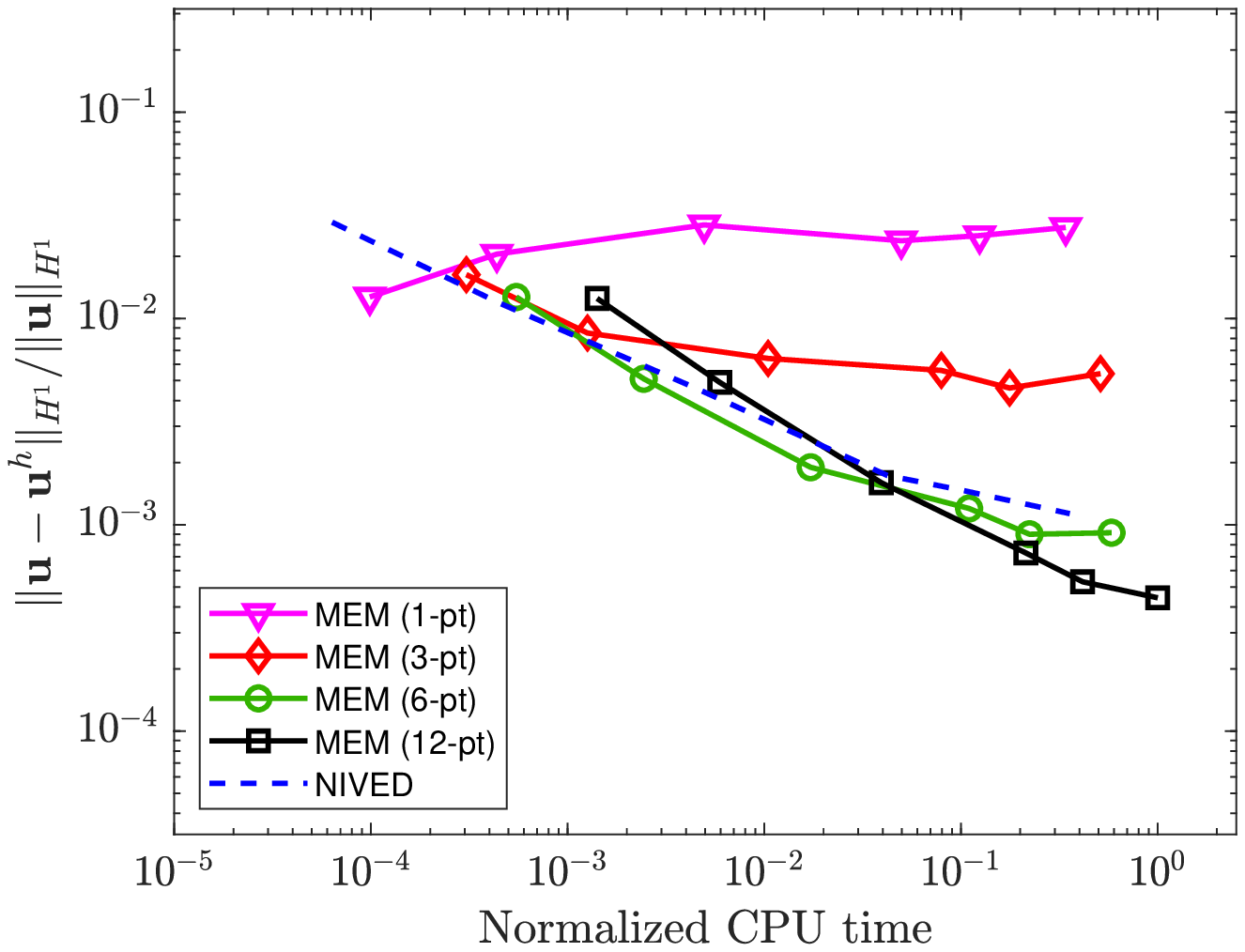, width = 0.5\textwidth}}
}
\subfigure[]{\label{fig:ex_manufactured_dynamic_CPU_cost_c}
\epsfig{file = 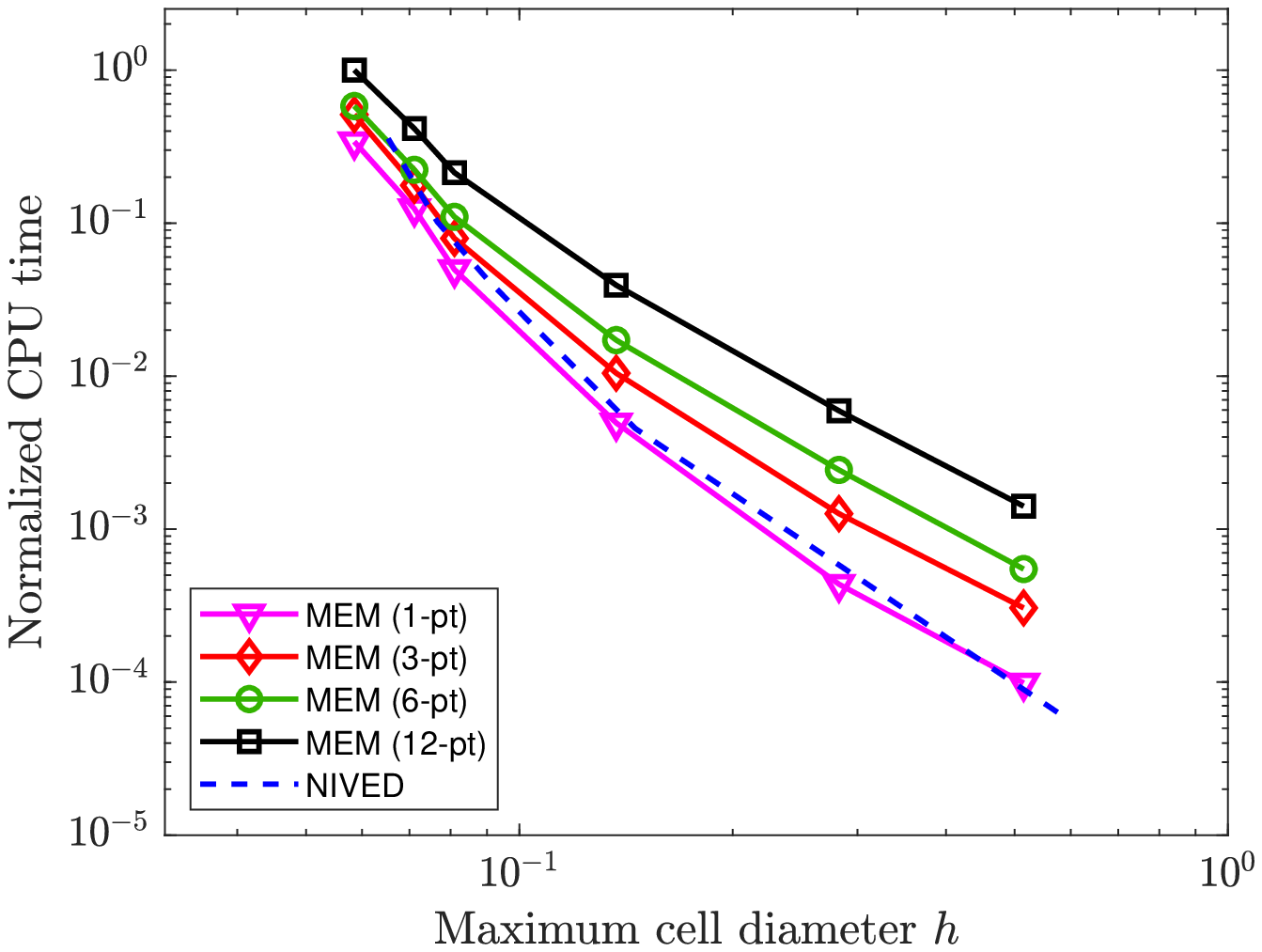, width = 0.5\textwidth}}
\caption{Computational cost of the MEM and NIVED methods for the manufactured elastodynamic problem. 
         (a)-(b) Computational cost of the accuracy and (c) computational cost
         of the cell refinement. The computational cost of the NIVED method 
         is similar to the computational cost of the MEM method using a 1-point 
         Gauss rule.}
\label{fig:ex_manufactured_dynamic_CPU_cost}
\end{figure}

Similar to the static case, \fref{fig:ex_manufactured_dynamic_CPU_cost} (also obtained
at $t=1$ s) reveals that in terms of cell refinement, the computational 
cost of a convergent MEM approach is greater than the NIVED counterpart 
thereby making the NIVED approach superior 
in performance. Moreover, in the dynamic regime the outperformance of the
NIVED over the MEM, in terms of cell refinement, is more pronounced 
as the computational cost of the NIVED approach is about the same as the 
computational cost of the MEM using
a 1-point Gauss rule (\fref{fig:ex_manufactured_dynamic_CPU_cost_c}).

\subsection{Thick-walled viscoelastic cylinder subjected to internal pressure}\label{sec:numexamples_thickwalled}

\begin{figure}[!tbhp]
\centering
\epsfig{file = 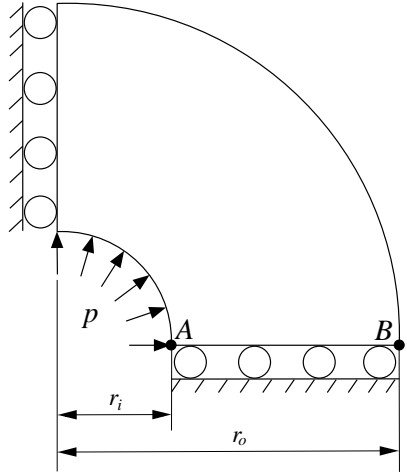, width = 0.35\textwidth}
\caption{Domain and boundary conditions for the thick-walled viscoelastic cylinder subjected to internal pressure.}
\label{fig:thickwalled_domain}
\end{figure}

We conclude this section with a demonstration of the nonlinear NIVED formulation for 
small-displacement and small-strain analysis. The test problem consists of a
thick-walled cylinder subjected to internal pressure whose material 
constitution is a linear isotropic solid described by the Generalized Maxwell 
viscoelastic model considered in Section~\ref{sec:nonlinear_nived}. This example 
has been used to assess the virtual element method for inelastic analysis in 
Reference~\cite{artioli:AOVEMIN:2017}. Due to symmetry,
only a quarter of the cross section of the cylinder is considered, as shown 
in~\fref{fig:thickwalled_domain}, where $r_i=2$ inch, $r_0=4$ inch and $p=10$ psi. 
The material parameters are set to $E=1000$ psi, $\nu=0.3$ and $\lambda_1=1$. 
The pressure is applied suddenly and the structural response is computed
through 20 unit time steps. 

In the first part of this study, the convergence of the radial displacement at 
control points $A$ and $B$ (\fref{fig:thickwalled_domain}) is assessed. For comparison
purposes, the numerical solution obtained with a highly refined mesh of 9-node FEM quadrilateral 
elements (FEM-Q9) is used. Three sets of choices are considered for the 
material parameters $\mu_0$ and $\mu_1$: $(\mu_0,\mu_1)=\{(0.7,0.3),\,(0.3,0.7),\,(0.01,0.99)\}$.
Using $E$ and $\nu$ to calculate the material's bulk modulus $K$ and~\eref{eq:relaxmodulus} 
to calculate $G(t=20)$, the effective Poisson's ratio is calculated 
as $\nu_{\mathrm{eff}}=\frac{3K-2G(20)}{2(3K+G(20))}$. For the three preceding 
sets of parameters this gives at time $t=20$ the values
$\nu_{\mathrm{eff}}=\{0.3542,\, 0.4338,\, 0.4977\}$, respectively. 
Thus, for the last set, the material is near-incompressible.

\begin{figure}[!tbhp]
\centering
\mbox{
\subfigure[]{\label{fig:thickwalled_mesh_a} \epsfig{file = 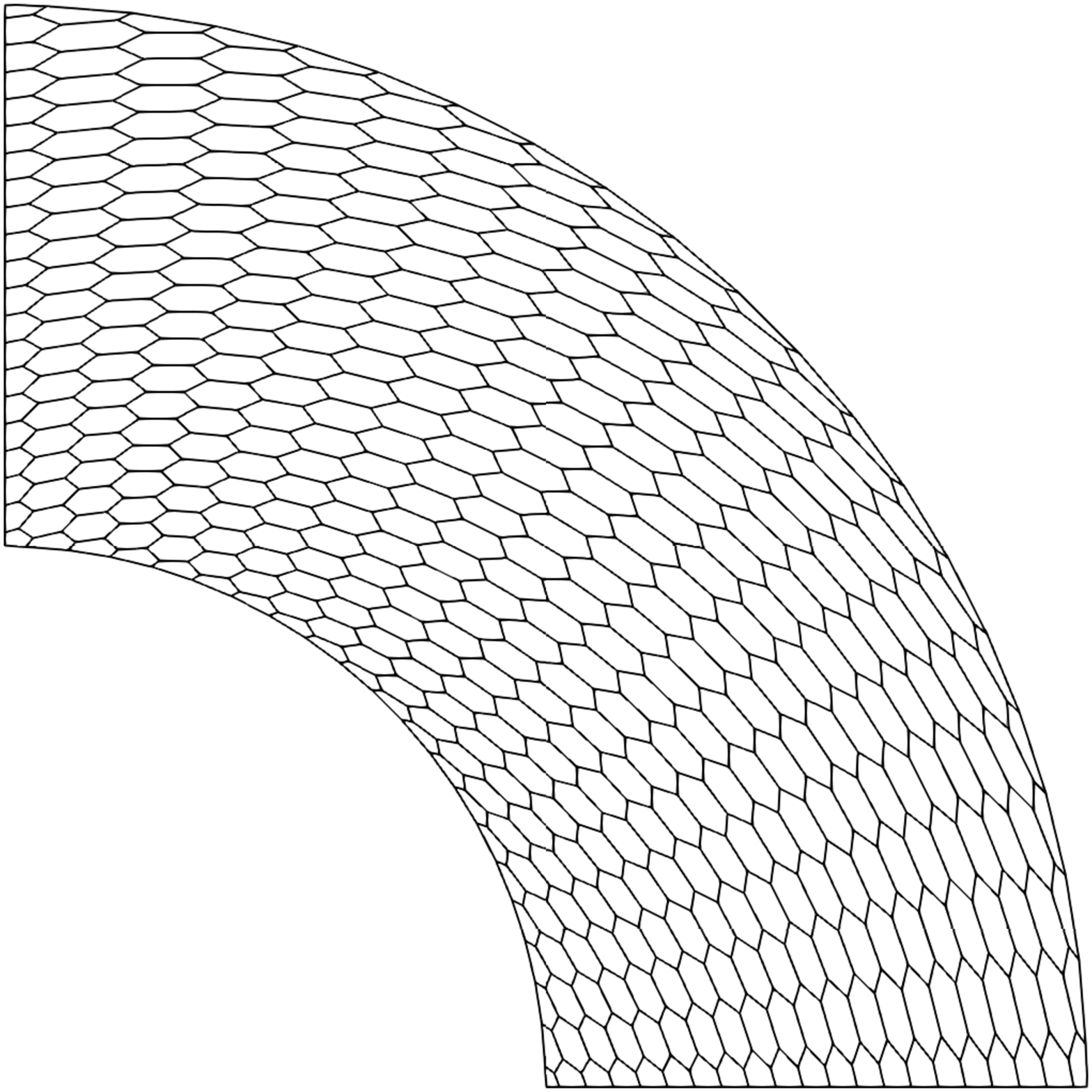, width = 0.185\textwidth}}
\subfigure[]{\label{fig:thickwalled_mesh_b} \epsfig{file = 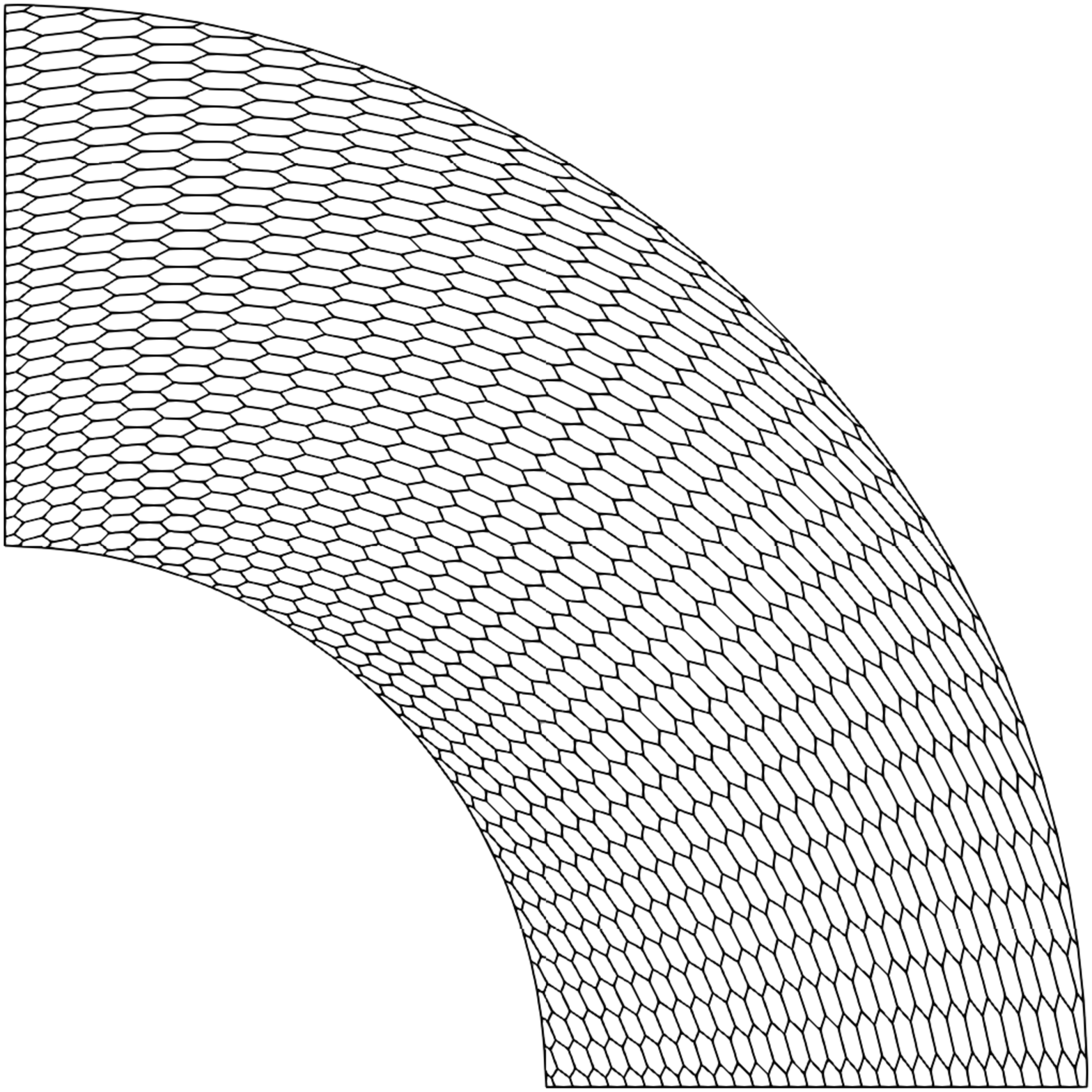, width = 0.185\textwidth}}
\subfigure[]{\label{fig:thickwalled_mesh_c} \epsfig{file = 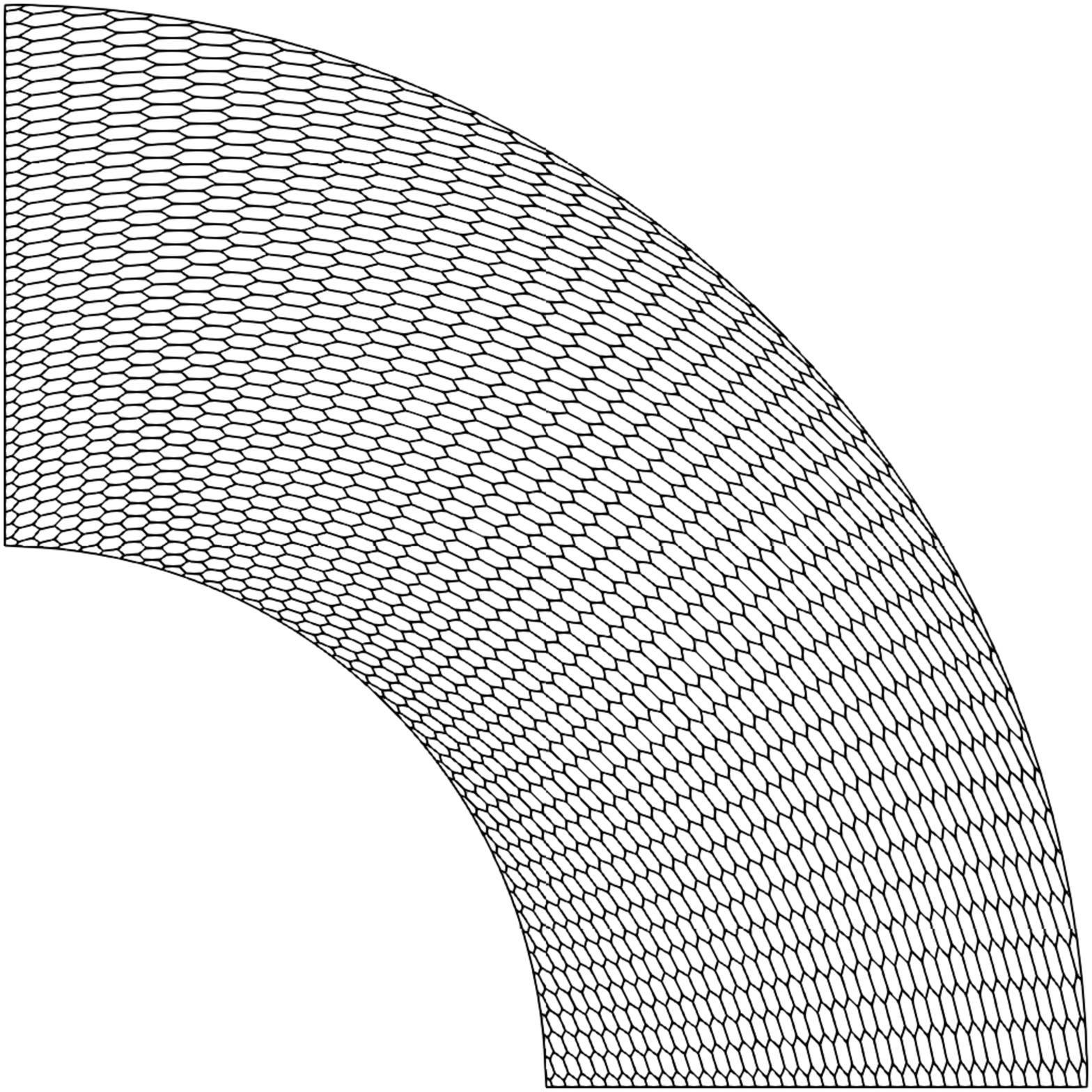, width = 0.185\textwidth}}
\subfigure[]{\label{fig:thickwalled_mesh_d} \epsfig{file = 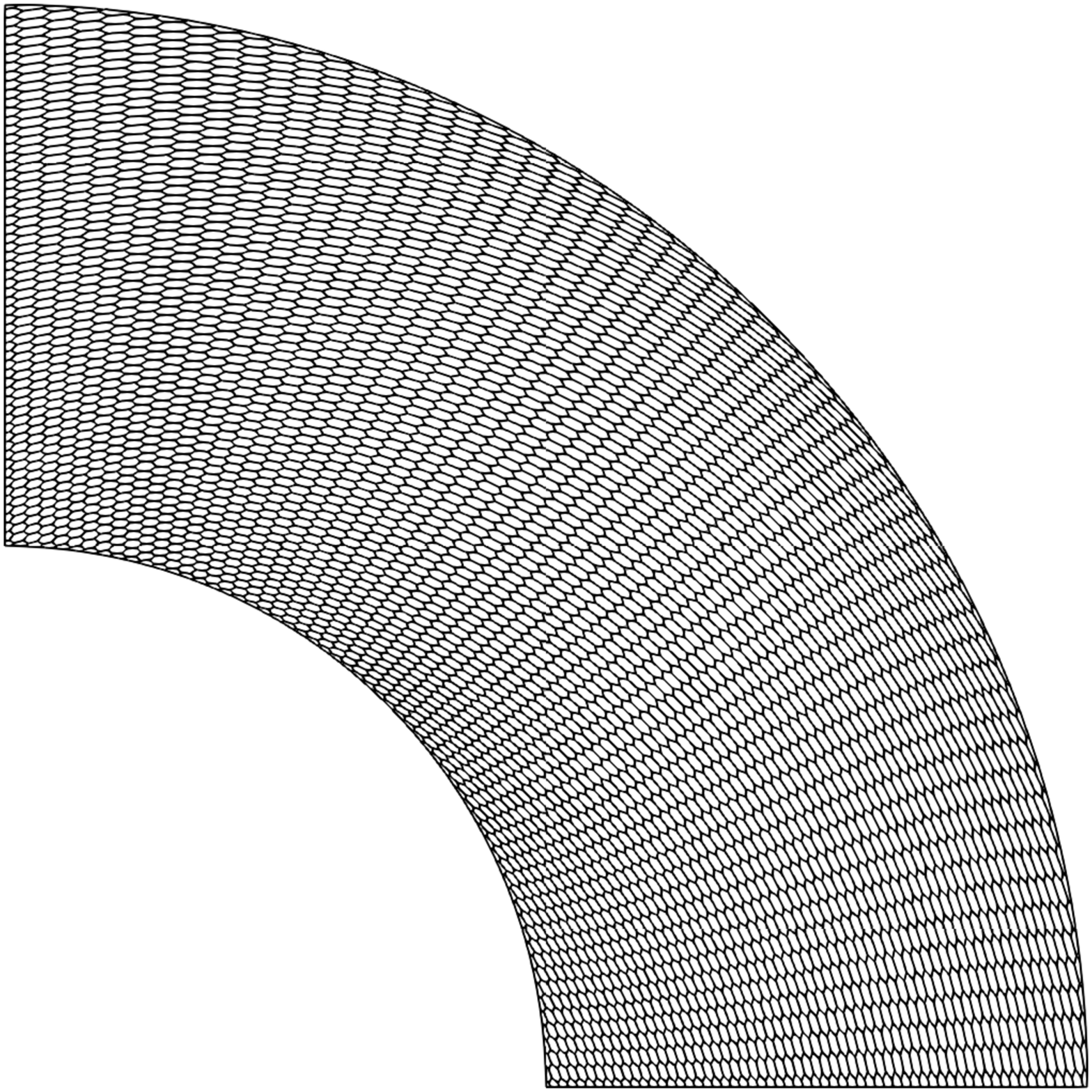, width = 0.185\textwidth}}
\subfigure[]{\label{fig:thickwalled_mesh_e} \epsfig{file = 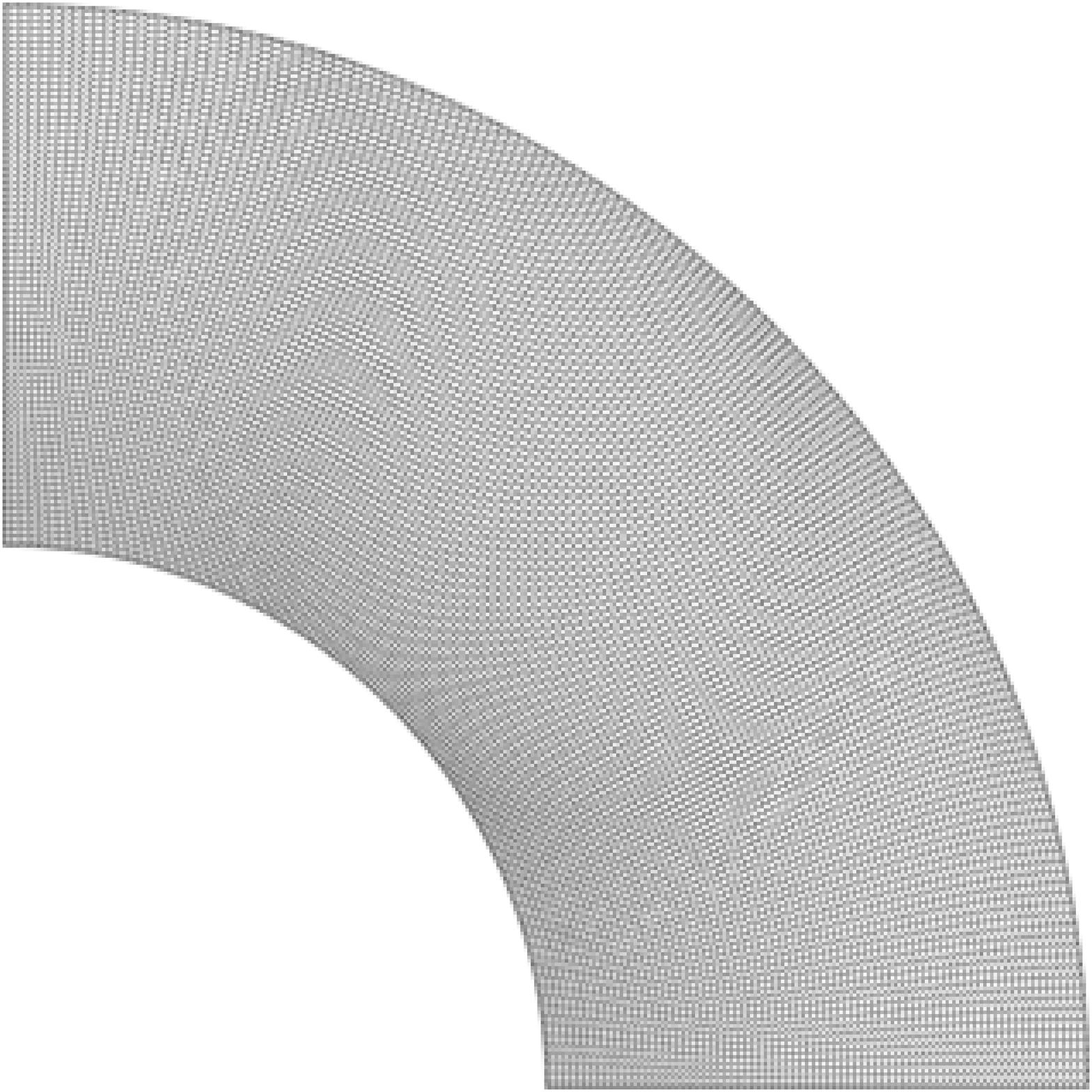, width = 0.185\textwidth}}
}
\caption{Sequence of NIVED background integration meshes used for the analysis of the thick-walled viscoelastic 
         cylinder subjected to internal pressure.}
\label{fig:thickwalled_mesh}
\end{figure}

The convergence of the NIVED radial displacement at control points $A$ and $B$ is summarized
in Table~\ref{table:thickwalled_disp_convergence} for the sequence of background
integration meshes depicted in~\fref{fig:thickwalled_mesh}. Good agreement
with the FEM-Q9 reference solution is observed for the three sets of material parameters.

\begin{table}[!tbhp]
\setlength{\arrayrulewidth}{.1em}
\begin{center}
\small\addtolength{\tabcolsep}{-1pt}
\caption{NIVED radial displacements at points $A$ and $B$ for the thick-walled viscoelastic cylinder problem at time $t=20$.}
\vspace*{0pt}         
\renewcommand{\arraystretch}{1.5}
\begin{tabular}{cccccccc}
\hline
  \multicolumn{2}{r}{$(\mu_0,\mu_1)=$}   & \multicolumn{2}{c}{$(0.7,0.3)$} & \multicolumn{2}{c}{$(0.3,0.7)$} & \multicolumn{2}{c}{$(0.01,0.99)$} \\ 
\hline
mesh                          & \# dofs & $u_A$       & $u_B$       & $u_A$       & $u_B$       & $u_A$       & $u_B$       \\ 
\hline
\fref{fig:thickwalled_mesh_a} & 882     & 0.053012    & 0.031904    & 0.118830    & 0.065142    & 0.646690    & 0.330050    \\
\fref{fig:thickwalled_mesh_b} & 1922    & 0.053077    & 0.031950    & 0.118970    & 0.065219    & 0.647820    & 0.330170    \\
\fref{fig:thickwalled_mesh_c} & 3362    & 0.053101    & 0.031966    & 0.119020    & 0.065246    & 0.648160    & 0.330210    \\
\fref{fig:thickwalled_mesh_d} & 7442    & 0.053118    & 0.031977    & 0.119050    & 0.065264    & 0.648370    & 0.330240    \\
\fref{fig:thickwalled_mesh_e} & 26082   & 0.053131    & 0.031983    & 0.119070    & 0.065272    & 0.648440    & 0.330210    \\ 
\hline
Ref. (FEM-Q9)                & 79242   & 0.053584    & 0.032425    & 0.119970    & 0.066177    & 0.652880    & 0.334810    \\ 
\hline
\end{tabular}
\label{table:thickwalled_disp_convergence}
\end{center}
\end{table}

\begin{figure}[!tbhp]
\centering
\mbox{
\subfigure[]{\label{fig:thickwalled_disp_curves_a}
\epsfig{file = 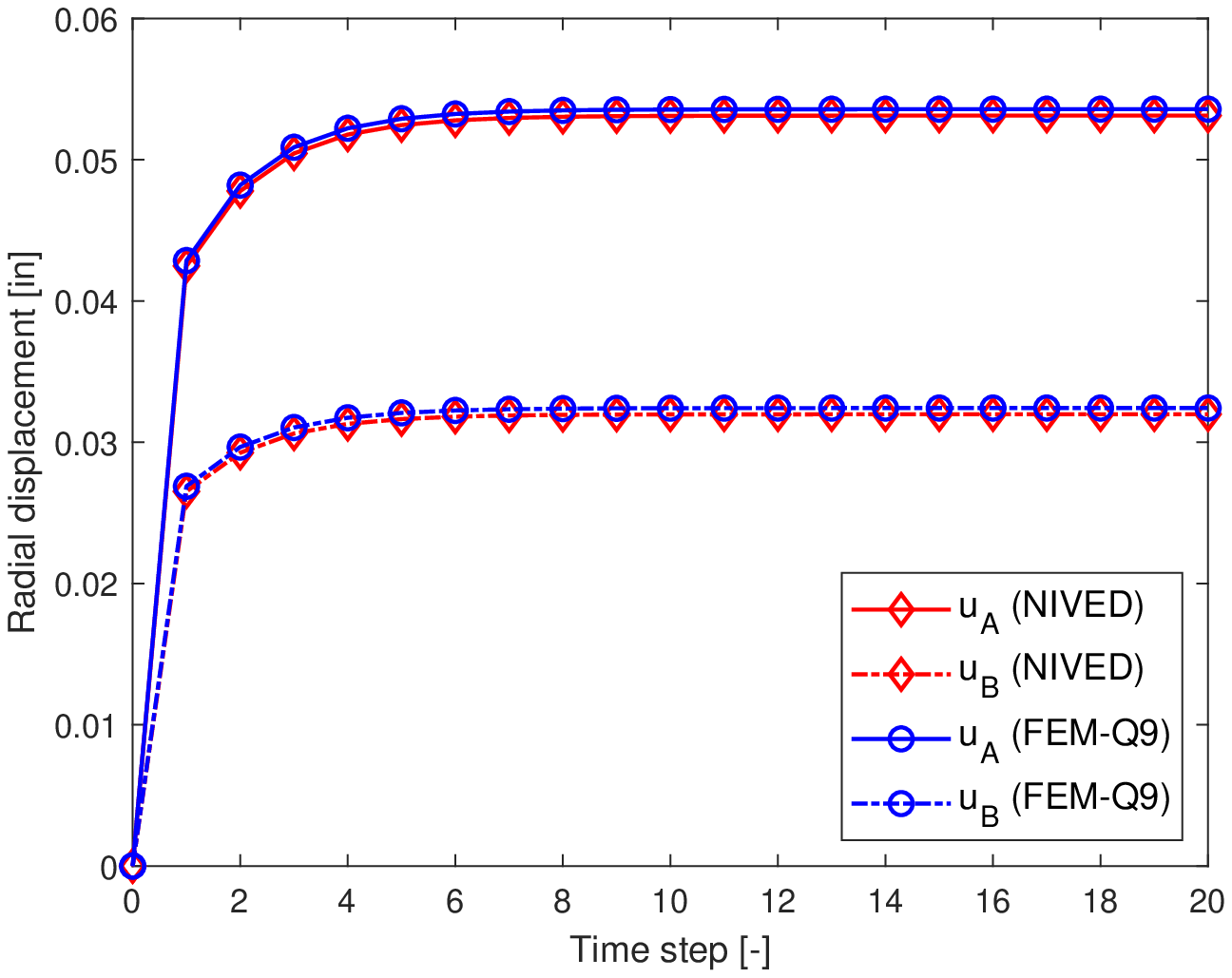, width = 0.5\textwidth}}
\subfigure[]{\label{fig:thickwalled_disp_curves_b}
\epsfig{file = 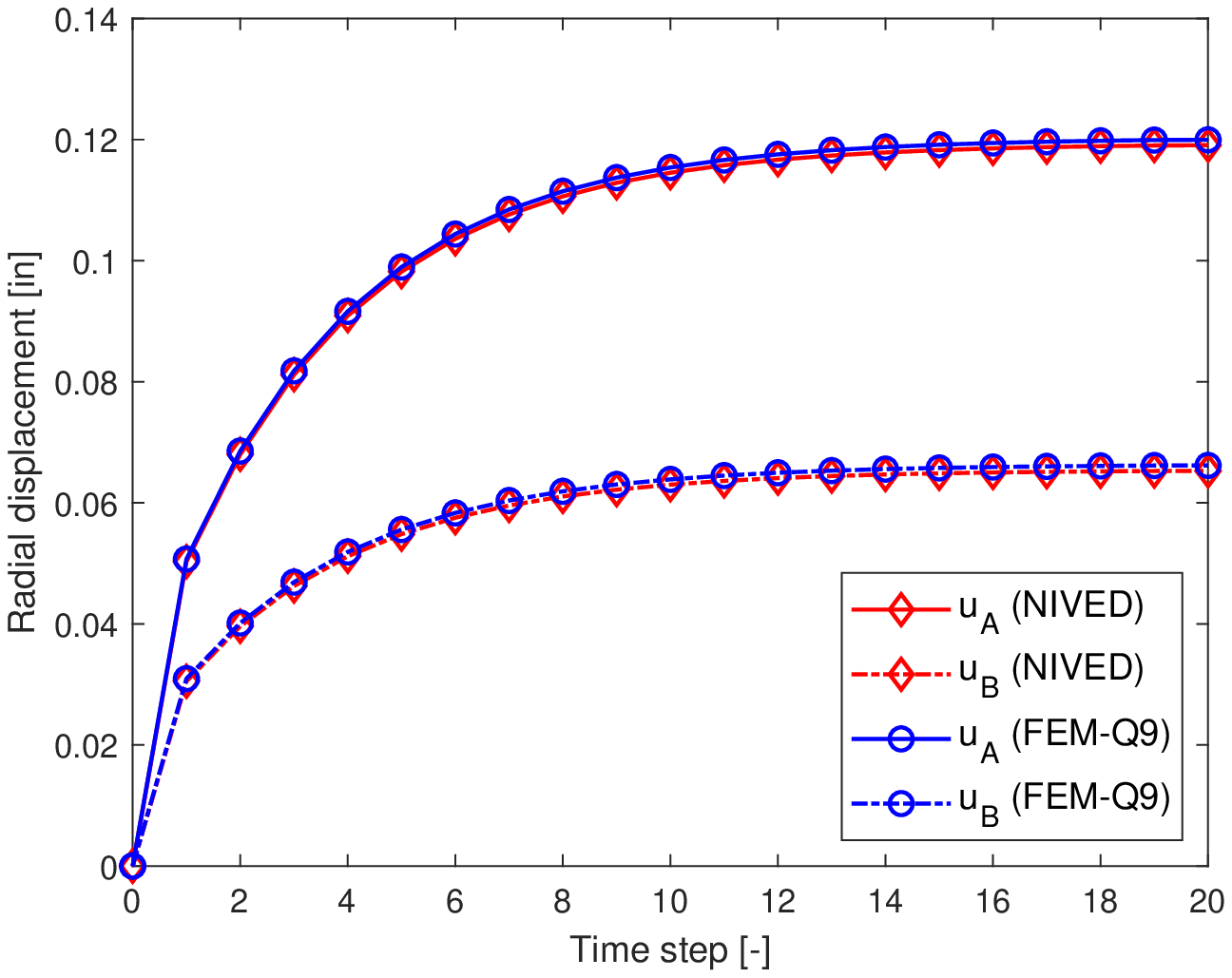, width = 0.5\textwidth}}
}
\subfigure[]{\label{fig:thickwalled_disp_curves_c}
\epsfig{file = 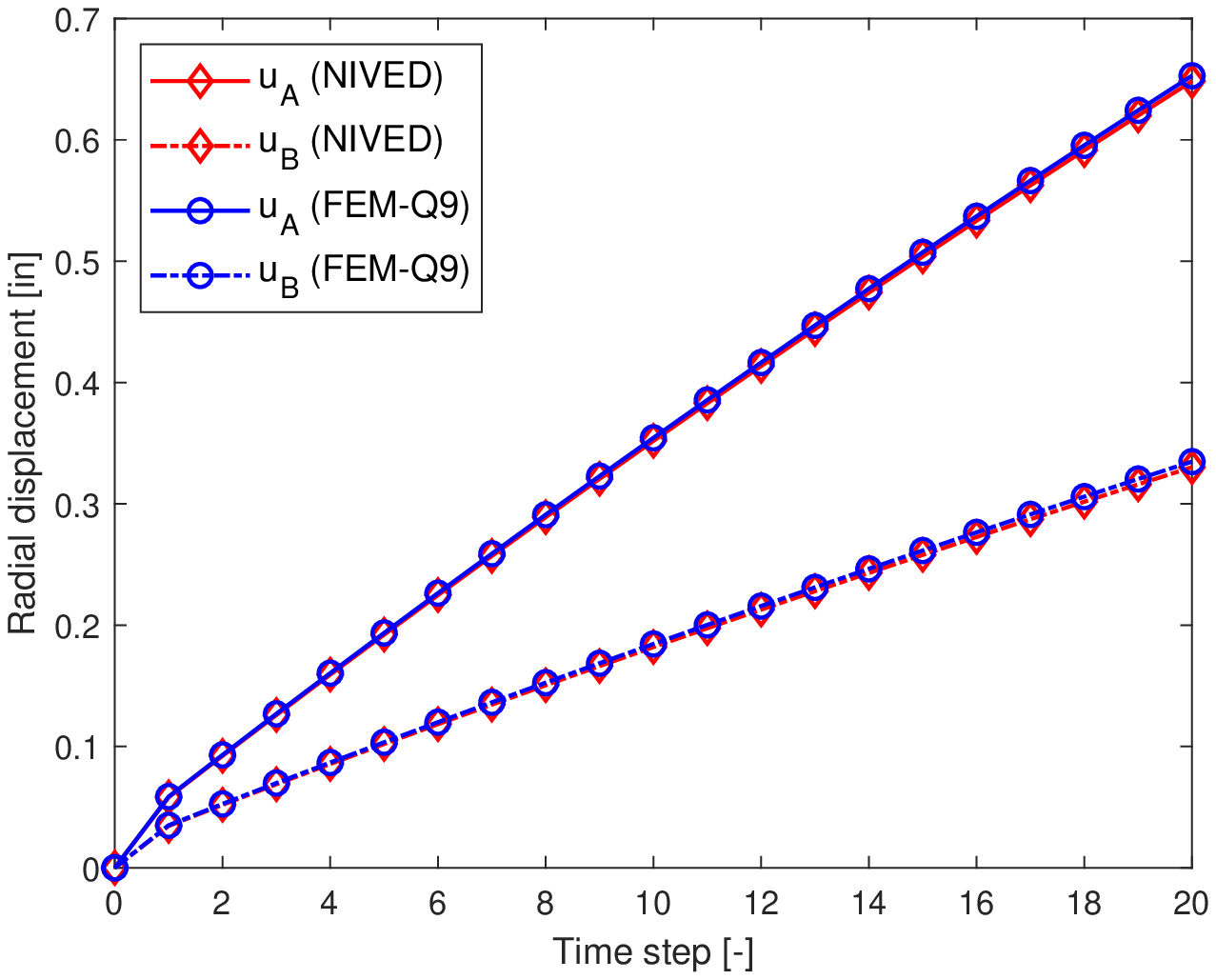, width = 0.5\textwidth}}
\caption{Thick-walled viscoelastic cylinder subjected to internal pressure.
         Time history of the radial displacement measured at control 
         points $A$ (upper curves) and $B$ (lower curves) for three sets of viscoelastic
         parameters: (a) $(\mu_0,\mu_1)=(0.7,0.3)$, (b) $(\mu_0,\mu_1)=(0.3,0.7)$ and
         (c) $(\mu_0,\mu_1)=(0.01,0.99)$.}
\label{fig:thickwalled_disp_curves}
\end{figure}

Curves depicting the time history of the radial displacement measured at control 
points $A$ and $B$ are shown in~\fref{fig:thickwalled_disp_curves} for the three 
sets of material parameters. The response curves for the NIVED scheme is 
plotted for the mesh shown in~\fref{fig:thickwalled_mesh_e}. Good agreement 
with the FEM-Q9 reference solutions is observed.

\begin{figure}[!tbhp]
\centering
\subfigure[]{\label{fig:thickwalled_quad_mesh_a} \epsfig{file = 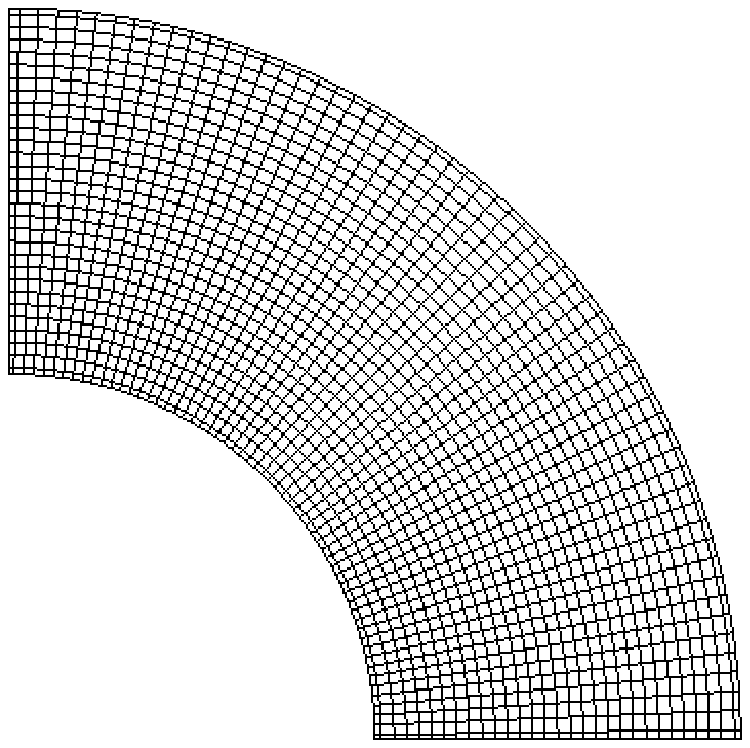, width = 0.35\textwidth}}
\subfigure[]{\label{fig:thickwalled_quad_mesh_b} \epsfig{file = 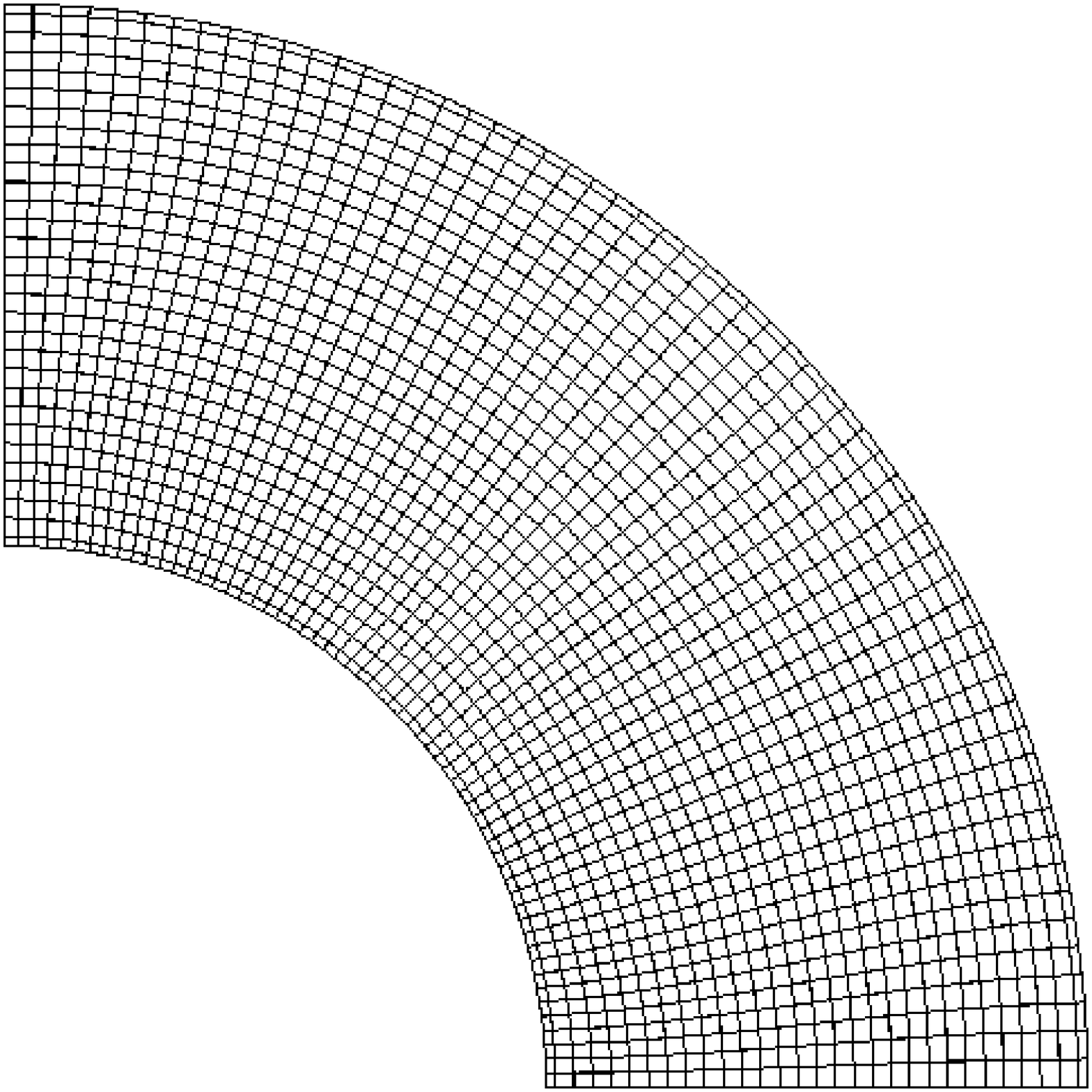, width = 0.35\textwidth}}
\caption{Integration meshes to study the performance of the numerical solutions of the 
         thick-walled viscoelastic cylinder problem in the near-incompressible limit. 
         Meshes for (a) NIVED scheme, and (b) FEM-Q4 and FEM-Q9 approaches.}
\label{fig:thickwalled_quad_mesh}
\end{figure}

\begin{figure}[!tbhp]
\centering
\mbox{
\subfigure[]{\label{fig:thickwalled_incomp_a} \epsfig{file = 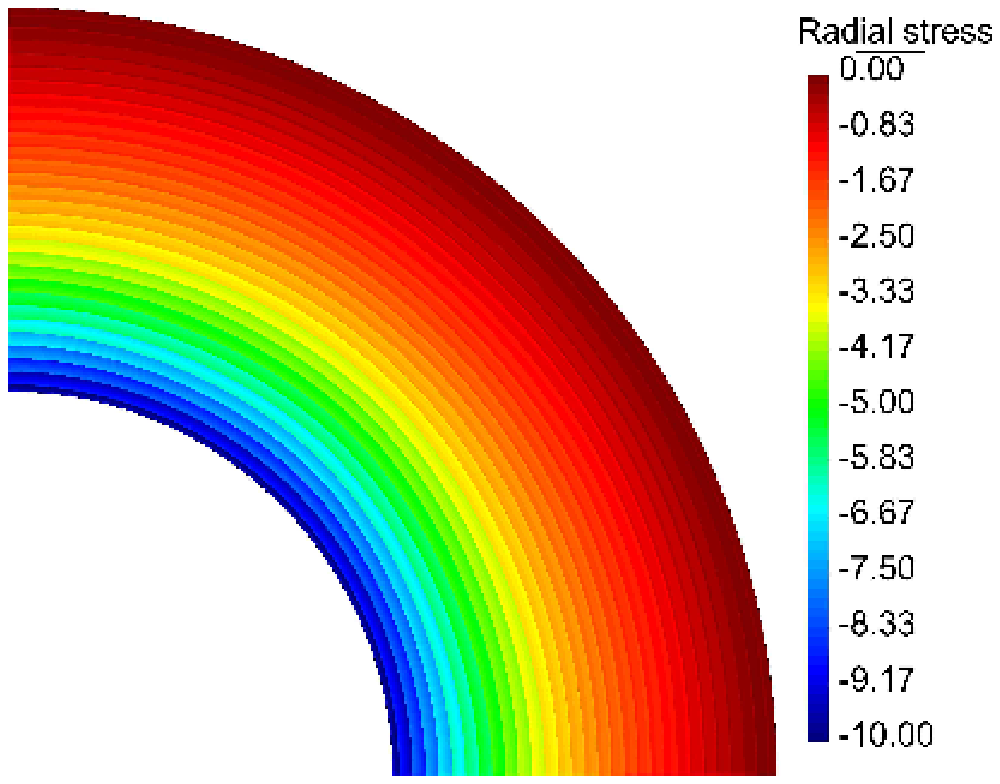, width = 0.32\textwidth}}
\subfigure[]{\label{fig:thickwalled_incomp_b} \epsfig{file = 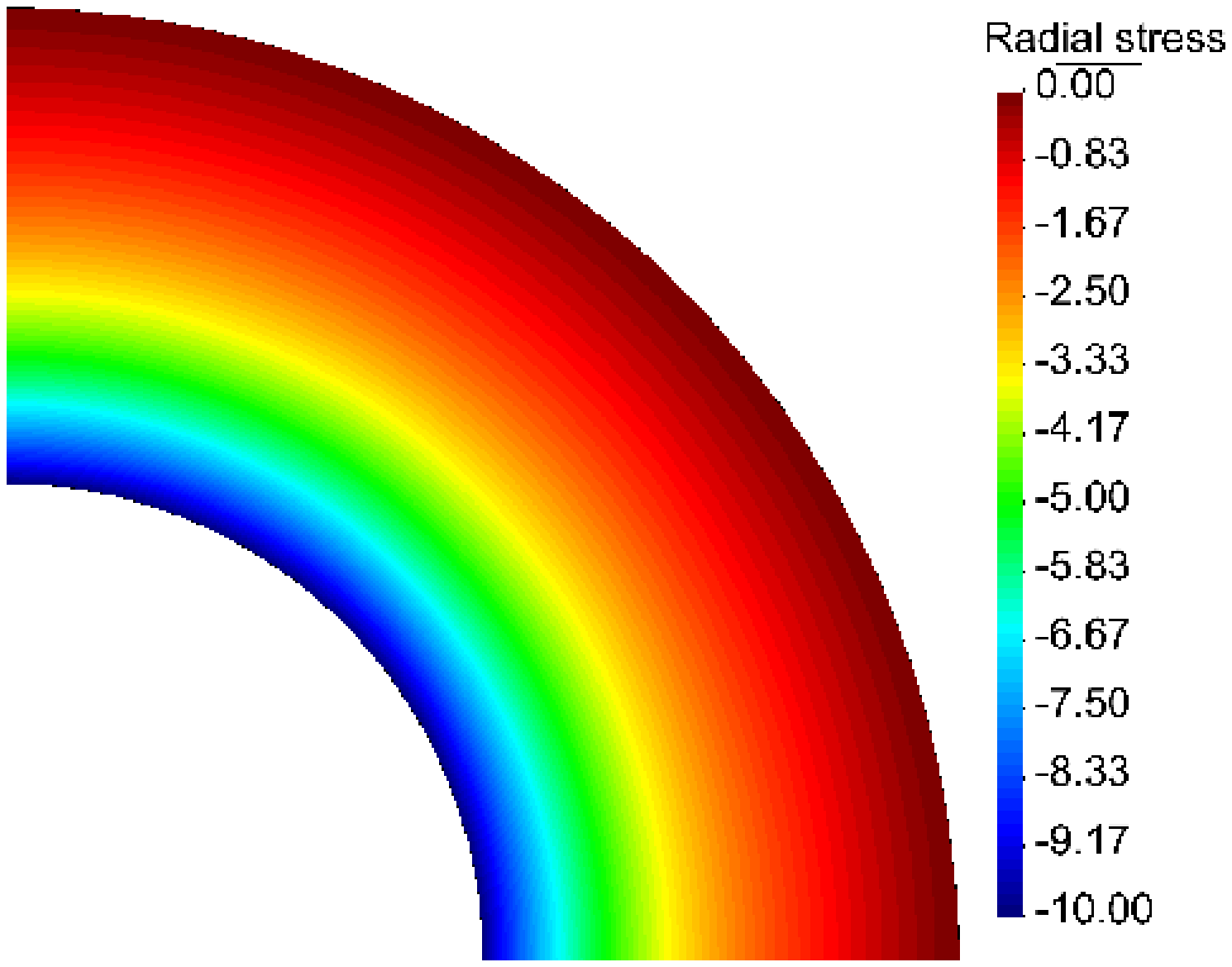, width = 0.32\textwidth}}
\subfigure[]{\label{fig:thickwalled_incomp_c} \epsfig{file = 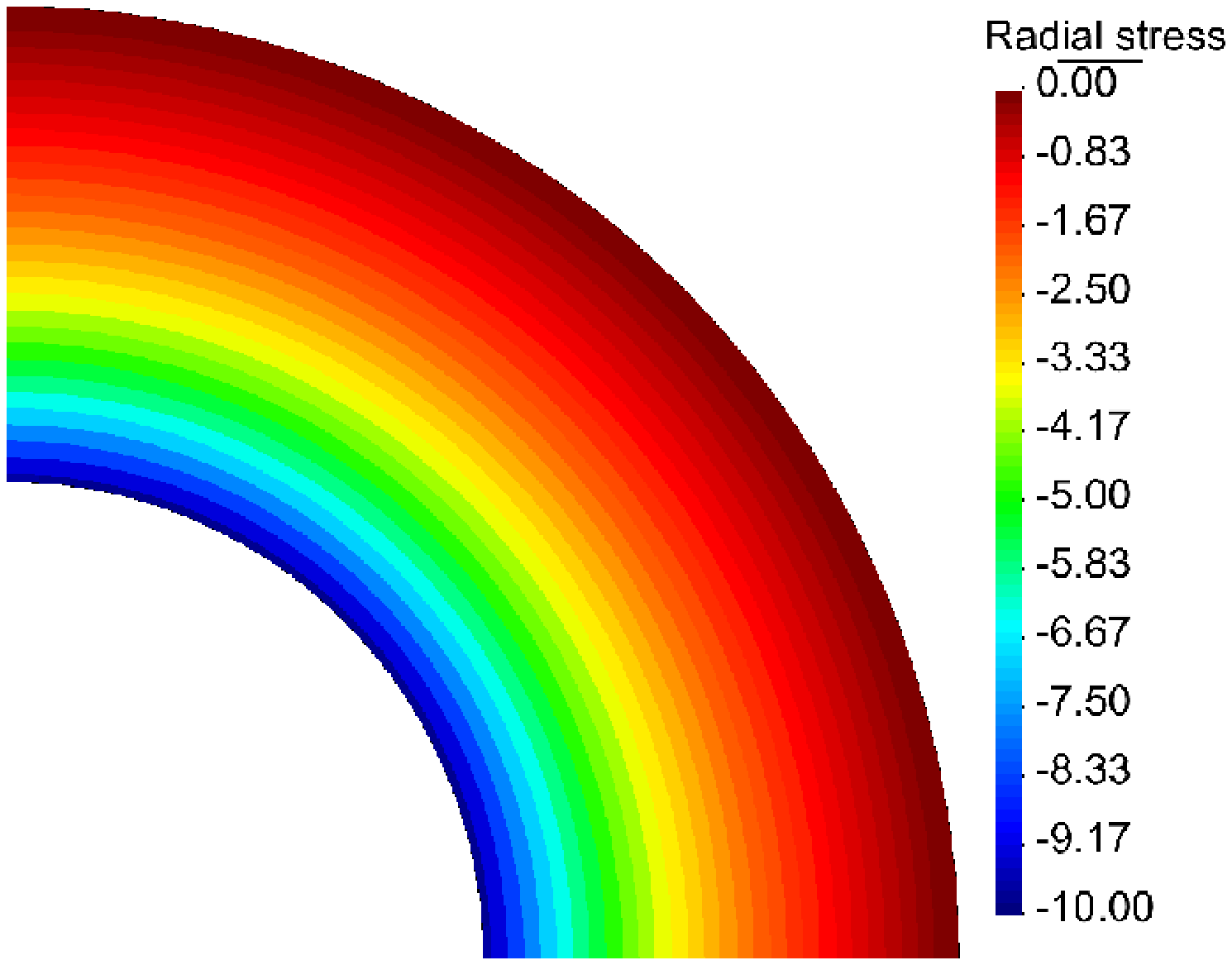, width = 0.32\textwidth}}
}
\mbox{
\subfigure[]{\label{fig:thickwalled_incomp_d} \epsfig{file = 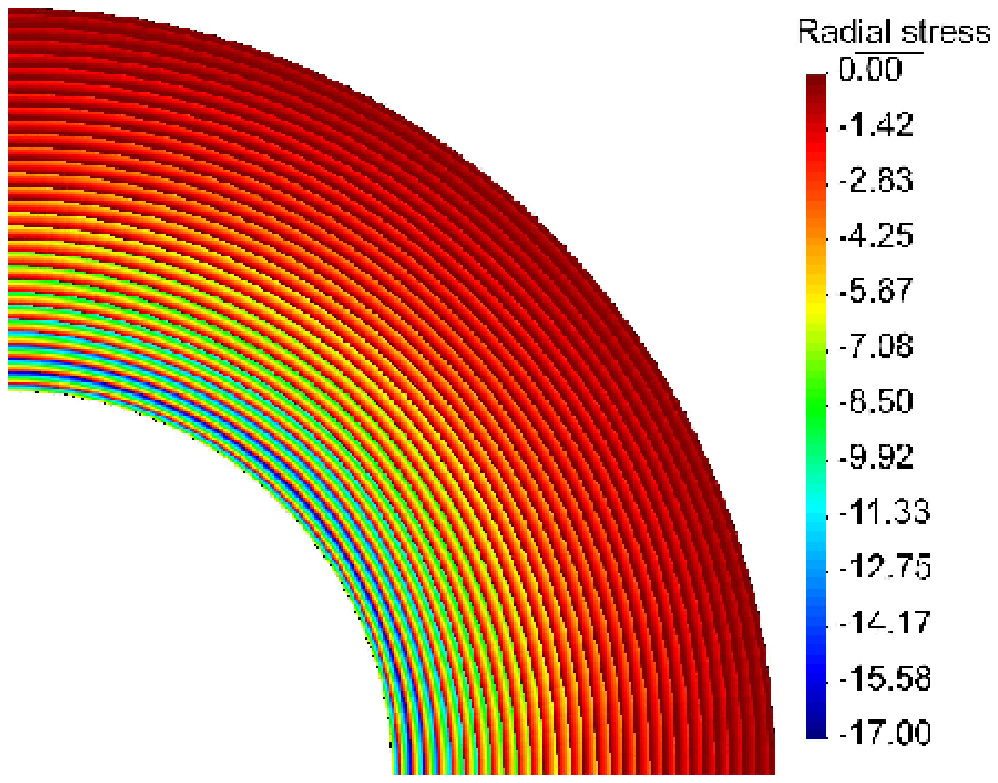, width = 0.32\textwidth}}
\subfigure[]{\label{fig:thickwalled_incomp_e} \epsfig{file = 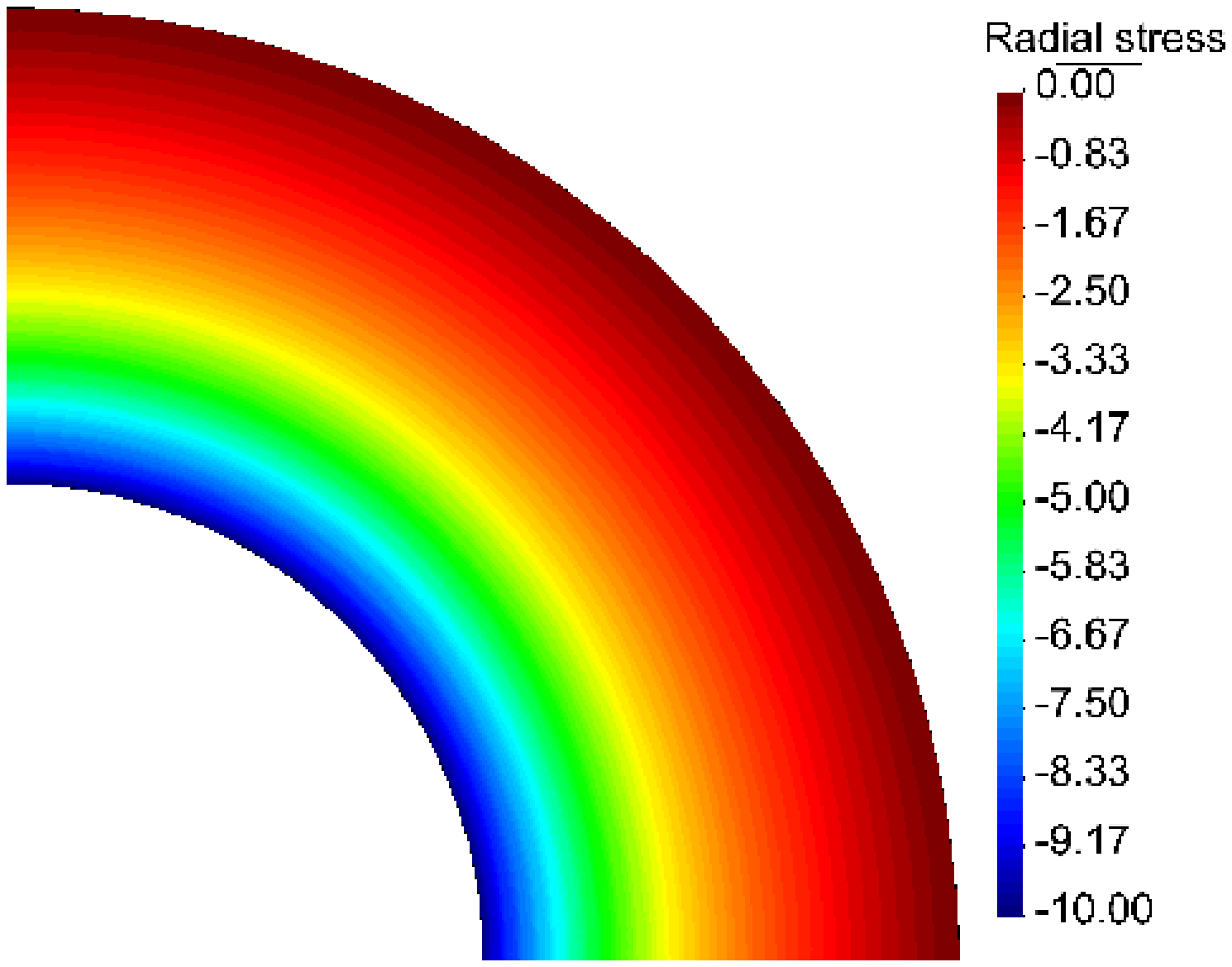, width = 0.32\textwidth}}
\subfigure[]{\label{fig:thickwalled_incomp_f} \epsfig{file = 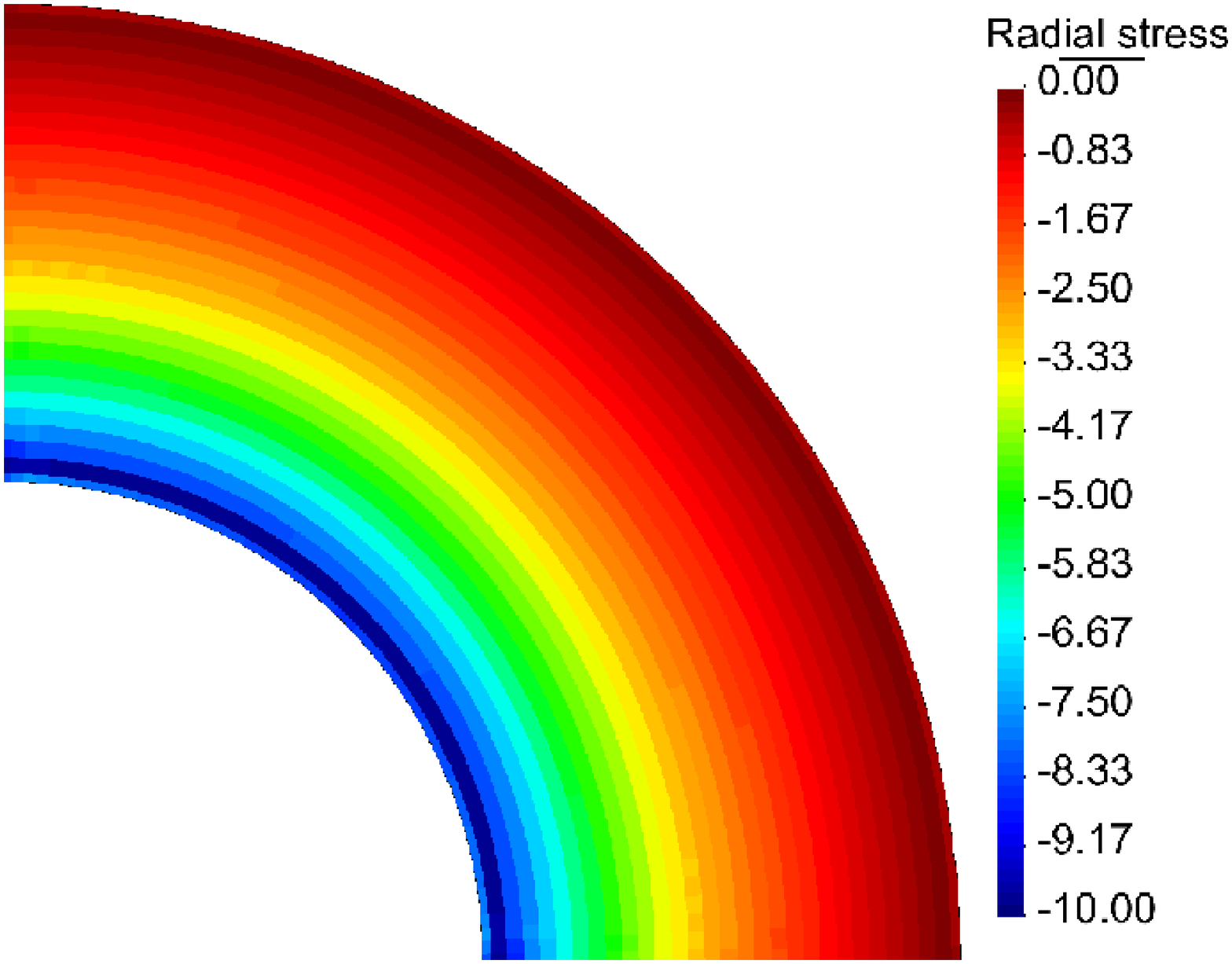, width = 0.32\textwidth}}
}
\caption{Thick-walled viscoelastic cylinder subjected 
         to internal pressure. Comparison 
         of radial stress contour plots using FEM and NIVED schemes.
         For viscoelastic parameters $(\mu_0,\mu_1)=(0.7,0.3)$:
         (a) FEM-Q4, (b) FEM-Q9, and (c) NIVED. 
         For viscoelastic parameters $(\mu_0,\mu_1)=(0.01,0.99)$:
         (d) FEM-Q4, (e) FEM-Q9, and (f) NIVED.}
\label{fig:thickwalled_incomp}
\end{figure}

Lastly, the performance of the NIVED method for
near-incompressible material behavior is studied. The NIVED approach is tested on the background 
mesh shown in~\fref{fig:thickwalled_quad_mesh_a} and its radial stress field solution
is compared with the corresponding FEM-Q9 solution on the 
mesh shown in~\fref{fig:thickwalled_quad_mesh_b}. Also on this mesh, the radial
stress field is computed with 4-node FEM quadrilateral elements (FEM-Q4). The comparisons
are summarized in~\fref{fig:thickwalled_incomp}. The radial stress field corresponding
to the compressible case $(\mu_0,\mu_1)=(0.7,0.3)$ ($\nu_{\mathrm{eff}}=0.3542$) is shown 
in Figures~\ref{fig:thickwalled_incomp_a}--\subref{fig:thickwalled_incomp_c} for
FEM-Q4, FEM-Q9 and NIVED methods, respectively. It is observed that the three methods
perform well in this case. Of course, because of its higher approximation order,
the FEM-Q9 radial stress field is smoother than the FEM-Q4 and NIVED radial stress fields.

The radial stress field corresponding to the near-incompressible case 
$(\mu_0,\mu_1)=(0.01,0.99)$ ($\nu_{\mathrm{eff}}=0.4977$) is
depicted in Figures~\ref{fig:thickwalled_incomp_d}--\subref{fig:thickwalled_incomp_f} for
FEM-Q4, FEM-Q9 and NIVED methods, respectively. In this case, the FEM-Q9 radial stress 
field looks smooth and no sign of locking is observed. However, locking behavior
is exhibited in the FEM-Q4 case as evidenced by the oscillatory radial stress field.
The NIVED scheme performs better than the FEM-Q4 in the incompressible limit, but
some small oscillations do occur in the radial stress field and thus it is not completely 
free of locking. This is not surprising since use of the virtual element framework 
in the NIVED approach ensures consistency (patch test satisfaction) and stability 
for a displacement-based formulation. To render the NIVED scheme locking-free in the 
near-incompressible limit would require introducing, for instance, a form of 
a nodal averaging technique~\cite{dohrmann:2000:NBU,puso:2006:ASN}.

\section{CONCLUDING REMARKS}\label{sec:conclusions}

In this paper, we proposed a novel nodal integration scheme for meshfree
Galerkin methods that is devised from the virtual 
element decomposition~\cite{BeiraoDaVeiga-Brezzi-Cangiani-Manzini-Marini-Russo:2013}, 
wherein the bilinear form is decomposed into a consistency part and 
a stability part that ensure that the method is consistent and stable. 
Linear maximum-entropy meshfree basis functions were adopted, but the formulation 
is applicable to any other linear meshfree approximant. We referred to this 
new nodal integration scheme as NIVED. As in any nodal integration
method, a nodal representative cell is needed in the NIVED approach.
The nodal cell can be constructed from a Voronoi diagram
or a Delaunay triangulation by connecting the centroids 
of triangles surrounding a node.

Numerical tests in linear elastostatic and linear elastodynamic boundary-value
problems were conducted, where the focus was on comparing 
the performance of the NIVED method and the 
maximum-entropy meshfree method (MEM) using standard Gauss integration.
We also provided the extension of the NIVED approach to
nonlinear analysis where the material constitution is a linear isotropic 
solid described by the Generalized Maxwell viscoelastic model.
Comparisons with reference solutions obtained with 
4-node and 9-node FEM quadrilateral elements
were conducted. In the NIVED approach, the weak 
form integrals are integrated by sampling them at the nodes and the
integration at a node is performed using a 1-point Gauss rule over 
the edges of its nodal representative cell. This task only involves 
the evaluation of basis functions (no derivatives are needed). 
In the MEM approach, the integration of the weak form integrals
is performed using standard Gauss integration on the interior 
of triangular cells, which requires the evaluation of
basis functions derivatives. 

Our main findings from this work are as follows. 
The NIVED scheme passes the patch test to machine 
precision, whereas the MEM does not achieve this level of 
accuracy due to integration errors. The numerical stability
test showed that both the NIVED and MEM methods deliver the three zero-energy 
rigid-body modes, but instability is exhibited only by the MEM method
as evidenced by the presence of the nonsmooth eigenmodes that follow 
the three rigid-body modes. Convergence was assessed through 
several numerical experiments, which included a
cantilever beam subjected to a parabolic end load, an infinite
plate with a circular hole and manufactured (elastostatic and
elastodynamic) problems. The convergence studies revealed that the NIVED method 
delivers the optimal rate of convergence in the $L^2$ norm
and the $H^1$ seminorm. On the other hand, the convergence of the MEM
was dependent on the number of Gauss points used inside
the triangular integration cell. In the tests that were conducted, 
the MEM required a 3-point Gauss rule for optimal convergence 
in the $L^2$ norm and, depending on the problem, a 6-point or 
a 12-point Gauss rule for optimal convergence in the $H^1$ seminorm.
In terms of computational cost, it was shown 
that for the same number of degrees of freedom
the proposed NIVED scheme is faster but less accurate 
than a convergent MEM approach. 
We also tested the NIVED and MEM methods in an $L$-shaped domain
under traction load, where the presence of a re-entrant corner introduces 
a stress singularity that results in a nonsmooth solution.
In this case, a study to assess the convergence to a reference value of 
the strain energy was conducted. The results showed that 
only the NIVED scheme uniformly approaches the reference 
strain energy. The performance of the nonlinear NIVED formulation
was assessed by solving the problem of a thick-walled viscoelastic 
cylinder subjected to internal pressure. The radial displacement 
solution was found to be in good agreement with a reference solution 
obtained using quadratic (9-node) quadrilateral finite elements.
We also showed that for a near-incompressible viscoelastic
material, the NIVED scheme is not completely devoid of locking.

In closing, we mention that a desirable feature that
can be offered by the (nodally integrated) NIVED approach 
is that state variables such as strains, stresses and other internal 
variables in nonlinear computations can be stored at 
the nodes, which is attractive for Lagrangian large deformation meshfree simulations. 
This extension of the NIVED method to large deformations is appealing, and is planned 
as part of future work.

\section*{ACKNOWLEDGEMENTS}
RSV and AOB acknowledge the support provided by the Chilean National Fund for Scientific
and Technological Development (FONDECYT) through grant CONICYT/FONDECYT No. 1181192.
NS acknowledges the research support of Sandia National Laboratories.
EA gratefully acknowledges the partial financial support of the University of Rome
Tor Vergata Mission Sustainability Programme through project SPY-E81I18000540005; 
and the partial financial support of PRIN 2017 project ``3D PRINTING: A BRIDGE TO 
THE FUTURE (3DP\_Future). Computational methods, innovative applications, experimental
validations of new materials and technologies,'' grant 2017L7X3CS\_004.
NHK is grateful for the support provided by the Chilean National Fund for Scientific
and Technological Development (FONDECYT) through grant CONICYT/FONDECYT No. 1181506.

\renewcommand{\refname}{REFERENCES}
\bibliographystyle{unsrt}

\end{document}